%% file: LT.tex
\documentclass{memo-l}

\usepackage{amsxtra,amssymb,amsmath,url,amsrefs, pdfsync, mathrsfs,enumerate}
\usepackage{array}
\usepackage[utf8]{inputenc}
\usepackage{scrtime}
\usepackage{xcolor}
\definecolor{cite}{rgb}{0.30,0.60,1.00}
\definecolor{url}{rgb}{1.00,0.10,0.80}
\definecolor{link}{rgb}{0.00,0.00,1.00}
\usepackage[colorlinks,linkcolor=link,urlcolor=url,citecolor=cite,breaklinks]{hyperref}
\hypersetup{
pdfstartpage=1,
pdfstartview=FitH}

\usepackage{charter}

\allowdisplaybreaks

\usepackage{multirow}
\usepackage{makecell}
\usepackage{booktabs}
\usepackage{float}

\usepackage{scalerel,stackengine}
\usepackage{pict2e}
\usepackage{longtable}
\stackMath

\newcommand\widecheck[1]{%
\savestack{\tmpbox}{\stretchto{%
  \scaleto{%
    \scalerel*[\widthof{\ensuremath{#1}}]{\kern-.4pt\bigwedge\kern-.4pt}%
    {\rule[-\textheight/2]{1ex}{\textheight}}
  }{\textheight}%
}{0.5ex}}%
\stackon[2pt]{#1}{\scalebox{-1}{\tmpbox}}%
}

\numberwithin{section}{chapter}
\numberwithin{equation}{chapter}

\usepackage{color}

\newcommand{\ud}{\mathrm{d}}
\newcommand{\ue}{\mathrm{e}}

\newcommand{\tr}{\mathrm{tr}}
\newcommand{\Tr}{\mathrm{Tr}}

\newcommand{\Frob}{\mathrm{Frob}}

\newcommand{\sgn}{\mathrm{sgn}}

\newcommand{\GL}{\mathrm{GL}}
\newcommand{\Gal}{\mathrm{Gal}}
\newcommand{\ord}{\mathrm{ord}}
\newcommand{\Aut}{\mathrm{Aut}}
\newcommand{\End}{\mathrm{End}}
\newcommand{\GSp}{\mathrm{GSp}}

\DeclareMathOperator{\Mod}{mod}

\renewcommand{\bmod}[1]{\,(\Mod{ #1})}

\newcommand{\F}{\mathbb{F}}
\newcommand{\N}{\mathbb{N}}
\newcommand{\Q}{\mathbb{Q}}
\newcommand{\R}{\mathbb{R}}
\newcommand{\Z}{\mathbb{Z}}
\newcommand{\C}{\mathbb{C}}

\newcommand{\cA}{\mathcal{A}}
\newcommand{\cB}{\mathcal{B}}
\newcommand{\cC}{\mathcal{C}}

\newcommand{\cE}{\mathcal{E}}

\newcommand{\cG}{\mathcal{G}}

\newcommand{\cJ}{\mathcal{J}}

\newcommand{\cL}{\mathcal{L}}

\newcommand{\cN}{\mathcal{N}}
\newcommand{\cP}{\mathcal{P}}
\newcommand{\cQ}{\mathcal{Q}}
\newcommand{\cR}{\mathcal{R}}

\newcommand{\cT}{\mathcal{T}}

\newcommand{\fc}{\mathfrak{c}}
\newcommand{\ff}{\mathfrak{f}}
\newcommand{\fh}{\mathfrak{h}}
\newcommand{\fm}{\mathfrak{m}}

\newcommand{\fp}{\mathfrak{p}}

\newcommand{\fA}{\mathfrak{A}}

\newcommand{\fC}{\mathfrak{C}}
\newcommand{\fD}{\mathfrak{D}}
\newcommand{\fF}{\mathfrak{F}}
\newcommand{\fG}{\mathfrak{G}}

\newcommand{\fS}{\mathfrak{S}}
\newcommand{\fU}{\mathfrak{U}}

\theoremstyle{plain}
\newtheorem{theorem}{Theorem}[chapter]
\newtheorem*{theorem*}{Theorem}
\newtheorem{lemma}[theorem]{Lemma}
\newtheorem{corollary}[theorem]{Corollary}
\newtheorem{conjecture}[theorem]{Conjecture}

\newtheorem{proposition}[theorem]{Proposition}

  \newenvironment{theorembis}[1]
  {\renewcommand{\thetheorem}{\ref{#1}$'$}%
   \addtocounter{theorem}{-1}%
   \begin{theorem}}
  {\end{theorem}}

\theoremstyle{remark}

\theoremstyle{definition}

\newtheorem{definition}[theorem]{Definition}

\newtheorem{remark}[theorem]{Remark}

\usepackage{tikz}
\newcommand\cube{\begin{tikzpicture}[scale=1.8]
    \coordinate (A1) at (0, 0);
    \coordinate (A2) at (0, 0.1);
    \coordinate (A3) at (0.1, 0.1);
    \coordinate (A4) at (0.1, 0);
    \coordinate (B1) at (0.03, 0.03);
    \coordinate (B2) at (0.03, 0.13);
    \coordinate (B3) at (0.13, 0.13);
    \coordinate (B4) at (0.13, 0.03);

    \draw (A1) -- (A2);
    \draw (A2) -- (A3);
    \draw (A3) -- (A4);
    \draw (A4) -- (A1);
    \draw[densely dotted] (A1) -- (B1);
    \draw[densely dotted] (B1) -- (B2);
    \draw (A2) -- (B2);
    \draw (B2) -- (B3);
    \draw (A3) -- (B3);
    \draw (A4) -- (B4);
    \draw (B4) -- (B3);
    \draw[densely dotted] (B1) -- (B4);
\end{tikzpicture}}

\begin{document}
\frontmatter

\title[Lang--Trotter Conjecture for CM Elliptic Curves]{Lang--Trotter Conjecture for CM Elliptic Curves}

 \author{Daqing Wan}

\address{Department of Mathematics, University of California at Irvine, 340 Rowland Hall, Irvine, CA, 92697, USA}

\email{dwan@math.uci.edu}

\author{Ping Xi}

\address{School of Mathematics and Statistics, Xi'an Jiaotong University, Xi'an 710049, P. R. China}

\email{ping.xi@xjtu.edu.cn}




\thanks{Ping XI is supported in part by NSFC (No. 12025106, No. 11971370) and
Shaanxi NSF (No. 2025JC-QYCX-002).}

\begin{abstract}

Given an elliptic curve $E$ over $\Q$ and non-zero integer $r$, the Lang--Trotter conjecture predicts a striking asymptotic formula for the number of good primes $p\leqslant x$, denoted by $\pi_{E,r}(x)$, such that the Frobenius trace of $E$ at $p$ is equal to the given integer $r$. We focus on the CM case in this memoir, and show how to realize the following two goals:
\begin{itemize}
\item to give an unconditional estimate for $\pi_{E,r}(x)$, which confirms the upper bound part of the conjecture up to a constant multiple;
\item to give a conditional explicit asymptotic formula for $\pi_{E,r}(x)$ based on the Hardy--Littlewood conjecture on primes presented by quadratic polynomials.
\end{itemize}

The merit of the above two results can be reflected as follows:
\begin{itemize}
\item We obtain an unconditional upper bound for $\pi_{E,r}(x)$, which is of the correct order of magnitude as suggested by the Lang--Trotter conjecture. This has not been achieved before for any special cases.
\item The comparison of our conditional asymptotic formula with the 
Lang--Trotter conjecture suggests a new constant comparison conjecture:  
our constant defined analytically should be equal to the Lang--Trotter constant defined algebraically in terms of Galois images. 
The significance of this comparison conjecture is that it gives a quick explicit 
way to understand and compute the much harder Lang--Trotter constant in terms of 
our constant, which is much easier and completely explicit. For instance, 
it allows us to explicitly classify when the constant is positive and 
when the constant is symmetric in terms of $r$.

\item Assuming the Hardy--Littlewood conjecture, we propose and prove several more conjectures on the distribution of rational primes that split in the imaginary quadratic field $\Q(\sqrt{-D})$ with additional arithmetic constraints. 
These new conjectures are the key in the proof of our version of the Lang--Trotter conjecture for CM elliptic curves. They could have potential further applications in arithmetics of CM elliptic curves. 
\end{itemize}

For completeness, we also summarize classical results on quadratic, cubic and quartic residues, as well as the corresponding reciprocity laws. This part should be of independent interests and could provide useful materials for more junior readers.
We also highlight some possible extensions of the arguments in this memoir that may work for other statistical problems of CM elliptic curves.

\vspace{\baselineskip}
\begin{flushright}\noindent
Irvine and Xi'an,\hfill {\it Daqing Wan}\\
October 2025\hfill {\it Ping Xi}\\
\end{flushright}

\end{abstract}

\maketitle

\setcounter{tocdepth}{1}

\tableofcontents

\input{Notation}

\mainmatter

\input{1.Background}
\input{2.HLupperbound}

\input{3.Residuesandreciprocity}
\input{4.D=1}

\input{5.D=3}

\input{6.otherCM}

\input{7.LTconstant}

\input{8.Historicalnotesperspectives}

\input{9.Numericaldata}


\include{appendix}

\backmatter


\end{document}

%% file: Notation.tex
\chapter*{Notation}

\begin{itemize}
\item $\Z=\{0,\pm1,\pm2,\pm3,\cdots\},$ $\Z^+=\{1,2,3,\cdots\},$ $\N=\Z^+\cup\{0\},$ $\Z^*=\Z\setminus\{0\}.$
\item $\C$ denotes the field of complex numbers.
\item $\ue(z):=\ue^{2\pi iz}$ for $z\in\C.$
\item $\varphi$ denotes the Euler totient function over $\Z$, i.e., $\varphi(n)=|(\Z/n\Z)^\times|$.
\item $(a,b)$ denotes the greatest common factor of $a,b\in\Z$ or a tuple given by $a$ and $b$.
\item $n^\infty$ denotes a suitably large power of $n$, so that $m\mid n^\infty$ means $m$ divides a certain large power of $n$.
\item $p$ is reserved for primes in $\Z$, and assumed to be positive if there is no particular emphasis.
\item $p^\nu\parallel q$ means $p^\nu\mid q$ but $p^{\nu+1}\nmid q$.
\item $\ord_p(n)$ denotes the exponent of $p$ in $n,$ i.e., $p^{v_p(n)}\parallel n.$
\item $\cP$ denotes the set of all primes in $\Z^+.$
\item $\pi$ appears in many occasions to denote the counting function of rational primes, or prime elements in rings of algebraic integers, or the transcendental number $3.1415926\cdots.$ One may find the correct meaning from the contexts.
\item $*$ in the summation indicates to sum over primitive elements.
\item $\Tr(\cdot)$ and $N(\cdot)$ denote the trace and norm of algebraic integers, respectively.
\item We also introduce some notation in asymptotic analysis:
\begin{itemize}
\item $A(x)\sim B(x)$ means $ A(x)/B(x)\rightarrow1$ as $x\rightarrow+\infty;$
\item $A(x)\ll B(x)$, or $A(x)=O(B(x))$, means $|A(x)|\leqslant cB(x)$ for some positive constant $c$ as $x\rightarrow+\infty;$
\item $A(x)\asymp B(x)$ means $ B(x)\ll A(x)\ll B(x).$
\end{itemize}
\end{itemize}

%% file: 1.Background.tex
\chapter{Backgrounds: history, conjectures and theorems}\label{ch:background}

Elliptic curves contribute as very active areas mixing algebra, geometry, analysis and also number theory.
As the ultimate refinement of the celebrated Sato--Tate conjecture, the Lang--Trotter conjecture aims to reveal statistical distributions of Frobenius traces with a fixed integer value instead of a fixed interval. In this chapter, we will give a general background on elliptic curves over rationals and finite fields, and recall several popular problems and known results around such topics. In particular, we will try to summarize the existing literature on the Lang--Trotter conjecture from different aspects, and formulate our main results in a suitably concise way. We hope we can motivate the readers sufficiently to open this memoir.

\section{Elliptic curves over rationals and finite fields}

Let $E$ be an elliptic curve over $\Q$ defined in its global minimal Weierstrass form \index{Weierstrass form}
\begin{align*}
E:y^2=x^3+ax+b
\end{align*}
with $a,b\in\Z.$ Denote by $N_E$ the conductor of $E$ that can be defined explicitly via the discriminant $\Delta_E=-16(4a^3 + 27b^2).$ We always assume $\Delta_E$ does not vanish.
For any prime
number $p\nmid N_E$, which we say a good prime, let $E_p$ denote the reduction of $E$ in $\F_p$. 
Let $$a_p(E):=p+1-|E_p(\F_p)|,$$
 the trace of Frobenius at $p$.\index{trace of Frobenius} More precisely, we may write
\begin{align*}
a_p(E)=-\sum_{x\in\F_p}\Big(\frac{x^3+ax+b}{p}\Big)
\end{align*}
for $p>3$, where $(\frac{\cdot}{p})$ denotes the Legendre symbol modulo $p$. For bad primes $p,$ define
\begin{align*}
a_p(E)=
\begin{cases}
1, &\text{ if } E \text{ has split multiplicative reduction at }p, \\
-1, &\text{ if } E \text{ has non-split multiplicative reduction at }p, \\
0, &\text{ if } E \text{ has additive reduction at }p.
\end{cases} 
\end{align*}
The importance and mysteriousness of $a_p(E)$ have been encoded in the behaviours of other objects concerning $E$. 
The Hasse--Weil $L$-function\index{Hasse--Weil $L$-function} then can be defined via the Euler product 
\begin{align*}
L(s, E) = \prod_{p\nmid N_E}({1-a_p p^{-s}} + p^{1-2s})^{-1}
\prod_{p\mid N_E}({1- a_p p^{-s}})^{-1}.
\end{align*}
In view of Hasse's bound $|a_p(E)|<2\sqrt{p}$\index{Hasse's bound}, it is clear that $L(s,E)$ absolutely converges in the half plane $\Re s>3/2.$
A lot of analytic information of $L(s,E)$ is now known thanks to the deep work of A. Wiles \cite{Wi95}, R. Taylor and A. Wiles \cite{TW95}, and C. Breuil, B. Conrad, F. Diamond and R. Taylor \cite{BCDT01} on Taniyama--Shimura--Weil conjecture\index{Taniyama--Shimura--Weil conjecture}, i,e., there exists a cusp form $f$ of weight $2$ and level $N_E$ such that $L(s,E)$ is also the $L$-function of $f$. Consequently, $L(s,E)$ admits an analytic continuation to the whole complex plane $\C$ and also satisfies a functional equation.
On the other hand, Hasse's bound becomes an immediate consequence of P. Deligne's proof \cite{De80} on Weil's conjecture for algebraic varieties over finite fields or Ramanujan's conjecture for holomorphic modular forms. 

A finer question is then to ask about statistical behaviours of $a_p(E)$ when $p$ or $E$ varies in a suitable family. The Sato--Tate conjecture\index{Sato--Tate conjecture} illustrates the equidistribution of $a_p(E)$ when $p$ varies over good primes, which is one of the motivations of Lang--Trotter conjecture\index{Lang--Trotter conjecture} considered in this article. We would like to explain the 
Sato--Tate conjecture and  Lang--Trotter conjecture in  next section. There are also many other celebrated conjectures (Birch and Swinnerton-Dyer conjecture for instance), but those are not the main concern here.

\section{Sato--Tate and Lang--Trotter conjectures}
We start this section with the following conjecture.\index{Sato--Tate conjecture}
\begin{conjecture}[Sato--Tate Conjecture]
Suppose $E/\Q$ is an elliptic curve as above.
For any fixed interval $I\subseteq[-2,2],$ we have
\[\lim_{x\rightarrow+\infty}\frac{|\{p\leqslant x: a_p(E)/\sqrt{p}\in I,~p\nmid N_E\}|}{\pi(x)}=\mu_{\mathrm{ST}}(I),\]
where $\mu_{\mathrm{ST}}$ is the Sato--Tate measure, and $\pi(x)$ denotes the number of primes $\leqslant x.$
\end{conjecture}

The Sato--Tate measure in $[-2,2]$ can be given explicitly according to $E$ has complex multiplication (CM for short) or not. More precisely,
\begin{align*}
\ud \mu_{\mathrm{ST}}=
\begin{cases}
\frac{1}{2}(\frac{\ud t}{\pi\sqrt{4-t^2}}+\delta_0), &\text{ if } E \text{ has CM}, \\\noalign{\vskip 2mm}
\frac{\sqrt{4-t^2}}{2\pi}\ud t, &\text{ if } E \text{ has no CM},
\end{cases} 
\end{align*}
where $\delta_0$ denotes the Dirac measure at $0$.
See some concise backgrounds on CM elliptic curves in Section \ref{sec:CMellipticcurve}. 
The CM case of Conjecture \ref{conj:LT} is much earlier proven by E. Hecke \cite{He18}, based on a criterion of M. Deuring (see Lemma \ref{lm:Deuring} in the next chapter) on the distribution of primes splitting in the relevant imaginary quadratic field. The non-CM case was not solved until around ten years ago thanks to a series of papers by L. Clozel, M. Harris and R. Taylor \cite{CHT08}, R. Taylor \cite{Ta08} and M. Harris, N. Shepherd-Barron and R. Taylor \cite{HST10}.
This was later generalized to non-CM holomorphic elliptic modular newforms by T. Barnet-Lamb, D. Geraghty, M. Harris and R. Taylor \cite{BGHT11}.

In 1976, S. Lang and H. Trotter \cite{LT76} formulated a very precise conjecture on the number of good primes $p$ such that 
$a_p(E)$ takes a fixed integer value\index{Lang--Trotter conjecture}. To this end, put
\begin{align}\label{eq:piE,r(x)-definition}
\pi_{E,r}(x):=|\{p\leqslant x:a_p(E)=r,~p\nmid N_E\}|
\end{align}
for all $r\in\Z$ and $x\geqslant2.$

\begin{conjecture}[Lang--Trotter Conjecture]\label{conj:LT}
Let $E$ be an elliptic curve over $\Q$ and assume additionally that $r\neq0$ if $E$ has complex multiplication.
Then
\begin{align*}
\pi_{E,r}(x)\sim \frac{c_{E,r}\cdot\sqrt{x}}{\log x},
\end{align*}
where $c_{E,r}\geqslant0$ is a constant that can be described 
in terms of the image of the associated Galois representation. 
\end{conjecture}

If $c_{E,r}=0$, we interpret the asymptotic formula in the above conjecture to mean that there are only finitely many primes $p$ for which $a_p=r.$ This happens when there is a congruence obstruction and then the conjecture becomes trivial. For instance, when $E$ has a $\Q$-rational 
point of order two, this implies that $a_p(E) \equiv p+1 \bmod 2$, and thus $r$ must be even in order for $a_p(E)=r$ with $p>2$; see N. Katz \cite{Ka09} for a precise 
discussion of congruence obstruction in a more general framework.   
In general, the constant $c_{E,r}$ is positive if and only if the associated Galois group has an element with trace $r$. But this condition is not easy to check in practice as one does not know the Galois group in general, even in the CM case.  

Before S. Lang and H. Trotter, B. Mazur \cite{Ma72} considered the primes $p$ such that $a_p\equiv1\bmod p$, which he called anomalous primes\index{anomalous prime}. In view of Hasse's bound, the above congruence is equivalent to $a_p=1$ for $p\geqslant5.$ Hence the problem of counting anomalous primes is a special case of the Lang--Trotter conjecture. When $E$ has CM, the case $r=0$ is excluded in Conjecture \ref{conj:LT} since a classical result of Deuring \cite{De41} implies that 
\begin{align*}
\pi_{E,0}(x)\sim \frac{1}{2}\frac{x}{\log x}.
\end{align*}
In view of Hasse's bound, the Sato--Tate conjecture implies that there is a positive proportion of primes such that $a_p(E)$ is of order of magnitude $\sqrt{p}$. Since $\pi(x)\sim x/\log x$ by Prime Number Theorem, a suitable probabilistic model predicts that the counting function of primes $p$ up to $x$, such that $a_p(E)$ is a fixed integer $r$, should be off by a factor $\sqrt{x}$ . This reveals the reasonableness of the term $\sqrt{x}/\log x$ in Conjecture \ref{conj:LT}. Moreover, the full strength of Conjecture \ref{conj:LT}, most notably the Lang--Trotter constant $c_{E,r}$\index{Lang--Trotter constant},  is also encoded in the 
Chebotarev density theorem\index{Chebotarev density theorem}, which reflects the subtle congruence contributions depending on the image of the associated Galois representation. 

Applying a similar 
probabilistic model, S. Lang and H. Trotter \cite{LT76} proposed another precise 
conjecture on the number of primes with a fixed trace in an imaginary quadratic field. 
We will show in Section \ref{sec:LTsecondconjecture} that this conjecture is equivalent to the Hardy--Littlewood conjecture. This should motivate us to believe the Hardy--Littlewood conjecture is a reasonable approach to illustrate the Lang--Trotter conjecture in the CM case.

The full strength of Conjecture \ref{conj:LT} seems out of reach at present.  In fact, no single non-trivial case (with $c_{E,r}>0$) of the Lang--Trotter conjecture has been proven, even assuming the Generalized Riemann hypothesis (GRH)\index{generalized Riemann Hypothesis}. 

In view of underlying difficulties, it is natural to ask weaker or alternative problems as approximations to the Lang--Trotter conjecture. Here are some instances:
\begin{enumerate}
\item Is there a non-trivial upper bound beating the trivial estimate
\begin{align*}
\pi_{E,r}(x)\ll \frac{x}{\log x}?
\end{align*}
\item Under suitable assumptions on $E$ and $r$, can one find a suitable increasing function $f(x)$, which could grow very slowly as $x\rightarrow+\infty$, such that
\begin{align*}
\pi_{E,r}(x)\gg f(x)?
\end{align*}
\item Are there any conjectures or hypothesis that can imply Conjecture \ref{conj:LT} or part of its full strength?
\item Once a conditional proof of Conjecture \ref{conj:LT} was found, how to interpret the resultant asymptotic formula? Are there any geometric consequences describing the trace of Frobenius or torsion points on $E$?
\item Instead of dealing with a single elliptic curve $E$ and a fixed integer $r,$ can one prove/disprove Conjecture \ref{conj:LT} when $E$ or $r$ runs over a suitable family?
\item In the spirit of vertical Sato-Tate distribution of elliptic curves and modular forms, can one prove a vertical analogue of Conjecture \ref{conj:LT}?
\end{enumerate}

The first part of our work is on the upper bound for $\pi_{E,r}(x)$. Here we would like to recall a lot of pioneer works before us. N. D. Elkies \cite{El87,El92} proved that for non-CM curve $E/\Q$ and $r=0$, 
$\pi_{E,0}(x)$ is unbounded (there exist infinitely many supersingular primes) \index{supersingular prime}
and 
\begin{align*}
\pi_{E,0}(x)\ll x^{\frac{3}{4}}. 
\end{align*}
This last upper bound was previously obtained by J. P. Serre \cite{Se81} under GRH for certain number fields\index{generalized Riemann Hypothesis}. The lower bound was later made effectively by \'E. Fouvry and M. R. Murty \cite{FM96}, i.e., 
\begin{align*}
\pi_{E,0}(x)\gg (\log\log\log x)(\log\log\log\log x)^{-\delta}
\end{align*}
for any $\delta>1.$ They also proved that $\pi_{E,0}(x)\gg \log\log x$ holds for infinitely many $x$, and this was previously known to Elkies \cite{Elthesis} and Murty for all sufficiently large $x$ if assuming the GRH for Dirichlet $L$-functions\index{generalized Riemann Hypothesis}.
For $\pi_{E,r}(x)$ with $r\not=0$, no single non-trivial lower bound 
is known, but there exist a lot of works focusing on the upper bound.

For $r\not=0$, all upper bounds so far are based on the effective Chebotarev density theorem\index{Chebotarev density theorem} and the images of the $\ell$-adic Galois representations attached to $E$, as initiated 
by Serre \cite{Se81} who proved that for any $\varepsilon>0$ and non-CM curve $E/\Q$, 
\begin{align*}
\pi_{E,r}(x)\ll_{\varepsilon}x(\log x)^{-\frac{5}{4}+\varepsilon}. 
\end{align*}

Serre's method and subsequent refinements by a number of authors have led to several 
improved upper bounds for $\pi_{E,r}(x)$.  
In the non-CM case, the first improvement was given 
by D. Wan \cite{Wa90} with the bound 
\begin{align*}
\pi_{E,r}(x)\ll_{\varepsilon}x(\log x)^{-2+\varepsilon}. 
\end{align*}
This was later improved by V. K. Murty \cite{Mu97} to
\begin{align*}
\pi_{E,r}(x)\ll \frac{x(\log\log x)^3}{(\log x)^2}.
\end{align*}
More recently, it was further sharpened by J. Thorner and A. Zaman \cite{TZ18} with the factor $(\log\log x)^2$ in place of $(\log\log x)^3$, which seems to be the best known unconditional upper bound for $\pi_{E,r}(x).$
Under GRH for Dedekind zeta functions, 
D. Zywina \cite{Zy15a} obtained
\begin{align*}
\pi_{E,r}(x)\ll x^{\frac{4}{5}}(\log x)^{-\frac{3}{5}},
\end{align*}
which seems to be the best known conditional upper bound. See  \cite{MMS88}  for an earlier result with a factor $(\log x)^{-\frac{1}{5}}$ in place of $(\log x)^{-\frac{3}{5}}.$ It is interesting to note that 
even assuming GRH, the estimate $\pi_{E,r}(x)\ll x^{\frac{1}{2}+\varepsilon}$ has not been proven yet.

For CM elliptic curves $E/\Q$, the effective Chebotarev \index{Chebotarev density theorem} and Galois image 
arguments also apply and give somewhat better bounds. An unconditional upper bound 
is stated in D. Wan \cite{Wa90} with 
\begin{align*}
\pi_{E,r}(x)\ll x\exp(-c_E\sqrt{\log x})
\end{align*}
for some positive constant $c_E$ depending on $E$. 
On the other hand, V. K. Murty \cite{Mu85} indicated that for any $\varepsilon>0$, 
\begin{align*}
\pi_{E,r}(x)\ll_{\varepsilon}x^{\frac{1}{2}+\varepsilon}
\end{align*}
holds if assuming the analytic continuation, functional equation, and GRH for all symmetric power L-functions of $E$\index{generalized Riemann Hypothesis}.

There are also some investigations of $\pi_{E,r}(x)$ on average and generalizations of Conjecture \ref{conj:LT} in different aspects. We will turn to such issues in the last chapter.

\section{Our main results}

We are now ready to state our first main result on the upper bound for $\pi_{E,r}(x).$

\begin{theorem}\label{thm:upperbound}
Let $r$ be a fixed non-zero integer and $E/\Q$ an elliptic curve with CM by the imaginary quadratic field $\Q(\sqrt{-D})$ with squarefree $D\geqslant1$. Let $\varepsilon>0.$ For all sufficiently large $x,$ we have
\begin{align*}
\pi_{E,r}(x)\leqslant \Big(\frac{712}{203}+\varepsilon\Big)\cdot h_{D,r}\cdot \frac{\sqrt{x}}{\log x}+O(1), 
\end{align*}
where $h_{D,r}$ is a constant depending only on $D$ and $r$, defined by 
\begin{align}\label{eq:constant-hD,r}
h_{D,r} = \xi(D,r)\cdot\frac{\sqrt{D}}{\varphi(D)}\prod_{p\nmid2r}\Big(1-\frac{(\frac{-D}{p})}{p-1}\Big)
\end{align}
with 
$\xi(D,r)$ taking values in $\{0, 1, 2\}$, explicitly given as in $\eqref{eq:xi(D,r)}$.
\end{theorem}

The proof of Theorem \ref{thm:upperbound} will be given in Chapter \ref{ch:HL&upperbound}. As we will see therein, the stronger inequality
\begin{align*}
\pi_{E,r}(x)+ \pi_{E,-r}(x) \leqslant \Big(\frac{712}{203}+\varepsilon\Big)\cdot h_{D,r}\cdot \frac{\sqrt{x}}{\log x}+O(1)
\end{align*}
also holds. This is the first time reaching the correct order of magnitude in the upper bound as conjectured by S. Lang and H. Trotter. Our new approach involves neither the Chebotarev density theorem 
nor images of Galois representations. Instead, in the 
CM case, we use the classical Deuring formula (Lemma \ref{lm:Deuring}) for $a_p(E)$ 
to transform the mixed counting function
$\pi_{E,r}(x)+\pi_{E,-r}(x)$ 
to the problem of presenting primes by quadratic polynomials. Roughly speaking, we expect any irreducible polynomial with integral coefficients can capture infinitely many primes if there are no local obstructions.
This can be illustrated by Schinzel Hypothesis \cite{SS58}\index{Schinzel Hypothesis} for a qualitative version or Bateman--Horn conjecture \cite{BH62}\index{Bateman--Horn conjecture} for a quantitative version. When $f$ is quadratic, one may refer to the Hardy--Littlewood conjecture \cite{HL23}\index{Hardy--Littlewood conjecture}; see also Conjecture \ref{conj:HL} in the next chapter. Therefore, after turning to the problem on the distribution of primes in quadratic progressions, many tools and ideas from analytic number theory can be employed to produce the upper bound for $\pi_{E,r}(x).$ We leave a remark that the scalar $\frac{712}{203}$ comes from applications of sieve methods, as will be shown in the next chapter.

After reaching the correct order of magnitude, it is natural to ask how far the upper bound in Theorem \ref{thm:upperbound} is away from the asymptotic formula as predicted in Conjecture \ref{conj:LT}.
Consequently, it is then necessary to compare the constants $h_{D,r}$ and $c_{E,r}:$ the former one comes from 
studying primes in quadratic progressions, and the latter one, we would like to call the Lang--Trotter constant\index{Lang--Trotter constant}, can be described in terms of Galois images\index{Galois representation}. However, an explicit description for $c_{E,r}$ seems to 
be subtle.  
N. Jones \cite{Jo09} initially interpreted $c_{E,r}$ by 
a formula of the shape 
\begin{equation}\label{eq:LTconstant-Galois}
\begin{split}
c_{E,r}&=\frac{m_E}{2\pi}\cdot\frac{|\Gal(K(E[m_E])/K)_r|}{|\Gal(K(E[m_E])/K)|}\cdot \prod_{p\mid r,~p\nmid m_E}\Big(1-\frac{(\frac{-D}{p})}{p}\Big)^{-1}\\
&\ \ \ \ \ \times\prod_{p\nmid rm_E}\Big(1-\frac{(\frac{-D}{p})}{(p-1)(p-(\frac{-D}{p}))}\Big), 
\end{split}
\end{equation}
see Section \ref{sec:Galoisimage} for more details. 
In principle, \eqref{eq:LTconstant-Galois} only gives a conceptual description on $c_{E,r}$, and it is usually difficult to compute its exact value when $E$ and $r$ are given directly. 
This phenomenon motivates the second part of this article. We would like   
to give an explicit analytic description of $c_{E,r}$ under the Hardy--Littlewood conjecture (Conjecture \ref{conj:HL}). This leads us to our second main result. 

\begin{theorem}\label{thm:formula}
Let $r$ be a fixed non-zero integer and $E/\Q$ an elliptic curve with CM by the imaginary quadratic field $\Q(\sqrt{-D})$ with squarefree $D\geqslant1$. 
The Hardy--Littlewood conjecture\index{Hardy--Littlewood conjecture} implies 
\begin{align*}
\pi_{E,r}(x)\sim\frac{\varpi_{E,r}\cdot\sqrt{x}}{\log x},
\end{align*}
where $\varpi_{E,r}\geqslant0$ is some explicit constant defined using analytic number theory. 
\end{theorem}
A key part of this theorem is a completely explicit description of 
the constant $\varpi_{E,r}$. 
The exact shape of $\varpi_{E,r}$, for instance in the CM case by $\Q(\sqrt{-3}),$ depends heavily on the sign and parity of coefficients in the equation defining $E$, as well as the residue class $r\bmod{24}.$ See Chapters \ref{ch:D=1}, \ref{ch:D=3} and \ref{ch:otherCM}  for concrete descriptions of the constant $\varpi_{E,r}$ . 
Theorem \ref{thm:formula} is also partially motivated by the works of Ji and Qin \cite{JQ14} and Qin \cite{Qi16,Qi21}.
In fact, Qin \cite{Qi16} gave, under the Hardy--Littlewood conjecture, a full description of 
the constant $\varpi_{E,r}$ for $r=\pm1$ and $E$ having CM by $\Q(\sqrt{-3}),$ in connection with Mazur's conjecture for anomalous primes.
This is a rather difficult case since there are six roots of unity in the ring of integers. More recently, Qin \cite{Qi21} also gave a full conditional description of 
the constant $\varpi_{E,r}$ for all $r\neq0$ and $E$ having CM by $\Q(\sqrt{-1}).$
We should emphasize that Qin is the first who worked out the details and gave explicit expressions for the constants in the above cases, although Mazur and Lang--Trotter should have some original ideas on the connections between the Mazur/CM Lang--Trotter conjecture and the Hardy-Littlewood conjecture. More importantly, Qin's work \cite{Qi16} shows that the values of $\varpi_{E,r}$ depend heavily on the coefficients of the equation defining $E$. This phenomenon gives a negative answer to ``a natural guess'' of Mazur who originally imagined that the values of $\varpi_{E,r}$ should be always the same in all possible cases; see the comments after Theorem \ref{thm:HL-LT:D=3} for more precise statements.

Theorem \ref{thm:formula} proves a conditional version of the 
Lang--Trotter conjecture. To show that this version agrees with the original Lang--Trotter conjecture, one is led to the following constant comparison conjecture. \index{comparison conjecture}

\begin{conjecture}[Comparison Conjecture]\label{conj:comparison}
Let $E/\Q$ be an elliptic curve with CM and $r\neq0$ a fixed integer. Then we have 
the equality
\begin{align*}
c_{E,r} = \varpi_{E,r}.
\end{align*}
\end{conjecture}

This conjecture is automatically true, if we assume both the Hardy--Littlewood conjecture and
the Lang--Trotter conjecture. But it should be much more realistic to prove than 
the later two conjectures. This comparison conjecture gives a quick explicit way 
to understand and compute the much harder Lang--Trotter constant $c_{E,r}$ 
in terms of the much easier and completely explicit constant $\varpi_{E,r}$.
For instance, we will give an explicit complete classification of the pair $(E,r)$ such that  the constant $\varpi_{E,r}$ is positive, 
see Chapter \ref{ch:LTconstant}. In particular, taking $r=1$, this gives an explicit 
complete classification of $E$ such that $\varpi_{E,1}$ is positive. 
This has application to Mazur's question \cite{Ma72} 
on when $E$ has infinitely many anomalous primes. 
The same classification 
is also true for the Lang--Trotter constant $c_{E,r}$ assuming the comparison conjecture.  
It would be interesting to give an unconditional proof of this classification for $c_{E,r}$. 
Recently, the Galois 
image (up to a conjugation) in the CM case is case by case explicitly classified in \'A. Lozano-Robledo \cite{LR22}. This makes it possible to give 
an explicit description of the Lang--Trotter constant $c_{E,r}$. 
This has been worked out in detail for some particular curves by A. Ray \cite{Ra23}.
On one hand, Ray developed a probabilistic model to explain Jones' constant $c_{E,r}$ as in \eqref{eq:LTconstant-Galois}; on the other hand, he evaluated $c_{E,r}$ for some $E$ and $r$ with the aid of \'A. Lozano-Robledo \cite{LR22} and F. Campagna and R. Pengo \cite{CP22}, and he has shown that the two constants 
$\varpi_{E,r}$ and $c_{E,r}$ coincide for such $E$ and $r$.

The case $r=2$ is of particular interests as the Lang--Trotter constant $c_{E,2}$ 
is known to be positive. This is because the Galois group always contains 
the $2\times 2$ identity matrix which has trace $2$ and so there is no congruence obstruction. 
Our positivity classification implies that $\varpi_{E,2}$ is also positive. 
Thus, the approximation 
\begin{align*}
\pi_{E,2}(x)\asymp\frac{\sqrt{x}}{\log x}
\end{align*}
holds for all CM elliptic curves $E$ under the Hardy--Littlewood conjecture. 
This is consistent with what the Lang--Trotter conjecture predicts. 
Note that the upper bound here is now unconditional (by Theorem \ref{thm:upperbound}), while the lower bound requires the Hardy--Littlewood conjecture\index{Hardy--Littlewood conjecture}, and the underlying constant 
can be made completely explicit, see Corollary \ref{coro:LTconstantpositive,r=2} in Chapter \ref{ch:LTconstant}. Even in the case $r=2$ where we know both $c_{E,2}>0$ and 
$\varpi_{E,2}>0$, we do not know an unconditional proof that they are equal.

To illustrate the shape of our main result in the second part, we present the conjectural asymptotic formula for $\pi_{E,r}(x)$ in several typical special cases.

\begin{corollary}\label{coro:HL-LT:D=1,case1}
Suppose $E/\Q$ is the elliptic curve defined by 
\begin{align*}
y^2=x^3-x.
\end{align*}
For each non-zero integer $r$, the Hardy--Littlewood conjecture $($Conjecture $\ref{conj:HL})$ implies 
\begin{align*}
\pi_{E,r}(x)\sim \frac{1}{2}\prod_{p\nmid r}\Big(1-\frac{(\frac{-1}{p})}{p-1}\Big)\cdot\frac{\sqrt{x}}{\log x}
\end{align*}
if $r\equiv2\bmod4,$ and $\pi_{E,r}(x)=O(1)$ otherwise.
\end{corollary}

\begin{corollary}\label{coro:HL-LT:D=1,case2}
Let $g$ be an odd integer. Suppose $E/\Q$ is the elliptic curve defined by 
\begin{align*}
y^2=x^3-2gx.
\end{align*}
For each integer $r\equiv2\bmod4$, the Hardy--Littlewood conjecture implies 
\begin{align*}
\pi_{E,r}(x)\sim \frac{1}{4}\prod_{p\nmid r}\Big(1-\frac{(\frac{-1}{p})}{p-1}\Big)\cdot\frac{\sqrt{x}}{\log x}.
\end{align*}
\end{corollary}

\begin{corollary}\label{coro:HL-LT:D=1,case3}
Let $g\in\Z^+$ be an odd prime. Suppose $E/\Q$ is the elliptic curve defined by 
\begin{align*}
y^2=x^3-gx.
\end{align*}
For each even integer $r$ with $g\nmid r$, the Hardy--Littlewood conjecture implies 
\begin{align*}
\pi_{E,r}(x)\sim \frac{1}{4}\cdot\Big\{1+\frac{\delta_{g,r}}{g-1-(\frac{-1}{g})}\Big\}\cdot\prod_{p\nmid r}\Big(1-\frac{(\frac{-1}{p})}{p-1}\Big)\cdot\frac{\sqrt{x}}{\log x},
\end{align*}
where
\begin{align*}
\delta_{g,r}
&=\begin{cases}
1,\ \ &\text{for~}  4\mid r,\\
-1,\ \ &\text{for~} 2\parallel r\text{~and~}g\equiv1\bmod4,\\
-1-2(-1)^{\frac{r-2}{4}+\frac{g+1}{4}},\ \ &\text{for~} 2\parallel r\text{~and~}g\equiv3\bmod4.
\end{cases}
\end{align*}
\end{corollary}

\begin{corollary}\label{coro:HL-LT:D=3,1}
Suppose $E/\Q$ is the elliptic curve defined by 
\begin{align*}
y^2=x^3+2.
\end{align*}
For each non-zero integer $r$, the Hardy--Littlewood conjecture implies 
\begin{align*}
\pi_{E,r}(x)\sim \frac{h_{3,r}}{6}\cdot\eta_r\cdot\frac{\sqrt{x}}{\log x},
\end{align*}
where $h_{3,r}\geqslant0$ is defined by $\eqref{eq:constant-hD,r}$
and
\begin{align*}
\eta_r&=\begin{cases}
2,\ \ &\text{for~} r\equiv16,20\bmod{24},\\ 
1,\ \ &\text{for~} r\equiv\pm1,\pm2,\pm5\bmod{12},\\
0,\ \ &\text{otherwise}.
\end{cases}
\end{align*}
\end{corollary}

To formulate the last asymptotic formula in this chapter, we introduce the following arithmetic function
\begin{align}\label{eq:Omegafunction}
\Omega_j(D;q,r)=\prod_{\substack{p^\nu\parallel q,~ p\nmid r\\j\nmid \nu}}\frac{-1}{p-1-(\frac{-D}{p})}
\end{align}
for all $D,q,r\in\Z$ and $j\in\Z^+$ with $2\nmid q.$

\begin{corollary}\label{thm:HL-LT:D=163}
Let $g\equiv1\bmod4$ be a positive integer coprime to $163$. Suppose $E/\Q$ is the elliptic curve defined by 
\begin{align*}
y^2=4x^3-53360\cdot 163g^2x+185801\cdot163^2g^3.
\end{align*}

For each non-zero integer $r$ coprime to $163,$ the Hardy--Littlewood conjecture implies 
\begin{align*}
\pi_{E,r}(x)
&\sim\frac{h_{163,r}}{2}\cdot\Big\{1-\Big(\frac{gr}{163}\Big)\Omega_2(163;g,r)\Big\}\cdot \frac{\sqrt{x}}{\log x},
\end{align*}
where $h_{163,r}\geqslant0$ is defined by $\eqref{eq:constant-hD,r}.$
\end{corollary}

Before concluding this chapter, we would like to mention that
the potential connection between CM Lang--Trotter conjecture \index{Lang--Trotter conjecture} and Hardy--Littlewood conjecture \index{Hardy--Littlewood conjecture} was realized 
long ago, for instance already by B. Mazur  
\cite{Ma72} and S. Lang and H. Trotter \cite{LT76}, at least in some special cases. 
The precise relationship, however, was not so clear due to several subtleties. 
The main difficulty lies in the understanding of the Lang--Trotter constant defined 
in term of Galois image.  
This constant is necessarily 
quite delicate because its conjectural answer (by the Comparison Conjecture) 
is already quite delicate as will be seen later.   

For a general CM elliptic curve $E$, one cannot reduce the Lang--Trotter 
conjecture for $\pi_{E,r}(x)$ to the Hardy--Littlewood conjecture 
for a single quadratic polynomial. One would need many 
quadratic polynomials in some sense. This part is doable but can be 
cumbersome to work out everything explicitly, as it may depend on a long list 
of case by case computations of higher reciprocity laws\index{reciprocity law}. This is the main input of Chapters \ref{ch:D=1}, \ref{ch:D=3} and \ref{ch:otherCM}. Our main idea here is to use 
the orthogonality of characters and analytic methods to keep all cases 
remain reasonably compact.  
See also Q. Ji and H. Qin \cite{JQ14} and H. Qin \cite{Qi16} for earlier partial results with more explicit polynomials via a different method. In particular, H. Qin \cite{Qi16} considered the curve $y^2=x^3+d$, which has CM by $\Q(\sqrt{-3}).$ He obtained an asymptotic formula for $\pi_{E,r}(x)$ with $r=\pm1$ under the Hardy--Littlewood conjecture,  and the case $r=1$ corresponds to anomalous primes \index{anomalous prime} initiated by B. Mazur \cite{Ma72}. Quite recently, H. Qin \cite{Qi21} turned his attention to the curve $y^2=x^3+gx$ which has CM by $\Q(\sqrt{-1})$, and was able to treat all non-zero $r$. It is worthwhile to compare our current approach with Qin's.
First of all, to determine the shape of the constant $\varpi_{E,r}$, we are both conditional 
under the Hardy--Littlewood conjecture, which is used to predict the distribution of primes in quadratic progressions (see Chapter \ref{ch:HL&upperbound} for details); in both approaches, laws of reciprocity are utilized. But we are different in at least three aspects:

\begin{itemize}
\item In both papers \cite{Qi16} and \cite{Qi21}, Qin started with combinatorial congruences to determine the Frobenius trace $a_p(E)$; see for instance \cite{Qi16} the details which appeal to B. C. Berndt, R. J. Evans and K. S. Williams \cite[Theorem 6.4.3]{BEW98}. We instead appeal to the original work of Gau{\ss} and the subsequent works of H. M. Stark and A. R. Rajwade, expressing $a_p(E)$
in terms of (quadratic, cubic and quartic) residue symbols; see Lemmas \ref{lm:Frobenius:D=1}, \ref{lm:Frobenius:D=3}, \ref{lm:Frobenius:D=2} and \ref{lm:Frobenius:D>=7} for details.

\item After applications of laws of reciprocity, both approaches need to evaluate sums of residue symbols to deal with distributions of prime elements in rings of integers (or equivalently, rational primes in quadratic progressions).
Qin computed sums of residue symbols by solving systems of linear equations (see \cite[Page 425]{Qi16} for instance). The introduction of Gau{\ss} sums in our approach allows one to study the relevant sums of (quadratic, cubic and quartic) residue symbols
in a unified way, and this is one key point that we can obtain a clean expression for $\varpi_{E,r}$ in Theorem \ref{thm:formula}. We believe that Gau{\ss} sums introduced here can serve as essential tools to study the Lang--Trotter conjecture in the current approach, and the evaluations of such Gau{\ss} sums in Chapter \ref{ch:Residuesandreciprocity} also have interests in themselves.
We also expect that such Gau{\ss} sums will find further applications in many other problems involving quadratic, cubic and quartic residues.
\item As one may see from Theorems \ref{thm:HL-LT:D=1} and \ref{thm:HL-LT:D=3}, the exact shape of $\varpi_{E,r}$ should depend heavily on the prime factorization of $r$, as well as the coefficients in the equation defining $E$. We have now treated all non-zero $r$ even in the most difficult case $D=3$, not just $r=\pm1$ as in the pioneer work of Qin \cite{Qi16}.
\end{itemize}

Both Qin and us have given explicit shapes of the constant in the asymptotic formula for $\pi_{E,r}(x)$ under the Hardy--Littlewood conjecture\index{Hardy--Littlewood conjecture}. A further  
problem is to prove that these two constants are equal, which is of course a consequence of the Hardy--Littlewood conjecture. But it is probably provable directly, at least in some cases. As an example, we provide a little air 
of compatibility here. For $r=1$ and CM elliptic curve $E_d$ of the form $y^2 = x^3 +d$, 
Qin \cite{Qi16} classified all $d$ such that his constant is positive, 
addressing Mazur's conjecture \cite{Ma72} \index{Mazur's conjecture} on when this curve $E_d$ possesses infinitely 
many anomalous primes. As will be seen in Theorem \ref{thm:anomalousprimes:D=3}, 
our classification of $d$ such that our constant $\varpi_{E_d, 1}$ is positive coincides with Qin's (after including a missing case mentioned in \cite[Page 3]{Qi21}).
It would be interesting to prove that this same classification also agrees with 
the classification of $d$ such that 
the Lang--Trotter constant $c_{E_d,1}$ is positive. This follows from a weak special case 
of Conjecture \ref{conj:comparison}, but it would be nice to prove it unconditionally. 

In Chapter \ref{ch:LTconstant}, we will give a short discussion on the geometric background of \eqref{eq:LTconstant-Galois}. Moreover, we will also present some other features of $\varpi_{E,r}$, including a complete classification on its positivity and symmetry in $r$.

Chapter \ref{ch:Historicalnotesperspectives} consists of some historical notes and perspectives that might be helpful for further researches. We also include some numerical data in the last chapter to support our theoretical constants in Chapters \ref{ch:D=1}, \ref{ch:D=3} and \ref{ch:otherCM}.

Overall, our work employs a lot of analytic number theory, which is becoming a more powerful and active area in recent years. We believe a lot of deep works will be realized by enthusiastic and ambitious scholars combining insights from analytic number theory, arithmetic geometry, and many other fields.

%% file: 2.HLupperbound.tex
\chapter{The Hardy--Littlewood conjecture and upper bound for $\pi_{E,r}(x)$}\label{ch:HL&upperbound}

In this chapter, we first formulate the Hardy--Littlewood conjecture \index{Hardy--Littlewood conjecture} for representations of primes by quadratic polynomials. After recalling some basic facts on CM elliptic curves, we will show how to relate the Lang--Trotter conjecture to distributions of primes in quadratic progressions. The upper bound in Theorem \ref{thm:upperbound} is then proven by using sieve methods.

\section{The Hardy--Littlewood conjecture and sieve methods}
The problem of presenting primes by quadratic polynomials can date back to Euler, and the precise conjecture with above asymptotic prediction were formulated by G. H. Hardy and J. E. Littlewood (\cite[Conjecture F]{HL23}) based on their circle method. \index{Hardy--Littlewood conjecture} 

\begin{conjecture}[Hardy--Littlewood Conjecture]\label{conj:HL}
Suppose $a,b,c\in\Z$ such that $a>0,$ $b^2-4ac$ is not a square and $(a,b,c)=1,$ and also assume that $a+b$ and $c$ are not simultaneously even. Then there exist infinitely many primes of the form $am^2+bm+c,m\in\Z^+.$ 

More precisely, denote by $\pi_{a,b,c}(x)$ the number of primes up to $x$ of the form $am^2+bm+c,m\in\Z^+$ with given $a,b,c$ as above. Then
\begin{align*}
\pi_{a,b,c}(x)\sim \frac{(2,a+b)}{\sqrt{a}}\frac{\delta}{\varphi(\delta)}\prod_{p\nmid 2a}\Big(1-\frac{(\frac{b^2-4ac}{p})}{p-1}\Big)\cdot \frac{\sqrt{x}}{\log x},
\end{align*}
where $\delta$ denotes the odd part of $(a,b)$, i.e., $\delta=(a,b)/(a,b,2^\infty)$.
\end{conjecture}

Of course, Conjecture \ref{conj:HL} is still out of reach at present, and on the other hand, H. Iwaniec \cite{Iw78} proved that there are infinitely many $m$ such that $am^2+bm+c$
has at most 2 prime factors, where $a,b,c$ are given as in Conjecture \ref{conj:HL}.
By virtue of Lemma \ref{lm:Deuring}, for upper bound purpose, we may transfer from Lang--Trotter conjecture \index{Lang--Trotter conjecture} in the CM case to Hardy--Littlewood conjecture with certain explicit quadratic polynomial, and a linear sieve applies to bound the number of primes captured by such a polynomial, which then yields the upper bound part of the Lang--Trotter conjecture as shown in Theorem \ref{thm:upperbound}.

We may modify Conjecture \ref{conj:HL} by mixing restrictions from both of arithmetic  and quadratic progressions. To do so, define
\begin{align*}
\pi_{a,b,c}(x;q,u)=|\{p\leqslant x: p=am^2+bm+c,~m\equiv u\bmod q,~m\in\Z^+\}|,
\end{align*}
where  $q$ is a positive integer and $(u,q)=1,$ and $a,b,c$ are the same as in Conjecture \ref{conj:HL}.
Note that $\pi_{a,b,c}(x;q,u)$ is not a new generalization of $\pi_{a,b,c}(x)$ indeed, because it may be  transformed to $\pi_{a',b',c'}(x)$ for suitable $a',b',c'\in\Z$ by elementary means. The counting function of the shape $\pi_{a,b,c}(x;q,u)$ will be used in Chapters \ref{ch:D=1}, \ref{ch:D=3} and \ref{ch:otherCM} to characterize the Lang--Trotter constant.

For $m\equiv u\bmod q$, we write $m=kq+u$ for some $k\in\Z.$ Hence
\begin{align*}
am^2+bm+c=a(kq+u)^2+b(kq+u)+c=Ak^2+Bk+C,
\end{align*}
where
\begin{align*}
A&=aq^2,\\
B&=(2au+b)q,\\
C&=au^2+bu+c.
\end{align*}
Note that
\begin{align*}
A+B&\equiv(a+b)q\bmod2,\\
C&\equiv (a+b)u+c\bmod2,\\
(A,B)&=q(aq,2au+b),\\
(A,B,C)&=(q(aq,2au+b),au^2+bu+c),
\end{align*}
and
\begin{align*}
B^2-4AC=(2au+b)^2q^2-4aq^2(au^2+bu+c)=q^2(b^2-4ac).
\end{align*}
Therefore, $B^2-4AC$ is a square if and only if $b^2-4ac$ is a square, and $(A,B,C)=1$ if and only if
$(q,au^2+bu+c)=(a,b,c)=1.$ 

We are now ready to formulate the modified version of Conjecture \ref{conj:HL}.\index{Hardy--Littlewood conjecture} 
\begin{conjecture}[Modified Hardy--Littlewood Conjecture]\label{conj:HL-modified}
Suppose $a,b,c,q,u\in\Z$ such that $a>0,$ $q>0,$ $b^2-4ac$ is not a square and $(q,au^2+bu+c)=(a,b,c)=1,$ and also assume that $(a+b)q$ and $(a+b)u+c$ are not simultaneously even. Then there exist infinitely many primes of the form $am^2+bm+c$ with $m\equiv u\bmod q$ and $m\in\Z^+.$ 

More precisely, we have
\begin{align*}
\pi_{a,b,c}(x;q,u)\sim \frac{(2,(a+b)q)}{q\sqrt{a}}\frac{\Delta}{\varphi(\Delta)}\prod_{p\nmid 2qa}\Big(1-\frac{(\frac{b^2-4ac}{p})}{p-1}\Big)\cdot \frac{\sqrt{x}}{\log x},
\end{align*}
where $\Delta$ denotes the odd part of $q(aq,2au+b)$. In particular, if $q$ is even, then
\begin{align*}
\pi_{a,b,c}(x;q,u)\sim \frac{2}{q\sqrt{a}}\frac{q^\sharp}{\varphi(q^\sharp)}\prod_{p\nmid qa}\Big(1-\frac{(\frac{b^2-4ac}{p})}{p-1}\Big)\cdot \frac{\sqrt{x}}{\log x}
\end{align*}
for $(a,b)=1,$ and correspondingly,
\begin{align*}
\pi_{a,b,c}(x;q,u)\sim \frac{2}{q\sqrt{a}}\frac{(qa)^\sharp}{\varphi((qa)^\sharp)}\prod_{p\nmid qa}\Big(1-\frac{(\frac{b^2-4ac}{p})}{p-1}\Big)\cdot \frac{\sqrt{x}}{\log x}
\end{align*}
for $a\mid b,$
where $n^\sharp$ denotes the odd part of $n$.
\end{conjecture}

In the remaining part of this section, we would like to introduce some backgrounds on the (linear) sieve method that will be used later.\index{sieve method} 
Generally speaking, the sieve method works with a non-negative sequence $\cA=(a_n)$ of integers, and the original aim is to capture primes or almost primes along this sequence. To this end, define the sifting function
\begin{align*}S(\cA,\cP,z)=\sum_{(n,P(z))=1}a_n
\end{align*}
with $P(z)=\prod_{p<z}p.$
Given a squarefree number $d\geqslant1$, consider the congruence sum 
\begin{align*}A_d=\sum_{n\equiv 0\bmod d}a_n.\end{align*}
In applications to sieve methods, it is expected that $\cA$ equidistributes in the special arithmetic progression 
$n\equiv0\bmod d$: There exists an appropriate approximation $\fS$ to $A_1$ and 
a multiplicative function $g$ supported on squarefree numbers with
\begin{align}\label{eq:density-g}
0\leqslant g(p)<1
\end{align}
such that the remainder
\begin{align}\label{eq:sieveremainder}
r_d(\cA):=A_d-g(d)\fS
\end{align}
is small on average over $d\mid P(z)$, and there exists a constant $L>1$ such that
\begin{align}\label{eq:sieve-V-linearsieve}
\frac{V(z_1)}{V(z_2)}
\leqslant\frac{\log z_2}{\log z_1}\Big(1+\frac{L}{\log z_1}\Big)
\end{align}
for $2\leqslant z_1< z_2$, where
\begin{align*}
V(z)=\prod_{p<z}(1-g(p)).
\end{align*}

Let $\fF$ and $\ff$ be the continuous solutions to the system
\begin{align}\label{eq:sieve-F}
\begin{cases}
s\fF(s) =2\mathrm{e}^\gamma & \text{for $\,0< s\leqslant 2$},
\\
s\ff(s)=0                                    & \text{for $\,0<s\leqslant 2$},
\\
(s\fF(s))'=f(s-1)                         & \text{for $\,s>2$},
\\
(s\ff(s))'=F(s-1)                         & \text{for $\,s>2$},
\end{cases}
\end{align}
where $\gamma$ is the Euler constant.

 An arithmetic function $\lambda(d)$ is of {\it level}
$Q$ if
\begin{align*}
\lambda(d)=0
\quad(d>Q)
\qquad\text{and}\qquad
|\lambda(d)|\leqslant1
\quad
(d\leqslant Q).
\end{align*}
Let $r\geqslant2$ be a positive integer. We say that $\lambda$ is {\it well-factorable}, if for every decomposition $Q = Q_1Q_2$ with $Q_1, Q_2\geqslant1,$ there exist $\lambda_1, \lambda_2$ such that
\begin{align*}
\lambda=\lambda_1*\lambda_2\end{align*}
with each $\lambda_j$ of level $Q_j$.

We now state the following fundamental result of H. Iwaniec \cite{Iw80}.

\begin{lemma}\label{lm:Iwaniec}
Let $0<\varepsilon<\tfrac{1}{8}$.
Under the above hypothesis, we have
\begin{align*}
S(\cA,\cP,z)
\leqslant \fS V(z)\Big\{\fF\Big(\frac{\log Q}{\log z}\Big)+E\Big\}
+\sum_{t\leqslant T}\sum_{d\mid P(z)}\lambda_t(d) r_d(\cA)
\end{align*}
for all $z\geqslant 2$,
where $\fF$ is given by the system \eqref{eq:sieve-F}, 
$T$ depends only on $\varepsilon$, 
$\lambda_t$ is well-factorizable of level $Q,$ 
and $E\ll\varepsilon +\varepsilon^{-8} \mathrm{e}^L(\log Q)^{-1/3}.$
\end{lemma}

As one can show, $\fF(s)$ decays rapidly to $1$ as $s\rightarrow+\infty.$ In applications, we would like to explore the admissible $Q$ as large as possible such that the total remainder
\begin{align*}
\sum_{t\leqslant T}\sum_{d\mid P(z)}\lambda_t(d) r_d(\cA)
\end{align*}
is still under control. In principle, a larger $Q$ would produce a better upper bound for $S(\cA,\cP,z)$ upon a suitable choice for $z.$ The merit of Lemma \ref{lm:Iwaniec} lies in the well-factorization of $\lambda_t$, so that one can estimate the resultant bilinear forms more effectively to produce a better admissible $Q$ (compared with earlier sieve methods before Iwaniec).

To prove Theorem \ref{thm:upperbound}, we need to specialize $\cA$ to sequences defined by quadratic polynomials. The details will be given in the next section.

\smallskip

\section{Elliptic curves with complex multiplication}\label{sec:CMellipticcurve}

In this section, we present some backgrounds on elliptic curves with complex multiplication\index{complex multiplication}, and reduce the mixed counting function $\pi_{E,r}(x)+\pi_{E,-r}(x)$ to bounding the number of primes in quadratic progressions.
As mentioned before, connections between two of them have been known for a long time in various special cases. For our more modest upper bound purpose, we shall give a simple uniform (and elementary) treatment in all situations. 

We assume that $E/\Q $ has CM
by an imaginary quadratic field $K=\Q(\sqrt{-D})$, where $D$ is a positive squarefree number. The endomorphism ring $R={\rm End}(E\otimes \C)$ is then an order of $K$. The theory of complex 
multiplication implies that the $j$-invariant $j(E)$ is in $\Q$ 
if and only if $K$ has class number one. Since our elliptic curve 
$E$ is defined over $\Q$, this means that $K$ is one of the 9 
imaginary quadratic fields with class number one. That is, $K=\Q (\sqrt{-D})$, where 
$$D\in\{1, 2, 3, 7, 11, 19, 43, 67, 163\}$$
as proven independently by A. Baker \cite{Ba66} and H. Stark \cite{St67}.\index{imaginary quadratic field}

Let $R_D$ denote 
the ring of integers in $K=\Q (\sqrt{-D})$. 
The order $R$ (in fact, every order in $K$) is of the form 
$R=\Z +fR_D$ for a unique positive integer $f$, called the conductor 
of the order $R$. The maximal order $R_D$ has conductor $1$. 
For an elliptic curve $E$ defined over $\Q$ with CM
by $K=\Q(\sqrt{-D})$, the conductor $f$ of its endomorphism ring $R$ 
is $1$, except that it can also be $2$ if $D=1,3,7$, and also $3$ if $D=3$. 
This gives in total $13$ different $\overline{\Q}$-isomorphism classes 
of elliptic curves defined over $\Q$ with CM. 
Each class is a family of infinitely many twists of 
one elliptic curve over $\Q$. 
A list of $13$ representatives from the $13$ classes with minimal 
discriminants and conductors is given in Appendix A in 
\cite{Si94}. If $E$ is an elliptic curve over $\Q$ with 
CM by the order $R={\rm End}(E\otimes \C)$ in $K$, 
there is another elliptic curve $E'$ over $\Q$ which is $\Q$-isogenous 
to $E/\Q$ such that ${\rm End}(E'\otimes \C)=R_D$ (the maximal order), 
see Ex 2.12 in \cite{Si94}.  
Since two isogenous elliptic curves over $\Q$ have the same 
Frobenius traces for all large primes $p$, we can restrict to 
CM elliptic curves with endomorphism rings having conductor $1$. 
In this way, it suffices to 
study the 9 $\overline{\Q}$-isomorphism classes of 
CM elliptic curves over $\Q$ with endomorphism ring $R$ 
equal to the ring $R_D$ of integers in one of the above $9$ imaginary quadratic fields $K=\Q(\sqrt{-D})$.

Let us recall the following classical criterion from Deuring \cite{De41}.  \index{Deuring's criterion}

\begin{lemma}\label{lm:Deuring}
Let $E/{\Q}$ be an elliptic curve with complex 
multiplication by the ring $R_D$ of integers in an imaginary quadratic field $K=\Q (\sqrt{-D})$. Let $p$ be a prime 
where $E$ has good reduction. If $p$ is inert in $R_D$, then $a_p(E)=0$. If $p=\pi\bar{\pi}$ splits in $R_D$, then 
$a_p(E)=\pi+\bar{\pi}$, where $\pi\in R_D$ is chosen so that the 
endomorphism $[\pi]$ induces the $p$-th power Frobenius automorphism of 
the reduction of $E$ modulo $p$. 
\end{lemma} 

We are now ready to relate the CM Lang--Trotter conjecture to primes presented by quadratic polynomials. We consider the cases $D\equiv 1,2 \bmod 4$ and $D\equiv 3 \bmod 4$ separately.

{\bf Case I.} Assume $D\equiv 1,2 \bmod 4$. Then, $R_D= \Z + \Z \sqrt{-D}$. 
In this case, if $p=\pi\bar{\pi}$ with $\pi=m+n\sqrt{-D} \in R_D$ 
such that the endomorphism $[\pi]$ induces the $p$-th power Frobenius automorphism of the reduction of $E$ modulo $p$, then the choice of $m$ and $n$ depends on $p$.
Now we have 
$$p=m^2 + Dn^2, \ a_p(E)=\pi+\bar{\pi}=2m, \ m,n\in \Z.$$
It then follows that $a_p(E)$ is an even integer and $p$ is of the form 
$$p=\Big(\frac{a_p(E)}{2}\Big)^2 + Dn^2, \ n\in \Z.$$
If $a_p(E)=r$, then we have 
\begin{align}\label{eq:p=(r/2)^2 + Dn^2}
p=\frac{1}{4}r^2+ Dn^2.
\end{align}
For large $p$, this equation determines  
$r$ uniquely up to a sign as $n$ has to be large. 
This sign can also be determined, but we do not need 
this more technical result for our purpose.  
Note that $\pi_{E,r}(x)$ vanishes unless $r$ is even, in which case we find
$$\pi_{E,r}(x)+ \pi_{E,-r}(x) \leqslant|\{p\leqslant x: p=\tfrac{1}{4}r^2+ Dn^2, \ n\in \Z\}|+O(1).$$
The inequality holds since the choice of $n$ in \eqref{eq:p=(r/2)^2 + Dn^2} may depend on $p$ and $E$, and we now allow $n$ to run over all integers freely to obtain an upper bound.
The term $O(1)$ comes from a finite number of small primes $p$ and the possible contributions with $(\tfrac{r}{2}, D)>1$. In particular, $\pi_{E,r}(x)=O(1)$ if $(\frac{r}{2}, D)>1$.

{\bf Case II.} Assume $D\equiv 3\bmod 4$. Then, $R_D =\Z +\Z \frac{1+\sqrt{-D}}{2}$. 
In this case, if $p=\pi\bar{\pi}$ with $\pi=m+n\frac{1+\sqrt{-D}}{2}
\in R_D$ such that the endomorphism $[\pi]$ induces the $p$-th power Frobenius automorphism of the reduction of $E$ modulo $p$, then we have 
$$p=\Big(m +\frac{n}{2}\Big)^2 + D\Big(\frac{n}{2}\Big)^2, \ a_p(E)=\pi+\bar{\pi}= 2m+n, \ m, n \in \Z.$$
If $a_p(E)=r$, then $2m+n=r$ and
$$p=\frac{1}{4}(D+1)r^2-Drm+Dm^2, \ m \in \Z.$$
It follows that
$$\pi_{E,r}(x)+ \pi_{E,-r}(x) \leqslant|\{p\leqslant x:p=\tfrac{1}{4}(D+1)r^2-Drm+Dm^2, \ m \in \Z\}|+O(1).$$
Note that $\pi_{E,r}(x)=O(1)$ if $(r, D)>1$.

We leave here a remark that the above $m$ can be restricted to be positive since $\tfrac{1}{4}(D+1)r^2-Drm+Dm^2$ is stable under the map $m\mapsto r-m.$ Moreover, since there are only two roots of unity in $R_D$ when $D\neq1,3,$ the above inequalities for $\pi_{E,r}(x)+ \pi_{E,-r}(x)$ then become identities up to an error $O(1)$.
In summary, we have proved 

\begin{lemma} \label{lm:LT-quadraticpolynomial}
Let $E/\Q$ be an elliptic curve with CM by the imaginary quadratic field 
$\mathbb{Q}(\sqrt{-D})$. For any integer $r\neq0$, 
there exists a quadratic polynomial $\fh_{D,r}\in \Z[X]$ depending only on $D$ and $r$ such that 
\begin{align}\label{eq:LT-HL}
\pi_{E,r}(x)+ \pi_{E,-r}(x) \leqslant|\{p\leqslant x:p=\fh_{D,r}(n), \ n\in \Z^+\}|+O(1).
\end{align}
Precisely, we may take 
\begin{align*}
\fh_{D,r}(n)=
\begin{cases}
Dn^2+\frac{1}{4}r^2,\ \ & \text{if~}D\equiv1,2\bmod4,~2\mid r,\\
Dn^2-Drn+\frac{1}{4}(D+1)r^2, & \text{if~}D\equiv3\bmod4,\\
0, &\text{otherwise}.
 \end{cases}
\end{align*}

Moreover, $``\leqslant"$ in $\eqref{eq:LT-HL}$ can be replaced by $``="$ if $D\neq1,3.$
\end{lemma}

\begin{remark}

For some of the $13$ elliptic curves $E/\Q$ 
in the table in \cite{Si94}, Q. Ji and H. Qin \cite{JQ14} constructed 
an explicit quadratic polynomial $\fh_{E,r}(x)$ such that 
$$\pi_{E,r}(x) = |\{p\leqslant x: p=\fh_{E,r}(n), \ n\in \Z\}| +O(1).$$
This polynomial $\fh_{E,r}$ depends more information on $E$, not just on $D$. 
Such results are more precise than Lemma \ref{lm:LT-quadraticpolynomial} for those special
elliptic curves. Since we are working with upper bounds only, 
the above simpler lemma, which works for all CM elliptic curves, 
suffices for our purpose. 

\end{remark}

As an illustration, we review the most classical example, 
where $E$ is the CM elliptic curve over $\Q$ defined by $y^2=x^3 +4x$. 
It has CM by $\Q(\sqrt{-1})$. 
This curve $E$ is the non-singular model of the celebrated 
lemnscate curve $x^2+y^2+x^2y^2=1$ studied in the last entry  
in the mathematical diary of Gau{\ss}\index{last entry of Gau{\ss}' diary}; see Chapter \ref{ch:Historicalnotesperspectives} for more history. For each prime $p\equiv 3 \bmod 4$, 
$E$ is supersingular at $p$ and thus $a_p(E)=0$. 
For each prime $p\equiv 1 \bmod 4$, 
Gau{\ss} wrote down an explicit formula which is equivalent to 
$$|E_p(\F_p)|=(m-1)^2 +n^2,$$
where $p=m^2 +n^2$, $m$ odd, $n$ even, 
and $\pi=m+ni \equiv 1 \bmod {2+2i}$ in $\Z[i]$. 
It then follows immediately that $a_p(E)=2m$. 
In particular, $\pi_{E,r}(x)=0$ if $r$ is either odd or divisible by $4$. 
If $r\equiv 2 \bmod 8$, one deduces that 
$a_p(E)=r$ if and only if $p=(r/2)^2 + 16n^2$ for some $n\in \Z$. 
If $r\equiv -2\bmod 8$, one has
$a_p(E)=r$ if and only if $p=(r/2)^2 + 4(2n+1)^2$ for some $n\in \Z$. 
For $r \equiv \pm2\bmod 8$, we thus obtain the mixed  
statement that $a_p(E)=\pm r$ 
if and only $p=(r/2)^2 +4n^2$ for some $n\in \Z$. It follows that 
the Lang--Trotter conjecture for this special curve $E$ and 
$r\equiv 2\bmod 8$ should be equivalent to the classical Hardy--Littlewood conjecture  
which asserts that 
$$|\{p\leqslant x:p=(r/2)^2+16n^2, \ n\in \Z\}|
\sim \frac{c(r)\cdot\sqrt{x}}{\log x},$$
where $c(r)>0$ is the Hardy--Littlewood constant given in Conjecture \ref{conj:HL}. We do not know if the Lang--Trotter constant $c_{E,r}$ and the Hardy--Littlewood constant $c(r)$  
have been checked to match in this special case. 
Similar type of equivalences have been studied  for several more 
CM elliptic curves over $\Q$ mentioned above in \cite{JQ14} and \cite{Qi16}. 
This simple example should already explain the subtlety and difficulty for the 
Lang--Trotter conjecture in the CM case. The non-CM case 
would apparently be more difficult, and no such classical interpretation exists! 

Another interesting example is the curve $E_d:y^2=x^3+d$, where $d$ is neither a square nor a cube 
in $\Q(\sqrt{-3})$. This curve has CM 
by $\Q(\sqrt{-3})$. As mentioned in Chapter \ref{ch:background}, B. Mazur \cite{Ma72} conjectured, in his work on Iwasawa theory and 
anomalous primes\index{anomalous prime}, that there is a 
positive constant $c_d$ such that 
\begin{align*}
\pi_{E_d,1}(x)\sim \frac{c_d\cdot\sqrt{x}}{\log x}. 
\end{align*}
Mazur did not describe the constant $c_d$ explicitly. \index{Mazur's conjecture}
H. Qin \cite{Qi16} gave a complete description of the constant $c_d$ under the Hardy--Littlewood conjecture, 
in particular showing that the constant $c_d$ can be zero in some cases.  
Again, it is not proven that the two constants $c_d$ and $c_{E_d,1}$ 
are the same. Of course, they should match exactly if one assumes both the 
Lang--Trotter conjecture and Hardy--Littlewood conjecture.\index{Hardy--Littlewood conjecture}

\section{Proof of Theorem \ref{thm:upperbound}}
We now complete the proof of Theorem \ref{thm:upperbound}. Throughout this section, we fix $r\neq0$ and let $\fh$ denote the quadratic polynomial $\fh_{D,r}$ with leading coefficient $D$, given explicitly in Lemma \ref{lm:LT-quadraticpolynomial}. If $E$ has CM by $\mathbb{Q}(\sqrt{-D})$, 
then we have  
\begin{align*}
\pi_{E,r}(x)+ \pi_{E,-r}(x) \leqslant |\{p\leqslant x:p=\fh(n), \ n\in \Z^+\}|+O(1).
\end{align*}
Switching the roles of $p$ and $n$, we may write
\begin{align*}
\pi_{E,r}(x)+ \pi_{E,-r}(x) \leqslant |\{n\leqslant X:\fh(n)\in\cP\}|+O(1)
\end{align*}
for some $X$ with $X\sim\sqrt{x/D}.$ Moreover, we define a constant $\delta(D,r)$ by
\begin{align*}
\delta(D,r)=
\begin{cases}
1,\ \ &\text{if~}D\equiv1\bmod4,2\mid r,(D,r)=1,\\
1,&\text{if~}D\equiv2\bmod4,2\mid r,(D,r/2)=1,\\
1,&\text{if~}D\equiv3\bmod4,(D,r)=1\text{~and~}2\mid r,\\
1, &\text{if } D\equiv3\bmod 8,(D,r)=1\text{~and~}2\nmid r,\\
0,&\text{otherwise}.
\end{cases}
\end{align*}
The introduction of $\delta(D,r)$ is to guarantee that $\fh$ is irreducible and has no fixed prime factors as long as $\delta(D,r)=1.$
Note that for $D=2,$ the condition $\delta(D,r)=1$ is equivalent to $r\equiv2\bmod4.$

The following inequality can be derived by sieve methods. 
\begin{proposition}\label{prop:LT-upperbound-1}\index{sieve method}
Let $\varepsilon>0$. Suppose $\fh$ is the quadratic polynomial 
$\fh_{D,r}$ in Lemma $\ref{lm:LT-quadraticpolynomial}$, where $D\in\{1, 2, 3, 7, 11, 19, 43, 67, 163\}.$ Then for all sufficiently large $x,$ we have
\begin{align}\label{eq:LT-upperbound-1}
\pi_{E,r}(x)+ \pi_{E,-r}(x) \leqslant \Big(\frac{712}{203}+\varepsilon\Big)\cdot\delta(D,r)\cdot \Gamma_\fh\cdot \frac{\sqrt{x/D}}{\log x}+O(1),
\end{align}
where $\delta(D,r)$ is defined as above and
\begin{align}\label{eq:Gamma-fh}
\Gamma_\fh=\prod_p\Big(1-\frac{\varrho(\fh,p)}{p}\Big)\Big(1-\frac{1}{p}\Big)^{-1}
\end{align}
with $\varrho(\fh,p)=|\{\nu\bmod p:\fh(\nu)\equiv0\bmod p\}|.$
\end{proposition}

In the remaining part of this section, we always work in the conditions such that $\delta(D,r)=1$.
Let $Y=X^{1-\varepsilon}$ and we introduce a smooth function $\phi$ supported on $[Y,X]$, satisfying
\begin{align*}
\begin{cases}
\phi(x)=1\ \ &\text{for } x\in[Y+1,X-1],\\
\phi(x)\in[0,1]\ \ &\text{for } x\in\R,\\
\phi^{(j)}\ll_j Y^{-j} &\text{for all } j\geqslant0.
\end{cases}
\end{align*}
Therefore,
\begin{align*}
\pi_{E,r}(x)+\pi_{E,-r}(x)
&\leqslant\sum_{\substack{m\in\Z\\\fh(m)\in\cP}}\phi(m)+O(X^{1-\varepsilon})\\
&\leqslant\sum_{\substack{m\in\Z\\(\fh(m),\Pi)=1}}\phi(m)+O(X^{1-\varepsilon}),
\end{align*}
where $\Pi=\prod_{p<\sqrt{x}}p$.
We would like to apply Lemma \ref{lm:Iwaniec} with
\begin{align*}
a_n=\sum_{\substack{m\in\Z\\\fh(m)=n}}\phi(m)
\end{align*}
and $z=\sqrt{x},$ in which situation we have
\begin{align*}
\fS=\sum_{m\in\Z}\phi(m)=X+o(X)
\end{align*}
and $g(d)=\varrho(\fh,d)/d$. It also requires to verify the conditions given by \eqref{eq:density-g} and \eqref{eq:sieve-V-linearsieve}. The former one is clear by Lagrange's theorem for $p>2$ and trivial for $p=2$.
More precisely, we have
\begin{align*}
\varrho(\fh,p)=
\begin{cases}
0,\ \ &\text{if } p\mid D,\\
1+(\frac{-Dr^2}{p}), &\text{if } p\nmid2D,
\end{cases}
\end{align*}
and
\begin{align*}
\varrho(\fh,2)=
\begin{cases}
1,\ \ &\text{if } D\equiv1\bmod 4,\\
1, &\text{if } D\equiv3\bmod 4\text{~and~}2\mid r,\\
0, &\text{if } D\equiv3\bmod 8\text{~and~}2\nmid r.\end{cases}
\end{align*}
The latter condition is guaranteed by 
\begin{align*}
\prod_{p<z}\Big(1-\frac{\varrho(\fh,p)}{p}\Big)
=\{1+o(1)\}\frac{\ue^{-\gamma}\Gamma_\fh}{\log z}
\end{align*}
for all sufficiently large $z$, where $\Gamma_{\fh}$ is defined by \eqref{eq:Gamma-fh}.
Therefore, from Lemma \ref{lm:Iwaniec} it follows that
\begin{align}
\pi_{E,r}(x)+ \pi_{E,-r}(x)
& \leqslant 2\ue^{-\gamma}\Gamma_\fh\cdot\{\fF(\vartheta)+\varepsilon\}\cdot\frac{X}{\log x}\label{eq:initialupperbound}
\end{align}
as long as 
\begin{align}\label{eq:level-D}
\sum_d\lambda(d)\Big(\sum_{\substack{n\in\Z\\\fh(n)\equiv0\bmod d}}\phi(n)-\frac{\varrho(\fh,d)}{d}\widehat{\phi}(0)\Big)\ll X(\log X)^{-A}
\end{align}
holds for any well-factorable function $\lambda$ of level $X^\vartheta$ and any $A>0$. To do so, we appeal to the recent work of
R. de la Bret\`eche and S. Drappeau \cite{BD20}.
\begin{proposition}[de la Bret\`eche--Drappeau]\label{prop:dlBD}
Let $$0<\theta<1+\frac{25}{178}.$$ Then $\eqref{eq:level-D}$ holds for any well-factorable function $\lambda$ of level $X^\theta$ and any $A>0$. 
\end{proposition}

Taking $\vartheta=1+\frac{25}{178}-\eta$ for some $\eta>0$ in \eqref{eq:initialupperbound}, and noting that $$2\ue^{-\gamma} \fF\Big(1+\frac{25}{178}\Big)=\frac{4}{1+\frac{25}{178}},$$ we arrive at \eqref{eq:LT-upperbound-1} in view of $X\sim\sqrt{x/D}.$

By virtue of evaluations for $\varrho(\fh,p)$, we find
\begin{align*}
\delta(D,r)\cdot\Gamma_\fh
&=\xi(D,r)\frac{D}{\varphi(D)}\prod_{p\nmid2D}\Big(1-\frac{1+(\frac{-Dr^2}{p}) }{p}\Big)\Big(1-\frac{1}{p}\Big)^{-1}\\
&=\xi(D,r)\frac{D}{\varphi(D)}\prod_{p\nmid2r}\Big(1-\frac{(\frac{-D}{p})}{p-1}\Big),
\end{align*}
where
\begin{align}\label{eq:xi(D,r)}
\xi(D,r)=
\begin{cases}
1,\ \ & \text{if~} D\equiv1\bmod4,2\mid r\text{~and~}(D,r)=1,\\
1,\ \ & \text{if~} D\equiv2\bmod4,2\mid r\text{~and~}(D,r/2)=1,\\
1, &\text{if } D\equiv3\bmod 4,2\mid r\text{~and~}(D,r)=1,\\
2, &\text{if } D\equiv3\bmod 8,2\nmid r\text{~and~}(D,r)=1,\\
0, &\text{otherwise}.
\end{cases}
\end{align}
This proves Theorem \ref{thm:upperbound} by combining Proposition \ref{prop:LT-upperbound-1}.

\begin{remark}
The strong level of distribution in Proposition \ref{prop:dlBD} is not necessary if one is not sensitive about the scalar in Theorem \ref{thm:upperbound}. A level $X^\vartheta$ with an absolute $\vartheta>0$ is sufficient to produce the correct order of magnitude in the upper bound for $\pi_{E,r}(x)+\pi_{E,-r}(x).$ In fact, H. Iwaniec \cite{Iw78} obtained a weaker level $X^{1+\frac{1}{15}-\eta}$, which enabled him to prove the infinitude of almost primes (with at most two prime factors) captured by fixed quadratic polynomials with no local obstructions. A direct application of Fourier analysis would lead to the level $X^{1-\eta}$. Iwaniec's level relies on Weil's bound for individual Kloosterman sums\index{Kloosterman sum}, and the work of R. de la Bret\`eche and S. Drappeau \cite{BD20} then borrows deep results from spectral theory of automorphic forms\index{spectral theory of automorphic forms}.
\end{remark}

\smallskip

\section{Fixed trace in imaginary quadratic fields}\label{sec:LTsecondconjecture}

In this section, we consider another conjecture of 
Lang--Trotter \cite{LT76} which counts primes with a fixed trace in 
an arbitrary imaginary quadratic field $\mathbb{Q}(\sqrt{-D})$, 
not necessarily with class number $1$, where $D$ is a positive squarefree integer. As before, let $R_D$ denote the ring of integers in $\mathbb{Q}(\sqrt{-D})$. \index{imaginary quadratic field}

Denote by $\Pi_D$ the set of those elements $\pi\in R_D$ 
such that the norm $N(\pi)=\pi\bar{\pi}$ is a rational prime number 
$p\in\Z^+$. Note that such an element $\pi$ is just a prime 
element in $R_D$ whose norm is a  prime number, that is, the principal ideal $(\pi)$ in $R_D$ 
is a prime ideal of $R_D$ with a prime residue field.

Given an element $\pi \in \Pi_D$, its trace 
${\rm Tr}(\pi)=\pi +\bar{\pi}$ 
is a rational integer. Note that if $\pi \in \Pi_D$, then 
its conjugate $\bar{\pi}$ is also in $\Pi_D$ with the same trace and norm. If $u$ is a unit in $R_D$ and $\pi\in \Pi_D$, then $u\pi \in \Pi_D$ with the same norm $p$ but generally different trace. 

For each non-zero integer $r$,  define  
\begin{align}\label{eq:piD,r(x)-definition}
2\pi_{D,r}(x):=|\{\pi \in \Pi_D: ~N(\pi) \leqslant x,~{\rm Tr}(\pi)=r\}|. 
\end{align}
Thus, $\pi_{D,r}(x)$ counts the number of un-ordered pairs 
$(\pi, \bar{\pi})$ such that $\pi\in \Pi_D$, $N(\pi)=p\leqslant x$ 
and ${\rm Tr}(\pi)=r$. As we shall see, we can alternatively write
\begin{align*}
\pi_{D,r}(x)=|\{p=N(\pi)\leqslant x:\pi \in \Pi_D,~{\rm Tr}(\pi)=r\}|. 
\end{align*}
Replacing $\pi$ by $-\pi$, one sees that 
$\pi_{D,r}(x)=\pi_{D,-r}(x)$. One can thus assume that $r$ is positive if 
one wishes. 

Based on a similar probabilistic model as in the original  
conjecture for $\pi_{E,r}(x)$, S. Lang and H. Trotter \cite[Page 80]{LT76} proposed 
the following precise conjecture for $\pi_{D,r}(x)$. \index{Lang--Trotter second conjecture}
\begin{conjecture}[Lang--Trotter second Conjecture]\label{conj:LT-2} 
Keep the same notation and convention as above.
For each non-zero integer $r,$ we have 
\begin{align*}
\pi_{D,r}(x)\sim \frac{c_{D,r}\cdot\sqrt{x}}{\log x},
\end{align*}
where $c_{D,r}\geqslant0$ is an explicit constant described below.  
\end{conjecture}

\begin{remark}
Conjecture \ref{conj:LT-2} is also beyond our current approach, and even a positive lower bound for $\pi_{D,r}(x)$ is unknown yet. On the other hand, one can do much better if summing over $r$ in a suitable range. For instance, G. Harman and P. A. Lewis \cite{HL01} proved that 
\begin{align}\label{eq:HarmanLewis}
\frac{1}{R}\sum_{|r|\leqslant R}\pi_{D,r}(x)\gg\frac{\sqrt{x}}{\log x}
\end{align}
for $D=1$ with any $R\gg x^{0.119}$ (More earlier references can also be found in \cite{HL01}).
We would remark that Conjecture \ref{conj:LT-2} has not been solved even under a suitable GRH. 
Assuming the Riemann Hypothesis for Hecke $L$-functions with Gr\"ossencharacters over $\Q(\sqrt{-1})$\index{generalized Riemann Hypothesis}, N. C. Ankeny \cite{An52} and J. Kubilius \cite{Ku55} have shown independently, also for $D=1$, that \eqref{eq:HarmanLewis} holds for any $R\gg \log x.$ All the above works might be extended to general $D$, at least for those $\Q(\sqrt{-D})$ with class number $1.$
\end{remark}

We now describe the explicit constant $c_{D,r}$. Put
\begin{align*}
\cA_p(D)=
\begin{cases}
0,\ \ &\text{if } p\mid D,\\
(1-\frac{1}{p})^{-1},\ \ &\text{if } (\frac{-D}{p})=1,\\
(1+\frac{1}{p})^{-1},\ \ &\text{if } (\frac{-D}{p})=-1,
\end{cases}
\end{align*}
\begin{align*}
\cB_p(D)=
\begin{cases}
(1-\frac{1}{p})^{-1},\ \ &\text{if } p\mid D,\\
(1-\frac{2}{p})(1-\frac{1}{p})^{-2},\ \ &\text{if } (\frac{-D}{p})=1,\\
(1-\frac{1}{p^2})^{-1},\ \ &\text{if } (\frac{-D}{p})=-1.
\end{cases}
\end{align*}
Following Lang and Trotter \cite[Theorem 2.3]{LT76}, we have
\begin{align*}
c_{D,r}&=\frac{w_D}{2\pi h(D)}F_4(D,r)\cdot\prod_{\substack{p>2\\ p\mid r}}\cA_p(D)\cdot\prod_{\substack{p>2\\ p\nmid r}}\cB_p(D),
\end{align*}
where $h(D)$
denotes the class number of $\Q(\sqrt{-D})$, 
$w_D$ denotes the number of roots of unity in $\Q(\sqrt{-D})$, and $F_4(D,r)$  with $r\bmod4$ is given explicitly by the following table.
\begin{table}[htbp]
\renewcommand{\arraystretch}{1.6}
\renewcommand{\tabcolsep}{1.9mm}
\begin{tabular}{|c|c|c|c|c|}
\hline
& $F_4(D,0)$ & $F_4(D,1)$ & $F_4(D,2)$ & $F_4(D,3)$\\ \hline
$D\equiv1\bmod4$ & $2$ & $0$ & $2$ & $0$ \\  \hline
$D\equiv2\bmod4$ & $0$ & $0$ & $4$ & $0$ \\  \hline
$D\equiv3\bmod8$ & $\frac{2}{3}$ & $\frac{4}{3}$ & $\frac{2}{3}$ & $\frac{4}{3}$ \\  \hline
$D\equiv7\bmod8$ & $2$ & $0$ & $2$ & $0$ \\  \hline
\end{tabular}
\\~\\
\vskip 2mm
\caption{Values of $F_4(D,r)$ with $r\bmod4$}
\end{table}

Note that there is an error in the original statement of Lang and Trotter \cite[Theorem 2.3]{LT76}, in which $\cA_p(D)$, with a different notation, is given by $(1-\frac{1}{p})^{-1}$ for odd prime $p\mid D.$ 
The above corrected value for $\cA_p(D)$ is zero for odd prime $p\mid D$. 

Lang and Trotter suggested that the above conjecture is a consequence 
of the Hardy--Littlewood conjecture \index{Hardy--Littlewood conjecture} with matching constants. They checked that this is indeed the case 
when $D=1$ and $r=1$; see \cite[Page 81]{LT76}. In the following, 
we check that this is also true for general $D$ and $r$. Furthermore, 
we can show, in a similar sieving process, that the upper bound part of Conjecture \ref{conj:LT-2} 
holds up to a constant factor. 

It is obvious that $c_{D,r}$ vanishes if the odd part of $(D,r)$ is larger than $1.$ We henceforth assume $(D,r)$ has no odd prime factors. Therefore,
\begin{align*}
\prod_{\substack{p>2\\ p\mid r}}\cA_p(D)
&=\prod_{\substack{p>2\\ p\mid r}}\Big(1-\frac{(\frac{-D}{p})}{p}\Big)^{-1}
\end{align*}
and
\begin{align*}
\prod_{\substack{p>2\\ p\nmid r}}\cB_p(D)
&=\prod_{\substack{p>2\\ p\mid D}}\Big(1-\frac{1}{p}\Big)^{-1}\cdot\prod_{\substack{p\nmid 2r\\(\frac{-D}{p})=1}}\Big(1-\frac{2}{p}\Big)\Big(1-\frac{1}{p}\Big)^{-2}\cdot\prod_{\substack{p\nmid 2r\\(\frac{-D}{p})=-1}}\Big(1-\frac{1}{p^2}\Big)^{-1}\\
&=\frac{D^\sharp}{\varphi(D^\sharp)}\prod_{p\nmid 2r}\Big(1-\frac{(\frac{-D}{p})}{(p-1)(p-(\frac{-D}{p}))}\Big),
\end{align*}
where $D^\sharp$ denotes the odd part of $D$. Combining the above evaluations, we infer
\begin{equation}\label{eq:cD,r-initial}
\begin{split}
c_{D,r}&=\frac{w_D}{2\pi  h(D)}F_4(D,r)\frac{D^\sharp}{\varphi(D^\sharp)}\prod_{\substack{p>2\\p\mid r}}\Big(1-\frac{(\frac{-D}{p})}{p}\Big)^{-1}\\
&\ \ \ \ \times\prod_{p\nmid 2r}\Big(1-\frac{(\frac{-D}{p})}{(p-1)(p-(\frac{-D}{p}))}\Big)
\end{split}
\end{equation}
as long as $(D,r)$ has no odd prime factors, in which case we further derive that
\begin{align*}
c_{D,r}&=\frac{w_D}{2\pi  h(D)}F_4(D,r)\frac{D^\sharp}{\varphi(D^\sharp)}\prod_{p>2}\Big(1-\frac{(\frac{-D}{p})}{p}\Big)^{-1}\prod_{p\nmid 2r}\Big(1-\frac{(\frac{-D}{p})}{p-1}\Big)\\
&=\frac{w_D}{2\pi  h(D)}F_4(D,r)\frac{D^\sharp}{\varphi(D^\sharp)}\eta_D L(1,\chi_D)\prod_{p\nmid 2r}\Big(1-\frac{(\frac{-D}{p})}{p-1}\Big),
\end{align*}
where $D^\sharp$ denotes the odd part of $D$, $\chi_D=(\frac{-D}{\cdot})$ or $(\frac{-4D}{\cdot})$ denotes the Kronecker symbol of $\Q(\sqrt{-D})$, and $\eta_D=1$ if $D\equiv1,2\bmod4,$ and $\eta_D=1-\frac{1}{2}(\frac{-D}{2})$ otherwise; i.e.,
\begin{align*}
\eta_D=\begin{cases}
1,\ \ &\text{if~}D\equiv1,2\bmod4,\\
\frac{3}{2}, &\text{if~}D\equiv3\bmod8,\\
\frac{1}{2}, &\text{if~}D\equiv7\bmod8.
\end{cases}
\end{align*}
Let $\Delta$ denote the discriminant of $\Q(\sqrt{-D})$. Then the class number formula  gives
\begin{align*}
w_D\sqrt{\Delta}L(1,\chi_D)=2\pi  h(D).
\end{align*}
Hence
\begin{align*}
\frac{w_DL(1,\chi_D)}{h(D)}=\begin{cases}
\pi D^{-\frac{1}{2}},\ \ &\text{if~}D\equiv1,2\bmod4,\\
2\pi D^{-\frac{1}{2}}, &\text{if~}D\equiv3\bmod4,
\end{cases}
\end{align*}
from which we infer that $\frac{w_D}{2\pi h(D)}\eta_D L(1,\chi_D)$ is $\frac{3}{2\sqrt{D}}$ for $D\equiv3\bmod8,$ and is
$\frac{1}{2\sqrt{D}}$ otherwise. These yield
\begin{align*}
\frac{w_D}{2\pi h(D)}\frac{D^\sharp}{\varphi(D^\sharp)}\eta_D L(1,\chi_D)
&=\frac{\sqrt{D}}{\varphi(D)}\times\begin{cases}
\frac{1}{4},\ \ &\text{if~}D\equiv2\bmod4,\\
\frac{1}{2},\ \ &\text{if~}D\equiv1,5,7\bmod8,\\
\frac{3}{2}, &\text{if~}D\equiv3\bmod8.
\end{cases}
\end{align*}
We then arrive at
\begin{align}\label{eq:c(D,r)=h(D,r)}
c_{D,r}&=h_{D,r}=\xi(D,r)\cdot\frac{\sqrt{D}}{\varphi(D)}\prod_{p\nmid2r}\Big(1-\frac{(\frac{-D}{p})}{p-1}\Big),
\end{align}
where $\xi(D,r)$ is given by \eqref{eq:xi(D,r)}. 
This shows that the second Lang--Trotter constant $c_{D,r}$ agrees with 
the Hardy--Littlewood constant $h_{D,r}$. 

We are now ready to prove the following 
upper bound for $\pi_{D,r}(x)$. 
\begin{theorem}\label{thm:LT-2} 
For any non-zero integer $r$, we have 
\begin{align*}
\pi_{D,r}(x) \leqslant \Big(\frac{712}{203}+\varepsilon\Big)\cdot h_{D,r}\cdot \frac{\sqrt{x}}{\log x}+O(1).
\end{align*}
\end{theorem}

The proof of this theorem  
is essentially the same as the proof of Theorem \ref{thm:upperbound}.  
We reduce it to Proposition \ref{prop:LT-upperbound-1}. 
Fix a non-zero integer $r$. Let $p=\pi\bar{\pi}$ be a prime which splits in $R_D$. We claim that 
for $p$ sufficiently large compared to $r$, there is at most one unit $u$ in $R_D$ such that ${\rm Tr}(u\pi)=r$. 

In the case that $D\not= 1,3$, there are only two units $\pm 1$.   
If ${\rm Tr}(\pi)=r\not=0$, then clearly ${\rm Tr}(-\pi)=-r\not\not=r$. 
The uniqueness follows. 

Suppose $D=1$. Let $\pi=m+ni \in R_D=\mathbb{Z}[i]$ 
such that $p=\pi\bar{\pi}=m^2+n^2$ with ${\rm Tr}(\pi)=2m=r$. 
There are four units $\{\pm1, \pm i\}$. As indicated above, 
${\rm Tr}(-\pi)=-r\not\not=r$. One checks that ${\rm Tr}(\pm i\pi)=\pm n\not=r$ since we assumed that $p$ is large and hence $|n|>|r|$. 
The uniqueness follows. 

It remains to consider the case $D=3$. Let $\pi=m+n\sqrt{-3} \in R_D=\mathbb{Z}+\mathbb{Z}\frac{1+\sqrt{-3}}{2}$  
such that $p=\pi\bar{\pi}=m^2+3n^2$ with ${\rm Tr}(\pi)=2m=r$.
There are now six units in $R_D$. If $u$ is one of the four 
units in $\{ \frac{\pm1 \pm\sqrt{-3}}{2}\}$, then one checks that 
$${\rm Tr}(u\pi) = (\pm m \pm 3n)= \pm\Big(\frac{r}{2}\pm 3n\Big).$$
The last number cannot be $r$ since we have assumed that $p$ is large 
and thus $|n|$ is large. The uniqueness again follows, completing the proof of the claim.  

Let  $\fh_{D,r}$ be the quadratic polynomial in Lemma \ref{lm:LT-quadraticpolynomial}, 
which depends only on $D$ and $r$. As shown in the proof of Lemma \ref{lm:LT-quadraticpolynomial}, 
if $\pi\in \Pi_D$ with ${\rm Tr}(\pi)=r$, then $p=\pi\bar{\pi}$ is of the 
form $\fh_{D,r}(n)$. Conversely, if $p=\fh_{D,r}(n)$ is a prime, 
then $p=\pi\bar{\pi}$ splits in $R_D$. One checks that for some unit $u$ in $R_D$, $u\pi$ has trace $r$. By 
the above uniqueness claim, we conclude that 
\begin{align*}
\pi_{D,r}(x)=|\{p\leqslant x: p=\fh_{D,r}(n), ~n\in \Z^+\}| +O(1).
\end{align*}
Theorem \ref{thm:LT-2} then follows from Proposition \ref{prop:LT-upperbound-1}.

%% file: 3.Residuesandreciprocity.tex
\chapter{Power residues and laws of reciprocity}\label{ch:Residuesandreciprocity}

A fundamental problem in number theory is to solve the congruence equation
\begin{align}\label{eq:congruenceequation}
f(x)\equiv0\bmod m
\end{align}
in $x\in\Z$, where $m$ is a positive integer and $f(X)\in\Z[X]$.  Take $f$ to be a quadratic polynomial for instance. After completing squares, the problem reduces to solving the diagonal equation of the shape
\begin{align*}
x^2\equiv a\bmod m.
\end{align*}
In the generic case $(a,m)=1$, we say $a$ is a quadratic residue modulo $m$ if this equation is solvable, and correspondingly, a quadratic non-residue modulo $m$ otherwise.
In a more general setting, for a given positive integer $k\geqslant2$ and also $(a,m)=1$, we say $a$ is a $k$-th residue modulo $m$ if the equation 
\begin{align}\label{eq:k-th-residueequation}
x^k\equiv a\bmod m
\end{align}
is solvable. The solvability of \eqref{eq:k-th-residueequation} is then reduced to consider the case of prime modulus $m$ by the Chinese remainder theorem and Hensel's lift ($p$-adic analysis).

The theory of quadratic residues dates back to Euler, Legendre, Gau{\ss} and Eisenstein, and subsequent generalizations 
to higher powers (cubic, quartic, quintic, sextic, etc) are also very fruitful. These also lead to the class field theory and automorphic representations in the 20th century. The heart of the theory lies in the laws of reciprocity, and one may refer to 
F. Lemmermeyer \cite{Le00} for the history and classical literature. We also suggest \cite{IR90} as a well-organized material with reasonable details and results.

Let $K$ be a number field with $O_K$ as the ring of integers. Let $n$ be a positive integer and denote by $\mu_n$ the set of $n$-th roots of unity in $\C.$ Suppose $\mu_n\subseteq K.$ For a prime ideal $\fp$ in $O_K$ with $n\not\in\fp,$ the $n$-th power residue symbol $(\frac{\cdot}{\fp}):O_K\rightarrow\mu_n\cup\{0\}$ can be defined 
via
\begin{align*}
\Big(\frac{\alpha}{\fp}\Big)\equiv \alpha^{(N\fp-1)/n}\bmod\fp
\end{align*}
for all $\alpha\in O_K,$ where $N\fp=|O_K/\fp|$ denotes the norm of $\fp$. The condition $N\fp\equiv1\bmod n$ can be verified via the bijection
\begin{align*}
\mu_n&\rightarrow\mu_n(O_K/\fp),\\
\zeta&\mapsto \zeta\bmod \fp.
\end{align*}
The above power residue symbol is well-defined since $X^n\equiv1\bmod\fp$ has exactly $n$ distinct roots and $\alpha^{\frac{N\fp-1}{n}}\bmod\fp$ gives one of them, provided that $\alpha\not\in\fp$.

The power residue symbol is truly a generalization of Legendre symbol over $\Q$ (see next section for details), and it admits a lot of beautiful properties.
For instance, for any $\alpha\in O_K$ and $\fp\nmid(\alpha,n),$ one has $(\frac{\alpha}{\fp})=1$ if and only if $\alpha$ becomes a $n$-th power in $O_K$. The multiplicativity 
\begin{align*}
\Big(\frac{\alpha\beta}{\fp}\Big)=\Big(\frac{\alpha}{\fp}\Big)\Big(\frac{\beta}{\fp}\Big)
\end{align*} 
holds for any $\alpha,\beta\in O_K.$ It is trivial that $(\frac{\alpha}{\fp})=(\frac{\beta}{\fp})$ if $\alpha\equiv\beta\bmod\fp.$ Hence $\alpha\mapsto(\frac{\alpha}{\fp})$ defines a (multiplicative) character on $(O_K/\fp)^\times.$

Let $\fm$ be an ideal in $O_K$ such that $(n,\fm)=1.$ According to the prime ideal decomposition $\fm=\prod_j \fp_j^{\nu_j},$ we may generalize the above power residue symbol modulo prime ideals by
\begin{align*}
\Big(\frac{\alpha}{\fm}\Big)=\prod_j\Big(\frac{\alpha}{\fp_j}\Big)^{\nu_j},
\end{align*}
following the same manner in generalizing Legendre symbol to Jacobi symbol over $\Q.$

The laws of reciprocity \index{reciprocity law} for $(\frac{\cdot}{\fp})$ and $(\frac{\cdot}{\fm})$ are due to Gau{\ss}, Jacobi, Eisenstein, Kummer, Hilbert, Hasse, E. Artin, et al in different situations, which motivate the birth of class field theory. We do not intend to expand the details here, and alternatively,  many concrete instances will be discussed in the following sections on quadratic, quartic and cubic residue symbols, as well as evaluations on relevant Gau{\ss} sums. More precisely, this chapter aims to outline some basic and essentially elementary results on quadratic, quartic and cubic residues and corresponding reciprocity laws. We omit the proofs of reciprocity laws here since none is due to us (this is the strongest reason!) and they do not belong to the main concern of this article. However, some results will be proved in detail occasionally because they usually appear in sporadic references. For the purpose of the main body of this memoir, we will also evaluate two kinds of Gau{\ss} sums with quartic and cubic characters. Their appearance seems to be very natural and should have been studied in depth already, but unfortunately we did not succeed in finding helpful references for our study on the Lang--Trotter conjecture in later chapters!

\section{Quadratic residues and reciprocity}

Let $p\in\Z^+$ be an odd prime. The Legendre symbol $(\frac{\cdot}{p}):\Z/p\Z\rightarrow\{-1,0,1\}$ is defined as \index{Legendre symbol}
\begin{align*}
\Big(\frac{a}{p}\Big)=
\begin{cases}
1,\ \ \ &\text{if $a$ is a quadratic residue modulo $p$},\\
-1,\ \ \ &\text{if $a$ is a quadratic non-residue modulo $p$},\\
0,\ \ \ &\text{if $p\mid a$}.
\end{cases}
\end{align*}
Euler's criterion asserts that 
\begin{align}\label{eq:Eulercriterion}
\Big(\frac{a}{p}\Big)\equiv a^{\frac{p-1}{2}}\bmod p.
\end{align}

The domain of $(\frac{\cdot}{p})$ can be extended to $\Z$ due to the periodicity. Moreover, we also have the (complete) multiplicativity
\begin{align*}
\Big(\frac{ab}{p}\Big)=\Big(\frac{a}{p}\Big)\Big(\frac{b}{p}\Big),\ \ a,b\in\Z.
\end{align*}
This transforms the evaluation of $(\frac{a}{p})$ to those of $(\frac{-1}{p})$, $(\frac{2}{p})$ and $(\frac{q}{p})$, where $q\in\Z^+$ is another odd prime different from $p$.

\begin{theorem}[Quadratic reciprocity law]\label{thm:quadraticreciprocity}\index{quadratic reciprocity law}
Let $p$ and $q$ be two distinct {\rm(}positive{\rm)} odd primes. Then
\[\Big(\frac{q}{p}\Big)\Big(\frac{p}{q}\Big)=(-1)^\frac{(p-1)(q-1)}{4}.\]
We also have the supplementary laws 
\begin{align*}
\Big(\frac{-1}{p}\Big)=(-1)^\frac{p-1}{2},\ \ \Big(\frac{2}{p}\Big)=(-1)^\frac{p^2-1}{8}.\end{align*}
\end{theorem}

One may also extend the validity of Theorem \ref{thm:quadraticreciprocity} to all primes which are not necessarily positive:
\[\Big(\frac{q}{|p|}\Big)\Big(\frac{p}{|q|}\Big)=(-1)^\frac{(p-1)(q-1)}{4}(p,q)_\infty,\]
where $p,q$ are distinct odd primes and
\begin{align}\label{eq:Hilbertsymbol}
(a,b)_\infty=
\begin{cases}
-1,\ \ \ &\text{if $a<0,~b<0$},\\
1,\ \ \ &\text{otherwise}
\end{cases}
\end{align}
is the Hilbert symbol (at $\infty$).\index{Hilbert symbol}

The Jacobi symbol \index{Jacobi symbol} extends Legendre symbol to the set of all positive odd integers by multiplicativity. More precisely, for each positive odd integer $m$, define
\[\Big(\frac{a}{m}\Big)=\prod_{p^\nu\parallel m}\Big(\frac{a}{p}\Big)^\nu.\]
The symbol $(\frac{\cdot}{m})$ enjoys similar properties with $(\frac{\cdot}{p})$ (i.e., periodicity, multiplicativity and the quadratic reciprocity law with the supplementary ones), and we omit the details here.

We now turn to the solvability of \eqref{eq:congruenceequation} in the very beginning of this chapter. When $f$ is a quadratic polynomial, the problem can be settled perfectly using Legendre and Jacobi symbols. 
Furthermore, we may also examine the counting function
\begin{align*}
\cN_\pm(f,q)&=\Big|\Big\{t\bmod q:\Big(\frac{f(t)}{q}\Big)=\pm 1\Big\}\Big|.
\end{align*}
Here and in what follows, assume $q$ is a positive odd integer and $f(X)=aX^2+bX+c\in\Z[X]$ is a fixed primitive quadratic polynomial, i.e., $(a,b,c)=1$. We also write $\Delta=b^2-4ac$. 
Moreover, put
\begin{align*}
\cN_0(f,q)&=|\{t\bmod q:f(t)\equiv0\bmod q\}|,\\
\cN_1(f,q)&=|\{t\bmod q:(f(t),q)=1\}|,
\end{align*}
and
\begin{align*}
\cN_2(f,q)&=\sum_{t\bmod q}\Big(\frac{f(t)}{q}\Big),
\end{align*}
so that
\begin{align*}
2\cN_\pm(f,q)&=\cN_1(f,q)\pm \cN_2(f,q).
\end{align*}

The following lemma is an exercise in elementary number theory.

\begin{lemma}\label{lm:N0andN1-evaluations}
Let $p$ be an odd prime. Then
\begin{align*}
\cN_0(f,p)=
\begin{cases}
1+(\frac{\Delta}{p}),\ \ &\text{for~} p\nmid a,\\
1,\ \ &\text{for~} p\mid a\text{~and~}p\nmid b,\\
0,\ \ &\text{for~} p\mid (a,b)\text{~and~}p\nmid c,
\end{cases}
\end{align*}
and
\begin{align*}
\cN_2(f,p)&=\begin{cases}
p(\frac{c}{p}),\ \ &\text{for~} p\mid (a,b),\\
0,\ \ &\text{for~} p\mid a\text{~and~}p\nmid b,\\
(p-1)(\frac{a}{p}),\ \ &\text{for~} p\nmid a\text{~and~}p\mid \Delta,\\
-(\frac{a}{p}),&\text{for~} p\nmid a\Delta.
\end{cases}
\end{align*}
\end{lemma}

Note that the functions $q\mapsto \cN_0(f,q),$ $q\mapsto \cN_1(f,q)$ and $q\mapsto \cN_2(f,q)$
are all multiplicative for each fixed $f$. Hence it suffices to consider the case when $q$ is an odd prime power.

Let $p$ be an odd prime and $\nu\geqslant1.$ Then
\begin{align*}
\cN_1(f,p^\nu)&=p^\nu-p^{\nu-1}\cN_0(f,p)\\
&=p^\nu-p^{\nu-1}\times
\begin{cases}
1+(\frac{\Delta}{p}),\ \ &\text{for~} p\nmid a,\\
1,\ \ &\text{for~} p\mid a\text{~and~}p\nmid b,\\
0,\ \ &\text{for~} p\mid (a,b)\text{~and~}p\nmid c.
\end{cases}
\end{align*}
On the other hand,
\begin{align*}
\cN_2(f,p^\nu)&=\cN_1(f,p^\nu)
\end{align*}
if $\nu$ is even, and otherwise 
\begin{align*}
\cN_2(f,p^\nu)&=p^{\nu-1}\cN_2(f,p)=p^{\nu-1}\times 
\begin{cases}
p(\frac{c}{p}),\ \ &\text{for~} p\mid (a,b),\\
0,\ \ &\text{for~} p\mid a\text{~and~}p\nmid b,\\
(p-1)(\frac{a}{p}),\ \ &\text{for~} p\nmid a\text{~and~}p\mid \Delta,\\
-(\frac{a}{p}),&\text{for~} p\nmid a\Delta.
\end{cases}
\end{align*}

Therefore, we may obtain explicit expressions for $\cN_1(f,q)$ and $\cN_2(f,q)$ with general $q$, and these would be much simpler if $(a,q)=1.$ In summary,  we obtain for $(a,q)=1$ that
\begin{align*}
\cN_1(f,q)
&=\varphi(q)\prod_{p\mid q}\Big(1-\frac{(\frac{\Delta}{p})}{p-1}\Big)
\end{align*}
and
\begin{align*}
\cN_2(f,q)
&=\Big(\frac{a}{q}\Big)\cN_1(f,q)\cdot\prod_{\substack{p^\nu\parallel q,~p\nmid \Delta\\ 2\nmid \nu}}\frac{-1}{p-1-(\frac{\Delta}{p})}.
\end{align*}
We thus conclude the following lemma.

\begin{lemma}\label{lm:countingresidues}
Keep the above convention. If $(q,2a)=1,$ then we have 
\begin{align*}
\cN_\pm(f,q)&=\frac{\varphi(q)}{2}\prod_{p\mid q}\Big(1-\frac{(\frac{\Delta}{p})}{p-1}\Big)\cdot\Big\{1\pm\Big(\frac{a}{q}\Big)\prod_{\substack{p^\nu\parallel q,~p\nmid \Delta\\ 2\nmid \nu}}\frac{-1}{p-1-(\frac{\Delta}{p})}\Big\}.
\end{align*}
\end{lemma}

\section{Quartic residues and reciprocity}
We turn to the theory of quartic residues in this section, and we concentrate to the case of prime moduli $p$. An integer $a$ coprime to $p$ is said to be a quartic residue modulo $p$ if the equation
\begin{align}
x^4\equiv a\bmod p
\end{align}
has solutions in $x\in\Z.$ Motivated by the Legendre symbol for quadratic residues in $\Z/p\Z$, one may also introduce the quartic residue symbol by characterizing $a^{(p-1)/(4,p-1)}\bmod p.$
Similar to Legendre symbol, we hope this quartic residue symbol also satisfies many important properties. Unfortunately, however, Gau{\ss} was not able to establish a reciprocity law for quartic residues in $\Z$ or $\Z/p\Z.$ He then moved to 
$$\Z[i]:=\{a+bi:a,b\in\Z\},$$
the so-called ring of Gaussian integers\index{Gaussian integer}. 
He did obtain quartic (and also of course quadratic) reciprocity laws in $\Z[i]$ as we will state below.

There are four roots of unity: $\pm1,\pm i.$
It is well-known that a prime $p\in\Z^+$ ramifies or splits in $\Z[i]$ if and only if $p=2$ or $p\equiv1\bmod4.$ The primes in $\Z[i]$, usually called Gaussian primes\index{Gaussian prime}, must come from one of the following three cases up to unit multiples:
\begin{enumerate}[(i)]
\item $1+i;$
\item $\pi\in\Z[i]$ with $N(\pi)\equiv1\bmod4$ being a rational prime;
\item $p\in\Z^+$ with $p\equiv 3\bmod4.$
\end{enumerate}

\begin{definition}
A Gaussian integer $\alpha$ is called odd if $1+i\nmid\alpha;$ it is called primary if $\alpha\equiv1\bmod{2+2i}.$
\end{definition}

Suppose $\xi=a+bi$ with $a,b\in\Z.$ Then $\xi$ is primary in $\Z[i]$ if and only if
\begin{align*}
a+b\equiv 1\bmod4,\ \ 2\mid b.
\end{align*}

Let $\pi$ be an odd primary prime in $\Z[i].$ The quartic residue symbol\index{quartic residue symbol} $(\frac{\cdot}{\pi})_4$ takes values in $\{0,\pm1,\pm i\}$ and can be defined explicitly via
\[\Big(\frac{\alpha}{\pi}\Big)_4\equiv \alpha^{(N(\pi)-1)/4}\bmod \pi\]
for all $\alpha\in\Z[i]$.
This is an analogue of the Legendre symbol in $\Z$, and the corresponding Jacobi symbol can also be defined similarly by multiplicativity.
In what follows, all results on quartic residues are valid for odd Gaussian integers that are not necessarily primes.

\begin{lemma}\label{lm:quarticsymbol-basicproperty}
Suppose $\xi,\xi_1,\xi_2\in\Z[i]$ are odd and primary. For $\alpha,\alpha_1,\alpha_2\in\Z[i],$ we have
\begin{enumerate}[{\rm(i)}]
\item $(\frac{i}{\xi})_4=i^{(N(\xi)-1)/4}.$
\item $(\frac{\alpha_1\alpha_2}{\xi})_4=(\frac{\alpha_1}{\xi})_4(\frac{\alpha_2}{\xi})_4.$
\item $(\frac{\alpha}{\xi_1\xi_2})_4=(\frac{\alpha}{\xi_1})_4(\frac{\alpha}{\xi_2})_4.$
\item $\overline{(\frac{\alpha}{\xi})_4}=(\frac{\overline{\alpha}}{\overline{\xi}})_4.$
\end{enumerate}
\end{lemma}

\begin{lemma}[Quartic reciprocity law]\label{lm:quarticreciprocity}\index{quartic reciprocity law}
Suppose $\xi,\eta\in\Z[i]$ are primary. Then
\begin{align*}
\Big(\frac{\eta}{\xi}\Big)_4=\Big(\frac{\xi}{\eta}\Big)_4(-1)^{\frac{N(\xi)-1}{4}\frac{N(\eta)-1}{4}}.
\end{align*}

We also have the supplementary laws 
\begin{align*}
\Big(\frac{i}{\xi}\Big)_4=i^{\frac{N(\xi)-1}{4}},\ \ \Big(\frac{1+i}{\xi}\Big)_4=i^{\frac{a-b-b^2-1}{4}},
\end{align*}
where $\xi=a+bi$ is primary with $a+b-1\equiv2b\equiv0\bmod 4$.
\end{lemma}

Note that $2=i^3(1+i)^2.$ The following lemma is obvious thanks to the above two supplementary laws.
\begin{lemma}\label{lm:quarticsymbol,at2}
Suppose $\xi=a+bi\in\Z[i]$ is primary with $a,b\in\Z$ and $(N(\xi),2)=1$. Then
\begin{align*}
\Big(\frac{2}{\xi}\Big)_4=i^{\frac{3a^2+b^2+2a-2b-5}{4}}.
\end{align*} 
\end{lemma}

\begin{lemma}\label{lm:Sunlemma}
Suppose $\xi\in\Z[i]$ and $p$ is an odd rational prime. Then
\begin{align*}
\Big(\frac{\xi}{p}\Big)_4^2=\Big(\frac{N(\xi)}{p}\Big).
\end{align*}
\end{lemma}

\proof
This elegant formula can be found in \cite[Lemma 2.1]{Su01}, and we present Sun's proof here. Without loss of generality, assume $\xi$ and $p$ are coprime.
We first consider $p=\pi\overline{\pi}\equiv1\bmod 4$ with primary $\pi\in\Z[i]$. Then
\begin{align*}
\Big(\frac{\xi}{p}\Big)_4^2=\Big(\frac{\xi}{\pi}\Big)_4^2\Big(\frac{\xi}{\overline{\pi}}\Big)_4^2=\Big(\frac{\xi}{\pi}\Big)_4^2\Big(\frac{\overline{\xi}}{\pi}\Big)_4^2
\end{align*}
since $(\frac{\xi}{\overline{\pi}})_4^2$ is real. Hence
\begin{align*}
\Big(\frac{\xi}{p}\Big)_4^2=\Big(\frac{N(\xi)}{\pi}\Big)_4^2\equiv N(\xi)^{\frac{N(\pi)-1}{2}}\equiv N(\xi)^{\frac{p-1}{2}}\bmod\pi.
\end{align*}
The desired formula follows by noting that $(\frac{N(\xi)}{p})\equiv N(\xi)^{\frac{p-1}{2}}\bmod\pi$ as given by Euler's criterion \eqref{eq:Eulercriterion}.

The case of $p\equiv3\bmod4$ can be proven similarly. In fact, 
\begin{align*}
\Big(\frac{\xi}{p}\Big)_4^2\equiv \xi^{\frac{p^2-1}{2}}\equiv (a+bi)^{\frac{(p-1)(p+1)}{2}}\bmod p.
\end{align*}
Since $\xi^p\equiv \overline{\xi}\bmod p$, we then find
\begin{align*}
\xi^{\frac{p+1}{2}}\Big(\frac{\xi}{p}\Big)_4^2\equiv\overline{\xi}^{\frac{p+1}{2}}\bmod p,
\end{align*}
from which it follows that
\begin{align*}
N(\xi)^{\frac{p+1}{2}}\equiv\xi^{p+1}\Big(\frac{\xi}{p}\Big)_4^2\equiv N(\xi)\Big(\frac{\xi}{p}\Big)_4^2\bmod p.
\end{align*}
This proves the lemma in view of $N(\xi)^{\frac{p-1}{2}}\equiv(\frac{N(\xi)}{p})\bmod p$ by Euler's criterion again.
\endproof

\begin{lemma}\label{lm:quarticsymbol=1}
Let $q\in\Z^+$ be odd. For any rational integer $a$ coprime to $q,$ we have
\begin{align*}
\Big(\frac{a}{q}\Big)_4=1.
\end{align*}
\end{lemma}
\proof
Due to multiplicativity, it suffices to consider the case when $q$ is prime. If $q\equiv1\bmod4,$ then $q$ splits in $\Z[i]$, so that $q=\pi\overline{\pi},$ which gives
\begin{align*}
\Big(\frac{a}{q}\Big)_4=\Big(\frac{a}{\pi}\Big)_4\Big(\frac{a}{\overline{\pi}}\Big)_4=\Big(\frac{a}{\pi}\Big)_4\overline{\Big(\frac{a}{\pi}\Big)_4}=1.
\end{align*}
If $q\equiv3\bmod4,$ then
\begin{align*}
\Big(\frac{a}{q}\Big)_4\equiv a^{\frac{q^2-1}{4}}\equiv a^{(q-1)\cdot\frac{q+1}{4}}\equiv1\bmod q.
\end{align*}
The desired result then follows by Fermat little theorem or Euler's theorem.
\endproof

\section{An incomplete character sum and quartic Gau{\ss} sums}

The Legendre/Jacobi symbol can define various character sums when studying quadratic residues and non-residues in $\Z$.
The similar phenomenon can also happen to the quartic residue symbol in $\Z[i]$.

Let $q\in\Z^+$ be an odd integer, $t\in\Z$ and $\beta\in\{\pm1,\pm i\}.$ Put
\begin{align}\label{eq:quarticcharactersum-incomplete}
\cQ_\beta(q,t):=\sum_{\substack{z\in\Z[i]/(q)\\ \Tr(\beta z)\equiv2t\bmod q}}\Big(\frac{z}{q}\Big)_4.
\end{align}
This is an incomplete character sum due to the restraint $\Tr(\beta z)\equiv2t\bmod q.$ The main task of this section is to explore some basic properties of $\cQ_\beta(q,t),$ and in particular, to give an explicit formula for evaluations.
To this end, we also introduce a complete character sum
\begin{align}\label{eq:quarticGausssum}
\fG(q):=\sum_{z\in\Z[i]/(q)}\Big(\frac{z}{q}\Big)_4\ue\Big(\frac{\Tr(z)}{q}\Big).
\end{align}
We shall call it a quartic Gau{\ss}  sum\index{quartic Gau{\ss}  sum}. We may find two identities between $\cQ_\beta(q,t)$ and $\fG(q)$, which allow us to succeed in finishing the exact evaluations.

Making the change of variable $z\rightarrow z/\beta,$ we find 
\begin{align*}
\cQ_\beta(q,t)=\overline{\Big(\frac{\beta}{q}\Big)_4}\cQ(q,t)=\Big(\frac{\beta}{q}\Big)_4\cQ(q,t),
\end{align*}
where we write $\cQ_1(q,t)=\cQ(q,t)$ henceforth. Here a basic fact $(\frac{\beta}{q})_4\in\R$ has been used. Furthermore,
the Chinese remainder theorem yields
\begin{align*}
\cQ(q,t)=\prod_{p^\nu\parallel q}\cQ(p^\nu,t),
\end{align*}
so that it suffices to treat $\cQ(q,t)$ when $q=p^\nu$ is a prime power. We now proceed separately according to $\nu\equiv1,2,3,0\bmod4.$

{\bf Case I.} $\nu\equiv1\bmod4$.
\begin{align*}
\cQ(p^\nu,t)=\sum_{\substack{z\in\Z[i]/(p^\nu)\\ \Tr(z)\equiv2t\bmod{p^\nu}}}\Big(\frac{z}{p}\Big)_4=p^{\nu-1}\cQ(p,t).
\end{align*}

{\bf Case II.} $\nu\equiv2\bmod4$.
\begin{align*}
\cQ(p^\nu,t)&=\sum_{\substack{z\in\Z[i]/(p^\nu)\\ \Tr(z)\equiv2t\bmod{p^\nu}}}\Big(\frac{z}{p}\Big)_4^2=\sum_{\substack{z\in\Z[i]/(p^\nu)\\ \Tr(z)\equiv2t\bmod{p^\nu}}}\Big(\frac{N(z)}{p}\Big)\\
&=p^{\nu-1}\sum_{\substack{z\in\Z[i]/(p)\\ \Tr(z)\equiv2t\bmod p}}\Big(\frac{N(z)}{p}\Big)=p^{\nu-1}\sum_{b\bmod p}\Big(\frac{b^2+t^2}{p}\Big),
\end{align*}
from which it follows that
\begin{align*}
\cQ(p^\nu,t)=
\begin{cases}
-p^{\nu-1},\ \ & p\nmid t,\\
p^{\nu-1}(p-1),\ \ & p\mid t.\end{cases}
\end{align*}

{\bf Case III.} $\nu\equiv3\bmod4$.
\begin{align*}
\cQ(p^\nu,t)=\sum_{\substack{z\in\Z[i]/(p^\nu)\\ \Tr(z)\equiv2t\bmod{p^\nu}}}\Big(\frac{z}{p}\Big)_4^3=\sum_{\substack{z\in\Z[i]/(p^\nu)\\ \Tr(z)\equiv2t\bmod{p^\nu}}}\Big(\frac{\overline{z}}{p}\Big)_4=p^{\nu-1}\cQ(p,t).
\end{align*}

{\bf Case IV.} $\nu\equiv0\bmod4$ with $\nu>0.$
\begin{align*}
\cQ(p^\nu,t)
&=\sum_{\substack{b\bmod {p^\nu}\\t+bi,~p\text{ coprime}}}1
=\sum_{\substack{b\bmod {p^\nu}\\(b^2+t^2,p)=1}}1
=p^\nu-\sum_{\substack{b\bmod {p^\nu}\\b^2\equiv-t^2\bmod p}}1\\
&=p^\nu-p^{\nu-1}\Big(1+\Big(\frac{-t^2}{p}\Big)\Big).
\end{align*}
More explicitly, we may write
\begin{align*}
\cQ(p^\nu,t)=
\begin{cases}
p^{\nu-1}(p-2),\ \ & p\equiv1\bmod4, p\nmid t,\\
p^\nu,\ \ & p\equiv3\bmod4, p\nmid t,\\
p^{\nu-1}(p-1),\ \ & p\mid t.\end{cases}
\end{align*}

Therefore, the essential step is to evaluate $\cQ(p,t)$.
To do so, we write
\begin{align*}
\cQ(p,t)&=\frac{1}{p}\sum_{h\bmod p}\ue\Big(\frac{-2ht}{p}\Big)\sum_{z\in\Z[i]/(p)}\Big(\frac{z}{p}\Big)_4\ue\Big(\frac{h\Tr(z)}{p}\Big)\\
&=\frac{1}{p}~\sideset{}{^*}\sum_{h\bmod p}\ue\Big(\frac{-2ht}{p}\Big)\sum_{z\in\Z[i]/(p)}\Big(\frac{z}{p}\Big)_4\ue\Big(\frac{h\Tr(z)}{p}\Big).
\end{align*}
Note that $h\Tr(z)=\Tr(hz)$ for each $h\in\Z.$ Making the change of variable $z\rightarrow h^{-1}z,$ we further have
\begin{align*}
\cQ(p,t)
&=\frac{\fG(p)}{p}~\sideset{}{^*}\sum_{h\bmod p}\Big(\frac{h^{-1}}{p}\Big)_4\ue\Big(\frac{-2ht}{p}\Big).
\end{align*}
By Lemma \ref{lm:quarticsymbol=1}, we find $(\frac{h^{-1}}{p})_4=1$ for any $h\in(\Z/p\Z)^\times$. It then follows that
\begin{align}\label{eq:identity1:Q-G}
\cQ(p,t)
&=\frac{\delta_p(t)\fG(p)}{p},
\end{align}
where $\delta_p(t)$ is $-1$ if $p\nmid t,$ and is $p-1$ otherwise. One can easily show that $|\fG(p)|=p$, the square root of $N(p)$ in $\Z[i],$ but the explicit determination of $\fG(p)$ requires a few arguments. 

\begin{proposition}\label{prop:quarticgausssum}
Let $p\in\Z^+$ be an odd prime. Then we have
\begin{align*}
\fG(p)=p\cdot\Big(\frac{i}{p}\Big)_4.
\end{align*}
\end{proposition}

\proof
Note that
\begin{align*}
\fG(p)&=\mathop{\sum\sum}_{x,y\bmod p}\Big(\frac{x+iy}{p}\Big)_4\ue\Big(\frac{2x}{p}\Big)\\
&=\sum_{x\bmod p}\Big(\frac{x}{p}\Big)_4\ue\Big(\frac{2x}{p}\Big)+\sideset{}{^*}\sum_{y\bmod p}\sum_{x\bmod p}\Big(\frac{x+iy}{p}\Big)_4\ue\Big(\frac{2x}{p}\Big).
\end{align*}
For $p\nmid x,$ one has $(\frac{x}{p})_4=1,$ so that the first term in the last line is exactly $-1$. Regarding the second term, we may make the change of variable $x\mapsto xy,$ so that
\begin{align*}
\fG(p)&=\sum_{x\bmod p}\Big(\frac{x+i}{p}\Big)_4~\sideset{}{^*}\sum_{y\bmod p}\Big(\frac{y}{p}\Big)_4\ue\Big(\frac{2xy}{p}\Big)-1\\
&=\sum_{x\bmod p}\Big(\frac{x+i}{p}\Big)_4~\sideset{}{^*}\sum_{y\bmod p}\ue\Big(\frac{2xy}{p}\Big)-1.
\end{align*}
For $p\mid x,$ the inner sum over $y$ is exactly $p-1$, and for $p\nmid x,$ the $y$-sum turns to be $-1$. Hence
\begin{align*}
\fG(p)&=\sum_{x\bmod p}\Big(\frac{x+i}{p}\Big)_4~\sideset{}{^*}\sum_{y\bmod p}\Big(\frac{y}{p}\Big)_4\ue\Big(\frac{2xy}{p}\Big)-1\\
&=(p-1)\Big(\frac{i}{p}\Big)_4-\sideset{}{^*}\sum_{x\bmod p}\Big(\frac{x+i}{p}\Big)_4-1\\
&=p\Big(\frac{i}{p}\Big)_4-\sum_{x\bmod p}\Big(\frac{x+i}{p}\Big)_4-1\\
&=p\Big(\frac{i}{p}\Big)_4-\Big(\frac{i}{p}\Big)_4\cQ(p,1)-1.
\end{align*}

On the other hand, we have
\begin{align}\label{eq:identity2:Q-G}
\cQ(p,1)
&=-\frac{\fG(p)}{p}.
\end{align}
The proposition then follows by combining \eqref{eq:identity1:Q-G} and \eqref{eq:identity2:Q-G}.
\endproof

Collecting the above evaluations for $\cQ(p^\nu,t)$ in different situations, i.e., $2\nmid \nu,$ $\nu\equiv2\bmod4$ and $4\mid \nu$, we may obtain an explicit formula evaluating $\cQ(q,t)$ for an arbitrary odd positive integer $q$.

We now summarize the above results as the following proposition.

\begin{proposition}\label{prop:incompletecharactersum-quarticgauss}
Let $q$ be an odd positive integer, $t\in\Z$ and $\beta\in\{\pm1,\pm i\}.$ Then
\begin{enumerate}[{\rm(i)}]
\item $\cQ_\beta(q,t)=(\frac{\beta}{q})_4\cQ(q,t).$
\item For fixed $t$ and $\beta,$ the function $q\mapsto \cQ_\beta(q,t)$ is multiplicative, i.e., for each odd $q_1,q_2\geqslant1$ with $(q_1,q_2)=1,$ it holds that $\cQ_\beta(q_1q_2,t)=\cQ_\beta(q_1,t)\cQ_\beta(q_2,t)$.
\item In general, we have
\begin{align*}
\cQ_\beta(q,t)&=\varphi(q)\Big(\frac{i\beta}{q}\Big)_4\cdot\Omega_4(1;q,t)\cdot \prod_{\substack{p\mid q\\ p\nmid t}}\Big(1-\frac{(\frac{-1}{p})}{p-1}\Big),
\end{align*}
where we recall the definition in $\eqref{eq:Omegafunction}$
\begin{align*}
\Omega_j(1;q,r)=\prod_{\substack{p^\nu\parallel q,~ p\nmid r\\j\nmid \nu}}\frac{-1}{p-1-(\frac{-1}{p})}.
\end{align*}
\end{enumerate}
\end{proposition}

\section{Cubic residues and reciprocity}
Rather than in $\Z,$ the beautiful theory of cubic residues are built in $\Z[\omega]$, where $\omega=\ue^{2\pi i/3}=-\frac{1}{2}+\frac{\sqrt{3}}{2}i.$ The original story starts from the solvability of 
\begin{align*}
x^3\equiv a\bmod p
\end{align*}
in $x\in\Z,$ where $p\in\Z^+$ is a prime and $(a,p)=1.$ Similar to the case on quartic residues, no satisfactory laws of cubic reciprocity have been established in $\Z,$ but the situation becomes reasonable in $\Z[\omega]$.

There are six roots of unity in $\Z[\omega]:$ $\pm1,\pm\omega,\pm\omega^2.$ A prime $p\in\Z^+$ ramifies or remains prime in $\Z[\omega]$ if and only if $p=3$ or $p\equiv2\bmod3.$ In particular, a prime $p\in\Z^+$ with $p\equiv1\bmod3$ splits in $\Z[\omega],$ and there exists $\pi=a+b\omega$ with $a,b\in\Z$ such that 
\begin{align*}
p=N(\pi)=\pi\overline{\pi}=a^2-ab+b^2.
\end{align*}
Elements in $\Z[\omega]$ are usually called Eisenstein integers\index{Eisenstein integer}, and primes in $\Z[\omega]$ are called Eisenstein primes \index{Eisenstein prime} correspondingly. Eisenstein primes must come from one of the following three cases up to unit multiples:
\begin{enumerate}[(i)]
\item $1-\omega;$
\item $\pi\in\Z[\omega]$ with $N(\pi)\equiv1\bmod3$ being a rational prime;
\item $p\in\Z^+$ with $p\equiv 2\bmod3.$
\end{enumerate}
Alternatively speaking, among all primes $p\in\Z^+,$ $3=-\omega^2(1-\omega)^2$ ramifies in $\Z[\omega]$, $p\equiv1\bmod3$ splits in $\Z[\omega]$, and $p\equiv2\bmod3$ is then inert in $\Z[\omega]$.

\begin{definition}
An Eisenstein integer $\alpha$ is called primary if $\alpha\equiv2\bmod3.$
\end{definition}

Suppose $\xi=a+b\omega$ with $a,b\in\Z.$ Then $\xi$ is primary in $\Z[\omega]$ if and only if
\begin{align*}
a\equiv -1\bmod3,\ \ 3\mid b.
\end{align*}

Let $\pi$ be a primary Eisenstein prime in $\Z[\omega].$ The cubic residue symbol $(\frac{\cdot}{\pi})_3$ takes values in $\{0,1, \omega,\omega^2\}$ and can be defined explicitly via
\[\Big(\frac{\alpha}{\pi}\Big)_3\equiv \alpha^{(N(\pi)-1)/3}\bmod \pi\]
for all $\alpha\in\Z[\omega]$.
The corresponding Jacobi symbol \index{Jacobi symbol} over $\Z[\omega]$ can also be defined similarly by multiplicativity, and the subsequent results are served for Eisenstein integers that are not necessarily primes.

\begin{lemma}\label{lm:cubicsymbol-basicproperty}
Suppose $\xi,\xi_1,\xi_2$ are primary Eisenstein integers. For $\alpha,\alpha_1,\alpha_2\in\Z[\omega],$ we have
\begin{enumerate}[{\rm(i)}]
\item $(\frac{\pm1}{\xi})_3=1,$ $(\frac{\omega}{\xi})_3=\omega^{(N(\xi)-1)/3}.$
\item $(\frac{\alpha_1\alpha_2}{\xi})_3=(\frac{\alpha_1}{\xi})_3(\frac{\alpha_2}{\xi})_3.$
\item $(\frac{\alpha}{\xi_1\xi_2})_3=(\frac{\alpha}{\xi_1})_3(\frac{\alpha}{\xi_2})_3.$
\item $\overline{(\frac{\alpha}{\xi})_3}=(\frac{\alpha}{\xi})_3^2=(\frac{\overline{\alpha}}{\overline{\xi}})_3.$
\end{enumerate}
\end{lemma}

\begin{lemma}[Cubic reciprocity law]\label{lm:cubicreciprocity}\index{cubic reciprocity law}
Suppose $\xi,\eta\in\Z[\omega]$ are primary with $3\nmid N(\xi)N(\eta)$ and $ N(\xi)\neq N(\eta)$. Then
\begin{align*}
\Big(\frac{\eta}{\xi}\Big)_3=\Big(\frac{\xi}{\eta}\Big)_3.
\end{align*}

We also have the supplementary laws 
\begin{align*}
\Big(\frac{\omega}{\xi}\Big)_3=\omega^{\frac{N(\xi)-1}{3}},\ \ \Big(\frac{1-\omega}{\xi}\Big)_3=\omega^{\frac{2(a+1)}{3}},
\end{align*}
where $\xi=a+b\omega$ is primary with $a+1\equiv b\equiv 0\bmod 3$ and $a,b\in\Z$.
\end{lemma}

\begin{lemma}\label{lm:cubicsymbol-specialvalues}
Suppose $\xi=a+b\omega\in\Z[\omega]$ is primary with $(N(\xi),6)=1$. Then
\begin{align*}
\Big(\frac{2}{\xi}\Big)_3=\omega^{2+(-1)^\frac{a-2}{3}-(-1)^{\frac{b}{3}}},\ \ \Big(\frac{1+2\omega}{\xi}\Big)_3=\omega^{\frac{b}{3}},\ \ \Big(\frac{3}{\xi}\Big)_3=\omega^{\frac{2b}{3}}.
\end{align*}
\end{lemma}

\proof
The condition $(N(\xi),6)=1$ implies that $\xi$ is coprime to 2, and thus
\begin{align*}
\Big(\frac{2}{\xi}\Big)_3=\Big(\frac{\xi}{2}\Big)_3=\Big(\frac{a+b\omega}{2}\Big)_3=
\begin{cases}
1,\ \ &\text{if~}2\nmid a,2\mid b,\\
(\frac{\omega}{2})_3,\ \ &\text{if~}2\mid a,2\nmid b,\\
(\frac{1+\omega}{2})_3,\ \ &\text{if~}2\nmid ab.
\end{cases}
\end{align*}
From Lemma \ref{lm:cubicreciprocity} we find $(\frac{1+\omega}{2})_3=(\frac{-\omega^2}{2})_3=(\frac{\omega}{2})_3^2=\omega^2,$ and thus
\begin{align*}
\Big(\frac{2}{\xi}\Big)_3=\begin{cases}
1,\ \ &\text{if~}2\nmid a,2\mid b,\\
\omega,\ \ &\text{if~}2\mid a,2\nmid b,\\
\omega^2,\ \ &\text{if~}2\nmid ab
\end{cases}
\end{align*}
as desired.

On the other hand, the supplementary laws in Lemma \ref{lm:cubicreciprocity} 
yield
\begin{align*}
\Big(\frac{1+2\omega}{\xi}\Big)_3
=\Big(\frac{\omega-\omega^2}{\xi}\Big)_3=\Big(\frac{\omega}{\xi}\Big)_3\Big(\frac{1-\omega}{\xi}\Big)_3=\omega^{\frac{b}{3}}
\end{align*}
and
\begin{align*}
\Big(\frac{3}{\xi}\Big)_3
=\Big(\frac{-(1+2\omega)^2}{\xi}\Big)_3=\omega^{\frac{2b}{3}}.
\end{align*}
This completes the proof of the lemma.
\endproof

\begin{lemma}\label{lm:cubicsymbol=1}
Let $q\in\Z^+$ be primary. For any rational integer $a$ coprime to $q,$ we have
\begin{align*}
\Big(\frac{a}{q}\Big)_3=1.
\end{align*}
\end{lemma}
The proof of Lemma \ref{lm:cubicsymbol=1} is left to interested readers.

\section{An incomplete character sum and cubic Gau{\ss} sums}
Similar to $\cQ_\beta(q,t)$ given by \eqref{eq:quarticcharactersum-incomplete}, we now consider an incomplete character sum involving the cubic residue symbol. Let $q$ be a positive integer with $(q,6)=1$ and $t\in\Z.$ For any 
$\beta\in\{\pm1,\pm \omega,\pm\omega^2\}$ and $\kappa\in\{0,1\},$ put
\begin{align}\label{eq:cubiccharactersum-incomplete}
\cC_\beta(q,t;\kappa)&=\sum_{\substack{z\in\Z[\omega]/(q)\\ \Tr(\beta z)\equiv t\bmod q}}\Big(\frac{z}{q}\Big)_3\Big(\frac{N(z)}{q}\Big)^\kappa.
\end{align} 
As an analogue of $\fG(q),$ we introduce the following cubic (in fact sextic if $\kappa=1$) Gau{\ss} sum \index{cubic Gau{\ss} sum}
\begin{align}\label{eq:cubicGausssum}
\fC(q,\kappa):=\sum_{z\in\Z[i]/(q)}\Big(\frac{z}{q}\Big)_3\Big(\frac{N(z)}{q}\Big)^\kappa\ue\Big(\frac{\Tr(z)}{q}\Big)
\end{align}
for the purpose of evaluating $\cC_\beta(q,t;\kappa).$

Making the change of variable $z\rightarrow z/\beta,$ we find 
\begin{align*}
\cC_\beta(q,t;\kappa)=\Big(\frac{\overline{\beta}}{q}\Big)_3\cC(q,t;\kappa),
\end{align*}
where we write $\cC_1(q,t;\kappa)=\cC(q,t;\kappa)$ henceforth. Furthermore,
the Chinese remainder theorem yields
\begin{align*}
\cC(q,t;\kappa)=\prod_{p^\nu\parallel q}\cC(p^\nu,t;\kappa),
\end{align*}
so that it suffices to treat the case of prime power moduli $q=p^\nu$. We now proceed separately according to $\nu\equiv1,2,0\bmod3.$

{\bf Case I.} $\nu\equiv1\bmod3$.
\begin{align*}
\cC(p^\nu,t;\kappa)&=\sum_{\substack{z\in\Z[\omega]/(p^\nu)\\ \Tr(z)\equiv t\bmod {p^\nu}}}\Big(\frac{z}{p}\Big)_3\Big(\frac{N(z)}{p}\Big)^{\nu\kappa} =p^{\nu-1}\cC(p,t;\nu\kappa).
\end{align*}

{\bf Case II.} $\nu\equiv2\bmod3$.
\begin{align*}
\cC(p^\nu,t;\kappa)&=\sum_{\substack{z\in\Z[\omega]/(p^\nu)\\ \Tr(z)\equiv t\bmod {p^\nu}}}\Big(\frac{z}{p}\Big)_3^2\Big(\frac{N(z)}{p}\Big)^{\nu\kappa} =\sum_{\substack{z\in\Z[\omega]/(p^\nu)\\ \Tr(z)\equiv t\bmod {p^\nu}}}\Big(\frac{\overline{z}}{p}\Big)_3\Big(\frac{N(z)}{p}\Big)^{\nu\kappa},
\end{align*}
so that
\begin{align*}
\cC(p^\nu,t;\kappa)&=p^{\nu-1}\overline{\cC(p,t;\nu\kappa)}.
\end{align*}

{\bf Case III.} $\nu\equiv0\bmod3$.
\begin{align*}
\cC(p^\nu,t;\kappa)&=\sum_{\substack{z\in\Z[\omega]/(p^\nu)\\ \Tr(z)\equiv t\bmod {p^\nu}}}\Big(\frac{N(z)}{p}\Big)^{\nu\kappa} =p^{\nu-1}\sum_{a\bmod p}\Big(\frac{3a^2+3at+t^2}{p}\Big)^{\nu\kappa}.
\end{align*}
Making the change of variable $a\mapsto a-\overline{2}t$ we find
\begin{align*}
\cC(p^\nu,t;\kappa)&=p^{\nu-1}\sum_{a\bmod p}\Big(\frac{12a^2+t^2}{p}\Big)^{\nu\kappa}.
\end{align*}
If $2\mid \nu\kappa,$ then
\begin{align*}
\cC(p^\nu,t;\kappa)&=p^{\nu-1}\sum_{\substack{a\bmod p\\(12a^2+t^2,p)=1}}1=p^{\nu-1}\Big(p-1-\Big(\frac{-3t^2}{p}\Big)\Big).
\end{align*}
If $2\nmid \nu\kappa,$ then
\begin{align*}
\cC(p^\nu,t;\kappa)&=p^{\nu-1}\sum_{a\bmod p}\Big(\frac{12a^2+t^2}{p}\Big)=p^{\nu-1}\Big(\frac{3}{p}\Big)\times
\begin{cases}
(p-1),\ \ & \text{if~}p\mid t,\\
-1, &\text{if~} p\nmid t.
\end{cases}
\end{align*}

Collecting all above arguments, the evaluation of $\cC(q,t;\kappa)$ reduces to the case of prime moduli.
In fact,
\begin{align*}
\cC(p,t;\kappa)&=\frac{1}{p}\sum_{h\bmod p}\ue\Big(\frac{-ht}{p}\Big)\sum_{z\in\Z[\omega]/(p)}\Big(\frac{z}{p}\Big)_3\Big(\frac{N(z)}{p}\Big)^\kappa\ue\Big(\frac{h\Tr(z)}{p}\Big)\\
&=\frac{1}{p}~\sideset{}{^*}\sum_{h\bmod p}\ue\Big(\frac{-ht}{p}\Big)\sum_{z\in\Z[\omega]/(p)}\Big(\frac{z}{p}\Big)_3\Big(\frac{N(z)}{p}\Big)^\kappa\ue\Big(\frac{h\Tr(z)}{p}\Big),
\end{align*}
where we have used a fact that
\begin{align*}
\sum_{z\in\Z[\omega]/(p)}\Big(\frac{z}{p}\Big)_3\Big(\frac{N(z)}{p}\Big)^\kappa=0
\end{align*}
since $z\mapsto(\frac{z}{p})_3(\frac{N(z)}{p})^\kappa$ defines a non-trivial character on $(\Z[\omega]/(p))^\times$ for each $\kappa\in\{0,1\}.$

Making the change of variable $z\mapsto h^{-1}z$, we get
\begin{align*}
\cC(p,t;\kappa)
&=\frac{\fC(p,\kappa)}{p}~\sideset{}{^*}\sum_{h\bmod p}\ue\Big(\frac{-ht}{p}\Big)
\end{align*}
in view of the fact that $(\frac{h^{-1}}{p})_3=1$ for any $h\in(\Z/p\Z)^\times$. Therefore,
\begin{align}\label{eq:cubicgausssum-relation1}
\cC(p,t;\kappa)
&=\frac{\delta_p(t)\fC(p,\kappa)}{p},
\end{align}
where $\delta_p(t)$ is $-1$ if $p\nmid t,$ and is $p-1$ otherwise. 

We now turn to establish another identity between $\cC(p,t;\kappa)$ and $\fC(p,\kappa)$. By definition we may write
\begin{align*}
\fC(p,\kappa)&=\mathop{\sum\sum}_{x,y\bmod p}\Big(\frac{x+y\omega}{p}\Big)_3\Big(\frac{x^2-xy+y^2}{p}\Big)^\kappa\ue\Big(\frac{2x-y}{p}\Big)\\
&=-1+\mathop{\sideset{}{^*}\sum\sum}_{y,x\bmod p}\Big(\frac{x+y\omega}{p}\Big)_3\Big(\frac{x^2-xy+y^2}{p}\Big)^\kappa\ue\Big(\frac{2x-y}{p}\Big).
\end{align*}
Making the change of variable $x\mapsto yx$, we get
\begin{align*}
\fC(p,\kappa)
&=-1+\sum_{x\bmod p}\Big(\frac{x+\omega}{p}\Big)_3\Big(\frac{x^2-x+1}{p}\Big)^\kappa\sideset{}{^*}\sum_{y\bmod p}\Big(\frac{y}{p}\Big)_3\ue\Big(\frac{(2x-1)y}{p}\Big)\\
&=-1+\sum_{x\bmod p}\Big(\frac{x+\omega}{p}\Big)_3\Big(\frac{x^2-x+1}{p}\Big)^\kappa\sideset{}{^*}\sum_{y\bmod p}\ue\Big(\frac{(2x-1)y}{p}\Big).
\end{align*}
If $2x\equiv1\bmod p,$ the $y$-sum is exactly $p-1$, and should be $-1$ otherwise. Therefore
\begin{align*}
\fC(p,\kappa)
&=-1+p\Big(\frac{1+2\omega}{p}\Big)_3\Big(\frac{3}{p}\Big)^\kappa-\sum_{x\bmod p}\Big(\frac{x+\omega}{p}\Big)_3\Big(\frac{x^2-x+1}{p}\Big)^\kappa.
\end{align*}
If $p\equiv-1\bmod3$ is primary in $\Z[\omega],$ then
\begin{align*}
\Big(\frac{1+2\omega}{p}\Big)_3&=\Big(\frac{\omega}{p}\Big)_3\Big(\frac{1-\omega}{p}\Big)_3=\omega^{\frac{p^2-1}{3}}\cdot \omega^{\frac{2(p+1)}{3}}=\omega^{\frac{(p+1)^2}{3}}=1,
\end{align*}
and thus for any prime $p>3$, we always have
\begin{align*}
\Big(\frac{1+2\omega}{p}\Big)_3&=1,
\end{align*}
which gives
\begin{align*}
\fC(p,\kappa)
&=-1+p\Big(\frac{3}{p}\Big)^\kappa-\sum_{x\bmod p}\Big(\frac{x+\omega}{p}\Big)_3\Big(\frac{x^2-x+1}{p}\Big)^\kappa.
\end{align*}

On the other hand, we have
\begin{align*}
\cC(p,-1;\kappa)&=\sum_{a\bmod p}\Big(\frac{a+(2a+1)\omega}{p}\Big)_3\Big(\frac{3a^2+3a+1}{p}\Big)^{\kappa}.
\end{align*}
Making the change of variable $a/(2a+1)\mapsto x$ for $p\nmid 2a+1$, we get
\begin{align*}
\cC(p,-1;\kappa)&=\sum_{\substack{x\bmod p\\ 2x\not\equiv1\bmod p}}\Big(\frac{x+\omega}{p}\Big)_3\Big(\frac{x^2-x+1}{p}\Big)^{\kappa}+\sum_{\substack{a\bmod p\\ p\mid 2a+1}}\Big(\frac{a}{p}\Big)_3\Big(\frac{3a^2+3a+1}{p}\Big)^{\kappa}.
\end{align*}
Note that 
\begin{align*}
\sum_{\substack{x\bmod p\\ 2x\equiv1\bmod p}}\Big(\frac{x+\omega}{p}\Big)_3\Big(\frac{x^2-x+1}{p}\Big)^{\kappa}=\Big(\frac{1+2\omega}{p}\Big)_3\Big(\frac{3}{p}\Big)^\kappa=\Big(\frac{3}{p}\Big)^\kappa,
\end{align*}
and
\begin{align*}
\sum_{\substack{a\bmod p\\ p\mid 2a+1}}\Big(\frac{a}{p}\Big)_3\Big(\frac{3a^2+3a+1}{p}\Big)^{\kappa}=1,
\end{align*}
from which we conclude that
\begin{align*}
\cC(p,-1;\kappa)&=\sum_{x\bmod p}\Big(\frac{x+\omega}{p}\Big)_3\Big(\frac{x^2-x+1}{p}\Big)^{\kappa}-\Big(\frac{3}{p}\Big)^\kappa+1,
\end{align*}
and hence
\begin{align}\label{eq:cubicgausssum-relation2}
\fC(p,\kappa)=(p-1)\Big(\frac{3}{p}\Big)^\kappa-\cC(p,-1;\kappa).
\end{align}
Combining \eqref{eq:cubicgausssum-relation1} and \eqref{eq:cubicgausssum-relation2} we arrive at the identities
\begin{align*}
\fC(p,\kappa)=p\Big(\frac{3}{p}\Big)^\kappa
\end{align*}
and
\begin{align*}
\cC(p,t;\kappa)=\delta_p(t)\cdot\Big(\frac{3}{p}\Big)^\kappa.
\end{align*}

The above arguments will lead to the following proposition.

\begin{proposition}\label{prop:incompletecharactersum-cubicgauss}
Let $q\in\Z^+$ with $(q,6)=1,$ $t\in\Z,$ $\beta\in\{\pm1,\pm \omega,\pm\omega^2\}$ and $\kappa\in\{0,1\}.$ Then
\begin{enumerate}[{\rm(i)}]
\item $\cC_\beta(q,t;\kappa)=(\frac{\overline{\beta}}{q})_3\cC(q,t;\kappa).$
\item For fixed $t,\kappa$ and $\beta,$ the function $q\mapsto \cC_\beta(q,t;\kappa)$ is multiplicative, i.e., for each $q_1,q_2\geqslant1$ with $(q_1,q_2)=(q_1q_2,6)=1,$ it holds that $\cC_\beta(q_1q_2,t;\kappa)=\cC_\beta(q_1,t;\kappa)\cC_\beta(q_2,t;\kappa)$.
\item In general, we have
\begin{align*}
\cC_\beta(q,t;\kappa)&=\varphi(q)\Big(\frac{\overline{\beta}}{q}\Big)_3\Big(\frac{3}{q}\Big)^\kappa\cdot \prod_{\substack{p\mid q\\ p\nmid t}}\Big(1-\frac{(\frac{-3}{p})}{p-1}\Big)\cdot\prod_{\substack{p^\nu\parallel q\\ p\nmid t\\ 3\nmid \nu\text{~or~}2\nmid \nu\kappa}}\frac{-1}{p-1-(\frac{-3}{p})}.
\end{align*}
\end{enumerate}
\end{proposition}

Taking $\kappa=0$ or $1$ separately, we arrive that the following consequence.
\begin{corollary}\label{coro:incompletecharactersum-cubicgauss}
Let $q\in\Z^+$ with $(q,6)=1,$ $t\in\Z,$ $\beta\in\{\pm1,\pm \omega,\pm\omega^2\}.$ Then
\begin{align*}
\cC_\beta(q,t;0)&=\varphi(q)\Big(\frac{\overline{\beta}}{q}\Big)_3\cdot \Omega_3(3;q,t)\cdot\prod_{\substack{p\mid q\\ p\nmid t}}\Big(1-\frac{(\frac{-3}{p})}{p-1}\Big)
\end{align*}
and
\begin{align*}
\cC_\beta(q,t;1)&=\varphi(q)\Big(\frac{\overline{\beta}}{q}\Big)_3\Big(\frac{3}{q}\Big)\cdot \Omega_6(3;q,t)\cdot\prod_{\substack{p\mid q\\ p\nmid t}}\Big(1-\frac{(\frac{-3}{p})}{p-1}\Big),
\end{align*}
where we recall the definition in $\eqref{eq:Omegafunction}$
\begin{align*}
\Omega_j(3;q,r)=\prod_{\substack{p^\nu\parallel q,~ p\nmid r\\j\nmid \nu}}\frac{-1}{p-1-(\frac{-3}{p})}.
\end{align*}
\end{corollary}

\section{Remarks}

Gau{\ss} sums occupy important positions in both of classical and modern number theory. Gau{\ss} himself introduced the sums
\[\sum_{x\in \Z/p\Z}\ue\Big(\frac{ax^2}{p}\Big)\]
and
\[\sum_{x\in \Z/p\Z}\Big(\frac{x}{p}\Big)\ue\Big(\frac{ax}{p}\Big),\]
both of which are now known as Gau{\ss} sums, to give his sixth proof of quadratic reciprocity law.
The above two sums are both defined over $\Z/p\Z.$ Moving the definition to $R/\pi R$ with $R=\Z[i]$ or $\Z[\omega]$, where $\pi$ is a prime element in $R$ such that $N(\pi)=p\in\Z^+$ is a prime, it is also possible to define Gau{\ss} sums. In fact, $R/\pi R$ is a finite field of $p$ elements, thus the relevant Gau{\ss} sum can be defined as
\begin{align*}
G(a,\chi_\pi):=\sum_{x\in \Z/p\Z}\chi_\pi(x)\ue\Big(\frac{ax}{p}\Big),
\end{align*}
where $\chi_\pi$ denotes the quartic and cubic residue symbols, respectively. This was originally utilized by Eisenstein to prove quartic and cubic reciprocity laws. Replacing $\chi_\pi$ by its powers, one may obtain a sequence of new Gau{\ss} sums, and a linear combination of them can produce the new Gau{\ss} sum \index{Gau{\ss} sum}
\[G_k(a,p):=\sum_{x\in \Z/p\Z}\ue\Big(\frac{ax^k}{p}\Big)\]
with $k=3,4.$

To characterize the Lang--Trotter constant in later chapters, we will employ the Gau{\ss} sums $\cQ_\beta(q,t)$ and $\cC_\beta(q,t;\kappa)$. They are defined in different means rather than $G(a,\chi_\pi)$. We expect these two should have many other potential applications.

Note that $-\omega^2$ is a sixth root of unity contained in $\Z[\omega]$, in which Eisenstein also established a sextic reciprocity law\index{sextic reciprocity law}. The sextic residue symbol \index{sextic residue symbol} $(\frac{\cdot}{\pi})_6$ can be defined in a similar fashion via
\[\Big(\frac{\alpha}{\pi}\Big)_6\equiv \alpha^{(N(\pi)-1)/6}\bmod \pi\]
for primary prime $\pi\in\Z[\omega],$
and it also enjoys very similar properties just as the quartic and cubic residue symbols. To state his sextic reciprocity law explicitly, Eisenstein himself introduced the concept ``E-primary" that is not quite easy to formulate. In this memoir, we avoid the use of sextic reciprocity law by observing
\begin{align*}
\Big(\frac{\alpha}{\pi}\Big)_6^3=\Big(\frac{\alpha}{N(\pi)}\Big)
\end{align*}
if $\alpha\in\Z,$ in which case one has
\begin{align*}
\overline{\Big(\frac{\alpha}{\pi}\Big)_6}=\Big(\frac{\alpha}{\pi}\Big)_6^5=\Big(\frac{\alpha}{\pi}\Big)_6^3\Big(\frac{\alpha}{\pi}\Big)_6^2=\Big(\frac{\alpha}{N(\pi)}\Big)\Big(\frac{\alpha}{\pi}\Big)_3.
\end{align*}
All necessary arguments then utilize the theory of quadratic and cubic residues.

%% file: 4.D=1.tex
\chapter{Complex multiplication by $\Q(\sqrt{-1})$}\label{ch:D=1}

Throughout this chapter, we assume $E/\Q$ is defined by
\begin{align}\label{eq:E:D=1}
y^2=x^3-gx
\end{align}
with $g\in\Z^*.$ This is the family of elliptic curves with CM by $\Q(\sqrt{-1})$. Under the Hardy--Littlewood conjecture\index{Hardy--Littlewood conjecture}, we would like to prove that
\begin{align}\label{eq:piE,r(x)-varpiE,r,D=1}
\pi_{E,r}(x)\sim\frac{\varpi_{E,r}\cdot\sqrt{x}}{\log x}
\end{align}
with some explicit constant $\varpi_{E,r}\geqslant0.$
To describe this constant, we write
\begin{align*}
g=(-1)^\delta2^\lambda g_1
\end{align*}
with
\begin{itemize}
\item $\delta\in\{0,1\},$
\item $\lambda\in\N,$
\item $2\nmid g_1\in\Z^+$.
\end{itemize}

Recall the definition \eqref{eq:Omegafunction}:
\begin{align*}
\Omega_j(D;q,r)=\prod_{\substack{p^\nu\parallel q,~ p\nmid r\\j\nmid \nu}}\frac{-1}{p-1-(\frac{-D}{p})}.
\end{align*}
The convention that $\Omega_j=\Omega_j(1;g_1,r)$ for $j\in\{2,4\}$ works throughout this chapter.
Denote by $\fU_1=\{\pm1,\pm i\}$, the group of roots of unity in $\Z[\sqrt{-1}]$.

\section{Asymptotics for $\pi_{E,r}(x)$: conditional results}

\begin{theorem}\label{thm:HL-LT:D=1}\index{Lang--Trotter conjecture}
Suppose $E/\Q$ is the elliptic curve defined by $\eqref{eq:E:D=1}$ with the above convention.
For each non-zero integer $r$, the Hardy--Littlewood conjecture \index{Hardy--Littlewood conjecture} implies that $\eqref{eq:piE,r(x)-varpiE,r,D=1}$ holds with
\begin{align}\label{eq:varpiE,r,D=1}
\varpi_{E,r}=\frac{h_{1,r}}{4}\cdot \kappa(g,r)
\end{align}
for $2\mid r,$
where $h_{1,r}\geqslant0$ is the constant defined by $\eqref{eq:constant-hD,r}$ and
\begin{align}\label{eq:LTconstant,D=1}
\kappa(g,r)
&=\begin{cases}
1-(-1)^{\frac{\lambda r}{4}}\Omega_2+(-1)^{\frac{r}{4}(\delta+\frac{g_1-1}{2})}(1-(-1)^{\frac{\lambda r}{4}})\Re(i^{1+\frac{\lambda r}{4}})\Omega_4,\ \ &\text{for~}  4\mid r,\\
1+\Omega_2+(-1)^{\frac{r-2}{4}+\frac{g_1^2-1}{8}}(1-(-1)^{\delta+\frac{\lambda+g_1-1}{2}})\Omega_4,\ \ &\text{for~} 2\parallel r,2\mid\lambda,\\
1,\ \ &\text{for~} 2\parallel r,2\nmid\lambda.
\end{cases}
\end{align}
Unconditionally, $\pi_{E,r}(x)$ is bounded by $O(1)$ if $2\nmid r$.
\end{theorem}

To evaluate $\pi_{E,r}(x)$ asymptotically for $E$ given by \eqref{eq:E:D=1}, we do not need the full strength of the Hardy--Littlewood conjecture. As we will see later, only Conjecture \ref{conj:Gd(x;r,alpha,beta)} below is sufficient!

\section{Initial transformations}

We start with the following explicit formula for Frobenius traces \index{trace of Frobenius} $a_p(E)$. See \cite[Chapter 18, Theorem 5]{IR90} for details.
\begin{lemma}\label{lm:Frobenius:D=1}
Suppose $p\in\Z^+$ with $p\nmid 2g.$ Then $a_p(E)=0$ if $p\equiv3\bmod4,$ and
\begin{align*}
a_p(E)=\Big(\frac{g}{\pi}\Big)_4\overline{\pi}+\Big(\frac{g}{\overline{\pi}}\Big)_4\pi=\Tr\Big(\Big(\frac{g}{\overline{\pi}}\Big)_4\pi\Big)
\end{align*}
if $p\equiv1\bmod4$ with $p=\pi\overline{\pi}$ such that $\pi\in\Z[i]$ is primary.
\end{lemma}

For $r\neq0,$ we may write
\begin{align*}
2\pi_{E,r}(x)
&=\sum_{\beta\in \fU_1}\sum_{\substack{\pi\in\Pi_1(x)\\ \Tr(\beta\pi)=r,~(\frac{g}{\pi})_4=\bar{\beta}\\ \pi\text{~primary}}}1+O(1),
\end{align*}
where the term $O(1)$ comes from bad primes, and the constant $2$ appears since $\pi$ and $\overline{\pi}$ should be taken into account or not simultaneously. For each primary $\pi=a+bi\in\Pi_1(x)$ with $a,b\in\Z$, we have
\begin{align*}
\Big(\frac{g}{\pi}\Big)_4
&=\Big(\frac{-1}{\pi}\Big)_4^\delta\Big(\frac{2}{\pi}\Big)_4^{\lambda}\Big(\frac{g_1}{\pi}\Big)_4
=i^{\varepsilon(a,b)}\Big(\frac{g_1}{\pi}\Big)_4,
\end{align*}
where $\varepsilon(a,b)=\frac{1}{2}\delta(a^2+b^2-1)-\frac{1}{2}\lambda b.$
According to the four choices of $\beta$ in $\fU_1,$
\begin{align*}
2\pi_{E,r}(x)
&=\sum_{0\leqslant d\leqslant 3}\sum_{\substack{\pi\in\Pi_1(x)\\ \Tr(i^d\pi)=r,~(\frac{g_1}{\pi})_4=i^{-d-\varepsilon(a,b)}\\ \pi=a+bi\text{~primary}}}1+O(1).
\end{align*}

To make further transformations on the inner sum over $\pi$, we introduce a new counting function
\begin{align}\label{eq:Gd-countingfunction}
\cG_d(x;r,\alpha,\beta,\boldsymbol\gamma)
&:=\Big|\Big\{\pi\in\Pi_1(x):
\begin{array}{c}
\Tr(\beta\pi)=r,~(\frac{\alpha}{\pi})_4=i^d,\\
\pi\equiv2\gamma_1+1+2\gamma_2 i\bmod{16}\end{array}
\Big\}\Big|
\end{align}
for $\boldsymbol\gamma=(\gamma_1,\gamma_2)$ with $\gamma_1\equiv\gamma_2\bmod2$ and
\begin{align*}
r,\alpha\in\Z^*,\ \ \beta\in \fU_1,\ \ 0\leqslant d\leqslant3,\ \ 0\leqslant\gamma_1,\gamma_2\leqslant7.
\end{align*}

By virtue of \eqref{eq:Gd-countingfunction} we arrive at
\begin{align}\label{eq:connection:pi-Gd}
2\pi_{E,r}(x)
&=\sum_{0\leqslant d\leqslant 3}\mathop{\sum\sum}_{\substack{0\leqslant \gamma_1,\gamma_2\leqslant7\\ \gamma_1\equiv\gamma_2\bmod2\\ \boldsymbol\gamma=(\gamma_1,\gamma_2)}}\cG_{-d+(2\delta+\lambda)\gamma_2}(x;r,g_1,i^d,\boldsymbol\gamma)+O(1).
\end{align}
It remains to evaluate $\cG_d(x;r,\alpha,\beta,\boldsymbol\gamma)$ in all possible cases.

\section{Asymptotic evaluation of $\cG_d(x;r,\alpha,\beta,\text{\boldmath$\gamma$})$}
It is clear that
\begin{align*}
\sum_{0\leqslant d\leqslant3}\sum_{\beta\in\fU_1}\mathop{\sum\sum}_{\substack{0\leqslant \gamma_1,\gamma_2\leqslant7\\ \gamma_1\equiv\gamma_2\bmod2\\ \boldsymbol\gamma=(\gamma_1,\gamma_2)}}\cG_d(x;r,\alpha,\beta,\boldsymbol\gamma)=2\pi_{1,r}(x)+O(1),
\end{align*}
where $\pi_{1,r}(x)$ is defined by \eqref{eq:piD,r(x)-definition}.
We would like to formulate the following conjecture on the asymptotic proportion of $\cG_d(x;r,\alpha,\beta,\boldsymbol\gamma)$ in $\pi_{1,r}(x)$.
\begin{conjecture}\label{conj:Gd(x;r,alpha,beta)}
Let $\alpha\equiv(-1)^\rho\bmod4$ be a positive integer with $\rho\in\{0,1\}$ and $\beta\in\{1,i\}$. For all large $x,$ we have
\begin{align*}
\cG_d(x;r,\alpha,\beta,\boldsymbol\gamma)&\sim\frac{h_{1,r}}{16}\cdot \Delta_G(r,\beta,\boldsymbol\gamma)\cdot G_d(r,\alpha,\beta,\boldsymbol\gamma)\cdot\frac{\sqrt{x}}{\log x},
\end{align*}
where $h_{1,r}\geqslant0$ is the constant defined by $\eqref{eq:constant-hD,r},$
\begin{align}\label{eq:DeltaG(r,beta,gamma)}
\Delta_G(r,\beta,\boldsymbol\gamma)=
\begin{cases}
1,\ \ & \text{for~}r\equiv4\gamma_1+2\bmod{32}\text{~and~}\beta=1,\\
1,\ \ & \text{for~}r\equiv-4\gamma_2\bmod{32}\text{~and~}\beta=i,\\
0,\ \ &\text{otherwise},
\end{cases}
\end{align}
and
\begin{align}\label{eq:Gd(r,alpha,beta,gamma)}
G_d(r,\alpha,\beta,\boldsymbol\gamma)
&=1+2\Re(i^d)(-1)^{\rho\gamma_2}\Big(\frac{i\beta}{\alpha}\Big)_4\Omega_4(1;\alpha,r)+(-1)^d\Omega_2(1;\alpha,r).
\end{align}
\end{conjecture}

\begin{remark}
In fact, Conjecture \ref{conj:Gd(x;r,alpha,beta)} deals with all possible $\beta\in\fU_1$ since
\begin{align*}
\cG_d(x;r,\alpha,\beta,\boldsymbol\gamma)=\cG_d(x;-r,\alpha,-\beta,\boldsymbol\gamma).
\end{align*}
\end{remark}

\begin{theorem}\label{thm:HL-Gd(x;r,alpha,beta)}
The Hardy--Littlewood conjecture \index{Hardy--Littlewood conjecture} implies Conjecture $\ref{conj:Gd(x;r,alpha,beta)}$.
\end{theorem}

We now present the proof of Theorem \ref{thm:HL-Gd(x;r,alpha,beta)}.
From the orthogonality, it follows that
\begin{align*}
\cG_d(x;r,\alpha,\beta,\boldsymbol\gamma)
&=\frac{1}{4}\sum_{0\leqslant j\leqslant 3}i^{-jd}\sum_{\substack{\pi\in\Pi_1(x)\\\Tr(\beta\pi)=r,~\pi\equiv2\gamma_1+1+2\gamma_2i\bmod{16}}}\Big(\frac{\alpha^j}{\pi}\Big)_4.
\end{align*}
Write $\alpha=(-1)^\rho\alpha_0$ so that $\alpha_0$ is primary in $\Z[i]$. From the quartic reciprocity law \index{quartic reciprocity law} (Lemma \ref{lm:quarticreciprocity}) it follows, for primary $\pi$, that
\begin{align*}
\Big(\frac{\alpha^j}{\pi}\Big)_4
&=\Big(\frac{\alpha_0^j}{\pi}\Big)_4\Big(\frac{-1}{\pi}\Big)_4^{-j\rho}
=\Big(\frac{\pi}{\alpha_0^j}\Big)_4(-1)^{\frac{N(\pi)-1}{4}\frac{N(\alpha_0)^j-1}{4}+\frac{(N(\pi)-1)j\rho}{4}}\\
&=\Big(\frac{\pi}{\alpha^j}\Big)_4(-1)^{\frac{N(\pi)-1}{4}(\frac{N(\alpha)^j-1}{4}+j\rho)}.
\end{align*}
This yields
\begin{align*}
\cG_d(x;r,\alpha,\beta,\boldsymbol\gamma)
&=\frac{1}{4}\sum_{0\leqslant j\leqslant 3}i^{-jd}\sum_{\substack{\pi\in\Pi_1(x)\\\Tr(\beta\pi)=r\\ \pi\equiv2\gamma_1+1+2\gamma_2i\bmod{16}}}\Big(\frac{\pi}{\alpha^j}\Big)_4(-1)^{\frac{N(\pi)-1}{4}(\frac{N(\alpha)^j-1}{4}+j\rho)}\\
&=\frac{1}{4}\sum_{0\leqslant j\leqslant 3}(-1)^{\gamma_2(\frac{N(\alpha)^j-1}{4}+j\rho)}i^{-jd}\cL_j,
\end{align*}
where
\begin{align*}
\cL_{j}:=\sum_{\delta\bmod{\alpha^j}}\Big(\frac{\delta}{\alpha^j}\Big)_4\sum_{\substack{\pi\in\Pi_1(x)\\\Tr(\beta\pi)=r,~\pi\equiv\delta\bmod{\alpha^j}\\ \pi\equiv2\gamma_1+1+2\gamma_2i\bmod{16}}}1.
\end{align*}

We would like to evaluate $\cL_j$ asymptotically assuming the (modified) Hardy--Littlewood conjecture.
To display the underlying idea, we now assume $\beta=1.$ For $\pi=m+ni$ with $\pi\equiv2\gamma_1+1+2\gamma_2i\bmod{16}$ and $\Tr(\pi)=r$, we find
\begin{align*}
2m=r,\ \ m\equiv2\gamma_1+1\bmod{16},\ \ n\equiv2\gamma_2\bmod{16}.
\end{align*}
Hence $\cL_j$ vanishes unless $r$ is even and $r/2\equiv2\gamma_1+1\bmod{16},$ in which case we find
\begin{align*}
\cL_j&=\sum_{h\bmod{\alpha^j}}\Big(\frac{\frac{r}{2}+hi}{\alpha^j}\Big)_4\sum_{\substack{n\in\Z\\ (r/2)^2+n^2\in\cP\cap[1,x]\\
n\equiv h\bmod{\alpha^j}\\ n\equiv 2\gamma_2\bmod{16}}}1.
\end{align*}
Following the convention in Conjecture \ref{conj:HL-modified}, the above $n$-sum is exactly $2\pi_{a,b,c}(x;q,u)$
with
\[a=1,~b=0,~c=(r/2)^2,~q=16\alpha^j\]
and a suitable $u\in\Z$ such that $u\equiv h\bmod{\alpha^j}$ and $u\equiv 2\gamma_2\bmod{16}$ by the Chinese remainder theorem\index{Chinese remainder theorem}. It is not difficult to verify all prerequisites in Conjecture \ref{conj:HL-modified}, which yields that the above $n$-sum should satisfy
\begin{align*}
\sum_n&\sim\frac{1}{4\varphi(\alpha^j)}\prod_{p\nmid 2\alpha^j r}\Big(1-\frac{(\frac{-1}{p})}{p-1}\Big)\cdot\frac{\sqrt{x}}{\log x},
\end{align*}
and thus implies an asymptotic formula for $\cL_j.$ 

Applying the above idea to general $\beta\in\fU_1,$ we may obtain
\begin{align}\label{eq:Lj-asymptotic}
\cL_j
\sim\frac{\cQ_\beta(\alpha^j,\frac{r}{2})}{4\varphi(\alpha^j)}\prod_{p\nmid 2\alpha^j r}\Big(1-\frac{(\frac{-1}{p})}{p-1}\Big)\cdot\frac{\sqrt{x}}{\log x},
\end{align}
where $\cQ_\beta(\alpha^j,\frac{r}{2})$ is an incomplete character sum defined by \eqref{eq:quarticcharactersum-incomplete},
and the necessary congruence conditions for $r$ should be $2\mid r$ and
\begin{align*}
r/2\equiv
\begin{cases}
\pm(2\gamma_1+1)\bmod{16},\ \ & \text{for~}\beta=\pm1,\\
\mp2\gamma_2\bmod{16},\ \ & \text{for~}\beta=\pm i.
\end{cases}
\end{align*}
Back to $\cG_d(x;r,\alpha,\beta,\boldsymbol\gamma),$ we derive that
\begin{align*}
\cG_d(x;r,\alpha,\beta,\boldsymbol\gamma)
&\sim\sum_{0\leqslant j\leqslant 3}(-1)^{\gamma_2(\frac{\alpha^{2j}-1}{4}+j\rho)}i^{-jd}\frac{\cQ_\beta(\alpha^j,\frac{r}{2})}{16\varphi(\alpha^j)}\prod_{p\nmid 2\alpha^j r}\Big(1-\frac{(\frac{-1}{p})}{p-1}\Big)\cdot\frac{\sqrt{x}}{\log x}.
\end{align*}

In order to conclude the asymptotic formula in Conjecture \ref{conj:Gd(x;r,alpha,beta)}, it remains to check
\begin{align}\label{eq:Gd(r,alpha,beta,gamma)-tobecheck}
\sum_{0\leqslant j\leqslant 3}(-1)^{\gamma_2(\frac{\alpha^{2j}-1}{4}+j\rho)}i^{-jd}\frac{\cQ_\beta(\alpha^j,\frac{r}{2})}{\varphi(\alpha^j)}\prod_{\substack{p\mid\alpha^j\\ p\nmid 2r}}\Big(1-\frac{(\frac{-1}{p})}{p-1}\Big)^{-1}=G_d(r,\alpha,\beta,\boldsymbol\gamma)
\end{align}
as given by \eqref{eq:Gd(r,alpha,beta,gamma)}. In fact, we infer from
Proposition \ref{prop:incompletecharactersum-quarticgauss} that
\begin{align*}
\frac{\cQ_\beta(\alpha^2,\frac{r}{2})}{\varphi(\alpha^2)}
&=\Big(\frac{i\beta}{\alpha}\Big)_4^2\cdot \Omega_2(1;\alpha,r)\cdot \prod_{\substack{p\mid \alpha\\ p\nmid r}}\Big(1-\frac{(\frac{-1}{p})}{p-1}\Big)\\
&=\Omega_2(1;\alpha,r)\cdot \prod_{\substack{p\mid \alpha\\ p\nmid r}}\Big(1-\frac{(\frac{-1}{p})}{p-1}\Big)
\end{align*}
and
\begin{align*}
\frac{\cQ_\beta(\alpha^3,\frac{r}{2})}{\varphi(\alpha^3)}
&=\Big(\frac{\beta}{\alpha}\Big)_4^2\frac{\cQ_\beta(\alpha,\frac{r}{2})}{\varphi(\alpha)}=\frac{\cQ_\beta(\alpha,\frac{r}{2})}{\varphi(\alpha)}\\
&=\Big(\frac{i\beta}{\alpha}\Big)_4\Omega_4(1;\alpha,r)\cdot \prod_{\substack{p\mid \alpha\\ p\nmid r}}\Big(1-\frac{(\frac{-1}{p})}{p-1}\Big),
\end{align*}
where $\Omega_2(\cdots)$ and $\Omega_4(\cdots)$ are defined by \eqref{eq:Omegafunction}. Inserting these into the LHS of \eqref{eq:Gd(r,alpha,beta,gamma)-tobecheck} we finish the checking and thus complete the proof of Theorem \ref{thm:HL-Gd(x;r,alpha,beta)}.

\section{From Conjecture \ref{conj:Gd(x;r,alpha,beta)} to Lang--Trotter}

Recall the connection \eqref{eq:connection:pi-Gd} between $\pi_{E,r}(x)$ and $\cG_d(x;r,\alpha,\beta,\boldsymbol\gamma)$:
\begin{align*}
2\pi_{E,r}(x)
&=\sum_{0\leqslant d\leqslant 3}\mathop{\sum\sum}_{\substack{0\leqslant \gamma_1,\gamma_2\leqslant7\\ \gamma_1\equiv\gamma_2\bmod2\\ \boldsymbol\gamma=(\gamma_1,\gamma_2)}}\cG_{-d+(2\delta+\lambda)\gamma_2}(x;r,g_1,i^d,\boldsymbol\gamma)+O(1),
\end{align*}
and define $\rho\in\{0,1\}$ via $g_1\equiv(-1)^\rho\bmod4$ in this section.

We first consider the case $2\parallel r$. Only $d=0,2$ contribute non-trivially. Thus
\begin{align*}
2\pi_{E,r}(x)
&=\mathop{\sum\sum}_{\substack{0\leqslant \gamma_1,\gamma_2\leqslant7\\ \gamma_1\equiv\gamma_2\bmod2\\ \boldsymbol\gamma=(\gamma_1,\gamma_2)}}\{\cG^{(1)}+\cG^{(2)}\}+O(1)
\end{align*}
with
\begin{align*}
\cG^{(1)}&=\cG_{(2\delta+\lambda)\gamma_2}(x;r,g_1,1,\boldsymbol\gamma),\\
\cG^{(2)}&=\cG_{2+(2\delta+\lambda)\gamma_2}(x;r,g_1,-1,\boldsymbol\gamma)
=\cG_{2+(2\delta+\lambda)\gamma_2}(x;-r,g_1,1,\boldsymbol\gamma).
\end{align*}

For $r\equiv4\gamma_1+2\bmod{32},$ Conjecture \ref{conj:Gd(x;r,alpha,beta)} yields
\begin{align*}
\cG^{(1)}&\sim\frac{h_{1,r}}{16}\Big(1+2\Re(i^{(2\delta+\lambda)\gamma_2})(-1)^{\rho\gamma_2}\Big(\frac{i}{g_1}\Big)_4\Omega_4+(-1)^{\lambda\gamma_2}\Omega_2\Big)\frac{\sqrt{x}}{\log x},
\end{align*}
and for $r\equiv-(4\gamma_1+2)\bmod{32},$ Conjecture \ref{conj:Gd(x;r,alpha,beta)}  yields
\begin{align*}
\cG^{(2)}&\sim\frac{h_{1,r}}{16}\Big(1-2\Re(i^{(2\delta+\lambda)\gamma_2})(-1)^{\rho\gamma_2}\Big(\frac{i}{g_1}\Big)_4\Omega_4+(-1)^{\lambda\gamma_2}\Omega_2\Big)\frac{\sqrt{x}}{\log x}.
\end{align*}
We then derive that
\begin{align*}
\pi_{E,r}(x)
&\sim\frac{h_{1,r}}{32}\sum_\pm\mathop{\sum\sum}_{\substack{0\leqslant \gamma_1,\gamma_2\leqslant7\\ \gamma_1\equiv\gamma_2\bmod2\\ \gamma_1\equiv \frac{\pm r-2}{4}\bmod8}}\Big(1\pm2\Re(i^{(2\delta+\lambda)\gamma_2})(-1)^{\rho\gamma_2}\Big(\frac{i}{g_1}\Big)_4\Omega_4\\
&\ \ \ \ \ \ \ \ \ \ \ \ \ \ \ \ \ \ \ \ \ \ \ \ \ \ \ \ \ \ \ \ \ \ \ \ \ \ \ \ \ +(-1)^{\lambda\gamma_2}\Omega_2\Big)\cdot \frac{\sqrt{x}}{\log x}\\
&=\frac{h_{1,r}}{8}\{2+((-1)^{\frac{\lambda(r-2)}{4}}+(-1)^{\frac{\lambda(r+2)}{4}})\Omega_2\}\cdot \frac{\sqrt{x}}{\log x}+\frac{c_{1,r}}{16}\Big(\frac{i}{g_1}\Big)_4\Omega_4\cdot \frac{\sqrt{x}}{\log x}\\
&\ \ \ \ \times\sum_{\pm}\Big(\mp(-1)^{\frac{\rho(r\pm2)}{4}}\sum_{\substack{0\leqslant\gamma_2\leqslant7\\ \gamma_2\equiv\frac{r\pm2}{4}\bmod2}}\Re(i^{(2\delta+\lambda)\gamma_2})\Big).
\end{align*}
Note that
\begin{align*}
\sum_{\substack{0\leqslant\gamma_2\leqslant7\\ \gamma_2\equiv\frac{r\pm2}{4}\bmod2}}\Re(i^{(2\delta+\lambda)\gamma_2})
&=\begin{cases}
4(-1)^{\frac{(r\pm2)(\delta+\frac{\lambda}{2})}{4}},\ \ &\text{for~}2\mid \lambda,\\
0,\ \ &\text{for~}2\nmid \lambda
\end{cases}
\end{align*}
for either choice of $\pm.$ We now conclude that
\begin{align*}
\pi_{E,r}(x)
&\sim\frac{h_{1,r}}{4}\Big\{1+\Omega_2+(-1)^{\frac{(\delta+\rho+\frac{\lambda}{2})(r-2)}{4}}(1-(-1)^{\delta+\rho+\frac{\lambda}{2}})\Big(\frac{i}{g_1}\Big)_4\Omega_4\Big\}\cdot \frac{\sqrt{x}}{\log x}
\end{align*}
if $2\mid \lambda$, and
\begin{align*}
\pi_{E,r}(x)
&\sim\frac{h_{1,r}}{4}\frac{\sqrt{x}}{\log x}
\end{align*}
otherwise. This proves Theorem \ref{thm:HL-LT:D=1} in the case of $2\parallel r$ by noting that
$$\Big(\frac{i}{g_1}\Big)_4=(-1)^{\frac{g_1^2-1}{8}}, \ \ (-1)^\rho=(-1)^{\frac{g_1-1}{2}}.$$

We now turn to the case $4\mid  r$. Only $d=1,3$ contribute non-trivially. Thus
\begin{align*}
2\pi_{E,r}(x)
&=\mathop{\sum\sum}_{\substack{0\leqslant \gamma_1,\gamma_2\leqslant7\\ \gamma_1\equiv\gamma_2\bmod2\\ \boldsymbol\gamma=(\gamma_1,\gamma_2)}}\{\cG^{(3)}+\cG^{(4)}\}+O(1)
\end{align*}
with
\begin{align*}
\cG^{(3)}&=\cG_{3+(2\delta+\lambda)\gamma_2}(x;r,g_1,i,\boldsymbol\gamma),\\
\cG^{(4)}&=\cG_{1+(2\delta+\lambda)\gamma_2}(x;r,g_1,-i,\boldsymbol\gamma)
=\cG_{1+(2\delta+\lambda)\gamma_2}(x;-r,g_1,i,\boldsymbol\gamma).
\end{align*}

For $r\equiv-4\gamma_2\bmod{32},$ Conjecture \ref{conj:Gd(x;r,alpha,beta)} yields
\begin{align*}
\cG^{(3)}&\sim\frac{h_{1,r}}{16}\Big(1-2\Re(i^{1+(2\delta+\lambda)\gamma_2})(-1)^{\rho\gamma_2}\Omega_4-(-1)^{\lambda\gamma_2}\Omega_2\Big)\frac{\sqrt{x}}{\log x},
\end{align*}
and for $r\equiv4\gamma_2\bmod{32},$ Conjecture \ref{conj:Gd(x;r,alpha,beta)} yields
\begin{align*}
\cG^{(4)}&\sim\frac{h_{1,r}}{16}\Big(1+2\Re(i^{1+(2\delta+\lambda)\gamma_2})(-1)^{\rho\gamma_2}\Omega_4-(-1)^{\lambda\gamma_2}\Omega_2\Big)\frac{\sqrt{x}}{\log x}.
\end{align*}
We then derive that
\begin{align*}
\pi_{E,r}(x)
&\sim\frac{h_{1,r}}{32}\sum_\pm\mathop{\sum\sum}_{\substack{0\leqslant \gamma_1,\gamma_2\leqslant7\\ \gamma_1\equiv\gamma_2\bmod2\\ \gamma_2\equiv \pm\frac{r}{4}\bmod8}}\Big(1\pm2\Re(i^{1+(2\delta+\lambda)\gamma_2})(-1)^{\rho\gamma_2}\Omega_4-(-1)^{\lambda\gamma_2}\Omega_2\Big)\frac{\sqrt{x}}{\log x}\\
&=\frac{h_{1,r}}{4}\{1+(-1)^{\frac{r}{4}(\delta+\frac{g_1-1}{2})}(1-(-1)^{\frac{\lambda r}{4}})\Re(i^{1+\frac{\lambda r}{4}})\Omega_4-(-1)^{\frac{\lambda r}{4}}\Omega_2\}\cdot \frac{\sqrt{x}}{\log x}.
\end{align*}
This proves Theorem \ref{thm:HL-LT:D=1} in the case of $4\mid r.$

%% file: 5.D=3.tex
\chapter{Complex multiplication by $\Q(\sqrt{-3})$}\label{ch:D=3}

Throughout this chapter, we assume $E/\Q$ is an elliptic curve with CM by $\Q(\sqrt{-3})$, defined by
\begin{align}\label{eq:E:D=3}
y^2=x^3+g
\end{align}
with $g\in\Z^*.$ As in Chapter \ref{ch:D=1}, the asymptotic formula
\begin{align}\label{eq:piE,r(x)-varpiE,r,D=3}
\pi_{E,r}(x)\sim\frac{\varpi_{E,r}\cdot\sqrt{x}}{\log x}
\end{align}
will also be proven, under the Hardy--Littlewood conjecture\index{Hardy--Littlewood conjecture}, with some explicit constant $\varpi_{E,r}\geqslant0.$ To set up, we write
\begin{align*}
g=(-1)^\delta2^\lambda3^\mu g_1
\end{align*}
with
\begin{itemize}
\item $\delta\in\{0,1\},$
\item $\lambda,\mu\in\N,$
\item $g_1\in\Z^+$ and $(6,g_1)=1$.
\end{itemize}

Recall the definition \eqref{eq:Omegafunction}:
\begin{align*}
\Omega_j(D;q,r)=\prod_{\substack{p^\nu\parallel q,~ p\nmid r\\j\nmid \nu}}\frac{-1}{p-1-(\frac{-D}{p})}.
\end{align*}
We adopt the convention $\Omega_j=\Omega_j(3;g_1,r)$ for $j\in\{2,3,6\}$ throughout this chapter.

Recall that $\omega=\ue^{2\pi i/3}=(-1+\sqrt{-3})/2$ is a primitive cubic root of unity. Denote by $\fU_3=\{\pm1,\pm \omega,\pm\omega^2\}$, the group of roots of unity in $\Z[\omega]$.

\section{Asymptotics for $\pi_{E,r}(x)$: conditional results}

\begin{theorem}\label{thm:HL-LT:D=3}\index{Lang--Trotter conjecture}
Suppose $E/\Q$ is the elliptic curve defined by $\eqref{eq:E:D=3}$ with the above convention.
For each non-zero integer $r$, the Hardy--Littlewood conjecture \index{Hardy--Littlewood conjecture} implies that $\eqref{eq:piE,r(x)-varpiE,r,D=3}$ holds with
\begin{align}\label{eq:varpiE,r,D=3}
\varpi_{E,r}=\frac{h_{3,r}}{6}\Big(1+\varsigma_2(g,r)\Big(\frac{3}{g_1}\Big)\Omega_2+\widehat{\mu}\varsigma_1(g,r)\Big(\Omega_3+\varsigma_2(g,r)\Big(\frac{3}{g_1}\Big)\Omega_6\Big)\Big)
\end{align}
for $3\nmid r,$ where $h_{3,r}\geqslant0$ is the constant defined by $\eqref{eq:constant-hD,r},$ $\widehat{\mu}$ vanishes unless $3\mid \mu,$ in which case it is $\frac{2}{3},$ and
\begin{align*}
\varsigma_1(g,r)&=\begin{cases}
1+2\Re(\omega^{\lambda+\frac{g_1^2-1}{3}}),\ \ &\text{for~} r\equiv\pm2\bmod6,\\ 
\Re(\omega^{1-\lambda}+\omega^{1+\frac{g_1^2-1}{3}}+\omega^{2-\lambda+\frac{g_1^2-1}{3}}),\ \ &\text{for~} r\equiv\pm5\bmod{6},\\
\end{cases}
\end{align*}
and
\begin{align*}
\varsigma_2(g,r)&=\begin{cases}
\pm(-1)^{\delta+\lambda+\mu+\frac{g_1-1}{2}},\ \ &\text{for~} r\equiv\pm8\bmod{24},\\ 
\pm(-1)^{\delta+\mu+\frac{g_1-1}{2}},\ \ &\text{for~} r\equiv\pm20\bmod{24},\\ 
\pm\frac{1+(-1)^{\lambda}}{2},\ \ &\text{for~} r\equiv\pm2\bmod{12},\\ 
\pm\frac{(1+(-1)^\lambda)(1+(-1)^{\delta+\mu+\frac{g_1-1}{2}})}{4},\ \ &\text{for~} r\equiv\pm5\bmod{6}.
\end{cases}
\end{align*}
Unconditionally, $\pi_{E,r}(x)$ is bounded by $O(1)$ if $3\mid r$.
\end{theorem}

The proof of Theorem \ref{thm:HL-LT:D=3} will be given in Section \ref{sec:proofofTheorem{thm:HL-LT:D=3}}. The full strength of the Hardy--Littlewood conjecture is not needed here, and its comparative power in Theorem \ref{thm:HL-LT:D=3} will be reflected by Conjecture \ref{conj:Ed(x;r,alpha,beta,gamma)}, which can be regarded as a special case of the original Hardy--Littlewood conjecture.

The constants $\varsigma_1(g,r)$ and $\varsigma_2(g,r)$ look a bit complicated at first sight, but they should be very easy to compute when $g$ and $r$ are given. Note that $\varsigma_1(g,r)$ depends on the exponents modulo 3, and $\varsigma_2(g,r)$ then depends on those modulo 2. We should also remark that the constant $\varpi_{E,r}$ is not always $\frac{h_{3,r}}{6}$, which was considered as a natural guess and imagined to be true by Mazur \cite{Ma72} in the case $r=1$.
Based on his pioneer work on the Mazur conjecture, Qin \cite{Qi16}, also assuming the Hardy--Littlewood conjecture, gave a negative answer to this observation of Mazur, and gave the complete evaluation of $\varpi_{E,r}$ for $r=\pm1.$

For some special curve $E$ and $r$, the constants in Theorem \ref{thm:HL-LT:D=3} can be significantly simplified.

\begin{corollary}\label{coro:HL-LT:D=3,g=-432}
Suppose $E/\Q$ is the elliptic curve defined by 
\begin{align*}
y^2=x^3-432.
\end{align*}
For each non-zero integer $r$, the Hardy--Littlewood conjecture implies 
\begin{align*}
\pi_{E,r}(x)\sim \frac{h_{3,r}}{3}\cdot\varsigma_r\cdot\frac{\sqrt{x}}{\log x},
\end{align*}
where $h_{3,r}\geqslant0$ is a constant defined by $\eqref{eq:constant-hD,r},$ and
\begin{align*}
\varsigma_r&=\begin{cases}
1,\ \ &\text{for~} r\equiv2\bmod3,\\ 
0,\ \ &\text{for~} r\equiv1\bmod3.
\end{cases}
\end{align*}
\end{corollary}

The curve treated in Corollary \ref{coro:HL-LT:D=3,g=-432} was originally considered by Gau{\ss} in his last diary entry. See Chapter \ref{ch:Historicalnotesperspectives} for more details. The following consequence gives a precise characterization of the Mazur conjecture \index{Mazur's conjecture} for anomalous primes \index{anomalous prime} (see \cite[Page 187]{Ma72}), and we assume $g$ is squarefree for simplicity. As one may see from below, the constant in the ``asymptotic formula'' can be zero in some cases, which yields there are at most finitely many anomalous primes. We will come back to this issue in Section \ref{sec:anomalousprimes} for any $g\in\Z^*.$

\begin{corollary}\label{coro:HL-LT:D=3,gsquarefree}
Suppose $E/\Q$ is the elliptic curve defined by 
\begin{align*}
y^2=x^3+g,
\end{align*}
where $g$ is a positive squarefree number with $2\nmid g$ and $3\mid g.$
The Hardy--Littlewood conjecture implies 
\begin{align*}
\pi_{E,\pm1}(x)\sim \kappa\Big(1\pm\Omega(g)\Big(\frac{g_1}{3}\Big)\varsigma(g)\Big)\cdot\frac{\sqrt{x}}{\log x},
\end{align*}
where 
\begin{align*}
\kappa&=\frac{\sqrt{3}}{6}\prod_{p>2}\Big(1-\frac{(\frac{-3}{p})}{p-1}\Big),\ \ \ \Omega(g)=\prod_{\substack{p\mid g\\ p>3}}\frac{-1}{p-1-(\frac{-3}{p})},
\end{align*}
and $\varsigma(g)$ vanishes unless $g\equiv1\bmod4,$ in which case it is $1.$
\end{corollary}

\section{Initial transformations}
Similar to Lemma \ref{lm:Frobenius:D=1}, we have the following explicit formula for Frobenius traces \index{trace of Frobenius} $a_p(E)$ in the case of CM by $\Q(\sqrt{-3})$. See \cite[Chapter 18, Theorem 4]{IR90} for details.
\begin{lemma}\label{lm:Frobenius:D=3}
Suppose $p\in\Z^+$ with $p\nmid 6g.$ Then $a_p(E)=0$ if $p\equiv2\bmod3,$ and
\begin{align*}
a_p(E)=-\Big(\frac{4g}{\pi}\Big)_6\overline{\pi}-\Big(\frac{4g}{\overline{\pi}}\Big)_6\pi=-\Tr\Big(\Big(\frac{4g}{\overline{\pi}}\Big)_6\pi\Big)
\end{align*}
if $p\equiv1\bmod3$ with $p=\pi\overline{\pi}$ such that $\pi\in\Z[\omega]$ is primary.
\end{lemma}

By Lemma \ref{lm:Frobenius:D=3}, we decompose $\pi_{E,r}(x)$ as
\begin{align*}
2\pi_{E,r}(x)
&=\sum_{\beta\in \fU_3}\sum_{\substack{\pi\in\Pi_3(x)\\ \Tr(\beta\pi)=-r,~\overline{(\frac{4g}{\pi})_6}=\beta\\ \pi\text{~primary}}}1+O(1).
\end{align*}
Note that
\begin{align*}
\overline{\Big(\frac{4g}{\pi}\Big)_6}=\Big(\frac{4g}{\pi}\Big)_6^5=\Big(\frac{4g}{\pi}\Big)_6^3\Big(\frac{4g}{\pi}\Big)_6^2=\Big(\frac{4g}{N(\pi)}\Big)\Big(\frac{4g}{\pi}\Big)_3
\end{align*}
thanks to \cite[Chapter 18, Section 7, Lemma 2]{IR90}. 
Therefore, for any primary $\pi=a+b\omega\in\Pi_3(x)$ with $N(\pi)\geqslant5,$ we have
\begin{align*}
\overline{\Big(\frac{4g}{\pi}\Big)_6}
&=\Big(\frac{(-1)^\delta2^\lambda3^\mu g_1}{N(\pi)}\Big)\Big(\frac{2}{\pi}\Big)_3^{2+\lambda}\Big(\frac{3}{\pi}\Big)_3^{\mu}\Big(\frac{g_1}{\pi}\Big)_3\\
&=\Big(\frac{(-1)^\delta2^\lambda3^\mu g_1}{N(\pi)}\Big)\Big(\frac{g_1}{\pi}\Big)_3\omega^{\varepsilon(a,b)}
\end{align*}
with
$\varepsilon(a,b)=(2+\lambda)(2+(-1)^\frac{a-2}{3}-(-1)^{\frac{b}{3}})+\frac{2\mu b}{3}.$
We then find
\begin{align*}
2\pi_{E,r}(x)
&=\sum_{\pm}\sum_{\beta\in \fU_3}\sum_{\substack{\pi=a+b\omega\in\Pi_3(x)\\ \Tr(\beta\pi)=-r,~(\frac{(-1)^\delta2^\lambda3^\mu g_1}{N(\pi)})=\pm1\\
\omega^{\varepsilon(a,b)}(\frac{g_1}{\pi})_3=\pm\beta\\ \pi\text{~primary}}}1+O(1).
\end{align*}

Note that the condition $(\frac{(-1)^\delta2^\lambda3^\mu}{N(\pi)})=\pm1$ is equivalent to
\begin{align*}
(-1)^{\frac{\lambda(N(\pi)^2-1)}{8}+\frac{(\delta+\mu)(N(\pi)-1)}{2}}=\pm1.
\end{align*}
Therefore,
\begin{align}\label{eq:initialtransformation:D=3}
2\pi_{E,r}(x)
&=\sum_{0\leqslant a\leqslant 1}\sum_{0\leqslant b\leqslant 1}\sum_{0\leqslant k\leqslant 3}\sum_{\pm}\sum_{\beta\in \fU_3}~\sideset{}{^\ddagger}\sum_{\pi=a+b\omega\in\Pi_3(x)}1+O(1),
\end{align}
where $\ddagger$ yields the following five simultaneous constraints:
\begin{align*}
\Tr(\beta\pi)=-r,\ \ \Big(\frac{g_1}{N(\pi)}\Big)=\pm(-1)^{\frac{\lambda k(k+1)}{2}+(\delta+\mu)k},\\
\omega^{(2+\lambda)(2+(-1)^\frac{a-2}{3}-(-1)^{\frac{b}{3}})+\frac{2\mu b}{3}}\Big(\frac{g_1}{\pi}\Big)_3=\pm\beta,\quad\quad\\
\pi\equiv -1\bmod3,\ \ \ N(\pi)\equiv2k+1\bmod8.\quad\quad
\end{align*}
To characterize the exponent $(2+\lambda)(2+(-1)^\frac{a-2}{3}-(-1)^{\frac{b}{3}})+\frac{2\mu b}{3},$ we need to consider the residue classes $a\bmod{6}$ and $b\bmod{18}.$

The above arguments lead to the following counting function
\begin{align*}
\cE_{d,k,\epsilon}(x;r,\alpha,\beta,\gamma,\eta)&=\Bigg|\Bigg\{\pi\in\Pi_3(x):\begin{array}{c}
\Tr(\beta\pi)=r,~(\frac{\alpha}{\pi})_3=\omega^d,~(\frac{\alpha}{N(\pi)})=\epsilon,\\
N(\pi)\equiv2k+1\bmod8,\\
\pi\equiv\gamma\bmod2,~\pi\equiv2+3\eta\bmod9\end{array}\Bigg\}\Bigg|
\end{align*}
for
\begin{align*}
r,\alpha\in\Z^*,\ \ \beta\in \fU_3,\ \ 0\leqslant d\leqslant2,\ \ 0\leqslant k\leqslant3
\end{align*}
and
\begin{align*}
\epsilon\in\{1,-1\},\ \ \gamma\in(\Z[\omega]/(2))^\times,\ \ \eta\in\Z[\omega]/(3).
\end{align*}
It is clear that
\begin{align*}
\sum_{\epsilon\in\{1,-1\}}\mathop{\sum\sum}_{0\leqslant d,l\leqslant2}\sum_{0\leqslant k\leqslant3}\sum_{\beta\in\fU_3}\sum_{\gamma\in(\Z[\omega]/(2))^\times}&\sum_{\eta\in\Z[\omega]/(3)}\cE_{d,k,\epsilon}(x;r,\alpha,\beta,\gamma,\eta)\\
&=2\pi_{3,r}(x)+O(1),
\end{align*}
where $\pi_{3,r}(x)$ is defined by \eqref{eq:piD,r(x)-definition}. The next section aims to predict the asymptotic proportion of $\cE_{d,k,\epsilon}(x;r,\alpha,\beta,\gamma,\eta)$
in $\pi_{3,r}(x)$, which in turn implies a conjectural characterization of $\pi_{E,r}(x)$ as stated in Theorem \ref{thm:HL-LT:D=3}.

\section{Asymptotic evaluation of $\cE_{d,k,\epsilon}(x;r,\alpha,\beta,\gamma,\eta)$}
In view of 
\begin{align*}\cE_{d,k,\epsilon}(x;r,\alpha,\beta,\gamma,\eta)=\cE_{d,k,\epsilon}(x;-r,\alpha,-\beta,\gamma,\eta),
\end{align*}
we only focus on the case $\beta\in\{1,\omega,\omega^2\}$.

\begin{conjecture}\label{conj:Ed(x;r,alpha,beta,gamma)}
Let $\alpha$ be a positive integer with $(\alpha,6)=1$ and $r$ a fixed non-zero integer. 
Let $\beta\in\{1,\omega,\omega^2\}$. For all large $x,$ we have
\begin{align*}
\cE_{d,k,\epsilon}(x;r,\alpha,\beta,\gamma,\eta)
&\sim \frac{h_{3,r}}{216}\cdot \Delta(r,k,\gamma,\beta,\eta)\cdot E_{d,k,\epsilon}(r,\alpha,\beta)\cdot\frac{\sqrt{x}}{\log x},
\end{align*}
where $h_{3,r}\geqslant0$ is the constant defined by $\eqref{eq:constant-hD,r},$ and
\begin{align}\label{eq:DeltaE(r,beta)1}
\Delta(r,k,\gamma,\beta,\eta)=
\begin{cases}
8,\ \ &\text{for~} r\equiv r_1(\beta,\eta)\bmod{72},k\equiv1\bmod4\text{~and~}\gamma\equiv\beta^2\bmod2,\\
8,\ \ &\text{for~} r\equiv r_2(\beta,\eta)\bmod{72},k\equiv3\bmod4\text{~and~}\gamma\equiv\beta^2\bmod2,\\
4,\ \ &\text{for~} r\equiv r_3(\beta,\eta)\bmod{36},2\mid k\text{~and~}\gamma\equiv\beta^2\bmod2,\\
1,\ \ &\text{for~} r\equiv r_4(\beta,\eta)\bmod{18}\text{~and~}\gamma\equiv\beta^2\omega,\beta^2(1+\omega)\bmod2,\\
0, &\text{otherwise}
\end{cases}
\end{align}
with
\begin{align*}
r_j(1,\eta)=
\begin{cases}
40-24\Tr(\beta\eta),\ \ &\text{for~} j=1,\\
4-24\Tr(\beta\eta),\ \ &\text{for~} j=2,\\
22-24\Tr(\beta\eta),\ \ &\text{for~} j=3,\\
13-24\Tr(\beta\eta),\ \ &\text{for~} j=4,
\end{cases}
\end{align*}
\begin{align*}
r_j(\omega,\eta)=r_j(\omega^2,\eta)=
\begin{cases}
16-24\Tr(\beta\eta),\ \ &\text{for~} j=1,\\
52-24\Tr(\beta\eta),\ \ &\text{for~} j=2,\\
34-24\Tr(\beta\eta),\ \ &\text{for~} j=3,\\
7-24\Tr(\beta\eta),\ \ &\text{for~} j=4,
\end{cases}
\end{align*}
and
\begin{equation}\label{eq:Ed,k,epsilon(r,alpha,beta)}
\begin{split}
E_{d,k,\epsilon}(r,\alpha,\beta)&=1+\epsilon(-1)^{\frac{k(\alpha-1)}{2}}\Big(\frac{3}{\alpha}\Big)\Omega_2(3;\alpha,r)\\
&\ \ \ \ +2\Re\Big(\omega^d\Big(\frac{\beta}{\alpha}\Big)_3\Big)\Big\{\Omega_3(3;\alpha,r)+\epsilon(-1)^{\frac{k(\alpha-1)}{2}}\Big(\frac{3}{\alpha}\Big)\Omega_6(3;\alpha,r)\Big\}.
\end{split}
\end{equation}
\end{conjecture}

Note that for $\eta=\eta_1+\eta_2\omega$ with $\eta_1,\eta_2\in\Z,$ we have
\begin{align*}
\Tr(\eta)=2\eta_1-\eta_2,\ \ \Tr(\beta\eta)=-(\eta_1+\eta_2),\ \ \Tr(\omega^2\eta)=2\eta_2-\eta_1.
\end{align*}

\begin{theorem}\label{thm:HL-Ed(x;r,alpha,beta,gamma)}
The Hardy--Littlewood conjecture implies Conjecture $\ref{conj:Ed(x;r,alpha,beta,gamma)}$.
\end{theorem}

We now give the proof of Theorem \ref{thm:HL-Ed(x;r,alpha,beta,gamma)}.
\subsection{Orthogonality and reciprocity}
It follows from orthogonality that
\begin{align*}
\cE_{d,k,\epsilon}(x;r,\alpha,\beta,\gamma,\eta)
&=\frac{1}{3}\sum_{0\leqslant j\leqslant 2}\omega^{-jd}\sum_{\substack{\pi\in\Pi_3(x)\\\ \Tr(\beta\pi)=r,~N(\pi)\equiv2k+1\bmod8\\
~\pi\equiv2+3\eta\bmod9,~\pi\equiv\gamma\bmod2\\(\frac{\alpha}{N(\pi)})=\epsilon}}\Big(\frac{\alpha}{\pi}\Big)_3^j.
\end{align*}

Define $\rho\in\{0,1\}$ such that $\alpha_0:=(-1)^\rho\alpha$ is primary in $\Z[\omega].$ 
A basic starting point is to transform $(\frac{\alpha}{\pi})_3$ using the cubic reciprocity law\index{cubic reciprocity law}. More precisely, for each primary $\pi\in\Z[\omega]$ coprime to $\alpha$ with $N(\pi)\nmid 6\alpha,$ we derive from Lemma \ref{lm:cubicreciprocity} that
\begin{align*}
\Big(\frac{\alpha}{\pi}\Big)_3
&=\Big(\frac{\alpha_0}{\pi}\Big)_3\Big(\frac{-1}{\pi}\Big)_3^{-\rho}
=\Big(\frac{\pi}{\alpha_0}\Big)_3=\Big(\frac{\pi}{\alpha}\Big)_3
\end{align*}
since $(\frac{-1}{\pi})_3=1$ by Lemma \ref{lm:cubicsymbol-basicproperty}.
It then follows that
\begin{align*}
\cE_{d,k,\epsilon}(x;r,\alpha,\beta,\gamma,\eta)
&=\frac{1}{3}\sum_{0\leqslant j\leqslant 2}\omega^{-dj}\sum_{\substack{\pi\in\Pi_3(x)\\ \Tr(\beta\pi)=r,~N(\pi)\equiv2k+1\bmod8\\
~\pi\equiv2+3\eta\bmod9,~\pi\equiv\gamma\bmod2\\ (\frac{\alpha}{N(\pi)})=\epsilon}}\Big(\frac{\pi}{\alpha}\Big)_3^j+O(1).
\end{align*}
The quadratic reciprocity law \index{quadratic reciprocity law} (for Jacobi symbols) yields
\begin{align*}
\Big(\frac{\alpha}{N(\pi)}\Big)=(-1)^{\frac{(\alpha-1)(N(\pi)-1)}{4}}\Big(\frac{N(\pi)}{\alpha}\Big)=(-1)^{\frac{k(\alpha-1)}{2}}\Big(\frac{N(\pi)}{\alpha}\Big)
\end{align*}
for all odd $\alpha,N(\pi)\geqslant3$ with $N(\pi)\equiv2k+1\bmod8$,
so that
\begin{align}\label{eq:Ed-decomposition-Rj}
\cE_{d,k,\epsilon}(x;r,\alpha,\beta,\gamma,\eta)
&=\frac{1}{3}\sum_{0\leqslant j\leqslant 2}\omega^{-dj}\cR_j+O(1),
\end{align}
where
\begin{align*}
\cR_j
&=\sum_{\substack{\theta\bmod \alpha\\ (\frac{\theta}{\alpha})=\epsilon(-1)^{\frac{k(\alpha-1)}{2}}}}\sum_{\substack{\pi\in\Pi_3(x)\\\Tr(\beta\pi)=r,~\pi\equiv\gamma_1+\gamma_2\omega\bmod2\\
N(\pi)\equiv\theta\bmod \alpha,~N(\pi)\equiv2k+1\bmod8\\ \pi\equiv2+3\eta\bmod9}}\Big(\frac{\pi}{\alpha}\Big)_3^j.
\end{align*}
In what follows, we write $\gamma=\gamma_1+\gamma_2\omega$ and $\eta=\eta_1+\eta_2\omega$ with $\gamma_1,\gamma_2,\eta_1,\eta_2\in\Z$. We would like to evaluate $\cR_j$ as consequences of the Hardy--Littlewood conjecture.

\subsection{Consequences of Hardy--Littlewood}
To show underlying ideas, we only give details for $\beta=1,$ which we henceforth assume. 
Note that $\cR_j$ is bounded by $O(1)$ unless
\begin{align}\label{eq:D=3restrictions:beta=1}
r\equiv\gamma_2\bmod2,\ \ r\equiv 4+3(2\eta_1-\eta_2)\bmod9,
\end{align}
in which case we further have
\begin{align*}
\cR_0
&=\sum_{\substack{\theta\bmod\alpha\\ (\frac{\theta}{\alpha})=\epsilon(-1)^{\frac{k(\alpha-1)}{2}}}}\sum_{\substack{m\in\Z\\ 3m^2-3mr+r^2\in\cP\cap[1,x]\\\ m\equiv\gamma_1\bmod2,~m\equiv 2+3\eta_1\bmod9\\3m^2-3mr+r^2\equiv \theta\bmod\alpha\\  3m^2-3mr+r^2\equiv 2k+1\bmod8}}1\\
&=\sum_{t\in\cT(\gamma_1,r,k)}\sum_{\substack{a\bmod\alpha\\ (\frac{3a^2-3ar+r^2}{\alpha})=\epsilon(-1)^{\frac{k(\alpha-1)}{2}}}}\sum_{\substack{m\in\Z\\ 3m^2-3mr+r^2\in\cP\cap[1,x]\\\ m\equiv\gamma_1+2t\bmod8,~m\equiv 2+3\eta_1\bmod9\\
m\equiv a\bmod{\alpha}}}1,
\end{align*}
where
\begin{align}\label{eq:T-countingt}
\cT(\gamma_1,r,k)
&:=\Big\{t\bmod4:\begin{array}{c}
4t^2+(4\gamma_1+2r) t+3\gamma_1^2-3r\gamma_1+r^2\\
\equiv2k+1\bmod 8\end{array}\Big\}
\end{align}
Following the convention in Conjecture \ref{conj:HL-modified}, the last $m$-sum is exactly 
$$\pi_{3,-3r,r^2}(x;q,u)+\pi_{3,3r,r^2}(x;q,-u)$$
with $q=72\alpha$ and some $u\in\Z$ such that $u\equiv\gamma_1+2t\bmod8,$ $u\equiv 2+3\eta_1\bmod9$ and $u\equiv a\bmod{\alpha}$ by the Chinese remainder theorem. \index{Chinese remainder theorem} It is not difficult to verify all prerequisites in Conjecture \ref{conj:HL-modified}, which yields that the above $m$-sum should satisfy
\begin{align}\label{eq:D=3:m-sum}
\sum_m&\sim\frac{1}{12\sqrt{3}\varphi(\alpha)}\prod_{p\nmid 6\alpha r}\Big(1-\frac{(\frac{-3}{p})}{p-1}\Big)\cdot\frac{\sqrt{x}}{\log x},
\end{align}
and thus implies
\begin{align*}
\cR_0
&\sim\frac{|\cT(\gamma_1,r,k)|}{12\sqrt{3}\varphi(\alpha)}\prod_{p\nmid 6\alpha r}\Big(1-\frac{(\frac{-3}{p})}{p-1}\Big)\cdot\frac{\sqrt{x}}{\log x}\cdot\sum_{\substack{a\bmod\alpha\\ (\frac{3a^2-3ar+r^2}{\alpha})=\epsilon(-1)^{\frac{k(\alpha-1)}{2}}}}1.
\end{align*}
By Lemma \ref{lm:countingresidues}, the last sum over $a$ is exactly
\begin{align*}
\frac{1}{2}\varphi(\alpha)\prod_{\substack{p\mid\alpha\\ p\nmid r}}\Big(1-\frac{(\frac{-3}{p})}{p-1}\Big)\Big\{1+\epsilon(-1)^{\frac{k(\alpha-1)}{2}}\Big(\frac{3}{\alpha}\Big)\Omega_2(3;\alpha,r)\Big\},
\end{align*}
which yields
\begin{align*}
\cR_0
&\sim\frac{|\cT(\gamma_1,r,k)|}{24\sqrt{3}}\prod_{p\nmid 2r}\Big(1-\frac{(\frac{-3}{p})}{p-1}\Big)\Big\{1+\epsilon(-1)^{\frac{k(\alpha-1)}{2}}\Big(\frac{3}{\alpha}\Big)\Omega_2(3;\alpha,r)\Big\}\cdot\frac{\sqrt{x}}{\log x}.
\end{align*}

Observe that $\cR_1=\overline{\cR_2}$, thus it suffices to deal with $\cR_1$. Following similar ideas as above, we find
\begin{align*}
\cR_1
&=\sum_{t\in\cT(\gamma_1,r,k)}\sum_{\substack{a\bmod\alpha\\ (\frac{3a^2-3ar+r^2}{\alpha})=\epsilon(-1)^{\frac{k(\alpha-1)}{2}}}}\Big(\frac{a+(2a-r)\omega}{\alpha}\Big)_3\\
&\ \ \ \ \times\sum_{\substack{m\in\Z\\ 3m^2-3mr+r^2\in\cP\cap[1,x]\\\ m\equiv\gamma_1+2t\bmod8,~m\equiv 2+3l\bmod9\\
m\equiv a\bmod\alpha}}1.
\end{align*}
The innermost sum over $m$ can be evaluated asymptotically as in \eqref{eq:D=3:m-sum}, and one important feature is that the main term does not depend on $a\bmod\alpha$. This gives
\begin{equation}\label{eq:R1-incompletecharactersum}
\begin{split}\cR_1
&\sim \frac{|\cT(\gamma_1,r,k)|}{12\sqrt{3}\varphi(\alpha)}\prod_{p\nmid 6\alpha r}\Big(1-\frac{(\frac{-3}{p})}{p-1}\Big)\cdot\frac{\sqrt{x}}{\log x}\\
&\ \ \ \ \times\sum_{\substack{a\bmod\alpha\\ (\frac{3a^2-3ar+r^2}{\alpha})=\epsilon(-1)^{\frac{k(\alpha-1)}{2}}}}\Big(\frac{a+(2a-r)\omega}{\alpha}\Big)_3.
\end{split}
\end{equation}
To evaluate the above sum over $a$, we appeal to Proposition \ref{prop:incompletecharactersum-cubicgauss}. In fact,
for each choice of $\pm,$ we have
\begin{align*}
\Sigma_{\pm}:&=\sum_{\substack{a\bmod{\alpha}\\ (\frac{3a^2-3ar+r^2}{\alpha})=\pm1}}\Big(\frac{a+(2a-r)\omega}{\alpha}\Big)_3\\
&=\frac{1}{2}\sum_{a\bmod{\alpha}}\Big(\frac{a+(2a-r)\omega}{\alpha}\Big)_3\Big(1\pm\Big(\frac{3a^2-3ar+r^2}{\alpha}\Big)\Big)\\
&=\frac{\cC_1(\alpha,r;0)\pm\cC_1(\alpha,r;1)}{2},
\end{align*}
where $\cC_1(\alpha,r;*)$ is defined as in \eqref{eq:cubiccharactersum-incomplete}.
Basic properties and evaluations of $\cC_1(\alpha,r;*)$ can be found in Proposition \ref{prop:incompletecharactersum-cubicgauss} and Corollary \ref{coro:incompletecharactersum-cubicgauss}. Consequently, we have
\begin{align*}
\Sigma_\pm
&=\frac{\varphi(\alpha)}{2}\prod_{\substack{p\mid \alpha\\ p\nmid r}}\Big(1-\frac{(\frac{-3}{p})}{p-1}\Big)\Big\{\Omega_3(3;\alpha,r)\pm\Big(\frac{3}{\alpha}\Big)\Omega_6(3;\alpha,r)\Big\},
\end{align*}
from which and \eqref{eq:R1-incompletecharactersum} we derive that
\begin{align*}
\cR_1
&\sim\frac{|\cT(\gamma_1,r,k)|}{24\sqrt{3}}\prod_{p\nmid6r}\Big(1-\frac{(\frac{-3}{p})}{p-1}\Big)\cdot\frac{\sqrt{x}}{\log x}\\
&\ \ \ \ \times \Big\{\Omega_3(3;\alpha,r)+\epsilon(-1)^{\frac{k(\alpha-1)}{2}}\Big(\frac{3}{\alpha}\Big)\Omega_6(3;\alpha,r)\Big\}.
\end{align*}

Collecting the above evaluations for $\cR_j,$ we derive from \eqref{eq:Ed-decomposition-Rj} a conjectural asymptotic formula for
$\cE_{d,k,\epsilon}(x;r,\alpha,\beta,\gamma)$ as long as $\beta=1$, in which case we have
\begin{align*}
\cE_{d,k,\epsilon}(x;r,\alpha,\beta,\gamma)
&\sim \frac{|\cT(\gamma_1,r,k)|\cdot \Delta(r,\beta,\eta)
\cdot E_{d,k,\epsilon}(r,\alpha,\beta)}{72\sqrt{3}}\\
&\ \ \ \ \times \prod_{p\nmid6r}\Big(1-\frac{(\frac{-3}{p})}{p-1}\Big)\cdot \frac{\sqrt{x}}{\log x},
\end{align*}
where for $\beta=1,$
\begin{align*}
E_{d,k,\epsilon}(r,\alpha,\beta)&=1+\epsilon(-1)^{\frac{k(\alpha-1)}{2}}\Big(\frac{3}{\alpha}\Big)\Omega_2(3;\alpha,r)\\
&\ \ \ \ +2\Re(\omega^d)\cdot\Big\{\Omega_3(3;\alpha,r)+\epsilon(-1)^{\frac{k(\alpha-1)}{2}}\Big(\frac{3}{\alpha}\Big)\Omega_6(3;\alpha,r)\Big\}
\end{align*}
and
\begin{align*}
\Delta(r,\beta,\eta)=
\begin{cases}
1,\ \ &\text{for~} r\equiv 4+3(2\eta_1-\eta_2)\bmod9,r\equiv\gamma_2\bmod2,\\
0, &\text{otherwise}.
\end{cases}
\end{align*}

\subsection{Counting in $\cT(\gamma_1,r,k)$}
To conclude Conjecture \ref{conj:Ed(x;r,alpha,beta,gamma)} for $\beta=1,$ 
it remains to evaluate $|\cT(\gamma_1,r,k)|$ as defined in \eqref{eq:T-countingt}. Recall the equation therein
\begin{align}\label{eq:T(gamma1,r,k)-equation}
2t^2+(2\gamma_1 +r)t+\frac{3\gamma_1^2-3r\gamma_1+r^2-1}{2}-k\equiv0\bmod 4,\ \ t\bmod 4.
\end{align}

If $2\nmid r$, then $2\nmid\gamma_2$ in view of \eqref{eq:D=3restrictions:beta=1}, and $|\cT(\gamma_1,r,k)|=1$ for any $k.$

If $2\mid r,$ then $2\mid\gamma_2$, and $\gamma_1$ should be odd.
The equation \eqref{eq:T(gamma1,r,k)-equation} reduces to
\begin{align*}
2t^2+(2+r)t+\frac{3\gamma_1^2-3r\gamma_1+r^2-1}{2}-k\equiv0\bmod 4.
\end{align*}

(1) If $r\equiv2\bmod4,$ then 
\begin{align*}
|\cT(\gamma_1,r,k)|
&=
\begin{cases}
2,\ \ &\text{~if~}-\gamma_1^2+\gamma_1r+r^2\equiv 2k+1\bmod4,\\
0,\ \ &\text{~if~}-\gamma_1^2+\gamma_1r+r^2\not\equiv 2k+1\bmod4,
\end{cases}\\
&=
\begin{cases}
2,\ \ &\text{~if~}2\mid k,\\
0,\ \ &\text{~if~}2\nmid k.
\end{cases}
\end{align*}

(2) If $r\equiv0\bmod4,$ then 
\begin{align*}
|\cT(\gamma_1,r,k)|
&=
\begin{cases}
4,\ \ &\text{~if~} 3\gamma_1^2-3\gamma_1r\equiv -r^2+2k+1\bmod8,\\
0,\ \ &\text{~if~} 3\gamma_1^2-3\gamma_1r\not\equiv -r^2+2k+1\bmod8.
\end{cases}\\
&=\begin{cases}
4,\ \ &\text{~if~} k\equiv1+r/2\bmod4,\\
0,\ \ &\text{~if~} k\not\equiv1+r/2\bmod4.
\end{cases}
\end{align*}

In summary, we find
\begin{align*}
|\cT(\gamma_1,r,k)|=
\begin{cases}
4,\ \ & \text{for~}4\mid r,~\gamma\equiv1\bmod2,~k\equiv1+r/2\bmod4,\\
2,\ \ & \text{for~}2\parallel r,~\gamma\equiv1\bmod2,~2\mid k,\\
1,\ \ & \text{for~}2\nmid r,~\gamma\equiv\omega,1+\omega\bmod2,\\
0, &\text{otherwise}.
\end{cases}
\end{align*}

\subsection{Concluding Conjecture \ref{conj:Ed(x;r,alpha,beta,gamma)}}~

\underline{\bf Case I: $\beta=1.$}
As given by \eqref{eq:constant-hD,r}, we have
\begin{align*}
h_{3,r}&=\xi(3,r)\cdot\frac{\sqrt{3}}{2}\prod_{p\nmid2r}\Big(1-\frac{(\frac{-3}{p})}{p-1}\Big),
\end{align*}
where
\begin{align*}
\xi(3,r)=
\begin{cases}
1, &\text{if } 2\mid r\text{~and~}3\nmid r,\\
2, &\text{if } (6,r)=1,\\
0, &\text{if } 3\mid r.
\end{cases}
\end{align*}
To prove Conjecture \ref{conj:Ed(x;r,alpha,beta,gamma)}, it suffices to verify that
\begin{align*}
\Delta(r,k,\gamma,\beta,\eta)\cdot h_{3,r}=\sqrt{3}\cdot \Delta(r,\beta,\eta)|\cT(\gamma_1,r,k)|\cdot \prod_{p\nmid6r}\Big(1-\frac{(\frac{-3}{p})}{p-1}\Big),
\end{align*}
which can be checked via
\begin{align*}
\Delta(r,\beta,\eta)|\cT(\gamma_1,r,k)|=
\begin{cases}
4,\ \ &\text{for~} r\equiv40-24(2\eta_1-\eta_2)
\bmod{72},k\equiv1\bmod4\text{~and~}\gamma\equiv1\bmod2,\\
4,\ \ &\text{for~} r\equiv4-24(2\eta_1-\eta_2)
\bmod{72},k\equiv3\bmod4\text{~and~}\gamma\equiv1\bmod2,\\
2,\ \ &\text{for~} r\equiv22-24(2\eta_1-\eta_2)
\bmod{36},2\mid k\text{~and~}\gamma\equiv1\bmod2,\\
1,\ \ &\text{for~} r\equiv13-24(2\eta_1-\eta_2)
\bmod{18}\text{~and~}\gamma\equiv\omega,1+\omega\bmod2,\\
0, &\text{otherwise}
\end{cases}
\end{align*}
in the case $\beta=1.$

\bigskip

\underline{\bf Case II: $\beta=\omega.$}
Modifying all above arguments
we arrive at
\begin{align*}
\cE_{d,k,\epsilon}(x;r,\alpha,\beta,\gamma,\eta)
&\sim \frac{|\cT(\gamma_1,r,k)|\cdot \Delta(r,\beta,\eta)
\cdot E_{d,k,\epsilon}(r,\alpha,\beta)}{72\sqrt{3}}\\
&\ \ \ \ \times\prod_{p\nmid6r}\Big(1-\frac{(\frac{-3}{p})}{p-1}\Big)\cdot \frac{\sqrt{x}}{\log x}
\end{align*}
for $\beta=\omega,$ in which case
\begin{align*}
E_{d,k,\epsilon}(r,\alpha,\beta)&=1+\epsilon(-1)^{\frac{k(\alpha-1)}{2}}\Big(\frac{3}{\alpha}\Big)\Omega_2(3;\alpha,r)\\
&\ \ \ \ +2\Re\Big(\omega^d\Big(\frac{\omega}{\alpha}\Big)_3\Big)\Big\{\Omega_3(3;\alpha,r)+\epsilon(-1)^{\frac{k(\alpha-1)}{2}}\Big(\frac{3}{\alpha}\Big)\Omega_6(3;\alpha,r)\Big\}
\end{align*}
and
\begin{align*}
\Delta(r,\beta,\eta)=
\begin{cases}
1,\ \ &\text{for~} r\equiv-2-3(\eta_1+\eta_2)\bmod9,r\equiv\gamma_1+\gamma_2\bmod2,\\
0, &\text{otherwise}.
\end{cases}
\end{align*}

We now evaluate the cardinality $|\cT(\gamma_1,r,k)|$ with respect to the condition $r\equiv\gamma_1+\gamma_2\bmod2.$
Recall the equation \eqref{eq:T(gamma1,r,k)-equation} in $\cT(\gamma_1,r,k)$:
\begin{align*}
2t^2+(2\gamma_1 +r)t+\frac{3\gamma_1^2-3r\gamma_1+r^2-1}{2}-k\equiv0\bmod 4,\ \ t\bmod 4.
\end{align*}

If $2\nmid r$, then $\gamma_1\not\equiv\gamma_2\bmod2$, and $|\cT(\gamma_1,r,k)|=1$ for any $k.$

If $2\mid r,$ then $\gamma_1\equiv\gamma_2\bmod2$, so that $\gamma_1\gamma_2$ should be odd.
The equation \eqref{eq:T(gamma1,r,k)-equation} reduces to
\begin{align*}
2t^2+(2+r)t+\frac{3\gamma_1^2-3r\gamma_1+r^2-1}{2}-k\equiv0\bmod 4.
\end{align*}

(1) If $r\equiv2\bmod4,$ then 
\begin{align*}
|\cT(\gamma_1,r,k)|
&=
\begin{cases}
2,\ \ &\text{~if~}-\gamma_1^2+\gamma_1r+r^2\equiv 2k+1\bmod4,\\
0,\ \ &\text{~if~}-\gamma_1^2+\gamma_1r+r^2\not\equiv 2k+1\bmod4,
\end{cases}\\
&=
\begin{cases}
2,\ \ &\text{~if~}2\mid k,\\
0,\ \ &\text{~if~}2\nmid k.
\end{cases}
\end{align*}

(2) If $r\equiv0\bmod4,$ then 
\begin{align*}
|\cT(\gamma_1,r,k)|
&=
\begin{cases}
4,\ \ &\text{~if~} 3\gamma_1^2-3\gamma_1r\equiv -r^2+2k+1\bmod8,\\
0,\ \ &\text{~if~} 3\gamma_1^2-3\gamma_1r\not\equiv -r^2+2k+1\bmod8.
\end{cases}\\
&=\begin{cases}
4,\ \ &\text{~if~} k\equiv1+r/2\bmod4,\\
0,\ \ &\text{~if~} k\not\equiv1+r/2\bmod4.
\end{cases}
\end{align*}

In summary, we find
\begin{align*}
|\cT(\gamma_1,r,k)|=
\begin{cases}
4,\ \ & \text{for~}4\mid r,~\gamma\equiv1+\omega\bmod2,~k\equiv1+r/2\bmod4,\\
2,\ \ & \text{for~}2\parallel r,~\gamma\equiv1+\omega\bmod2,~2\mid k,\\
1,\ \ & \text{for~}2\nmid r,~\gamma\equiv1,\omega\bmod2,\\
0, &\text{otherwise}.
\end{cases}
\end{align*}

It becomes a routine job to check that
\begin{align*}
\Delta(r,\beta,\eta)|\cT(\gamma_1,r,k)|=
\begin{cases}
4,\ \ &\text{for~} r\equiv16+24(\eta_1+\eta_2)
\bmod{72},k\equiv1\bmod4\text{~and~}\gamma\equiv1+\omega\bmod2,\\
4,\ \ &\text{for~} r\equiv52+24(\eta_1+\eta_2)
\bmod{72},k\equiv3\bmod4\text{~and~}\gamma\equiv1+\omega\bmod2,\\
2,\ \ &\text{for~} r\equiv34+24(\eta_1+\eta_2)
\bmod{36},2\mid k\text{~and~}\gamma\equiv1+\omega\bmod2,\\
1,\ \ &\text{for~} r\equiv7+24(\eta_1+\eta_2)
\bmod{18}\text{~and~}\gamma\equiv1,\omega\bmod2,\\
0, &\text{otherwise}
\end{cases}
\end{align*}
in the case $\beta=\omega$ as desired in Conjecture \ref{conj:Ed(x;r,alpha,beta,gamma)}.

\bigskip

\underline{\bf Case III: $\beta=\omega^2.$}
It remains to consider $\beta=\omega^2,$ for which similar arguments yield
\begin{align*}
\cE_{d,k,\epsilon}(x;r,\alpha,\beta,\gamma,\eta)
&\sim \frac{|\cT(\gamma_2,r,k)|\cdot \Delta(r,\beta,\eta)
\cdot E_{d,k,\epsilon}(r,\alpha,\beta)}{72\sqrt{3}}\\
&\ \ \ \ \times\prod_{p\nmid6r}\Big(1-\frac{(\frac{-3}{p})}{p-1}\Big)\cdot \frac{\sqrt{x}}{\log x},
\end{align*}
where
\begin{align*}
E_{d,k,\epsilon}(r,\alpha,\beta)&=1+\epsilon(-1)^{\frac{k(\alpha-1)}{2}}\Big(\frac{3}{\alpha}\Big)\Omega_2(3;\alpha,r)\\
&\ \ \ \ +2\Re\Big(\omega^d\Big(\frac{\omega^2}{\alpha}\Big)_3\Big)\Big\{\Omega_3(3;\alpha,r)+\epsilon(-1)^{\frac{k(\alpha-1)}{2}}\Big(\frac{3}{\alpha}\Big)\Omega_6(3;\alpha,r)\Big\}
\end{align*}
and
\begin{align*}
\Delta(r,\beta,\eta)=
\begin{cases}
1,\ \ &\text{for~} r\equiv-2+3(2\eta_2-\eta_1)\bmod9,r\equiv\gamma_1\bmod2,\\
0, &\text{otherwise}.
\end{cases}
\end{align*}

The treatment of $|\cT(\gamma_2,r,k)|$ is quite similar to that of $|\cT(\gamma_1,r,k)|$. In fact, we have
\begin{align*}
|\cT(\gamma_2,r,k)|=
\begin{cases}
4,\ \ & \text{for~}4\mid r,~\gamma\equiv\omega\bmod2,~k\equiv1+r/2\bmod4,\\
2,\ \ & \text{for~}2\parallel r,~\gamma\equiv\omega\bmod2,~2\mid k,\\
1,\ \ & \text{for~}2\nmid r,~\gamma\equiv1,1+\omega\bmod2,\\
0, &\text{otherwise},
\end{cases}
\end{align*}
which gives
\begin{align*}
\Delta(r,\beta,\eta)|\cT(\gamma_2,r,k)|=
\begin{cases}
4,\ \ &\text{for~} r\equiv16-24(2\eta_2-\eta_1)
\bmod{72},k\equiv1\bmod4\text{~and~}\gamma\equiv\omega\bmod2,\\
4,\ \ &\text{for~} r\equiv52-24(2\eta_2-\eta_1)
\bmod{72},k\equiv3\bmod4\text{~and~}\gamma\equiv\omega\bmod2,\\
2,\ \ &\text{for~} r\equiv34-24(2\eta_2-\eta_1)
\bmod{36},2\mid k\text{~and~}\gamma\equiv\omega\bmod2,\\
1,\ \ &\text{for~} r\equiv7-24(2\eta_2-\eta_1)
\bmod{18}\text{~and~}\gamma\equiv1,1+\omega\bmod2,\\
0, &\text{otherwise}.
\end{cases}
\end{align*}
Then the proof of Theorem \ref{thm:HL-Ed(x;r,alpha,beta,gamma)} is completed by collecting all above cases.

\section{From Conjecture \ref{conj:Ed(x;r,alpha,beta,gamma)} to Lang--Trotter}
\label{sec:proofofTheorem{thm:HL-LT:D=3}}

We now prove  Theorem \ref{thm:HL-LT:D=3}. To begin with, we derive from \eqref{eq:initialtransformation:D=3} that
\begin{equation}\label{eq:connection:pi-Ed}
\begin{split}
\pi_{E,r}(x)
\sim\frac{1}{2}&\mathop{\sum\sum}_{\substack{0\leqslant a,b\leqslant 1\\(a,b)\neq (0,0)}}\mathop{\sum\sum\sum}_{0\leqslant d,s,t\leqslant 2}\sum_{0\leqslant k\leqslant 3}\\
&\sum_{\pm}\cE_{d-d_0,k,\pm\epsilon}(x;\mp r,g_1,\omega^d,a+b\omega,s+t\omega),
\end{split}
\end{equation}
where
\begin{align*}
d_0=(2+\lambda)(2+(-1)^a-(-1)^b)-\mu t,\ \ \ \epsilon=(-1)^{\frac{\lambda k(k+1)}{2}+(\delta+\mu)k}.
\end{align*}
According to the choices of $a,b$ and $k$, we have
\begin{align*}
d_0=
\begin{cases}
-\mu t,\ \ &\text{for~} a=1\text{~and~}b=0,\\
2(2+\lambda)-\mu t,\ \ &\text{for~} a=1\text{~and~}b=1,\\
4(2+\lambda)-\mu t,\ \ &\text{for~} a=0\text{~and~}b=1,\end{cases}
\end{align*}
and also
\begin{align*}
\epsilon=\begin{cases}
1,\ \ &\text{for~} k=0,\\
(-1)^{\delta+\lambda+\mu},\ \ &\text{for~} k=1,\\
(-1)^{\lambda},\ \ &\text{for~} k=2,\\
(-1)^{\delta+\mu},\ \ &\text{for~} k=3.
\end{cases}
\end{align*}

Assuming Conjecture \ref{conj:Ed(x;r,alpha,beta,gamma)}, we find
\begin{equation}\label{eq:connection:Ed-asymptotic}
\begin{split}
&\cE_{d-d_0,k,\pm\epsilon}(x;\mp r,g_1,\omega^d,a+b\omega,s+t\omega)\\
&\ \ \ \sim\frac{h_{3,r}}{216}\cdot \Delta(\mp r,k,a+b\omega,\omega^d,s+t\omega)\cdot E_{d-d_0,k,\pm\epsilon}(r,g_1,\omega^d)\frac{\sqrt{x}}{\log x}.
\end{split}
\end{equation}
Hence it suffices to simplify
\begin{align*}
C&:=\sum_{0\leqslant d\leqslant2}C_d
\end{align*}
with
\begin{align*}
C_d
&=\mathop{\sum\sum}_{\substack{0\leqslant a,b\leqslant 1\\(a,b)\neq (0,0)}}\mathop{\sum\sum}_{0\leqslant s,t\leqslant 2}\sum_{0\leqslant k\leqslant 3}\sum_{\pm}\Delta(\mp r,k,a+b\omega,\omega^d,s+t\omega)\\
&\ \ \ \times E_{d-d_0,k,\pm\epsilon}(r,g_1,\omega^d).
\end{align*}

It is rather complicated to compute $C$ explicitly in all cases. To convince skeptical readers, we show some details for $r\equiv40\bmod{72}.$ To simply the presentation, we write $\Omega_j=\Omega_j(3;g_1,r).$
In fact,
\begin{align*}
C_0
&=\mathop{\sum\sum}_{\substack{0\leqslant s,t\leqslant 2\\2s-t\equiv0\bmod3}}\Delta(r,1,1,1,s+t\omega)E_{\mu t,1,-(-1)^{\delta+\lambda+\mu}}(r,g_1,1)\\
&=8\mathop{\sum\sum}_{\substack{0\leqslant s,t\leqslant 2\\2s-t\equiv0\bmod3}}\Big\{1-(-1)^{\delta+\lambda+\mu+\frac{g_1-1}{2}}\Big(\frac{3}{g_1}\Big)\Omega_2\\
&\ \ \ \ +2\Re(\omega^{\mu t})\cdot\Big(\Omega_3-(-1)^{\delta+\lambda+\mu+\frac{g_1-1}{2}}\Big(\frac{3}{g_1}\Big)\Omega_6\Big)\Big\}\\
&=24\Big\{1-(-1)^{\delta+\lambda+\mu+\frac{g_1-1}{2}}\Big(\frac{3}{g_1}\Big)\Omega_2\Big\}\\
&\ \ \ \ +16\sum_{0\leqslant t\leqslant 2}\Re(\omega^{\mu t})\cdot\Big\{\Omega_3-(-1)^{\delta+\lambda+\mu+\frac{g_1-1}{2}}\Big(\frac{3}{g_1}\Big)\Omega_6\Big\}.
\end{align*}
Correspondingly, we also have
\begin{align*}
C_1
&=\mathop{\sum\sum}_{\substack{0\leqslant s,t\leqslant 2\\s+t\equiv1\bmod3}}\Delta(r,1,1+\omega,\omega,s+t\omega)E_{1+\mu t,1,-(-1)^{\delta+\lambda+\mu}}(r,g_1,\omega)\\
&=24\Big\{1-(-1)^{\delta+\lambda+\mu+\frac{g_1-1}{2}}\Big(\frac{3}{g_1}\Big)\Omega_2\Big\}\\
&\ \ \ \ +16\sum_{0\leqslant t\leqslant 2}\Re\Big(\omega^{\lambda+\mu t}\Big(\frac{\omega}{g_1}\Big)_3\Big)\cdot\Big\{\Omega_3-(-1)^{\delta+\lambda+\mu+\frac{g_1-1}{2}}\Big(\frac{3}{g_1}\Big)\Omega_6\Big\},
\end{align*}
and
\begin{align*}
C_2
&=\mathop{\sum\sum}_{\substack{0\leqslant s,t\leqslant 2\\2t-s\equiv-1\bmod3}}\Delta(r,1,\omega,\omega^2,s+t\omega)E_{2+\mu t,1,-(-1)^{\delta+\lambda+\mu}}(r,g_1,\omega^2)\\
&=24\Big\{1-(-1)^{\delta+\lambda+\mu+\frac{g_1-1}{2}}\Big(\frac{3}{g_1}\Big)\Omega_2\Big\}\\
&\ \ \ \ +16\sum_{0\leqslant t\leqslant 2}\Re\Big(\omega^{-\lambda+\mu t}\Big(\frac{\omega^2}{g_1}\Big)_3\Big)\cdot\Big\{\Omega_3-(-1)^{\delta+\lambda+\mu+\frac{g_1-1}{2}}\Big(\frac{3}{g_1}\Big)\Omega_6\Big\}.
\end{align*}
It thus follows that
\begin{align*}
C&=72\Big\{1-(-1)^{\delta+\lambda+\mu+\frac{g_1-1}{2}}\Big(\frac{3}{g_1}\Big)\Omega_2\Big\}+16\Big(1+2\Re\Big(\omega^{\lambda+\frac{g_1^2-1}{3}}\Big)\Big)\\
&\ \ \ \ \times\sum_{0\leqslant t\leqslant 2}\Re(\omega^{\mu t})\cdot\Big\{\Omega_3-(-1)^{\delta+\lambda+\mu+\frac{g_1-1}{2}}\Big(\frac{3}{g_1}\Big)\Omega_6\Big\},
\end{align*}
which is in fact valid for all $r\equiv16\bmod{24}.$

In view of \eqref{eq:connection:pi-Ed} and \eqref{eq:connection:Ed-asymptotic}, it remains to check
\begin{align*}
\widehat{\mu}=\frac{2}{9}\sum_{0\leqslant t\leqslant 2}\Re(\omega^{\mu t})
\end{align*}
and
\begin{align*}
\frac{C}{72}=1+\widehat{\mu}\cdot\varsigma_1(g,r)\Omega_3+\varsigma_2(g,r)\Big(\frac{3}{g_1}\Big)(\Omega_2+\widehat{\mu}\cdot\varsigma_1(g,r)\Omega_6)
\end{align*}
as required in Theorem \ref{thm:HL-LT:D=3}. The above arguments can also be modified to give the desired value of $C$ for $r\equiv-16\bmod{24}.$

The cases $r\equiv\pm1\bmod6$ are more complicated. For instance, we have
\begin{align*}
C_0
&=\mathop{\sum\sum}_{\substack{0\leqslant s,t\leqslant 2\\2s\equiv t\bmod3}}
\big\{\Delta(r,0,\omega,1,s+t\omega)E_{1-\lambda+\mu t,0,-1}(r,g_1,1)\\
&\ \ \ \ \ +\Delta(r,1,\omega,1,s+t\omega)E_{1-\lambda+\mu t,1,-(-1)^{\delta+\lambda+\mu}}(r,g_1,1)\\
&\ \ \ \ \ +\Delta(r,2,\omega,1,s+t\omega)E_{1-\lambda+\mu t,2,-(-1)^{\lambda}}(r,g_1,1)\\
&\ \ \ \ \ +\Delta(r,3,\omega,1,s+t\omega)E_{1-\lambda+\mu t,3,-(-1)^{\delta+\mu}}(r,g_1,1)\\
&\ \ \ \ \ +\Delta(r,0,1+\omega,1,s+t\omega)E_{2+\lambda+\mu t,0,-1}(r,g_1,1)\\
&\ \ \ \ \ +\Delta(r,1,1+\omega,1,s+t\omega)E_{2+\lambda+\mu t,1,-(-1)^{\delta+\lambda+\mu}}(r,g_1,1)\\
&\ \ \ \ \ +\Delta(r,2,1+\omega,1,s+t\omega)E_{2+\lambda+\mu t,2,-(-1)^{\lambda}}(r,g_1,1)\\
&\ \ \ \ \ +\Delta(r,3,1+\omega,1,s+t\omega)E_{2+\lambda+\mu t,3,-(-1)^{\delta+\mu}}(r,g_1,1)\big\}
\end{align*}
for $r\equiv13\bmod{18},$ in which case we also have similar evaluations for $C_1$ and $C_2$.
The complicated details will be given in Appendix \ref{ch:appendix}.

%% file: 6.otherCM.tex
\chapter{Complex multiplication by remaining imaginary quadratic fields}\label{ch:otherCM}

This chapter deals with elliptic curves with other CM cases by $\Q(\sqrt{-D})$, i.e., $D\in\{2,7,11,19,43,67,163\}.$
We consider the CM elliptic curve 
\begin{align}\label{eq:E:otherD}
E/\Q:y^2=4x^3+ax+b,
\end{align}
where $a,b$ will be given in terms of an arbitrary integer $g$ as below such that the curve is non-singular:
\begin{itemize}
\item $D=2:$ $(a,b)=(-\frac{40}{3}g^2,-\frac{224}{27}g^3)$, \smallskip
\item $D=7:$ $(a,b)=(-\frac{35}{4}g^2,-\frac{49}{8}g^3)$,\smallskip
\item $D=11:$ $(a,b)=(-\frac{88}{3}g^2,\frac{847}{27}g^3)$,\smallskip
\item $D=19:$ $(a,b)=(-152g^2,361g^3)$,\smallskip
\item $D=43:$ $(a,b)=(-80\cdot 43g^2,21\cdot 43^2g^3)$,\smallskip
\item $D=67:$ $(a,b)=(-440\cdot 67g^2,217\cdot 67^2g^3)$,\smallskip
\item $D=163:$ $(a,b)=(-53360\cdot 163g^2,185801\cdot163^2g^3)$.
\end{itemize}
We keep all such choices of $a$ and $b$ for different $D$ throughout this chapter. There are many different choices of curves in an isogenous class, and choices of all above scalars make the expression of $a_p(E)$ to take a uniform shape when $D\geqslant 7.$ See Lemma \ref{lm:Frobenius:D>=7} below for details.

To characterize $\pi_{E,r}(x)$ under the Hardy--Littlewood conjecture, we write
\begin{align*}
g=(-1)^\delta2^\lambda D^\mu g_1
\end{align*}
with
\begin{itemize}
\item $\delta\in\{0,1\},$
\item $\lambda,\mu\in\N,$
\item $g_1\in\Z^+$ with $(2D,g_1)=1$,
\end{itemize}
where we always assume $\mu=0$ for $D=2.$

Recall the definition \eqref{eq:Omegafunction}:
\begin{align*}
\Omega_j(D;q,r)=\prod_{\substack{p^\nu\parallel q,~ p\nmid r\\j\nmid \nu}}\frac{-1}{p-1-(\frac{-D}{p})}.
\end{align*}
As in Chapter \ref{ch:HL&upperbound}, denote by $R_D$ the ring of algebraic integers in $\Q(\sqrt{-D})$, so that
$R_2=\Z[\sqrt{-2}]$ and $R_D=\Z[\frac{1+\sqrt{-D}}{2}]$ for $D\geqslant7.$
Moreover, there are only two roots of unity (i.e., $\pm1$) in $R_D$ for each of above $D$. Quadratic reciprocity law (amongst other arguments) is sufficient, so that the arguments here are much simpler than those for $D=1,3$.

\section{Asymptotics for $\pi_{E,r}(x)$: conditional results}
The primary task of this chapter is to prove the asymptotic formula
\begin{align}\label{eq:piE,r(x)-varpiE,r,otherD}
\pi_{E,r}(x)\sim\frac{\varpi_{E,r}\cdot\sqrt{x}}{\log x}
\end{align}
for all above CM elliptic curves $E$. The proof is also conditional on the Hardy--Littlewood conjecture\index{Hardy--Littlewood conjecture}, and the constant $\varpi_{E,r}$ can also be made explicitly.

\begin{theorem}\label{thm:HL-LT:D=2}\index{Lang--Trotter conjecture}
Let $D=2$. Suppose $E/\Q$ is the elliptic curve defined by $\eqref{eq:E:otherD}$ with the above convention.
For each non-zero integer $r\equiv2\bmod4,$ the Hardy--Littlewood conjecture \index{Hardy--Littlewood conjecture} implies that $\eqref{eq:piE,r(x)-varpiE,r,otherD}$ holds with
\begin{align}\label{eq:varpiE,r,D=2}
\varpi_{E,r}=\frac{h_{2,r}}{2}\Big(1+\frac{1}{2}(-1)^{\frac{(r-2)(r+10)}{32}+\delta+\lambda+\frac{g_1-1}{2}}\Big(\frac{2}{g_1}\Big)\Omega_2(2;g_1,r)\Big),
\end{align}
where $h_{2,r}\geqslant0$ is the constant defined by $\eqref{eq:constant-hD,r}.$
Unconditionally, $\pi_{E,r}(x)$ is bounded by $O(1)$ if $r\not\equiv2\bmod4$.
\end{theorem}

\begin{theorem}\label{thm:HL-LT:D>=7}\index{Lang--Trotter conjecture}
Let $D\in\{7,11, 19, 43, 67, 163\}.$ Suppose $E/\Q$ is the elliptic curve defined by $\eqref{eq:E:otherD}$ with the above convention.

For each non-zero integer $r,$ the Hardy--Littlewood conjecture \index{Hardy--Littlewood conjecture} implies that $\eqref{eq:piE,r(x)-varpiE,r,otherD}$ holds with
\begin{align}\label{eq:varpiE,r,otherD}
\varpi_{E,r}=\frac{h_{D,r}}{2}\cdot\Big(1+\varepsilon_D(g,r)\Big(\frac{2^{\lambda+1}g_1r}{D}\Big)\Omega_2(D;g_1,r)\Big)
\end{align}
for $D\nmid r,$
where $h_{D,r}\geqslant0$ is the constant defined by $\eqref{eq:constant-hD,r},$ and
\begin{align*}
\varepsilon_D(g,r)
&=\begin{cases}
\frac{1}{2}(-1)^{\frac{g_1-1}{2}}(1+(-1)^{\lambda}),\ \ &\text{for~}  2\parallel r,\\
(-1)^{\delta+\mu+\frac{\lambda r}{4}},\ \ &\text{for~} 4\mid r,\\
\frac{1}{4}(1+(-1)^{\lambda})((-1)^{\frac{g_1-1}{2}}+(-1)^{\delta+\mu}),\ \ &\text{for~} 2\nmid r.
\end{cases}
\end{align*}
Unconditionally, $\pi_{E,r}(x)$ is bounded by $O(1)$ if $D\mid r$.
\end{theorem}

As a consequence of Theorems \ref{thm:HL-LT:D=2} and \ref{thm:HL-LT:D>=7}, we may obtain the following asymptotic behaviour of the mixed sum $\pi_{E,r}(x)+\pi_{E,-r}(x)$ under the Hardy--Littlewood conjecture,  which coincides with Lemma \ref{lm:LT-quadraticpolynomial} for $D\neq1,3.$

\begin{corollary}\label{coro:HL-LT:D>=7}
 Suppose $E/\Q$ is the elliptic curve, with CM by $\Q(\sqrt{-D})$ as defined by $\eqref{eq:E:otherD}$ for $D\in\{7,11, 19, 43, 67, 163\}.$

For each non-zero integer $r,$ the Hardy--Littlewood conjecture \index{Hardy--Littlewood conjecture} implies 
\begin{align*}
\pi_{E,r}(x)+\pi_{E,-r}(x)
&\sim  \frac{h_{D,r}\cdot\sqrt{x}}{\log x}
\end{align*}
for $D\nmid r,$ where $h_{3,r}\geqslant0$ is the constant defined by $\eqref{eq:constant-hD,r}.$
Unconditionally, $\pi_{E,r}(x)$ is bounded by $O(1)$ if $D\mid r$.
\end{corollary}

Corollary \ref{coro:HL-LT:D>=7} follows by observing that $\varepsilon_D(g,r)=\varepsilon_D(g,-r)$ and $(\frac{-r}{D})=-(\frac{r}{D})$ since $D\equiv3\bmod4$ for $D\in\{7,11, 19, 43, 67, 163\}.$ In the case of $D=2,$ we have
\begin{align*}
(-1)^{\frac{(r-2)(r+10)}{32}}=-(-1)^{\frac{(-r-2)(-r+10)}{32}}
\end{align*}
for $r\equiv2\bmod4.$ This leads to the following consequence.

\begin{corollary}\label{coro:HL-LT:D=2}
 Suppose $E/\Q$ is the elliptic curve, with CM by $\Q(\sqrt{-D})$ as defined by $\eqref{eq:E:otherD}$ for $D=2.$

For each non-zero integer $r\equiv2\bmod4,$ the Hardy--Littlewood conjecture \index{Hardy--Littlewood conjecture} implies 
\begin{align*}
\pi_{E,r}(x)+\pi_{E,-r}(x)
&\sim  \frac{h_{2,r}\cdot\sqrt{x}}{\log x},
\end{align*}
where $h_{2,r}\geqslant0$ is the constant defined by $\eqref{eq:constant-hD,r}$.
Unconditionally, $\pi_{E,r}(x)$ is bounded by $O(1)$ if $r\not\equiv2\bmod4$.
\end{corollary}

\section{Initial transformations}

We start with the following explicit formulae for Frobenius traces \index{trace of Frobenius} $a_p(E)$. 
\begin{lemma}\label{lm:Frobenius:D=2}\index{trace of Frobenius}
Let $D=2.$ Suppose $p\in\Z^+$ with $p\nmid 2g.$ Then $a_p(E)=0$ if $p$ is inert in $\Z[\sqrt{-2}]$ and
\begin{align*}
a_p(E)=\Big(\frac{g}{p}\Big)(-1)^{[\frac{p}{8}]+\frac{\Tr(\pi)-2}{4}}\Tr(\pi)
\end{align*}
if $p$ splits in $\Z[\sqrt{-2}]$ with $p=\pi\overline{\pi}$ such that $\Tr(\pi)$ is positive.
\end{lemma}

\begin{lemma}\label{lm:Frobenius:D>=7}\index{trace of Frobenius}
Let $D\in\{7,11, 19, 43, 67, 163\}.$ Suppose $p\in\Z^+$ with $p\nmid 2g.$ Then $a_p(E)=0$ if $p$ is inert in $R_D,$ and
\begin{align*}
a_p(E)=\Big(\frac{g}{p}\Big)\Big(\frac{2\Tr(\pi)}{D}\Big)\Tr(\pi)
\end{align*}
if $p$ splits in $R_D$ with $p=\pi\overline{\pi}$ such that $\Tr(\pi)$ is positive.
\end{lemma}

Lemma \ref{lm:Frobenius:D>=7} is directly borrowed from H. M. Stark \cite{St96} with slight modifications. An alternative formulation of Lemma \ref{lm:Frobenius:D=2} can be found in A. R. Rajwade \cite{Ra68}. The original version, given as below, was stated for the 
curve $y^2=x(x^2-4gx+2g^2)$, which is isogenous to $y^2=4x^3-\frac{40}{3}g^2x-\frac{224}{27}g^3$ according to the morphism $(x,y)\mapsto (x+4g/3,y/2).$

\begin{lemma}\label{lm:Frobenius:D=2-original}\index{trace of Frobenius}
Let $D=2$ Suppose $p\in\Z^+$ with $p\nmid 2g.$ Then $a_p(E)=0$ if $p$ is inert in $\Z[\sqrt{-2}]$ and
\begin{align*}
a_p(E)=\Big(\frac{g}{p}\Big)(-1)^{\frac{p-1}{8}}\Tr(\pi)
\end{align*}
if $p\equiv1\bmod 8$ with $p=\pi\overline{\pi}$ such that $\pi\equiv1\bmod{2\sqrt{-2}},$ and
\begin{align*}
a_p(E)=\Big(\frac{g}{p}\Big)(-1)^{\frac{p-3}{8}}\Tr(\pi)
\end{align*}
if $p\equiv3\bmod 8$ with $p=\pi\overline{\pi}$ such that $\pi\equiv1+\sqrt{-2}\bmod{2\sqrt{-2}}.$
\end{lemma}

Before turning to characterize $\pi_{E,r}(x),$ we show how to obtain Lemma \ref{lm:Frobenius:D=2} via 
Lemma \ref{lm:Frobenius:D=2-original}.

Suppose $p=\pi\overline{\pi}\equiv1\bmod8,$ and $\pi=a+b\sqrt{-2}\equiv1\bmod{2\sqrt{-2}}.$ We then find
\[\Tr(\pi)=2a,\ \ a\equiv1\bmod4,\ \ 2\mid b.\]
Hence $\frac{\Tr(\pi)-2}{4}$ is even and $\frac{\Tr(-\pi)-2}{4}$ is odd, so that Lemma \ref{lm:Frobenius:D=2-original} yields
\begin{align*}
a_p(E)=\Big(\frac{g}{p}\Big)(-1)^{\frac{p-1}{8}}\Tr(\pi)=\Big(\frac{g}{p}\Big)(-1)^{[\frac{p}{8}]+\frac{\Tr(\pi)-2}{4}}\Tr(\pi)
\end{align*}
if $a>0,$ and
\begin{align*}
a_p(E)=\Big(\frac{g}{p}\Big)(-1)^{\frac{p-1}{8}}\Tr(\pi)=\Big(\frac{g}{p}\Big)(-1)^{[\frac{p}{8}]+\frac{\Tr(-\pi)-2}{4}}\Tr(-\pi)
\end{align*}
if $a<0$. This confirms Lemma \ref{lm:Frobenius:D=2} for $p\equiv1\bmod 8.$ The arguments for $p\equiv3\bmod8$ are quite similar.

We now come to evaluate $\pi_{E,r}(x)$ by virtue of Lemmas \ref{lm:Frobenius:D=2} and \ref{lm:Frobenius:D>=7}.
Consider the case $D=2$ firstly. Lemma  \ref{lm:Frobenius:D=2} yields
\begin{align*}
2\pi_{E,r}(x)
&=\sum_{\substack{\pi\in\Pi_2(x)\\ N(\pi)=p, ~\Tr(\pi)=u>0\\ (\frac{g}{p})(-1)^{[\frac{p}{8}]+\frac{u-2}{4}}u=r}}1+O(1).
\end{align*}
The quadratic reciprocity law \index{quadratic reciprocity law} then gives
\begin{align*}
\Big(\frac{g}{p}\Big)
=\Big(\frac{-1}{p}\Big)^\delta\Big(\frac{2}{p}\Big)^\lambda\Big(\frac{g_1}{p}\Big)
=(-1)^{\frac{\lambda(p^2-1)}{8}+\frac{(2\delta+g_1-1)(p-1)}{4}}\Big(\frac{p}{g_1}\Big),
\end{align*}
from which we derive that
\begin{align}\label{eq:initialtransformation:D=2}
\pi_{E,r}(x)
&\sim\frac{1}{2}~\sideset{}{^\ddagger}\sum_{\substack{\pi\in\Pi_2(x)\\ N(\pi)=p, ~\Tr(\pi)=|r|}}1,
\end{align}
where the superscript $\ddagger$ means
\begin{align}\label{eq:initialtransformation:D=2,ddagger}
\Big(\frac{p}{g_1}\Big)=(-1)^{[\frac{p}{8}]+\frac{r-2}{4}+\frac{\lambda(p^2-1)}{8}+\frac{(2\delta+g_1-1)(p-1)}{4}}.
\end{align}
Note that the LHS of \eqref{eq:initialtransformation:D=2,ddagger} depends on $p\bmod{g_1},$ and the RHS depends on $p\bmod{16}.$ Hence it reduces to consider the number of $\pi\in\Pi_2(x)$ with fixed trace $\Tr(\pi)$ such that the norm $N(\pi)$ falls into a certain residue class modulo $16g_1.$ 

Before going into more details for $D=2$, we consider $D\in\{7,11, 19, 43, 67, 163\}$ alternatively, in which case Lemma \ref{lm:Frobenius:D>=7} gives
\begin{align*}
2\pi_{E,r}(x)
&=\sum_{\substack{\pi\in\Pi_D(x)\\ N(\pi)=p, ~\Tr(\pi)=u>0\\(\frac{g}{p})(\frac{2u}{D})u=r}}1+O(1).
\end{align*}
From the quadratic reciprocity law, we obtain
\begin{align*}
\Big(\frac{g}{p}\Big)=\Big(\frac{-1}{p}\Big)^{\delta}\Big(\frac{2}{p}\Big)^\lambda\Big(\frac{D^\mu g_1}{p}\Big)=(-1)^{\frac{(\delta+\mu)(p-1)}{2}+\frac{\lambda(p^2-1)}{8}+\frac{(g_1-1)(p-1)}{4}}\Big(\frac{p}{D^\mu g_1}\Big),
\end{align*}
so that
\begin{align}\label{eq:initialtransformation:D>=7}
\pi_{E,r}(x)
&\sim\frac{1}{2}~\sideset{}{^\ddagger}\sum_{\substack{\pi\in\Pi_D(x)\\ N(\pi)=p, ~\Tr(\pi)=|r|}}1,
\end{align}
where $\ddagger$ means
\begin{align*}
(-1)^{\frac{(\delta+\mu)(p-1)}{2}+\frac{\lambda(p^2-1)}{8}+\frac{(g_1-1)(p-1)}{4}}\Big(\frac{p}{D^\mu g_1}\Big)\Big(\frac{2|r|}{D}\Big)=\sgn(r).
\end{align*}
This leads to another counting problem on the number of $\pi\in\Pi_D(x)$ with fixed  trace $\Tr(\pi)$ such that the norm $N(\pi)$ falls into a certain residue class modulo $8D^\mu g_1.$ 

To treat the cases $D=2$ and $D\in\{7,11, 19, 43, 67, 163\}$ in a uniform manner, we introduce the following counting function
\begin{align*}
\Pi_D(x;r,q,a):=|\{\pi\in\Pi_D(x):\Tr(\pi)=r,~N(\pi)\equiv a\bmod q\}|
\end{align*}
for all integers $r,q,a$ with $q\geqslant1$ and $(q,a)=1.$ It is clear that $\Pi_D(x;r,q,a)=\Pi_D(x;-r,q,a).$
In the next section, we give a conjectural asymptotic formula for $\Pi_D(x;r,q,a),$ which can be used to characterize $\pi_{E,r}(x)$ in view of \eqref{eq:initialtransformation:D=2} and \eqref{eq:initialtransformation:D>=7}.

To close this section, we give precise expressions of $\pi_{E,r}(x)$ in terms of $\Pi_D(x;r,q,a)$.

\underline{(I) $D=2.$} From \eqref{eq:initialtransformation:D=2} and the Chinese remainder theorem, it follows that
\begin{align}\label{eq:connection:pi-PiD(x),D=2}
\pi_{E,r}(x)
&\sim\frac{1}{2}\sum_{0\leqslant k\leqslant 7}\sum_{\substack{a\bmod{g_1}\\(\frac{a}{g_1})=(-1)^{\delta_k}}}\Pi_D(x;r,16g_1,\delta_{k,a}),
\end{align}
where 
\begin{align*}
\delta_k&=\Big[\frac{2k+1}{8}\Big]+\frac{r-2}{4}+\frac{\lambda k(k+1)}{2}+\frac{(2\delta+g_1-1)k}{2},\\
\delta_{k,a}&\equiv 2k+1\bmod{16},\\
\delta_{k,a}&\equiv a\bmod{g_1}.
\end{align*}

\underline{(II) $D\in\{7,11, 19, 43, 67, 163\}.$} Note that $(\frac{|r|}{D})=(\frac{r}{D})\sgn(r)$ since $D\equiv3\bmod4.$ From \eqref{eq:initialtransformation:D>=7} and the Chinese remainder theorem\index{Chinese remainder theorem}, it follows that
\begin{align}\label{eq:connection:pi-PiD(x),D>=7}
\pi_{E,r}(x)
&\sim\frac{1}{2}\sum_{0\leqslant k\leqslant 3}\sum_{\substack{a\bmod A\\ (\frac{a}{A})=(-1)^{\varepsilon_k}(\frac{2r}{D})}}\Pi_D(x;r,8A,\varepsilon_{k,a}),
\end{align}
where $A=D^\mu g_1,$
\begin{align*}
\varepsilon_k&=\Big(\delta+\mu+\frac{g_1-1}{2}\Big)k+\frac{\lambda k(k+1)}{2},\\
\varepsilon_{k,a}&\equiv2k+1\bmod8,\\
\varepsilon_{k,a}&\equiv a\bmod A.
\end{align*}

\section{Asymptotic evaluation of $\Pi_D(x;r,q,a)$}
In order to formulate the conjecture for $\Pi_D(x;r,q,a)$, we introduce the convention
\begin{align}\label{eq:theta(D,r)}
\theta(D,r)=
\begin{cases}
1,\ \ & \text{if~} D=2\text{~and~}r\equiv2\bmod4,\\
1, &\text{if } D\equiv3\bmod 4,2\mid r\text{~and~}(D,r)=1,\\
1, &\text{if } D\equiv3\bmod 8,2\nmid r\text{~and~}(D,r)=1,\\
0, &\text{otherwise},
\end{cases}
\end{align}
which can be compared with \eqref{eq:xi(D,r)}.

We expect $\Pi_D(x;r,q,a)$ to satisfy the following asymptotic behaviour.

\begin{conjecture}\label{conj:PiD(x;r,q,a)}
Let $D\in\{2,3, 7, 11, 19, 43, 67, 163\},$ $r\in\Z^*$ and $(q,a)=1$. For all large $x$, we have
\begin{align*}
\Pi_D(x;r,q,a)
&\sim \theta(D,r)\varrho(r,q,a)\cdot\frac{2\sqrt{2}}{q}\frac{q^\sharp}{\varphi(q^\sharp)}\prod_{p\nmid2qr}\Big(1-\frac{(\frac{-2}{p})}{p-1}\Big)\cdot\frac{\sqrt{x}}{\log x}
\end{align*}
if $D=2,$ and
\begin{align*}
\Pi_D(x;r,q,a)
&\sim \theta(D,r)\varrho_D(r,q,a)\cdot\frac{2(2,q(r+1))}{q}\frac{\sqrt{D}q^\sharp}{\varphi(Dq^\sharp)}\prod_{p\nmid2qr}\Big(1-\frac{(\frac{-D}{p})}{p-1}\Big)\cdot\frac{\sqrt{x}}{\log x}
\end{align*}
if $D>2,$
where $q^\sharp$ denotes the odd part of $q,$
\begin{align*}
\varrho(r,q,a)
&:=|\{t\bmod q:2t^2+\frac{r^2}{4}\equiv a\bmod q\}|,\\
\varrho_D(r,q,a)
&:=|\{t\bmod q:Dt^2-Drt+\frac{D+1}{4}r^2\equiv a\bmod q\}|.
\end{align*}
\end{conjecture}

\begin{theorem}\label{thm:HL-PiD(x;r,q,a)}
The Hardy--Littlewood conjecture \index{Hardy--Littlewood conjecture} implies Conjecture $\ref{conj:PiD(x;r,q,a)}$.
\end{theorem}

The proof of Theorem \ref{thm:HL-PiD(x;r,q,a)} goes separately for $D=2$ and $D\geqslant7$ since the ring $R_D$ of integers in $\Q(\sqrt{-D})$ is given by
\[R_2=\Z[\sqrt{-2}],\ \ R_D=\Z\Big[\frac{1+\sqrt{-D}}{2}\Big]\text{~for~}D\geqslant7.\]

We now assume $D=2$. Then each prime element $\pi\in\Pi_2$ takes the shape $\pi=m+n\sqrt{-2}$ with $m,n\in\Z.$ Hence
\begin{align*}
\Pi_D(x;r,q,a)
&=\sum_{\substack{m+n\sqrt{-2}\in\Pi_2(x)\\2m=r\\ m^2+2n^2\equiv a\bmod q}}1=\sum_{\substack{n\in\Z\\ 2n^2+(r/2)^2\in\cP\cap[1,x]\\ 2n^2+(r/2)^2\equiv a\bmod q}}1.
\end{align*}
The above sums are bounded unless $r\equiv2\bmod4$, which we henceforth assume (in the case of $D=2$).
Following the convention in Conjecture \ref{conj:HL-modified}, we find
\begin{align*}
\Pi_D(x;r,q,a)
&=2\sum_{\substack{u\bmod q\\2u^2+\frac{r^2}{4}\equiv a\bmod q }}\pi_{2,0,(r/2)^2}(x;q,u).
\end{align*}
All prerequisites in Conjecture \ref{conj:HL-modified} are satisfied since $r\equiv2\bmod4,$ so that
\begin{align*}
\Pi_D(x;r,q,a)
&\sim \varrho(r,q,a)\cdot\frac{2\sqrt{2}}{q}\frac{q^\sharp}{\varphi(q^\sharp)}\prod_{p\nmid2qr}\Big(1-\frac{(\frac{-2}{p})}{p-1}\Big)\cdot\frac{\sqrt{x}}{\log x}
\end{align*}
as implied by the modified Hardy--Littlewood conjecture. This proves Theorem \ref{thm:HL-PiD(x;r,q,a)} 
for $D=2$.

We turn to confirming  the remaining cases $D\in\{3, 7, 11, 19, 43, 67, 163\}$.
Firstly, we have
\begin{align*}
\Pi_D(x;r,q,a)
&=\sum_{\substack{m+n\frac{1+\sqrt{-D}}{2}\in\Pi_D(x)\\2m+n=r\\ m^2+mn+\frac{D+1}{4}n^2\equiv a\bmod q}}1
=\sum_{\substack{m\in\Z\\ Dm^2-Drm+\frac{D+1}{4}r^2\in\cP\cap[1,x]\\ Dm^2-Drm+\frac{D+1}{4}r^2\equiv a\bmod q}}1.
\end{align*}
Following the convention in Conjecture \ref{conj:HL-modified}, we may write
\begin{align*}
\Pi_D(x;r,q,a)
&=2\sum_{\substack{u\bmod q\\Du^2-Dru+\frac{D+1}{4}r^2\equiv a\bmod q}}\pi_{D,-Dr,\frac{D+1}{4}r^2}(x;q,u).
\end{align*}
It requires to verify all prerequisites in Conjecture \ref{conj:HL-modified} for the special quadratic polynomial $f(x)=Dx^2-Drx+\frac{D+1}{4}r^2.$ More precisely, we hope
\begin{align}
&(q,f(u))=1,\label{eq:PiD(x;r,q,a)-condition1}\\
&\Big(D,\frac{D+1}{4}r^2\Big)=1,\label{eq:PiD(x;r,q,a)-condition2}\\
&D(1-r)q\text{~and~}D(1-r)u+\frac{D+1}{4}r^2\text{~are not simultaneously even}.\label{eq:PiD(x;r,q,a)-condition3}
\end{align}
Note that \eqref{eq:PiD(x;r,q,a)-condition1} can be guaranteed by $(q,a)=1.$ Since $D\equiv3\bmod4,$ the condition \eqref{eq:PiD(x;r,q,a)-condition2} holds if and only if $(D,r)=1$. It remains to check \eqref{eq:PiD(x;r,q,a)-condition3}.
If $r$ is odd, then $D(1-r)q$ and $D(1-r)u$ are forced to be even, so that 
\begin{align*}
D(1-r)u+\frac{D+1}{4}r^2\text{~is odd}\Longleftrightarrow\frac{D+1}{4}\text{~is odd}\Longleftrightarrow D\equiv3\bmod8.
\end{align*}
If $r$ is even, \eqref{eq:PiD(x;r,q,a)-condition3} holds if and only if $D(1-r)q$ and $D(1-r)u$ are not simultaneously even, i.e., $q$ and $u$ are not simultaneously even, which can also be guaranteed by $(q,a)=1.$ In summary, we find 
all prerequisites in Conjecture \ref{conj:HL-modified} are satisfied if and only if $\theta(D,r)=1$, in which case 
Conjecture \ref{conj:HL-modified} implies
\begin{align*}
\Pi_D(x;r,q,a)
&\sim \varrho_D(r,q,a)\cdot\frac{2(2,q(r+1))}{q}\frac{\sqrt{D}q^\sharp}{\varphi(Dq^\sharp)}\prod_{p\nmid2qr}\Big(1-\frac{(\frac{-D}{p})}{p-1}\Big)\cdot\frac{\sqrt{x}}{\log x}
\end{align*}
as desired.

\section{From Conjecture \ref{conj:PiD(x;r,q,a)} to Lang--Trotter: $D=2$}

Assume $r\equiv2\bmod4$ and recall \eqref{eq:connection:pi-PiD(x),D=2}. From Conjecture \ref{conj:PiD(x;r,q,a)} it follows that
\begin{align*}
\Pi_D(x;r,16g_1,\delta_{k,a})\sim
\varrho(r,16g_1,\delta_{k,a})\cdot\frac{1}{4\sqrt{2}\varphi(g_1)}\prod_{p\nmid2g_1r}\Big(1-\frac{(\frac{-2}{p})}{p-1}\Big)\cdot\frac{\sqrt{x}}{\log x}.
\end{align*}
From the Chinese remainder theorem, 
\begin{align*}
\varrho(r,16g_1,\delta_{k,a})=\varrho(r,16,2k+1)\cdot\varrho(r,g_1,a).
\end{align*}
For $r\equiv2\bmod4$, we find
\begin{align*}
\varrho(r,16,2k+1)=
\begin{cases}
4,\ \ &\text{if~} 2k+1-r^2/4\equiv0,8\bmod{16},\\ 
8,\ \ &\text{if~} 2k+1-r^2/4\equiv2\bmod{16},\\ 
0,\ \ &\text{otherwise}.
\end{cases}
\end{align*}
Collecting the above evaluations, we derive from \eqref{eq:connection:pi-PiD(x),D=2} that
\begin{align}\label{eq:pi-Sigma1,Simga2,D=2}
\pi_{E,r}(x)
&\sim\frac{1}{4\sqrt{2}\varphi(g_1)}\prod_{p\nmid2g_1r}\Big(1-\frac{(\frac{-2}{p})}{p-1}\Big)\cdot\frac{\sqrt{x}}{\log x}\cdot \{4\Sigma_1+8\Sigma_2\},
\end{align}
where
\begin{align*}
\Sigma_1&=\sum_{\substack{0\leqslant k\leqslant 7\\k\equiv\frac{r^2-4}{8},\frac{r^2+28}{8}\bmod8}}\sum_{\substack{a\bmod{g_1}\\(\frac{a}{g_1})=(-1)^{\delta_k}}}\varrho(r,g_1,a),\\
\Sigma_2&=\sum_{\substack{0\leqslant k\leqslant 7\\ k\equiv\frac{r^2+4}{8}\bmod8}}\sum_{\substack{a\bmod{g_1}\\(\frac{a}{g_1})=(-1)^{\delta_k}}}\varrho(r,g_1,a).
\end{align*}

For each $\epsilon\in\{1,-1\},$ it follows from Lemma \ref{lm:countingresidues} that
\begin{align*}
\sum_{\substack{a\bmod{g_1}\\(\frac{a}{g_1})=\epsilon}}\varrho(r,g_1,a)
&=\sum_{\substack{a\bmod{g_1}\\(\frac{a}{g_1})=\epsilon}}\sum_{\substack{t\bmod{g_1}\\ 2t^2+\frac{r^2}{4}\equiv a\bmod {g_1}}}1
=\sum_{\substack{t\bmod{g_1}\\ (\frac{2t^2+\frac{r^2}{4}}{g_1})=\epsilon}}1\\
&=\frac{\varphi(g_1)}{2}\prod_{\substack{p\mid g_1\\ p\nmid r}}\Big(1-\frac{(\frac{-2}{p})}{p-1}\Big)\cdot\Big\{1+\epsilon\Big(\frac{2}{g_1}\Big)\Omega_2(2;g_1,r)\Big\}.
\end{align*}
Hence we obtain
\begin{align*}
\Sigma_1&=\frac{\varphi(g_1)}{2}\prod_{\substack{p\mid g_1\\ p\nmid r}}\Big(1-\frac{(\frac{-2}{p})}{p-1}\Big)\sum_{\substack{0\leqslant k\leqslant 7\\ k\equiv\frac{r^2-4}{8},\frac{r^2+28}{8}\bmod{8}}}\Big\{1+(-1)^{\delta_k}\Big(\frac{2}{g_1}\Big)\Omega_2(2;g_1,r)\Big\}\\
&=\frac{\varphi(g_1)}{2}\prod_{\substack{p\mid g_1\\ p\nmid r}}\Big(1-\frac{(\frac{-2}{p})}{p-1}\Big)\Big\{2+((-1)^{\eta_1}+(-1)^{\eta_2})\Big(\frac{2}{g_1}\Big)\Omega_2(2;g_1,r)\Big\}
\end{align*}
and
\begin{align*}
\Sigma_2&=\frac{\varphi(g_1)}{2}\prod_{\substack{p\mid g_1\\ p\nmid r}}\Big(1-\frac{(\frac{-2}{p})}{p-1}\Big)\sum_{\substack{0\leqslant k\leqslant 7\\ k\equiv\frac{r^2+4}{8}\bmod{8}}}\Big\{1+(-1)^{\delta_k}\Big(\frac{2}{g_1}\Big)\Omega_2(2;g_1,r)\Big\}\\
&=\frac{\varphi(g_1)}{2}\prod_{\substack{p\mid g_1\\ p\nmid r}}\Big(1-\frac{(\frac{-2}{p})}{p-1}\Big)\Big\{1+(-1)^{\eta_3}\Big(\frac{2}{g_1}\Big)\Omega_2(2;g_1,r)\Big\}
\end{align*}
with
\begin{align*}
\eta_1&=\Big[\frac{r^2}{32}\Big]+\frac{r-2}{4}+\frac{\lambda (r^4-16)}{128}+\frac{(2\delta+g_1-1)(r^2-4)}{16},\\
\eta_2&=\Big[\frac{r^2}{32}\Big]+1+\frac{r-2}{4}+\frac{\lambda (r^2+28)(r^2+36)}{128}+\frac{(2\delta+g_1-1)(r^2+28)}{16},\\
\eta_3&=\Big[\frac{r^2+8}{32}\Big]+\frac{r-2}{4}+\frac{\lambda(r^2+4)(r^2+12)}{128}+\frac{(2\delta+g_1-1)(r^2+4)}{16}.
\end{align*}
Recall the definition \eqref{eq:constant-hD,r} of $h_{2,r}.$ It then follows from \eqref{eq:pi-Sigma1,Simga2,D=2} that
\begin{align}\label{eq:pi-eta1,eta2,eta3,D=2}
\varpi_{E,r}=\frac{h_{2,r}}{2}\cdot \Big(1+\frac{(-1)^{\eta_1}+(-1)^{\eta_2}+2(-1)^{\eta_3}}{4}\Big(\frac{2}{g_1}\Big)\Omega_2(2;g_1,r)\Big).
\end{align}

In view of $2\parallel r,$ we find
\begin{align*}
\frac{r^4-16}{128}\equiv\frac{r^2-4}{16}\equiv\frac{r^2+28}{16}\equiv\frac{(r^2+28)(r^2+36)}{128}\equiv0\bmod2
\end{align*}
and
\begin{align*}
\frac{(r^2+4)(r^2+12)}{128}\equiv\frac{r^2+4}{8}\equiv1\bmod2.
\end{align*}
Moreover, 
\begin{align*}
\Big[\frac{r^2}{32}\Big]=\Big[\frac{r^2+8}{32}\Big]=\frac{r^2-4}{32}.
\end{align*}
Hence
\begin{align*}
\eta_1&\equiv\eta_2+1\equiv\frac{(r-2)(r+10)}{32}\bmod2,\\
\eta_3&\equiv\frac{(r-2)(r+10)}{32}+\delta+\lambda+\frac{g_1-1}{2}\bmod2,
\end{align*}
from which and \eqref{eq:pi-eta1,eta2,eta3,D=2}, Theorem \ref{thm:HL-LT:D=2} follows.

\section{From Conjecture \ref{conj:PiD(x;r,q,a)} to Lang--Trotter: $D\geqslant7$}

Recall \eqref{eq:connection:pi-PiD(x),D>=7}.
From Conjecture \ref{conj:PiD(x;r,q,a)} it follows that
\begin{align*}
\Pi_D(x;r,8A,\delta_{k,a})
&\sim\theta(D,r)\varrho_D(r,8A,\delta_{k,a})\cdot\frac{\sqrt{D}}{2\varphi(DA)}\prod_{p\nmid2g_1r}\Big(1-\frac{(\frac{-D}{p})}{p-1}\Big)\cdot\frac{\sqrt{x}}{\log x}.
\end{align*}
Since $\varrho_D(r,8A,\delta_{k,a})=\varrho_D(r,8,2k+1)\varrho_D(r,A,a)$ by the Chinese remainder theorem, we obtain from \eqref{eq:connection:pi-PiD(x),D>=7} that
\begin{equation}\label{eq:pi-nucounting,D>=7}
\begin{split}
\pi_{E,r}(x)
&\sim\frac{\theta(D,r)\cdot\sqrt{D}}{4\varphi(DA)}\prod_{p\nmid2g_1r}\Big(1-\frac{(\frac{-D}{p})}{p-1}\Big)\cdot\frac{\sqrt{x}}{\log x}\cdot \sum_{0\leqslant k\leqslant 3}\varrho_D(r,8,2k+1)\\
&\ \ \ \ \ \times\sum_{\substack{a\bmod A\\ (\frac{a}{A})=(-1)^{(\delta+\mu+\frac{g_1-1}{2})k+\frac{\lambda k(k+1)}{2}}(\frac{2r}{D})}}\varrho_D(r,A,a).
\end{split}
\end{equation}

Let $\epsilon\in\{1,-1\}.$ Then
\begin{align*}
\sum_{\substack{a\bmod A\\ (\frac{a}{A})=\epsilon}}\varrho_D(r,A,a)
&=\sum_{\substack{t\bmod A\\ (\frac{Dt^2-Drt+\frac{D+1}{4}r^2}{A})=\epsilon}}1
=D^\mu\sum_{\substack{t\bmod{g_1}\\ (\frac{Dt^2-Drt+\frac{D+1}{4}r^2}{g_1})=\epsilon}}1\\
&=\frac{D^\mu\varphi(g_1)}{2}\prod_{\substack{p\mid g_1\\p\nmid r}}\Big(1-\frac{(\frac{-D}{p})}{p-1}\Big)\cdot\Big\{1+\epsilon\Big(\frac{D}{g_1}\Big)\Omega_2(D;g_1,r)\Big\}
\end{align*}
by Lemma \ref{lm:countingresidues}.
We now derive from \eqref{eq:pi-nucounting,D>=7}
that
\begin{equation}\label{eq:pi-nucounting,second,D>=7}
\begin{split}
\pi_{E,r}(x)
&\sim\frac{\fc}{8}\cdot\frac{\sqrt{x}}{\log x}\cdot \sum_{0\leqslant k\leqslant 3}\varrho_D(r,8,2k+1)\\
&\ \ \ \ \ \times\Big\{1+(-1)^{(\delta+\mu+\frac{g_1-1}{2})k+\frac{\lambda k(k+1)}{2}}\Big(\frac{2r}{D}\Big)\Big(\frac{D}{g_1}\Big)\Omega_2(D;g_1,r)\Big\}
\end{split}
\end{equation}
with
\begin{align*}
\fc=\frac{\theta(D,r)\cdot \sqrt{D}}{\varphi(D)}\prod_{p\nmid2r}\Big(1-\frac{(\frac{-D}{p})}{p-1}\Big).
\end{align*}

It requires to evaluate $\varrho_D(r,8,2k+1)$ explicitly. We claim that
\begin{align}\label{eq:nuD(r,8,2k+1)-evaluation}
\varrho_D(r,8,2k+1)
&=\begin{cases}
2,\ \ &\text{if~}D\equiv3\bmod8\text{~and~}2\nmid r,\\
2,&\text{if~}k\equiv \frac{r^2-4}{8},\frac{r^2+12}{8}\bmod4\text{~and~}2\parallel r,\\
4, &\text{if~}k\equiv \frac{r^2}{8}+\frac{D-1}{2}\bmod4\text{~and~}4\mid r,\\
0, &\text{otherwise}.
\end{cases}
\end{align}
For $D\equiv3\bmod8$ and $2\nmid r,$ we derive from \eqref{eq:pi-nucounting,second,D>=7} and \eqref{eq:nuD(r,8,2k+1)-evaluation} that
\begin{align*}
\varpi_{E,r}=\fc\cdot\Big(1+\frac{(1+(-1)^{\lambda})(1+(-1)^{\delta+\mu+\frac{g_1-1}{2}})}{4}\Big(\frac{2r}{D}\Big)\Big(\frac{D}{g_1}\Big)\Omega_2(D;g_1,r)\Big).
\end{align*}
Similarly, we may also obtain
\begin{align*}
\varpi_{E,r}=\frac{\fc}{2}\cdot\Big(1+\frac{1+(-1)^{\lambda}}{2}\Big(\frac{2r}{D}\Big)\Big(\frac{D}{g_1}\Big)\Omega_2(D;g_1,r)\Big)
\end{align*}
for $2\parallel r,$ and
\begin{align*}
\varpi_{E,r}=\frac{\fc}{2}\cdot\Big(1+(-1)^{\delta+\mu+\frac{g_1-1}{2}+\frac{\lambda r}{4}}\Big(\frac{2^{\lambda+1}r}{D}\Big)\Big(\frac{D}{g_1}\Big)\Omega_2(D;g_1,r)\Big)
\end{align*}
for $4\mid r.$

It remains to prove the claim \eqref{eq:nuD(r,8,2k+1)-evaluation}. Recall that $\varrho_D(r,8,2k+1)$ counts solutions to
\begin{align*}
Dt^2-Drt+\frac{D+1}{4}r^2\equiv 2k+1\bmod 8
\end{align*}
in $t\bmod8.$ Suppose $2\nmid r$ and make the change of variable $t\rightarrow rt,$ we then find
\begin{align*}
\varrho_D(r,8,2k+1)=|\{t\bmod8:Dt^2-Dt+\frac{D+1}{4}\equiv 2k+1\bmod 8\}|.
\end{align*}
Note that $t^2-t\equiv a\bmod8$ has exactly two solutions in $t\bmod8$ for any even $a$, and has no solutions if $a$ is odd. Therefore, $\varrho_D(r,8,2k+1)$ is equal to $2$ if $D\equiv3\bmod8$, and vanishes otherwise.

We now assume $2\mid r$. Making the change of variable $t\rightarrow t+r/2,$ we find
\begin{align*}
\varrho_D(r,8,2k+1)=|\{t\bmod 8:Dt^2\equiv 2k+1-\frac{r^2}{4}\bmod 8\}|.
\end{align*}
Note that $t^2\equiv a\bmod8$ has exactly four solutions if $a\equiv1\bmod8$, and has two solutions if $a\equiv0,4\bmod8$, and has no solutions otherwise.
Therefore,
\begin{align*}
\varrho_D(r,8,2k+1)
&=\begin{cases}
2,\ \ &\text{if~}2k+1-r^2/4\equiv 0,4D\bmod8\text{~and~}2\parallel r,\\
4, &\text{if~}2k+1-r^2/4\equiv D\bmod8\text{~and~}4\mid r,\\
0, &\text{otherwise}.
\end{cases}
\end{align*}
This proves \eqref{eq:nuD(r,8,2k+1)-evaluation}, and thus proves Theorem \ref{thm:HL-LT:D>=7} by noting that $(\frac{D}{g_1})=(\frac{g_1}{D})(-1)^{\frac{g_1-1}{2}}$ since $D\equiv3\bmod4.$

%% file: 7.LTconstant.tex
\chapter{Further interpretations on the Lang--Trotter constant}\label{ch:LTconstant}

Put $$\fD:=\{1,2,3,7,11,19, 43, 67, 163\}.$$
For $D\in\fD$, we assume $E$ has CM by $\Q(\sqrt{-D})$ as defined by \eqref{eq:E:D=1}, \eqref{eq:E:D=3} and \eqref{eq:E:otherD}, respectively, as in the previous three chapters. Decompose the involved parameter $g$ by
\begin{align}\label{eq:g-factorization}
g=(-1)^\delta2^\lambda D^\mu g_1
\end{align}
with
\begin{itemize}
\item $\delta\in\{0,1\},$
\item $\lambda,\mu\in\N,$
\item $g_1\in\Z^+$ with $(2D,g_1)=1$,
\end{itemize}
where we always assume $\mu=0$ for $D=2.$ This convention follows Chapters \ref{ch:D=1}, \ref{ch:D=3} and \ref{ch:otherCM} exactly.

Recall that $\varpi_{E,r}\geqslant0$ is an explicit constant defined by \eqref{eq:varpiE,r,D=1}, \eqref{eq:varpiE,r,D=3}, \eqref{eq:varpiE,r,D=2} and \eqref{eq:varpiE,r,otherD} for different $D\in\fD$ and $r$ satisfying suitable congruence conditions. To extend the definition domain to all $r\in\Z,$ we put $\varpi_{E,r}=0$ if 
\begin{itemize}
\item $2\nmid r$ for $D=1;$
\item $3\mid r$ for $D=3;$
\item $r\not\equiv2\bmod4$ for $D=2;$
\item $D\mid r$ for $D\in\{7,11,19, 43, 67, 163\}.$
\end{itemize}
We have proven, under the Hardy--Littlewood conjecture\index{Hardy--Littlewood conjecture}, that
\begin{align*}
\pi_{E,r}(x)\sim\frac{\varpi_{E,r}\cdot\sqrt{x}}{\log x}
\end{align*}
as $x\rightarrow+\infty.$
Assuming the Hardy--Littlewood conjecture (Conjecture \ref{conj:HL}) and the original Lang--Trotter conjecture (Conjecture \ref{conj:LT}) simultaneously, the constant $\varpi_{E,r}$ should coincide with $c_{E,r}$. The former one 
is defined analytically and completely explicit, while the latter one is defined algebraically (see \eqref{eq:LTconstant-Galois}) and much harder to compute. This motivates us to formulate the following constant comparison conjecture \index{comparison conjecture} (Conjecture \ref{conj:comparison}).
\begin{conjecture}\label{conj:cE,r=varpiE,r}\index{Lang--Trotter constant}
Suppose $E/\Q$ is an elliptic curve defined by $\eqref{eq:E:D=1},$ $\eqref{eq:E:D=3}$ or $\eqref{eq:E:otherD}$. For each fixed non-zero integer $r,$ we have
\begin{align*}
c_{E,r}=\varpi_{E,r}.
\end{align*}
\end{conjecture}

The significance of this conjecture is that it gives a
simple explicit way to compute the Lang--Trotter constant following suggestions from the Hardy--Littlewood conjecture.
This chapter aims to explore more features of $\varpi_{E,r}$, which is expected to  be helpful to understand 
$c_{E,r}.$ For instance, we would like to determine when $\varpi_{E,r}$ is positive or not, and when it does not depend on the sign of $r$. The subsequent discussions will proceed case by case according to Theorems \ref{thm:HL-LT:D=1}, \ref{thm:HL-LT:D=3}, \ref{thm:HL-LT:D=2} and \ref{thm:HL-LT:D>=7} for $D\in\fD.$

Recall the arithmetic function
\begin{align*}
\Omega_j(D;q,r)=\prod_{\substack{p^\nu\parallel q,~ p\nmid r\\j\nmid \nu}}\frac{-1}{p-1-(\frac{-D}{p})}
\end{align*}
as defined by \eqref{eq:Omegafunction} for all $D\in\fD,q,r\in\Z$ and $j\in\Z^+$ with $2\nmid q.$

For $j\in\Z^+$, define an arithmetic function $\fA:\Z^+\rightarrow\Z^+:$
\begin{align*}
\fA_j(q)=\prod_{\substack{p^\nu\parallel q\\j\nmid \nu}}p.
\end{align*}
For each fixed positive integer $q$, it is clear that $\fA_i(q)\mid \fA_j(q)$ if $i\mid j.$ The condition $\fA_j(q)\mid r$ implies $\Omega_j(D;q,r)=1.$ Moreover, $\fA_j(q)=1$ implies that $q$ is a $j$-th power, including $1$.

\section{Positivity}

Define an auxiliary set 
\begin{align*}
\cJ_p(g,r):=\{j\in\Z^+:j\nmid \ord_p(g_1),\ \ p\nmid r,\ \ \fA_j(g_1)\mid pr\}
\end{align*}
with $p\in\cP,$ and $g,r\in\Z^*$, where $g_1$ is defined by \eqref{eq:g-factorization}. Clearly, if $j_1,j_2\in \cJ_p(g,r)$ for $(j_1,j_2)=1,$ then $j_1j_2\in \cJ_p(g,r)$.
The set of $\cJ_p(g,r)$ characterizes the order of prime factor $p$ of $g$. For instance, if $j\in\cJ_p(g,r)$ with $r=\pm1$ or $\pm2$, then $j\nmid \ord_p(g_1)\geqslant1,$ which yields $p\mid g_1$, so that $p\nmid r$ holds automatically. Moreover, the condition $\fA_j(g_1)\mid pr$ reduces to $\fA_j(g_1)\mid p,$ which means $p$ is the only prime factor $\ell$ of $g_1$ such that $j\nmid \ord_\ell(g_1).$ Hence we may write $g_1=pa^j$ for some $a\in\Z$ with $p\nmid a.$
When $r$ is given differently from $\pm1,\pm2,$ the condition $j\in\cJ_p(g,r)$ can also be translated to the exact information on prime decomposition of $g_1$, although a bit complicated.

We are now ready to state the precise results on positivity or vanishing of $\varpi_{E,r}$ for all CM elliptic curves.

\begin{theorem}\label{thm:LTconstant:D=1}\index{Lang--Trotter constant}
Suppose $r$ is a fixed non-zero integer with $2\mid r,$ and $E/\Q$ is the elliptic curve defined by $\eqref{eq:E:D=1}$ with CM by $\Q(\sqrt{-1})$.
Then $\varpi_{E,r}$ vanishes if and only if one of the following $10$ conditions holds:
\smallskip

$\mathrm{(I)}$ $\fA_2(g_1)\mid r:$

\begin{itemize}
\item $r\equiv0\bmod8;$\smallskip
\item $r\equiv4\bmod8,$ $2\mid\lambda.$\smallskip
\end{itemize}

$\mathrm{(II)}$ $\fA_4(g_1)\mid r:$

\begin{itemize}
\item $r\equiv4\bmod{16},$ $\lambda\equiv1\bmod4,$ $2\mid \delta+\frac{g_1-1}{2};$\smallskip
\item $r\equiv12\bmod{16},$ $\lambda\equiv3\bmod4,$ $2\mid \delta+\frac{g_1-1}{2};$\smallskip 
\item $r\equiv4\bmod{16},$ $\lambda\equiv3\bmod4,$ $2\nmid \delta+\frac{g_1-1}{2};$\smallskip
\item $r\equiv12\bmod{16},$ $\lambda\equiv1\bmod4,$ $2\nmid \delta+\frac{g_1-1}{2};$\smallskip
\item $r\equiv2\bmod8,$ $2\mid\lambda,$ $2\mid \delta+\frac{\lambda}{2},$ $g_1\equiv3\bmod8;$ \smallskip 
\item $r\equiv2\bmod8,$ $2\mid\lambda,$ $2\nmid \delta+\frac{\lambda}{2},$ $g_1\equiv5\bmod8;$ \smallskip
\item $r\equiv6\bmod8,$ $2\mid\lambda,$ $2\nmid \delta+\frac{\lambda}{2},$ $g_1\equiv1\bmod8;$ \smallskip
\item $r\equiv6\bmod8,$ $2\mid\lambda,$ $2\mid \delta+\frac{\lambda}{2},$ $g_1\equiv7\bmod8.$
\smallskip
\end{itemize}

$\mathrm{(III)}$ $2\in \cJ_3(g,r)\cup\cJ_5(g,r):$

\begin{itemize}
\item $r\equiv4\bmod{16},$ $\lambda\equiv1\bmod4,$ $2\nmid \delta+\frac{g_1-1}{2};$\smallskip
\item $r\equiv12\bmod{16},$ $\lambda\equiv3\bmod4,$ $2\nmid \delta+\frac{g_1-1}{2};$\smallskip 
\item $r\equiv4\bmod{16},$ $\lambda\equiv3\bmod4,$ $2\mid \delta+\frac{g_1-1}{2};$\smallskip
\item $r\equiv12\bmod{16},$ $\lambda\equiv1\bmod4,$ $2\mid \delta+\frac{g_1-1}{2};$\smallskip
\item $r\equiv2\bmod8,$ $2\mid\lambda,$ $2\nmid \delta+\frac{\lambda}{2},$ $g_1\equiv1\bmod8;$ \smallskip 
\item $r\equiv2\bmod8,$ $2\mid\lambda,$ $2\mid \delta+\frac{\lambda}{2},$ $g_1\equiv7\bmod8;$ \smallskip
\item $r\equiv6\bmod8,$ $2\mid\lambda,$ $2\mid \delta+\frac{\lambda}{2},$ $g_1\equiv3\bmod8;$ \smallskip
\item $r\equiv6\bmod8,$ $2\mid\lambda,$ $2\nmid \delta+\frac{\lambda}{2},$ $g_1\equiv5\bmod8.$
\end{itemize}
\end{theorem}

\begin{theorem}\label{thm:LTconstant:D=3}\index{Lang--Trotter constant}
Suppose $r$ is a fixed non-zero integer with $3\nmid r,$ $E/\Q$ is the elliptic curve defined by $\eqref{eq:E:D=3}$ with CM by $\Q(\sqrt{-3})$.
Then $\varpi_{E,r}$ vanishes if and only if one of the following $18$ conditions holds:

\smallskip
$\mathrm{(I)}$ $\{2,3,6\}\subseteq\cJ_{5}(g,r)\cup\cJ_{7}(g,r):$
\begin{itemize}
\item $r\equiv8\bmod{24},$ $3\mid(\mu,\lambda+\frac{g_1^2-1}{3}),$ $(-1)^{\delta+\lambda+\mu}=(\frac{g_1}{3});$ \smallskip
\item $r\equiv16\bmod{24},$ $3\mid(\mu,\lambda+\frac{g_1^2-1}{3}),$ $(-1)^{\delta+\lambda+\mu}=(\frac{-g_1}{3});$\smallskip
\item $r\equiv20\bmod{24},$ $3\mid(\mu,\lambda+\frac{g_1^2-1}{3}),$ $(-1)^{\delta+\mu}=(\frac{g_1}{3});$\smallskip
\item $r\equiv4\bmod{24},$ $3\mid(\mu,\lambda+\frac{g_1^2-1}{3}),$ $(-1)^{\delta+\mu}=(\frac{-g_1}{3});$ \smallskip
\item $r\equiv2\bmod{12},$ $g_1\equiv1,11\bmod{12},$ $3\mid(\mu,\lambda+\frac{g_1^2-1}{3}),$ $2\mid\lambda;$ \smallskip
\item $r\equiv10\bmod{12},$ $g_1\equiv5,7\bmod{12},$ $3\mid(\mu,\lambda+\frac{g_1^2-1}{3}),$ $2\mid\lambda;$ \smallskip
\item $r\equiv5\bmod{6},$ $ g_1\equiv1,11\bmod{12},$ $3\mid(\mu,2+\lambda,\frac{g_1^2+2}{3}),$ $2\mid(\lambda,\delta+\mu+\frac{g_1-1}{2});$ \smallskip
\item $r\equiv1\bmod{6},$ $ g_1\equiv5,7\bmod{12},$ $3\mid(\mu,2+\lambda,\frac{g_1^2+2}{3}),$ $2\mid(\lambda,\delta+\mu+\frac{g_1-1}{2}).$\smallskip
\end{itemize}

$\mathrm{(II)}$ $\fA_2(g_1)\mid r:$
\begin{itemize}
\item $r\equiv8\bmod{24},$ $(-1)^{\delta+\lambda+\mu}=(\frac{-g_1}{3});$ \smallskip
\item $r\equiv16\bmod{24},$ $(-1)^{\delta+\lambda+\mu}=(\frac{g_1}{3});$ \smallskip
\item $r\equiv20\bmod{24},$ $(-1)^{\delta+\mu}=(\frac{-g_1}{3});$ \smallskip
\item $r\equiv4\bmod{24},$ $(-1)^{\delta+\mu}=(\frac{g_1}{3});$ \smallskip
\item $r\equiv2\bmod{12},$ $g_1\equiv5,7\bmod{12},$ $2\mid\lambda;$ \smallskip
\item $r\equiv10\bmod{12},$ $g_1\equiv1,11\bmod{12},$ $2\mid\lambda;$ \smallskip
\item $r\equiv5\bmod{6},$ $g_1\equiv5,7\bmod{12},$ $2\mid(\lambda,\delta+\mu+\frac{g_1-1}{2});$\smallskip
\item $r\equiv1\bmod{6},$ $g_1\equiv1,11\bmod{12},$ $2\mid(\lambda,\delta+\mu+\frac{g_1-1}{2}).$ \smallskip
\end{itemize}

$\mathrm{(III)}$ $\fA_3(g_1)\mid r:$
\begin{itemize}
\item $r\equiv\pm1\bmod{6},$ $g_1\equiv1,8\bmod9,$ $3\mid(\lambda,\mu);$ \smallskip
\item $r\equiv\pm1\bmod{6},$ $g_1\equiv2,7\bmod9,$ $3\mid(\lambda+1,\mu).$
\end{itemize}
\end{theorem}

\begin{theorem}\label{thm:LTconstant:D=2}\index{Lang--Trotter constant}
Suppose $r$ is a fixed integer with $r\equiv2\bmod4,$ and $E/\Q$ is the elliptic curve defined by $\eqref{eq:E:otherD}$ with CM by $\Q(\sqrt{-2})$.
Then $\varpi_{E,r}$ is always positive.
\end{theorem}

\begin{theorem}\label{thm:LTconstant:D>=7}\index{Lang--Trotter constant}
Let $D\in\{7,11, 19, 43, 67, 163\}$ and $D\nmid r.$ Suppose $E/\Q$ is the elliptic curve defined by $\eqref{eq:E:otherD}$ with CM by $\Q(\sqrt{-D})$.
Then $\varpi_{E,r}$ vanishes if and only if one of the following $9$ conditions holds:
\begin{itemize}
\item $r\equiv2\bmod4,$ $(-1)^{\frac{g_1+1}{2}}=(\frac{2g_1r}{D}),$ $2\mid \lambda,$ $\fA_2(g_1)\mid r;$\smallskip
\item $r\equiv0\bmod4,$ $(-1)^{\delta+\mu}=-(\frac{2g_1r}{D}),$ $2\mid \lambda,$ $\fA_2(g_1)\mid r;$\smallskip
\item $r\equiv0\bmod4,$ $(-1)^{\delta+\mu+\frac{r}{4}}=-(\frac{g_1r}{D}),$ $2\nmid \lambda,$ $\fA_2(g_1)\mid r;$\smallskip
\item $2\nmid r,$ $(-1)^{\frac{g_1+1}{2}}=(\frac{2g_1r}{D}),$ $2\mid (\lambda,\delta+\mu+\frac{g_1-1}{2}),$ $\fA_2(g_1)\mid r;$\smallskip
\item $r\equiv2\bmod4,$ $(-1)^{\frac{g_1-1}{2}}=(\frac{2g_1r}{D}),$ $2\mid \lambda,$ $2\in\cJ_3(g,r),$ $D=11;$\smallskip
\item $r\equiv0\bmod4,$ $(-1)^{\delta+\mu}=(\frac{2g_1r}{D}),$ $2\mid \lambda,$ $2\in\cJ_3(g,r),$ $D=11;$\smallskip
\item $r\equiv0\bmod4,$ $(-1)^{\delta+\mu+\frac{r}{4}}=(\frac{g_1r}{D}),$ $2\nmid \lambda,$ $2\in\cJ_3(g,r),$ $D=11;$\smallskip
\item $2\nmid r,$ $(-1)^{\frac{g_1-1}{2}}=(\frac{2g_1r}{D}),$ $2\mid (\lambda,\delta+\mu+\frac{g_1-1}{2}),$ $2\in\cJ_3(g,r),$ $D=11;$\smallskip
\item $2\nmid r,$ $D=7.$
\end{itemize}
\end{theorem}

As an immediate consequence of the above criterion, we may obtain the positivity of $\varpi_{E,r}$ for $r=2.$
\begin{corollary}\label{coro:LTconstantpositive,r=2}\index{Lang--Trotter conjecture}
Suppose $E/\Q$ is the elliptic curve with CM as defined by $\eqref{eq:E:D=1},$ $\eqref{eq:E:D=3}$ or $\eqref{eq:E:otherD}.$
Then the Hardy--Littlewood conjecture \index{Hardy--Littlewood conjecture} implies
\begin{align*}
\pi_{E,2}(x)\gg \frac{\sqrt{x}}{\log x}.
\end{align*}
\end{corollary}

Corollary \ref{coro:LTconstantpositive,r=2} coincides with the definition \eqref{eq:LTconstant-Galois}, from which one may see that the identity element, which is of trace $2$, is always contained in $\Gal(K(E[m_E])/K)_r.$
The proof of Corollary \ref{coro:LTconstantpositive,r=2} can proceed according to $D\in\fD$ one by one. Write $\varpi_E=\varpi_{E,2}$
for brevity.

\begin{proof}[Proof of Corollary $\ref{coro:LTconstantpositive,r=2}$]
For $D=1$, $\varpi_E$ cannot vanish unless $\fA_4(g_1)\mid 2,$ which implies $\fA_4(g_1)=1,$ i.e., $g_1$ must be a fourth power in $\Z^+.$ This, however, is not compatible with $g_1\equiv3,5\bmod8.$

For $D=3$, $\varpi_E$ cannot vanish unless 
\begin{align*}
2\in \cJ_{5}(g,2)\cup\cJ_{7}(g,2),\ \ g_1\equiv\pm1\bmod{12},
\end{align*}
or
\begin{align*}
\fA_2(g_1)\mid 2, \ \ g_1\equiv\pm5\bmod{12}.
\end{align*}
The second possibility cannot happen due to the same reason as above. On the other hand, if $2\in \cJ_{5}(g,2)$, then $g_1=5a^2$ for some $a\in\Z$ with $(a,6)=1.$ This implies $g_1\equiv5\bmod{12}.$ Similarly, the condition $2\in \cJ_{7}(g,2)$ implies $g_1\equiv7\bmod{12}.$ Each of the above congruences contradicts to $g_1\equiv\pm1\bmod{12}$ as required. On the other hand, the condition $\fA_2(g_1)\mid 2$ yields $\fA_2(g_1)=1,$ so that $g_1$ is pushed to be a square integer, which contradicts to $g_1\equiv\pm5\bmod{12}$ as required. 

The case $D=2$ is obvious.

We now consider the remaining cases $D\in\{7,11, 19, 43, 67, 163\}.$ Note that 
$\varpi_E$ cannot vanish unless 
\begin{align*}
\fA_2(g_1)\mid 2, \ \ (-1)^{\frac{g_1+1}{2}}=\Big(\frac{g_1}{D}\Big)
\end{align*}
or
\begin{align*}
2\in \cJ_{3}(g,2),\ \ (-1)^{\frac{g_1-1}{2}}=\Big(\frac{g_1}{D}\Big),\ \ D=11.
\end{align*}
The condition  $\fA_2(g_1)\mid 2$ guarantees that $g_1$ is an odd square in $\Z^+,$ thus we have $g_1\equiv1\bmod4$, which implies
$(-1)^{\frac{g_1+1}{2}}=-1$ and $(\frac{g_1}{D})=1.$
On the other hand, $2\in \cJ_{3}(g,2)$ implies that $g_1=3h^2$ for some $h\in\Z$ with $11\nmid h.$ Hence
$g_1\equiv3\bmod4$ and $ (-1)^{\frac{g_1-1}{2}}=-1.$ However, $(\frac{g_1}{11})=(\frac{3}{11})=1.$
\end{proof}

\subsection{Preliminary analysis of $\Omega_j(D;q,r)$}

The proofs of Theorems \ref{thm:LTconstant:D=1}, \ref{thm:LTconstant:D=3}, \ref{thm:LTconstant:D=2} and \ref{thm:LTconstant:D>=7} involve computations of the arithmetic function $\Omega_j(D;q,r)$. Keep in mind the conditions
$D\in\fD,$ $q,r\in\Z$, $j\in\Z^+$ and $2\nmid q.$  We have $\Omega_j(D;q,r)=1$ if the involved product is empty, i.e., $\fA_j(q)\mid r$. We now assume $\fA_j(q)\nmid r$, and would like to show that $\Omega_j(D;q,r)$ is relatively small in most cases. To do so, write
\begin{align}\label{eq:sigmaDp}
\sigma_D(p):=p-1-\Big(\frac{-D}{p}\Big)
\end{align}
for $D\in\fD$ and $p\nmid 2D.$

For $D=1$, we find $(\frac{-D}{p})=-1$ if $p=3$, and $(\frac{-D}{p})=1$ if $p=5$. Hence $\sigma_D(p)=3$ for $p=3,5,$ and $\sigma_D(p)\geqslant5$ for any $p>5$.

For $D=2,$ we find $(\frac{-D}{p})=1$ if $p=3$, in which case $\sigma_D(p)=1$, and we also have 
$\sigma_D(p)\geqslant5$ for any $p>3$.

For $D=3,$ we find $(\frac{-D}{p})=-1$ if $p=5$, and $(\frac{-D}{p})=1$ if $p=7$. Hence 
$\sigma_D(p)=5$ for $p=5,7$, and $\sigma_D(p)\geqslant11$ for any $p>7$.

For the remaining $D\in\fD,$ we always have $\sigma_D(p)\geqslant3$ for any $p>3.$ It remains to treat the case $p=3$, which depends on the residue class $D\bmod3.$ More precisely, $\sigma_D(3)=1$ for $D=11$, and $\sigma_D(3)=3$ for $D\in\{7,19,43,67,163\}$.

In summary, we obtain
\begin{align*}
\sigma_D(p)
\begin{cases}
=3,\ \ &\text{if~}D=1\text{~and~}p=3,5,\\
=1,\ \ &\text{if~}D=2\text{~and~}p=3,\\
=5,\ \ &\text{if~}D=3\text{~and~}p=5,7,\\
=1,\ \ &\text{if~}D=11\text{~and~}p=3,\\
\geqslant5, &\text{if~}D=1\text{~and~}p>5,\\
\geqslant5, &\text{if~}D=2\text{~and~}p>3,\\
\geqslant11, &\text{if~}D=3\text{~and~}p>7,\\
\geqslant3, &\text{if~}D=11\text{~and~}p>3,\\
\geqslant3, &\text{if~}D\in\{7,19,43,67,163\}\text{~and~}p>2,\\
\end{cases}
\end{align*}
from which we may obtain some auxiliary bounds for $\Omega_j(D;q,r)$.
For instance, when $D=1$, we find $|\Omega_j(D;q,r)|\leqslant\frac{1}{3}$ if and only if $\fA_j(q)\nmid r$. 
For $D=2$, the crude bound $|\Omega_j(D;q,r)|\leqslant1$ will be sufficient. For $D=3$, $\Omega_j(D;q,r)=1$ can only happen when $\fA_j(q)\mid r,$ and otherwise $-\frac{1}{5}\leqslant\Omega_j(D;q,r)<\frac{1}{5};$ moreover, $\Omega_j(D;q,r)=-\frac{1}{5}$  if $5$ or $7$ appears in the factorization of $q$ to a suitable power.

For all positive integers $i,j$ with $(i,j)=1$, one has
\begin{align*}
\Omega_{ij}(D;q,r)=\Omega_{j}(D;q,r)
\end{align*}
if $\Omega_{i}(D;q,r)=1.$

\subsection{CM by $\Q(\sqrt{-1})$}
As in Theorem \ref{thm:HL-LT:D=1}, recall that
\begin{align*}
\varpi_{E,r}=\frac{h_{1,r}}{4}\cdot\kappa
\end{align*}
for $2\mid r,$
where (writing $\Omega_j=\Omega_j(1;g_1,r)$)
\begin{align*}
\kappa
&=\begin{cases}
1-(-1)^{\frac{\lambda r}{4}}\Omega_2+(-1)^{\frac{r}{4}(\delta+\frac{g_1-1}{2})}(1-(-1)^{\frac{\lambda r}{4}})\Re(i^{1+\frac{\lambda r}{4}})\Omega_4,\ \ &4\mid r,\\
1+\Omega_2+(-1)^{\frac{r-2}{4}+\frac{g_1^2-1}{8}}(1-(-1)^{\delta+\frac{\lambda+g_1-1}{2}})\Omega_4,\ \ &2\parallel r,2\mid\lambda,\\
1,\ \ &2\parallel r,2\nmid\lambda.
\end{cases}
\end{align*}
Since $h_{1,r}>0$ for all even $r\neq0,$ it suffices to give a condition equivalent to $\kappa=0.$

(I) Suppose $4\mid r$. Then
\begin{align*}
\kappa=0
&\Longleftrightarrow(-1)^{\frac{\lambda r}{4}}\Omega_2-(-1)^{\frac{r}{4}(\delta+\frac{g_1-1}{2})}(1-(-1)^{\frac{\lambda r}{4}})\Re(i^{1+\frac{\lambda r}{4}})\Omega_4=1.
\end{align*}
This leads to three possibilities:
\begin{itemize}
\item $\Omega_2=1,$ $8\mid \lambda r;$
\item $\Omega_2=\Omega_4=1,$ $8\nmid \lambda r,$ $(-1)^{\frac{r}{4}(\delta+\frac{g_1-1}{2})}\Re(i^{1+\frac{\lambda r}{4}})=-1;$
\item $\Omega_2=\Omega_4=-\frac{1}{3},$ $8\nmid \lambda r,$ $(-1)^{\frac{r}{4}(\delta+\frac{g_1-1}{2})}\Re(i^{1+\frac{\lambda r}{4}})=1.$
\end{itemize}

(II) Suppose $2\parallel r,2\mid\lambda$. Then $\kappa=0$ holds if and only if
\begin{align*}
\Omega_2+(-1)^{\frac{r-2}{4}+\frac{g_1^2-1}{8}}(1-(-1)^{\delta+\frac{\lambda+g_1-1}{2}})\Omega_4=-1.
\end{align*}
This leads to two possibilities:
\begin{itemize}
\item $\Omega_2=\Omega_4=1,$ $2\nmid(\frac{r-2}{4}+\frac{g_1^2-1}{8})(\delta+\frac{\lambda+g_1-1}{2});$
\item $\Omega_2=\Omega_4=-\frac{1}{3},$ $2\mid(\frac{r-2}{4}+\frac{g_1^2-1}{8}),$ $2\nmid (\delta+\frac{\lambda+g_1-1}{2}).$
\end{itemize}

Note that the condition $\Omega_2=\Omega_4=1$ is equivalent to $\fA_4(g_1)\mid r.$
On the other hand, $\Omega_2=\Omega_4=-\frac{1}{3}$ holds if and only if $2\in \cJ_3(g,r)\cup\cJ_5(g,r).$
Now Theorem \ref{thm:LTconstant:D=1} follows by analyzing the congruence classes of $r,\delta,g_1$ modulo powers of $2$.

\subsection{CM by $\Q(\sqrt{-3})$}
As in Theorem \ref{thm:HL-LT:D=3}, we have
\begin{align*}
\varpi_{E,r}=\frac{h_{3,r}}{6}\Big(1+\varsigma_2\Big(\frac{3}{g_1}\Big)\Omega_2+\widehat{\mu}\varsigma_1\Big(\Omega_3+\varsigma_2\Big(\frac{3}{g_1}\Big)\Omega_6\Big)\Big)
\end{align*}
for $3\nmid r$, 
where we write $\Omega_j=\Omega_j(3;g_1,r)$, $h_{3,r}\geqslant0$ is defined as in \eqref{eq:constant-hD,r},
\begin{align}\label{eq:varsigma1}
\varsigma_1&=\begin{cases}
1+2\Re(\omega^{\lambda+\frac{g_1^2-1}{3}}),\ \ &\text{for~} r\equiv\pm2\bmod6,\\ 
\Re(\omega^{1-\lambda}+\omega^{1+\frac{g_1^2-1}{3}}+\omega^{2-\lambda+\frac{g_1^2-1}{3}}),\ \ &\text{for~} r\equiv\pm5\bmod{6},\\
\end{cases}
\end{align}
and
\begin{align}\label{eq:varsigma2}
\varsigma_2&=\begin{cases}
\pm(-1)^{\delta+\lambda+\mu+\frac{g_1-1}{2}},\ \ &\text{for~} r\equiv\pm8\bmod{24},\\ 
\pm(-1)^{\delta+\mu+\frac{g_1-1}{2}},\ \ &\text{for~} r\equiv\pm20\bmod{24},\\ 
\pm\frac{1+(-1)^{\lambda}}{2},\ \ &\text{for~} r\equiv\pm2\bmod{12},\\ 
\pm\frac{(1+(-1)^\lambda)(1+(-1)^{\delta+\mu+\frac{g_1-1}{2}})}{4},\ \ &\text{for~} r\equiv\pm5\bmod{6}.
\end{cases}
\end{align}

Note that \begin{align*}
\varsigma_1\in\Big\{-\frac{3}{2},0,\frac{3}{2},3\Big\},\ \ \ \varsigma_2\in\{-1,0,1\},
\end{align*}
and $h_{3,r}\neq0$ since $3\nmid r.$ Therefore, $\varpi_{E,r}=0$ if and only if
\begin{align}\label{eq:LTvanish,D=3,1}
\varsigma_2\Big(\frac{3}{g_1}\Big)\Omega_2+\widehat{\mu}\varsigma_1\Big(\Omega_3+\varsigma_2\Big(\frac{3}{g_1}\Big)\Omega_6\Big)=-1.
\end{align}

If $\widehat{\mu}\varsigma_1=0,$ then
\eqref{eq:LTvanish,D=3,1} holds if and only if 
\begin{align*}
\varsigma_2=-\Big(\frac{3}{g_1}\Big),\ \ \Omega_2=1.
\end{align*}
We henceforth assume $\widehat{\mu}\varsigma_1\neq0,$ i.e., $3\mid \mu$ and $\varsigma_1\neq0,$ in which case \eqref{eq:LTvanish,D=3,1} becomes
\begin{align}\label{eq:LTvanish,D=3,2}
\varsigma_2\Big(\frac{3}{g_1}\Big)\Omega_2+\frac{2}{3}\varsigma_1\Big(\Omega_3+\varsigma_2\Big(\frac{3}{g_1}\Big)\Omega_6\Big)=-1.
\end{align}

(I) Suppose $\Omega_j\neq1$ for each $j=2,3,6.$ We then find $|\Omega_j|\leqslant\frac{1}{5}.$
Thus
\begin{align*}
\Big|\varsigma_2\Big(\frac{3}{g_1}\Big)\Omega_2\Big|\leqslant\frac{1}{5}
\end{align*}
and
\begin{align*}
\Big|\frac{2}{3}\varsigma_1\Big(\Omega_3+\varsigma_2\Big(\frac{3}{g_1}\Big)\Omega_6\Big)\Big|\leqslant\frac{4}{5}.
\end{align*}
Therefore, \eqref{eq:LTvanish,D=3,2} holds if and only if 
\begin{align*}
\varsigma_2\Big(\frac{3}{g_1}\Big)\Omega_2=-\frac{1}{5},\ \ \ \varsigma_1\Big(\Omega_3+\varsigma_2\Big(\frac{3}{g_1}\Big)\Omega_6\Big)=-\frac{6}{5},
\end{align*}
i.e.,
\begin{align*}
\varsigma_1=3,\ \ \varsigma_2\Big(\frac{3}{g_1}\Big)=1,\ \ \Omega_2=\Omega_3=\Omega_6=-\frac{1}{5}.
\end{align*}

(II) We now suppose $\Omega_2=1,$ so that $\Omega_3=\Omega_6,$ and \eqref{eq:LTvanish,D=3,2} becomes
\begin{align*}
\Big(1+\varsigma_2\Big(\frac{3}{g_1}\Big)\Big)\Big(1+\frac{2}{3}\varsigma_1\Omega_3\Big)=0.
\end{align*}
This holds if and only if one of the following two conditions holds:
\begin{itemize}
\item $\varsigma_2(\frac{3}{g_1})=-1;$
\item $\varsigma_1=-\frac{3}{2},$ $\Omega_3=1.$
\end{itemize}

(III) We assume $\Omega_3=1,$ so that $\Omega_2=\Omega_6,$ and \eqref{eq:LTvanish,D=3,2} becomes
\begin{align*}
\Big(1+\varsigma_2\Big(\frac{3}{g_1}\Big)\Omega_2\Big)\Big(1+\frac{2}{3}\varsigma_1\Big)=0.
\end{align*}
This holds if and only if one of the following two conditions holds:
\begin{itemize}
\item $\varsigma_1=-\frac{3}{2};$
\item $\varsigma_2(\frac{3}{g_1})=-1,$ $\Omega_2=1.$
\end{itemize}

(IV) We now suppose $\Omega_6=1,$ so that $\Omega_2=\Omega_3=1,$ and this can be covered by the last situation.

In summary, we find $\varpi_{E,r}=0$ if and only if one of the following three conditions holds:
\begin{itemize}
\item $3\mid \mu$, $\varsigma_1=3,$ $\varsigma_2(\frac{3}{g_1})=1,$ $\Omega_2=\Omega_3=\Omega_6=-\frac{1}{5};$
\item $\varsigma_2(\frac{3}{g_1})=-1,$ $\Omega_2=1;$
\item $3\mid \mu$, $\varsigma_1=-\frac{3}{2},$ $\Omega_3=1;$
\end{itemize}

The values of $\varsigma_1$ and $\varsigma_2$ can be determined according to modulo arithmetics of $r,\delta,\lambda,\mu$. For each $j\in\Z^+$, the condition $\Omega_j=1$ is equivalent to $\fA_j(g_1)\mid r.$
On the other hand, the condition $\Omega_j=-\frac{1}{5}$ is equivalent to
\begin{align*}
j\in \cJ_5(g,r)\cup\cJ_7(g,r).
\end{align*}
Theorem \ref{thm:LTconstant:D=3} then follows by routine calculations.

\subsection{CM by $\Q(\sqrt{-2})$}
For $r\equiv2\bmod4,$ we have
\begin{align*}
\varpi_{E,r}=\frac{h_{2,r}}{2}\Big(1+\frac{1}{2}(-1)^{\frac{(r-2)(r+10)}{32}+\delta+\lambda+\frac{g_1-1}{2}}\Big(\frac{2}{g_1}\Big)\Omega_2(2;g_1,r)\Big).
\end{align*}
Since $|\Omega_2(2;g_1,r)|\leqslant1,$ we find $\varpi_{E,r}=0$ if and only if $h_{2,r}=0,$ which can only happen when $r\not\equiv2\bmod4$ as shown by \eqref{eq:constant-hD,r}.

\subsection{CM by $\Q(\sqrt{-D})$ with remaining $D$}
As in Theorem \ref{thm:HL-LT:D>=7}, 
\begin{align*}
\varpi_{E,r}=\frac{h_{D,r}}{2}\cdot\Big(1+\varepsilon_D(g,r)\Big(\frac{2^{\lambda+1}g_1r}{D}\Big)\Omega_2(D;g_1,r)\Big)
\end{align*}
for $D\nmid r$, where
\begin{align*}
\varepsilon_D(g,r)
&=\begin{cases}
\frac{1}{2}(-1)^{\frac{g_1-1}{2}}(1+(-1)^{\lambda}),\ \ &\text{for~}  2\parallel r,\\
(-1)^{\delta+\mu+\frac{\lambda r}{4}},\ \ &\text{for~} 4\mid r,\\
\frac{1}{4}(1+(-1)^{\lambda})((-1)^{\frac{g_1-1}{2}}+(-1)^{\delta+\mu}),\ \ &\text{for~} 2\nmid r.
\end{cases}
\end{align*}
In view of \eqref{eq:constant-hD,r}, $h_{D,r}$ never vanishes if $D\equiv3\bmod8.$ On the other hand, for $D\equiv7\bmod8,$ $h_{D,r}$ vanishes if and only if $2\nmid r.$

We now assume $h_{D,r}\neq0,$  so that $\varpi_{E,r}=0$ if and only if 
\begin{align*}
\varepsilon_D(g,r)\Big(\frac{2^{\lambda+1}g_1r}{D}\Big)\Omega_2(D;g_1,r)=-1.
\end{align*}
There will be two possibilities:
\begin{itemize}
\item $\varepsilon_D(g,r)(\frac{2^{\lambda+1}g_1r}{D})=-1,$ $\Omega_2(D;g_1,r)=1;$
\item $\varepsilon_D(g,r)(\frac{2^{\lambda+1}g_1r}{D})=1,$ $\Omega_2(D;g_1,r)=-1.$
\end{itemize}
The condition $\Omega_2(D;g_1,r)=1$ holds if and only if $\fA_2(g_1)\mid r$. On the other hand, the condition 
$\Omega_2(D;g_1,r)=-1$ can only hold if $D=11$, in which case we also require $2\in\cJ_3(g,r)$.

\section{Anomalous primes}\label{sec:anomalousprimes}

Restricting to the case $r=1$, $\pi_{E,r}(x)$ then counts anomalous primes \index{anomalous prime} as defined by B. Mazur \cite{Ma72}. 
The following theorem determines when there are infinitely many anomalous primes for elliptic curves with CM by $\Q(\sqrt{-3}).$

\begin{theorem}\label{thm:anomalousprimes:D=3}\index{anomalous prime} 
Suppose $E/\Q$ is an elliptic curve with CM by $\Q(\sqrt{-3}),$ defined by
\begin{align*}
y^2=x^3+g
\end{align*}
with $g\in\Z^*.$ Assume the Hardy--Littlewood conjecture\index{Hardy--Littlewood conjecture}.
Then there are only finitely many anomalous primes if and only if
\begin{align}\label{eq:anomalousprimes,D=3}
g\in\{\square,~~-3\square,~~\cube,~~80\cube^2,~~-2160\cube^2,~~-268912 \cube^2\text{~~or}~~~7260624\cube^2\},
\end{align}
where $\square$ and $\cube$ denote a square and cube in $\Z,$ respectively, and thus $\cube^2$ denotes a sixth power in $\Z.$
\end{theorem}

\begin{remark}
Note that
\begin{align*}
80=2^4\cdot 5,~2160=2^4\cdot3^3\cdot 5,~268912=2^4\cdot 7^5,~7260624=2^4\cdot3^3\cdot 7^5.
\end{align*}
Moreover, one can compare Theorem \ref{thm:anomalousprimes:D=3} with \cite[Theorem 4.11]{Qi16}, in which the last two possibilities were missed, but this has been included in \cite{Qi21}. In fact, $-3=(1-2\omega)^2$ is a square in $\Z[\omega],$ hence an alternative formulation to \eqref{eq:anomalousprimes,D=3} should be
\begin{align*}
g\in\{\alpha^2,~~\beta^3,~~80\gamma^6,~~-268912\delta^6\}\cap\Z
\end{align*}
with $\alpha,\beta,\gamma,\delta\in\Z[\omega]$, as stated in \cite[Theorem 4.11]{Qi16} together with the subsequent remark after \cite[Theorem]{Qi21}.
\end{remark}

There are no anomalous primes if the elliptic curve has CM by $\Q(\sqrt{-D})$ $(D=1,2,7)$ due to congruence obstructions. We now turn to the remaining $D$.

\begin{theorem}\label{thm:anomalousprimes:D>=11}\index{anomalous prime} 
Let $D\in\{11, 19, 43, 67, 163\}.$ Suppose $E/\Q$ is the elliptic curve defined by $\eqref{eq:E:otherD}$ with CM by $\Q(\sqrt{-D})$.
Assume the Hardy--Littlewood conjecture\index{Hardy--Littlewood conjecture}.
Then there are only finitely many anomalous primes if and only if
\begin{align*}
g\in\{\square,~~~~~-D\square\}
\end{align*}
with two extra possibilities $33\square,$ $-3\square$
if $D=11,$
where $\square$ denotes a square in $\Z.$
\end{theorem}

\begin{remark}
It might be unconditional (and even trivial), in some special cases, to prove the finiteness of anomalous primes due to congruence obstructions. But the full criterion for sufficiency and necessity in Theorems \ref{thm:anomalousprimes:D=3} and \ref{thm:anomalousprimes:D>=11} requires the Hardy--Littlewood conjecture.
\end{remark}

We only give the proof of Theorem \ref{thm:anomalousprimes:D=3}, and Theorem \ref{thm:anomalousprimes:D>=11} can be proven similarly. Hence we will work with $E$ as defined by \eqref{eq:E:D=3} in what follows. 

According to Theorem \ref{thm:LTconstant:D=3}, the constant $\varpi_{E,1}$ vanishes if and only if one of the following conditions holds:
\begin{enumerate}[(I)]
\item $\{2,3,6\}\subseteq\cJ_{5}(g,1)\cup\cJ_{7}(g,1),$ $g_1\equiv5,7\bmod{12},$ $3\mid(\mu,2+\lambda,\frac{g_1^2+2}{3}),$ $2\mid(\lambda,\delta+\mu+\frac{g_1-1}{2});$ \smallskip
\item $\fA_2(g_1)=1,$ $g_1\equiv1,11\bmod{12},$ $2\mid(\lambda,\delta+\mu+\frac{g_1-1}{2});$ \smallskip
\item $\fA_3(g_1)=1,$ $g_1\equiv1,8\bmod9,$ $3\mid(\lambda,\mu);$ \smallskip
\item $\fA_3(g_1)=1,$ $g_1\equiv2,7\bmod9,$ $3\mid(\lambda+1,\mu).$
\end{enumerate}
We would like to analyze the four possibilities one by one.

Case (I): Since $6\in\cJ_{5}(g,1)\cup\cJ_{7}(g,1),$ we may find
\begin{align*}
g_1=5^{\alpha}a^6\text{~or~}7^{\beta}b^6,
\end{align*}
where $1\leqslant \alpha,\beta\leqslant5$ and $(ab,6)=1.$ The extra constraint $\{2,3\}\subseteq\cJ_{5}(g,1)\cup\cJ_{7}(g,1)$ yields the shape of $g_1$ more explicitly:
\begin{align*}
g_1=5a_1^6,~5^5a_2^6,~7b_1^6,\text{~or~}7^5b_2^6,
\end{align*}
where $(a_1a_2b_1b_2,6)=1.$ Note that $n^6\equiv1\bmod m$ for any $n\in\Z$ with $(n,6)=1$ and $m\in\{9,12\}.$ The condition $3\mid \frac{g_1^2+2}{3}$ eliminates two possibilities, so that
\begin{align*}
g_1=5a_1^6\text{~or~}7^5b_2^6.
\end{align*}
In summary, Case (I) happens if and only if one of the following two possibilities occurs:
\begin{itemize}
\item $g_1=5a^6$ with some $(a,6)=1,$ $3\mid(\mu,2+\lambda),$ $2\mid(\lambda,\delta+\mu);$ \smallskip
\item $g_1=7^5b^6$ with some $(b,6)=1,$ $3\mid(\mu,2+\lambda),$ $2\mid(\lambda,\delta+\mu+1).$
\end{itemize}

Case (II): Note that $\fA_2(g_1)=1$ is equivalent to the condition that $g_1$ is a perfect square in $\Z^+,$ so that $g_1\equiv1\bmod{12}.$

Cases (III) \& (IV): The condition $\fA_3(g_1)=1$ can be translated to that $g_1=c^3$, where $c>0$ and $(c,6)=1.$ Therefore,
\begin{align*}
g_1\equiv
\begin{cases}
1\bmod9,\ \ &\text{if~}c\equiv1\bmod6,\\
8\bmod9,\ \ &\text{if~}c\equiv5\bmod6.
\end{cases}
\end{align*}

In summary, we find the constant $\varpi_{E,1}$ vanishes if and only if one of the following conditions holds:

\begin{itemize}
\item $g_1=5a^6$ for some $a\in\Z$ with $(a,6)=1,$ $3\mid(\mu,2+\lambda),$ $2\mid(\lambda,\delta+\mu);$ \smallskip
\item $g_1=7^5b^6$ for some $b\in\Z$ with $(b,6)=1,$ $3\mid(\mu,2+\lambda),$ $2\mid(\lambda,\delta+\mu+1);$\smallskip
\item $g_1=c^2$ for some $c\in\Z$ with $(c,6)=1,$ $2\mid(\lambda,\delta+\mu);$ \smallskip
\item $g_1=d^3$ for some $d\in\Z$ with $(d,6)=1,$ $3\mid(\lambda,\mu).$
\end{itemize}
Then Theorem \ref{thm:anomalousprimes:D=3} follows by displaying $g_1$ explicitly in all above four cases.

\section{Symmetry}
According to Theorems \ref{thm:HL-LT:D=1}, \ref{thm:HL-LT:D=3}, \ref{thm:HL-LT:D=2} and \ref{thm:HL-LT:D>=7}, the shape of $\varpi_{E,r}$ depends heavily on the sign and modulo arithmetics of $r$. In this section, we would like to determine 
when $\varpi_{E,r}$ does not depend on the sign of $r$, i.e.,
\begin{align}\label{eq:LTconstant-symmetry}
\varpi_{E,r}=\varpi_{E,-r}.
\end{align}
Note that $\Omega_j(D;q,r)=\Omega_j(D;q,-r)$ and $h_{D,r}=h_{D,-r}$, which make the analysis of $\varpi_{E,r}$ much easier.

\begin{theorem}\label{thm:LTconstantsymmetry:D=1}\index{Lang--Trotter constant} 
Suppose $E/\Q$ is the elliptic curve defined by $\eqref{eq:E:D=1}$ with CM by $\Q(\sqrt{-1})$.
Then $\eqref{eq:LTconstant-symmetry}$ holds if and only if one of the following possibilities happens:
\begin{itemize}
\item $2\nmid r;$ \smallskip
\item $8\mid r;$\smallskip
\item $r\equiv4\bmod8,$ $2\mid\lambda;$\smallskip
\item $r\equiv2\bmod4,$ $2\nmid\lambda;$\smallskip
\item $r\equiv2\bmod4,$ $2\mid(\lambda,\delta+\frac{\lambda+g_1-1}{2}).$
\end{itemize}
\end{theorem}

\begin{theorem}\label{thm:LTconstantsymmetry:D=3}\index{Lang--Trotter constant} 
Suppose $E/\Q$ is the elliptic curve defined by $\eqref{eq:E:D=3}$ with CM by $\Q(\sqrt{-3})$.
Then $\eqref{eq:LTconstant-symmetry}$ holds if and only if one of the following possibilities happens:
\begin{itemize}
\item $3\mid r;$
\item $3\nmid r$ and $4\nmid r,$ $2\nmid\lambda;$ \smallskip
\item $r\equiv1,5\bmod6,$ $g_1\equiv1\bmod4,$ $2\nmid \delta+\mu;$\smallskip
\item $r\equiv1,5\bmod6,$ $g_1\equiv3\bmod4,$ $2\mid \delta+\mu;$\smallskip
\item $r\equiv1,5\bmod{6},$ $g_1\equiv1,8\bmod9,$ $3\mid \mu,$ $\lambda\not\equiv1\bmod3,$ $\fA_3(g_1^*)\mid r.$
\end{itemize}
Here $g_1^*$ denotes the largest square divisor of $g_1$.
\end{theorem}

\begin{theorem}\label{thm:LTconstantsymmetry:D=2}\index{Lang--Trotter constant} 
Suppose $E/\Q$ is the elliptic curve defined by $\eqref{eq:E:otherD}$ with CM by $\Q(\sqrt{-2})$.
Then $\eqref{eq:LTconstant-symmetry}$ holds if and only if 
$r\not\equiv2\bmod4,$ in which case we always have $\varpi_{E,r}=\varpi_{E,-r}=0.$
\end{theorem}

\begin{theorem}\label{thm:LTconstantsymmetry:D>=7}\index{Lang--Trotter constant} 
Let $D\in\{7,11, 19, 43, 67, 163\}.$ Suppose $E/\Q$ is the elliptic curve defined by $\eqref{eq:E:otherD}$ with CM by $\Q(\sqrt{-D})$.
Then $\eqref{eq:LTconstant-symmetry}$ holds if and only if one of the following possibilities happens:
\begin{itemize}
\item $D\mid r;$\smallskip
\item $D=7,$ $(r,14)=1;$\smallskip
\item $r\equiv2\bmod4,$ $2\nmid\lambda,$ $D\nmid r;$\smallskip
\item $2\nmid r\lambda,$ $D\nmid r;$\smallskip
\item $2\nmid r(\delta+\mu+\frac{g_1-1}{2}),$ $D\nmid r.$
\end{itemize}
\end{theorem}

\subsection{CM by $\Q(\sqrt{-1})$}
Note that $\varpi_{E,r}=0$ for all odd $r$, in which case \eqref{eq:LTconstant-symmetry} holds automatically (and trivially).
In view of \eqref{eq:LTconstant,D=1}, we need to examine \eqref{eq:LTconstant-symmetry} according to the parities of $r,\lambda$ one by one. If $4\mid r,$ we have $(-1)^{\frac{\lambda r}{4}}=(-1)^{\frac{-\lambda r}{4}}$, which yields that
\eqref{eq:LTconstant-symmetry} holds if and only if
\begin{align*}
(1-(-1)^{\frac{\lambda r}{4}})\Re(i^{1+\frac{\lambda r}{4}}-i^{1-\frac{\lambda r}{4}})=0,
\end{align*}
i.e., $8\mid \lambda r$.
If $2\parallel r$ and $2\nmid\lambda$, \eqref{eq:LTconstant-symmetry} does hold automatically. In the remaining case $2\parallel r,2\mid\lambda$, we have $(-1)^{\frac{r-2}{4}}\neq(-1)^{\frac{-r-2}{4}}$, so that we should expect $(1-(-1)^{\delta+\frac{\lambda+g_1-1}{2}})=0,$ i.e., $2\mid \delta+\frac{\lambda+g_1-1}{2}.$

\subsection{CM by $\Q(\sqrt{-3})$}
The case $3\mid r$ guarantees \eqref{eq:LTconstant-symmetry} trivially, and we henceforth assume $3\nmid r,$ in which case we have $h_{3,r}\neq0.$
Note that $\varsigma_1$ does not depend on the sign of $r$ as shown in \eqref{eq:varsigma1}.
Therefore, \eqref{eq:LTconstant-symmetry}  holds if and only if
\begin{align*}
\varsigma_2\cdot(\Omega_2+\widehat{\mu}\varsigma_1\Omega_6)=0,
\end{align*}
where we also adopt the convention $\Omega_j=\Omega_j(3;g_1,r)$ as before.

According to \eqref{eq:varsigma2}, $\varsigma_2=0$ can happen only if $r\equiv\pm2\bmod{12}$ or $r\equiv\pm1\bmod6.$
More precisely, $\varsigma_2=0$ if and only if $2\nmid\lambda$ for $r\equiv\pm2\bmod{12}$; and if and only if $2\nmid\lambda$ or $2\nmid \delta+\mu+\frac{g_1-1}{2}$ for $r\equiv\pm1\bmod6.$

We now analyze the possibilities with
\begin{align*}
\Omega_2+\widehat{\mu}\varsigma_1\Omega_6=0,
\end{align*}
which can never happen if $\widehat{\mu}=0,$ so that we henceforth assume $3\mid \mu,$ in which case the above equation becomes
\begin{align}\label{eq:LTconstantsymmetry,D=3,Omega26} 
\Omega_2+\frac{2}{3}\varsigma_1\Omega_6=0.
\end{align}
Note that $|\Omega_2|\geqslant|\Omega_6|$, then we should require $|\varsigma_1|\geqslant\frac{3}{2}.$ Thus the only possible cases are $\varsigma_1=3$ or $\pm\frac{3}{2}.$

If $\varsigma_1=3,$ then \eqref{eq:LTconstantsymmetry,D=3,Omega26} becomes $\Omega_2+2\Omega_6=0,$ which cannot happen since $\sigma_D(p)\neq2$ as defined by \eqref{eq:sigmaDp}. Similarly, we also have $\varsigma_1\neq\frac{3}{2}$ since \eqref{eq:LTconstantsymmetry,D=3,Omega26} requires $\Omega_2=-\Omega_6$ and $\sigma_D(p)=-1$ for some $p$, which is also impossible.

If $\varsigma_1=-\frac{3}{2},$ then \eqref{eq:LTconstantsymmetry,D=3,Omega26} becomes $\Omega_2=\Omega_6,$ which is equivalent to 
\begin{align}\label{eq:LTconstant,D=3,product=1}
\prod_{\substack{p^\nu\parallel g_1\\ p\nmid r\\ 3\nmid\nu,2\mid \nu}}\sigma_3(p)=1.
\end{align}
Since $\sigma_3(p)\geqslant5$ for any prime $p>3$, \eqref{eq:LTconstant,D=3,product=1} holds if and only if the involved product is empty, i.e., $\fA_3(g_1^*)\mid r,$ where $g_1^*$ denotes the largest square divisor of $g_1$.

In summary, we find \eqref{eq:LTconstant-symmetry}  holds if and only if
\begin{itemize}
\item $3\mid r;$
\item $r\equiv\pm2\bmod{12},$ $2\nmid\lambda;$ 
\item $r\equiv\pm1\bmod6,$ $2\nmid\lambda;$
\item $r\equiv\pm1\bmod6,$ $2\nmid \delta+\mu+\frac{g_1-1}{2};$
\item $3\mid \mu,$ $\varsigma_1=-\frac{3}{2},$ $\fA_3(g_1^*)\mid r.$
\end{itemize}
The value of $\varsigma_1$ can be analyzed via \eqref{eq:varsigma1}, and we omit the details here.

\subsection{CM by $\Q(\sqrt{-D})$ with remaining $D$}
The case $D=2$ is trivial since $$(-1)^{\frac{(r-2)(r+10)}{32}}\neq(-1)^{\frac{(-r-2)(-r+10)}{32}}$$ for $r\equiv2\bmod4$. 

We now assume $D\geqslant7$ and $D\nmid r.$
Hence \eqref{eq:LTconstant-symmetry}  holds if and only if $h_{D,r}=0$ or
\begin{align}\label{eq:LTconstantsymmetry,otherD,initial}
\varepsilon_D(g,r)\Big(\frac{2^{\lambda+1}g_1r}{D}\Big)=\varepsilon_D(g,-r)\Big(\frac{-2^{\lambda+1}g_1r}{D}\Big),
\end{align}
where $\varepsilon_D(g,r)$ is given as in Theorem \ref{thm:HL-LT:D>=7}, i.e.,
\begin{align*}
\varepsilon_D(g,r)
&=\begin{cases}
\frac{1}{2}(-1)^{\frac{g_1-1}{2}}(1+(-1)^{\lambda}),\ \ &\text{for~}  2\parallel r,\\
(-1)^{\delta+\mu+\frac{\lambda r}{4}},\ \ &\text{for~} 4\mid r,\\
\frac{1}{4}(1+(-1)^{\lambda})((-1)^{\frac{g_1-1}{2}}+(-1)^{\delta+\mu}),\ \ &\text{for~} 2\nmid r.
\end{cases}
\end{align*}
Since $(\frac{-1}{D})=-1$ for each $D\in\{7,11,19, 43, 67, 163\}$, then 
\eqref{eq:LTconstantsymmetry,otherD,initial}  holds if and only if
\begin{align*}
\varepsilon_D(g,r)=0,
\end{align*}
which can never happen for $4\mid r,$ and can probably happen for other $r$ according to parities of $\lambda$ and 
$\delta+\mu+\frac{g_1-1}{2}.$
This proves Theorem \ref{thm:LTconstantsymmetry:D>=7} by noting that $c_{D,r}$ vanishes if and only if $2\nmid r$ and $D=7$.

\section{Image of Galois representations}
\label{sec:Galoisimage}

Let $\overline{\Q}$ be an algebraic closure of $\Q.$ For each integer $m>1$, denote by $E[m]$ the $m$-torsion subgroup of $E(\overline{\Q}).$ The action of the absolute Galois group $\Gal(\overline{\Q}/\Q)$ on $E[m]$ induces
a Galois representation\index{Galois representation}
\begin{align*}
\rho_{E,m}: \Gal(\overline{\Q}/\Q)\rightarrow\Aut(E[m])\simeq\GL_2(\Z/m\Z),
\end{align*}
which is a continuous homomorphism. Note that $E[m]$ is a free $\Z/m\Z$-module of rank $2$, then the image of $\rho_{E,m}$, usually denoted by $G_m(E)$,  can be regarded as a subgroup of $\GL_2(\Z/m\Z)$.
Taking the inverse limit of the $\rho_{E,m}$ over positive integers $m$ (with bases chosen compatibly), one obtains a continuous group homomorphism
\begin{align*}
\rho_E: \Gal(\overline{\Q}/\Q)\rightarrow\GL_2(\widehat{\Z}):=\lim_{\leftarrow}\GL_2(\Z/m\Z).
\end{align*}
According to Serre’s open image theorem ([57, Section IV.3.2] and [58]), the image $\rho_E(\Gal(\overline{\Q}/\Q))$ has finite index in $\GL_2(\widehat{\Z})$ if $E$ does not have CM. 
In other words, there exists a positive integer level $m$ such that, if 
\begin{align*}
\varphi: \GL_2(\widehat{\Z})\rightarrow\GL_2(\Z/m\Z)
\end{align*}
is the canonical projection, one has
\begin{align*}
\rho_E(\Gal(\overline{\Q}/\Q))=\varphi^{-1}(G_m(E)).
\end{align*}
The smallest positive integer $m$, denoted by $m_E$, appears naturally in the definition of $c_{E,r}$ by Lang and Trotter \cite[Theorem 4.2]{LT76}.

The case of CM elliptic curves is more complicated  since $\rho_E(\Gal(\overline{\Q}/\Q))$ is not an open subgroup of $\GL_2(\widehat{\Z})$. But it has an index two subgroup which is abelian. We now follow the description by N. Jones \cite{Jo09}. 

Suppose $E$ is an elliptic curve with CM by a imaginary quadratic field $K=\Q(\sqrt{-D})$ with $D\in\fD$.
As before, $R_D$ denotes the ring of integers in $\Q(\sqrt{-D})$. Define the inverse limit $\widehat{R}_D=\displaystyle\lim_{\leftarrow}R_D/mR_D.$ Then the torsion part 
\begin{align*}
E_{\mathrm{tor}}:=\bigcup_{m\geqslant1}E[m]
\end{align*}
of $E$ is a one-dimensional $\widehat{R}_D$-module on which $\Gal(\overline{\Q}/K))$ acts,
preserving the $\widehat{R}_D$-action. Thus, the image of $\rho_E$ restricted to $\Gal(\overline{\Q}/K))$ maps into $(\widehat{R}_D)^*:$
\begin{align*}
\rho_E: \Gal(\overline{\Q}/K)\rightarrow (\widehat{R}_D)^*=\GL_1(\widehat{R}_D).
\end{align*}
Hence the Galois representation\index{Galois representation} on the torsion of a CM elliptic curve can be viewed as a $\GL_1$ analogue of the non-CM situation. 
The Galois image $\rho_E( \Gal(\overline{\Q}/K))$ is contained in $(\widehat{R}_D)^*$ with finite index according to Serre \cite[Section 4.5]{Se72}.

Denote by $K(E[m])$ the number field generated by the coordinates of the $m$-torsion points of $E(K)$.
We may regard $\Gal(K(E[m])/K)$ as a subgroup of $(R_D/mR_D)^*,$ and there is 
a positive integer $m_0$ such that for each positive integer $m$, we have
\begin{align*}
\Gal(K(E[m])/K)\simeq \pi^{-1}(\Gal(K(E[(m,m_0)])/K)),
\end{align*}
where $\pi:(R_D/mR_D)^*\rightarrow(R_D/(m,m_0)R_D)^*$ denotes the canonical projection. Denote by $m_E$ the smallest
positive integer $m_0$ with the above property and also dividing by $4p$ for any prime $p$ that ramifies in $R_D.$
For the odd unramified primes $p$ with $p\nmid m_E,$ we have
\begin{align*}
\Gal(K(E[p])/K)\simeq (R_D/pR_D)^*.
\end{align*}
Following heuristics for elliptic curves without CM by Lang and Trotter \cite{LT76}, Jones \cite[Section 2.2]{Jo09} interpreted the Lang--Trotter constant\index{Lang--Trotter constant}, in the CM situation, by
\begin{align*}
c_{E,r}&=\frac{m_E}{2\pi}\cdot\frac{|\Gal(K(E[m_E])/K)_r|}{|\Gal(K(E[m_E])/K)|}\cdot \prod_{p\mid r,~p\nmid m_E}\Big(1-\frac{(\frac{-D}{p})}{p}\Big)^{-1}\\
&\ \ \ \ \ \times\prod_{p\nmid rm_E}\Big(1-\frac{(\frac{-D}{p})}{(p-1)(p-(\frac{-D}{p}))}\Big)
\end{align*}
as given in \eqref{eq:LTconstant-Galois},
Note that for any subgroup $G\leqslant \GL_2(\Z/n\Z)$, $G_r$ is defined by
\[G_r:=\{g\in G: \tr ~g\equiv r\bmod n\}\]
since we may view $\GL_1(R_D/mR_D)$ as a subgroup of $\GL_2(\Z/m\Z)$ by fixing a $\Z/m\Z$-basis of $R_D/mR_D$. 

Anyway, the exact shape of $c_{E,r}$ depends heavily on computations of images of Galois representations. In the situation of non-CM curves $E/\Q$, there are quite a few works on computing Galois images. For instance, A. Reverter \& N. Vila \cite{RV01} computed the image for all $E/\Q$ with conductor up to 200 and for each prime $p$. A. Sutherland \cite{Su16} described two probabilistic algorithms (a Las Vegas algorithm and a Monte Carlo algorithm) to determine exceptional primes $\ell$ for which $G_\ell(E)$ is properly contained in $\GL_2(\Z/\ell\Z)$, and they are used to determine the groups $G_\ell(E)$ up to local conjugacy. One may also refer to \cite{RZ15, Zy15a, Lo17}. In the CM case, we would like to mention the recent work of \'A. Lozano-Robledo \cite{LR22},\index{Galois representation}
which gives an explicit classification of the possible $p$-adic Galois representations that are attached to elliptic curves $E$ with CM defined over $\Q(j(E))$.
As mentioned in Chapter \ref{ch:background}, A. Ray \cite{Ra23} employed the classification of Lozano-Robledo, as well as the work of F. Campagna and R. Pengo \cite{CP22} on the representations for Galois groups of division fields, 
to compute $c_{E,r}$ explicitly for 20 different CM elliptic curves $E$. As far as he has checked, the 
Lang--Trotter constant $c_{E,r}$, for such special $E$, is always equal to our constant $\varpi_{E,r}$ deduced from the Hardy--Litttlewood conjecture.
It should be very interesting and challenging to give a full description of $c_{E,r}$ based on \cite{LR22} and \cite{CP22}.
Conjecture \ref{conj:cE,r=varpiE,r}, together with many theorems in this chapter, suggests that analytic number theory may give some initial interpretations of $c_{E,r}$ via the
positivity and symmetry of $\varpi_{E,r}.$

%% file: 8.Historicalnotesperspectives.tex
\chapter{Historical notes and some perspectives}\label{ch:Historicalnotesperspectives}

The Lang--Trotter conjecture (Conjecture \ref{conj:LT}) predicts the distribution of primes constrained by fixed Frobenius traces related to elliptic curves. There are some other conjectures, also usually referred as Lang--Trotter conjecture, on the distribution of primes with Frobenius conditions. The determination of Frobenius traces can date back to Gau{\ss}, who essentially considered some typical elliptic curves with CM. In this chapter, we would like to mention a few historical notes and perspectives related to the abovementioned facts, as well as some selected generalizations, examples, and possible variants.

\section{Gau{\ss}' last diary entry and Disquisitiones Arithmeticae}

{\it Observatio per inductionem facta gravissima theoriam residuorum biquadraticorum cum functionibus lemniscaticis elegantissime nectens. Puta, si $a+bi$ est numerus primus, $a-1+bi$ per $2+2i$ divisibilis, multitudo omnium solutionum congruentiae $1=xx+yy+xxyy\bmod{a+bi}$ inclusis $x=\infty,y=\pm i,z=\pm i,y=\infty$ fit
$= (a-1)^2+bb.$}

Gau{\ss}' diary \index{last entry of Gau{\ss}' diary} recorded many mathematical discoveries of Gau{\ss} from 1796 to 1814. The above paragraph appeared in the last entry of this diary on 7 July 1814, and can be translated as follows:

{\it A most important observation made by induction which connects the theory of biquadratic residues most elegantly with the lemniscatic functions. Suppose, if $a+bi$ is a prime number, $a-1+bi$ divisible by $2+2i$, then the number of all solutions of the congruence
$1=xx+yy+xxyy\bmod{a+bi}$, including $x=\infty,y=\pm i,z=\pm i,y=\infty$, equals
$(a-1)^2+bb.$}

The text of the ``Tagebuch'' was rediscovered in 1897 and edited and published by Felix Klein \cite{Kl03}, with the Last Entry on page 33. A later publication appeared in \cite{GaX1}. For a brief history see \cite[Page 97]{Gr84}, and see \cite{Le00} and \cite{Oo16} for some comments from a modern perspective.

The equation $x^2+y^2+x^2y^2=1$ defines a non-singular affine curve over either $\Q$ or $\Q(\sqrt{-1})$. 
Put
\begin{align*}
z=\frac{1-x^2}{y},
\end{align*}
the original equation becomes
\begin{align*}
z^2=1-x^4.
\end{align*}
With further substitutions
\begin{align*}
u=\frac{1}{4}(v+2)^2z,\ \ x=\frac{v-2}{v+2},
\end{align*}
we arrive at the equation
\begin{align*}
u^2=v^3+4v,
\end{align*}
which gives exactly the same curve as defined by \eqref{eq:E:D=1} with $g=-4.$ In a modern language of elliptic curves, Gau{\ss}' observation can be reformulated as follows.
\begin{theorem}\label{thm:Gauss,lastentry}\index{trace of Frobenius}
Let $p\in\Z^+$ be a prime with $p\equiv1\bmod4.$ Suppose $E$ is an
elliptic curve given by $y^2=x^3+4x.$ Then
\begin{align*}
|E_p(\F_p)|=p+1-(\pi+\overline{\pi}),
\end{align*}
where $p=\pi\overline{\pi}\in\Z[i]$ with $\pi\equiv1\bmod{2+2i}.$
\end{theorem}

Theorem \ref{thm:Gauss,lastentry} is of course a special case of Lemma \ref{lm:Frobenius:D=1}, the latter of which is borrowed directly from Ireland and Rosen \cite[Chapter 18, Theorem 5]{IR90}. According to many references, G. Herglotz \cite{He21} gave the first proof for Theorem \ref{thm:Gauss,lastentry}, from which we learn that Herglotz was not aware of the original proof by Gau{\ss} himself contained in \cite[Page 89]{GaII}. In fact, Gau{\ss} directly evaluated the quadratic character sum
\begin{align*}
\sum_{a\bmod p}\Big(\frac{a^4+1}{p}\Big).
\end{align*}
We leave an exercise to interested readers relating this sum to
\begin{align*}
\sum_{a\bmod p}\Big(\frac{a^4-1}{p}\Big).
\end{align*}
S. Chowla \cite{Ch49} gave an elementary proof for Theorem \ref{thm:Gauss,lastentry}.

We now recall another celebrated result of Gau{\ss}.
\begin{theorem}\label{thm:Gauss,DA}
Let $p\in\Z^+$ be a prime with $p\equiv1\bmod3.$ Then $x^3-y^3\equiv1\bmod p$ has $p-2+a$ solutions in $x,y\bmod p,$ where $4p=a^2+27b^2$ with $a,b\in\Z$ and $a\equiv1\bmod3.$
\end{theorem}

Theorem \ref{thm:Gauss,DA} can be found in \S 358 of {\it Disquisitiones Arithmeticae}\index{Disquisitiones Arithmeticae} \cite{GaDA}. Here is a modern reformulation in the language of elliptic curves.
\begin{theorembis}{thm:Gauss,DA}\index{trace of Frobenius}
Let $p\in\Z^+$ be a prime with $p\equiv1\bmod3.$ Suppose $E$ is an
elliptic curve given by $y^2=x^3-432.$ Then
\begin{align*}
|E_p(\F_p)|=p+1+(\pi+\overline{\pi}),
\end{align*}
where $p=\pi\overline{\pi}\in\Z[\omega]$ with $\pi\equiv-1\bmod3.$
\end{theorembis}

We leave an exercise to readers relating the curve $x^3-y^3=1$ with $u^2=v^3-432.$ One can regard Theorem \ref{thm:Gauss,DA} or \ref{thm:Gauss,DA}$'$ as a special case of Lemma \ref{lm:Frobenius:D=3}, which is taken from \cite[Chapter 18, Theorem 4]{IR90}. See also H. Davenport and H. Hasse \cite[Section 7]{DH34}.  See S. Chowla \cite{Ch49} for an elementary proof.

In other words, elliptic curves with CM by $\Q(\sqrt{-1})$ and $\Q(\sqrt{-3})$ have a rich and long history. Starting from 1960's, A. R. Rajwade and his co-authors dealt with elliptic curves CM by $\Q(\sqrt{-D})$ with $D=2,7,11,19$ using cyclotomy and the theory of complex multiplication. In particular, Lemma \ref{lm:Frobenius:D=2} is due to Rajwade himself.
See \cite{Ra70,Ra77,RP81,RPR84} for other related works. For $D\in\{7,11,19,43,67,163\},$ B. H. Gross \cite[Section 24]{Gr80} in his thesis computed the trace of Frobenius for some representative curves for each of such $D$.
In the general case, we appeal to the work of H. M. Stark \cite{St96}
as shown by Lemma \ref{lm:Frobenius:D>=7}. 

All above explicit evaluations of Frobenius traces lie in the heart of our work, and they serve as the starting point to employ analytic number theory to characterize $\pi_{E,r}(x)$ under the Hardy--Littlewood conjecture. 
It is worth to mention that K. Rubin and A. Silverberg \cite{RS09} generalized Stark's work to the case of all CM elliptic curves over arbitrary number fields. See \cite{Si10} for a nice survey of their joint work with reasonable backgrounds.
There were some partial generalizations before that, for instance W. Miller \cite{Mi98} treated the case with CM by an imaginary quadratic field whose discriminant is even and not divisible by $3$.
 We expect the analytic number theory arguments in this memoir, together with the explicit formula for Frobenius trace in \cite{RS09}, can be correspondingly generalized to treat Lang--Trotter conjecture for all CM elliptic curves over arbitrary number fields.

\section{Variants and generalizations of Lang--Trotter conjecture}

Conjecture \ref{conj:LT} of Lang and Trotter admits a lot of variants and generalizations. We would like to mention some of them to motivate interested readers.

Consider an elliptic curve $E/\Q$ with conductor $N_E.$ For each good prime $p\nmid N_E,$ we have already defined the trace of Frobenius $a_p(E)=p+1-|E_p(\F_p)|$. Write $a_p=\pi_p+\overline{\pi}_p$, so that $\pi_p$ is the Frobenius endomorphism of $E$ at $p,$ and
a root of the polynomial $X^2-a_pX+p\in\Z[X].$ 
If $E$ has CM, we know from Deuring \cite{De41} that 
\begin{align*}
\Q(\pi_p)=\End_{\overline{\Q}}(E)\otimes_\Z\Q
\end{align*}
if $a_p\neq0,$
where $\End_{\overline{\Q}}(E)$ denotes the ring of endomorphisms of $E$ over the algebraic closure $\overline{\Q}$ of $\Q$. What will happen if $E$ is a non-CM elliptic curve? Lang and Trotter \cite{LT76} formulated a conjecture on the number of good primes such that $\Q(\pi_p)$ is fixed by a given imaginary quadratic field.
\begin{conjecture}\label{conj:LT-field}
Let $E$ be a non-CM elliptic curve over $\Q$ and $K$ a fixed imaginary quadratic field. Write
\begin{align*}
\pi_{E,K}(x)=|\{p\leqslant x:\Q(\pi_p)=K,~p\nmid N_E\}|.
\end{align*}
Then we have
\begin{align*}
\pi_{E,K}(x)\sim \frac{c_{E,K}\cdot\sqrt{x}}{\log x},
\end{align*}
where $c_{E,K}$ is a positive constant depending on $E$ and $K$.
\end{conjecture}

Conjecture \ref{conj:LT-field} is also known as the Lang--Trotter conjecture\index{Lang--Trotter conjecture}, and it also seems beyond the current methods.
Invoking the square sieve from analytic number theory, A. C. Cojocaru, \'E.  Fouvry and M. R. Murty \cite{CFM05} proved
\begin{align*}
\pi_{E,K}(x)\ll \frac{x}{\log x}\cdot \frac{(\log\log x)^{\frac{13}{12}}}{(\log x)^{\frac{1}{24}}}
\end{align*}
with an implied constant depending on $E$ and $K$, and they may sharp this upper bound with a power saving if a generalized Riemann Hypothesis \index{generalized Riemann Hypothesis} for Dedekind zeta functions is assumed. See \cite{Zy15b,TZ18} for subsequent progressions.

Besides individual upper bounds for $\pi_{E,r}(x)$ and $\pi_{E,K}(x)$ mentioned as above, both of Conjectures \ref{conj:LT} and \ref{conj:LT-field} 
have been verified on average. See \cite{DP99, JS11, JS13, Ja16} for some progressions.

Conjecture \ref{conj:LT} has been generalized to arbitrary number fields, and one may refer to \cite{JS13,Ja16} and references therein for many details. Function fields analogues can be found in Katz \cite{Ka09}.
Generalizations to abelian varieties have also received a lot of attentions in recent years.
We now mention the work of A. C. Cojocaru, R. Davis, A. Silverberg and K. E. Stange \cite{CDSS17}
very briefly. Let $A$ be an abelian variety over $\Q$ of dimension $g$ and conductor $N_A$. Denote by
\begin{align*}
\rho_A:\Gal(\overline{\Q}/\Q)\rightarrow\GSp_{2g}(\widehat{\Z})
\end{align*}
the absolute Galois representation defined by the inverse limit of the representations
\begin{align*}
\overline{\rho}_{A,m}:\Gal(\overline{\Q}/\Q)\rightarrow\GSp_{2g}(\Z/m\Z)
\end{align*}
of $\Gal(\overline{\Q}/\Q)$ on the $m$-torsion $A[m]$ for each positive integer $m$. For each prime $\ell$ we denote by
\begin{align*}
\overline{\rho}_{A,\ell}:\Gal(\overline{\Q}/\Q)\rightarrow\GSp_{2g}(\Z/m\Z)
\end{align*}
the $\ell$-adic representation, i.e., the representation of $\Gal(\overline{\Q}/\Q)$ on the $\ell$-adic Tate module 
$\displaystyle\lim_{\leftarrow} A[\ell^n]$. Similar to the situation of elliptic curves, define
$$a_p(A)=\tr(\rho_{A,\ell}(\Frob_p)).$$
Then the trace $a_p(A)$ is a sum of $2g$ Weil numbers of weight $1$, so that
$|a_p(A)|\leqslant2g\sqrt{p}.$ Put
\begin{align*}
\pi_{A,r}(x):=|\{p\leqslant x:a_p(A)=r,~p\nmid N_E\}|.
\end{align*}
A. C. Cojocaru, R. Davis, A. Silverberg and K. E. Stange \cite{CDSS17}  obtained various upper bounds for $\pi_{A,r}(x)$, conditionally and unconditionally. More importantly, they formulated the following conjecture in the spirit of Conjecture \ref{conj:LT}.

\begin{conjecture}\label{conj:CDSS}
Let $A/\Q$ be an Abelian variety \index{Abelian variety} of dimension $g$ such that the absolute Galois image $\rho_A(\Gal(\overline{\Q}/\Q))$ is open in $\GSp_{2g}(\widehat{\Z}).$ Suppose $a_p(A)/\sqrt{p}$ equidistributes in $[-2g,2g]$, with respect to some Harr measure, as $p$ runs over all good primes. Then for any fixed $r\in\Z^*,$
\begin{align*}
\pi_{A,r}(x):=|\{p\leqslant x:a_p(A)=r,~p\nmid N_E\}|\sim \frac{c_{A,r}\cdot\sqrt{x}}{\log x},
\end{align*}
where $c_{A,r}\geqslant0$ is a constant defined in terms of Galois images.
\end{conjecture}

Conjecture \ref{conj:CDSS} contains Conjecture \ref{conj:LT} as a special case when $A$ is an elliptic curve without CM and $r\neq0$ since the equidistribution is already known thanks to the Sato--Tate conjecture. We now move to Abelian surfaces $\Q$-isogenous to products of non-CM elliptic curves $E_1\times E_2$ defined over $\Q$, the equidistribution of which is also known to M. Harris \cite{Ha09}.
Quite recently, H. Chen, N. Jones and V. Serban \cite{CJS20} considered the case of abelian surfaces in which case the absolute Galois image is no longer open in $\GSp_{2g}(\widehat{\Z}).$ In particular, they concentrate their attentions to products of elliptic curves defined over $\Q$ and without CM. A precise conjecture on the relevant $\pi_{A,r}(x)$
is also formulated, as well as some upper bounds and explicit computations on the constant involved in the conjectural asymptotic formula.
We also mention that \'E. Fouvry and M. R. Murty \cite{FM95}, A. Akbary, C. David and R.  Juricevic \cite{ADJ04} and A. Akbary and J. Parks \cite{AP19} have worked with products of two elliptic curves such that both Frobenius traces are fixed by integers.

To be precise, we take two elliptic curves $E_1,E_2$ defined over $\Q$. Let $r_1,r_2$ be two fixed integers. Consider the counting function
\begin{align*}
\pi_{\mathbf{E},\mathbf{r}}(x):=|\{p\leqslant x: a_p(E_1)=r_1,~a_p(E_2)=r_2,~p\nmid N_{E_1}N_{E_2}\}|
\end{align*}
with $\mathbf{E}=(E_1,E_2),$ and $\mathbf{r}=(r_1,r_2).$ If $E_1$ and $E_2$ are not $\Q$-isogenous and without CM, 
A. Akbary and J. Parks \cite{AP19} conjectured that
\begin{align*}
\pi_{\bf{E},\bf{r}}(x)\sim c_{\mathbf{E},\mathbf{r}}\cdot \log\log x
\end{align*}
for some constant $c_{\mathbf{E},\mathbf{r}}.$ A similar phenomenon was treated in \cite{ADJ04} for $r_1=r_2$.
The case with $r_1=r_2=0$ was already considered by Lang and Trotter \cite[Page 37]{LT76}, and verified by \'E. Fouvry and M. R. Murty \cite{FM95} in an averaged sense.
The Frobenius restriction in H. Chen, N. Jones and V. Serban \cite{CJS20} is much weaker. In fact, they considered 
\begin{align*}
\sum_{r_1+r_2=r}|\{p\leqslant x: a_p(E_1)=r_1,~a_p(E_2)=r_2,~p\nmid N_{E_1}N_{E_2}\}|
\end{align*}
when $E_1,E_2$ are not $\Q$-isogenous and without CM. Here is their conjecture.
\begin{conjecture}\label{conj:CJS}
Let $E_1,E_2$ be two non-CM elliptic curves over $\Q$ that are not $\Q$-isogenous. Let $A/\Q$ be an abelian surface isogenous to $E_1\times E_2.$ For any $r\in\Z^*,$ we have
\begin{align*}
\pi_{E_1,E_2,r}(x):=|\{p\leqslant x:a_p(A)=r,~p\nmid N_E\}|\sim \frac{c_{E_1,E_2,r}\cdot\sqrt{x}}{\log x},
\end{align*}
where $c_{E_1,E_2,r}\geqslant0$ is a constant defined in terms of Galois images in {\rm\cite{CJS20}}.
\end{conjecture}

The CM analogue of Conjecture \ref{conj:CJS} is natural, but we did not succeed in finding any existing literature. Following the arguments in Chapter \ref{ch:HL&upperbound}, we will be led to counting primes which cannot be covered by the Hardy--Littlewood conjecture\index{Hardy--Littlewood conjecture}. For instance, suppose $E_1,E_2$ are two elliptic curves over $\Q$ with CM by $\Q(\sqrt{-1})$.
Then we need to capture primes $p$ of the shape
\begin{align*}
p=m_1^2+n^2=m_2^2+(r/2-n)^2,\ \ m_1,m_2,n\in\Z,
\end{align*}
which can be reduced to counting primes represented by quartic polynomials in two variables (by eliminating the variable $n$). Of course, this cannot be guaranteed by Hardy--Littlewood conjecture or Schinzel Hypothesis. The sieve methods should work for an upper bound of
$\pi_{E_1,E_2,r}(x)$ in this CM case, and circle methods of Hardy--Littlewood can also be employed to predict the precise asymptotic formula. We would like to come back to this topic in forthcoming works ($E_1,E_2$ are not necessarily with the same CM), and we hope analytic number theory, amongst other things, should work in many statistical problems for CM elliptic curves.

%% file: 9.Numericaldata.tex
\chapter{Numerical data}\label{ch:Numericaldata}

Recall that
\begin{align*}
\pi_{E,r}(x)&=|\{p\leqslant x: a_p(E)=r,\ \ p\nmid N_E\}|,\\
\pi_{D,r}(x)&=|\{p=N(\pi)\leqslant x: \pi\in \pi_D,~\Tr(\pi)=r\}|
\end{align*}
as defined by \eqref{eq:piE,r(x)-definition} and \eqref{eq:piD,r(x)-definition}.
We have proven in the previous chapters, under the Hardy--Littlewood conjecture,  that 
\begin{align*}
\pi_{E,r}(x)\sim \frac{\varpi_{E,r}\cdot \sqrt{x}}{\log x},
\end{align*}
where $\varpi_{E,r}$ is a non-negative constant which can be defined explicitly according the CM field of $E/\Q.$
The second conjecture of Lang and Trotter (Conjecture \ref{conj:LT-2}), which is also a special case of the Hardy--Littlewood conjecture \index{Hardy--Littlewood conjecture}, asserts that
\begin{align*}
\pi_{D,r}(x)\sim \frac{h_{D,r}\cdot \sqrt{x}}{\log x},
\end{align*}
where $h_{D,r}$ is a non-negative constant defined by \eqref{eq:constant-hD,r}.

The determination of $\varpi_{E,r}$ lies in the heart of the second part of our work. In order to illustrate the correctness of such determinations, we would like to do some numerical computations of $\pi_{E,r}(x)$ and $\pi_{D,r}(x)$ for suitably large $x$. More precisely, define
\begin{align*}
\kappa(x;E,r):=\frac{\pi_{E,r}(x)}{\pi_{D,r}(x)},\ \ \ \kappa(E,r):=\frac{\varpi_{E,r}}{h_{D,r}}
\end{align*}
if the denominators do not vanish. We expect that 
\begin{align*}
\lim_{x\rightarrow+\infty}\kappa(x;E,r)=\kappa(E,r).
\end{align*}
We will give some explicit computations for some selected elliptic curves $E$ with respect to different choices of $r$. The starting point would be the shape of Frobenius trace $a_p(E)$ as given by Lemmas \ref{lm:Frobenius:D=1}, \ref{lm:Frobenius:D=3}, \ref{lm:Frobenius:D=2} and \ref{lm:Frobenius:D>=7}.

\section{CM by $\Q(\sqrt{-1})$}

Note that $2=(1+i)(1-i)$ is a product of two even Gaussian primes in $\Z[\sqrt{-1}],$ and $3$ is inert in $\Z[\sqrt{-1}].$
Moreover, the evaluation of $(\frac{3}{\pi})_4$ requires the quartic reciprocity law.
Hence we take two typical elliptic curves $y^2=x^3-2x$ and $y^2=x^3-3x$. The numerical results will be shown in Tables \ref{table:y2=x3-2x} and \ref{table:y2=x3-3x}, respectively. In the former case, the computation by MATLAB finds no prime captured by $\pi_{E,r}(x)$ for $r=4,8$, and no primes can be captured by $\pi_{E,r}(x)$ for $r=12$ in the latter case. All of these coincide with the congruence obstruction.

\newpage

\begin{table}[H]
\renewcommand{\arraystretch}{1.1}
\renewcommand{\tabcolsep}{2.2mm}
\centering
\begin{tabular}{cl|cccc}
\toprule
\multicolumn{2}{l}{}            & $\pi_{E,r}(x)$ & $\pi_{1,r}(x)$& $\kappa(x;E,r)$ &  $\kappa(E,r)$ \\ \hline
\multirow{4}{*}{$r=2$} & $x=10^4$  & $4$ & $18$  & $0.2222222222$  &   \multirow{4}{*}{$\frac{1}{4}$} \\ \cline{2-5}
                          & $x=10^5$  &   $13$   & $50$  &  $0.26$ &   \\ \cline{2-5} 
                          & $x=10^6$  &  $32$    & $111$  &  $0.2882882883$  &   \\ \cline{2-5} 
                          & $x=10^7$  &   $82$   & $315$  &  $0.2603174603$ &     \\ \hline
\multirow{4}{*}{$r=6$} & $x=10^4$  & $2$ & $14$  & $0.1428571428$  &   \multirow{4}{*}{$\frac{1}{4}$} \\ \cline{2-5}
                          & $x=10^5$  &   $8$   & $32$  & $0.25$  &   \\ \cline{2-5} 
                          & $x=10^6$  &  $16$    & $80$  &  $0.2$  &   \\ \cline{2-5} 
                          & $x=10^7$  &   $48$   & $210$  & $0.2285714285$  &    \\ \hline
\multirow{4}{*}{$r=10$} & $x=10^4$  & $6$ & $24$  & $0.25$  &   \multirow{4}{*}{$\frac{1}{4}$} \\ \cline{2-5}
                          & $x=10^5$  &   $12$   & $62$  &  $0.1935483871$ &   \\ \cline{2-5} 
                          & $x=10^6$  &  $42$    & $162$  &  $0.2592592593$  &   \\ \cline{2-5} 
                          & $x=10^7$  &   $110$   & $415$  &  $0.2650602409$ &     \\ \hline
\multirow{4}{*}{$r=12$} & $x=10^4$  & $14$ & $14$  & $1$  &   \multirow{4}{*}{$1$} \\ \cline{2-5}
                          & $x=10^5$  &   $29$   & $29$   &  $1$  &   \\ \cline{2-5} 
                          & $x=10^6$  &  $77$    & $77$  &  $1$  &   \\ \cline{2-5} 
                          & $x=10^7$  &   $200$   &  $200$  & $1$   &  \\ \bottomrule
\end{tabular}
\caption{$E:y^2=x^3-2x.$}\label{table:y2=x3-2x}
\end{table}

\begin{table}[H]
\renewcommand{\arraystretch}{1.1}
\renewcommand{\tabcolsep}{2.2mm}
\centering
\begin{tabular}{cl|cccc}
\toprule
\multicolumn{2}{l}{}            & $\pi_{E,r}(x)$ & $\pi_{1,r}(x)$& $\kappa(x;E,r)$ &  $\kappa(E,r)$ \\ \hline
\multirow{4}{*}{$r=2$} & $x=10^4$  & $7$ & $18$  & $0.3888888889$  &   \multirow{4}{*}{$\frac{1}{3}$} \\ \cline{2-5}
                          & $x=10^5$  &   $19$   & $50$  &  $0.38$ &   \\ \cline{2-5} 
                          & $x=10^6$  &  $41$    & $111$  &  $0.3693693693$  &   \\ \cline{2-5} 
                          & $x=10^7$  &   $101$    & $315$  &  $0.3206349206$ &   \\ \hline
\multirow{4}{*}{$r=4$} & $x=10^4$  & $5$ & $21$  & $0.2380952381$  &   \multirow{4}{*}{$\frac{1}{3}$} \\ \cline{2-5}
                          & $x=10^5$  &   $15$   & $53$  & $0.2830188679$  &   \\ \cline{2-5} 
                          & $x=10^6$  &  $42$    & $125$  &  $0.336$  &   \\ \cline{2-5} 
                          & $x=10^7$  &   $110$    & $324$  & $0.3395061728$  &   \\ \hline
\multirow{4}{*}{$r=6$} & $x=10^4$  & $14$ & $14$  & $1$  &   \multirow{4}{*}{$1$} \\ \cline{2-5}
                          & $x=10^5$  &   $32$   & $32$  &  $1$ &   \\ \cline{2-5} 
                          & $x=10^6$  &  $80$    & $80$  &  $1$  &   \\ \cline{2-5} 
                          & $x=10^7$  &   $210$    & $210$  &  $1$ &   \\ \hline
\multirow{4}{*}{$r=8$} & $x=10^4$  & $5$ & $21$  & $0.2380952381$  &   \multirow{4}{*}{$\frac{1}{3}$} \\ \cline{2-5}
                          & $x=10^5$  &   $17$   & $48$   &  $0.3541666667$  &   \\ \cline{2-5} 
                          & $x=10^6$  &  $42$    & $126$  &  $0.3333333333$  &   \\ \cline{2-5} 
                          & $x=10^7$  &   $104$   &  $316$ & $0.3291139241$   &   \\ \hline
\multirow{4}{*}{$r=10$} & $x=10^4$  & $7$ & $24$  & $0.2916666667$  &   \multirow{4}{*}{$\frac{1}{3}$} \\ \cline{2-5}
                          & $x=10^5$  &   $19$   & $62$   &  $0.3064516129$  &   \\ \cline{2-5} 
                          & $x=10^6$  &  $51$    & $162$  &  $0.3148148148$  &   \\ \cline{2-5} 
                          & $x=10^7$  &    $137$   &  $415$  & $0.3301204819$   &   \\ \bottomrule
\end{tabular}
\caption{$E:y^2=x^3-3x.$}\label{table:y2=x3-3x}
\end{table}

\newpage

\section{CM by $\Q(\sqrt{-3})$}
For the case with CM by $\Q(\sqrt{-3})$, we take the curve $y^2=x^3-432$ which was previously considered by Gau{\ss} (as discussed in Chapter \ref{ch:Historicalnotesperspectives}).

\begin{table}[hp]
\renewcommand{\arraystretch}{1.1}
\renewcommand{\tabcolsep}{2.2mm}
\centering
\begin{tabular}{cl|cccc}
\toprule
\multicolumn{2}{l}{}            & $\pi_{E,r}(x)$ & $\pi_{1,r}(x)$& $\kappa(x;E,r)$ &  $\kappa(E,r)$ \\ \hline
\multirow{4}{*}{$r=-7$} & $x=10^4$  & $10$ & $31$  & $0.3225806452$  &   \multirow{4}{*}{$\frac{1}{3}$} \\ \cline{2-5}
                          & $x=10^5$  &   $27$   & $76$  &  $0.3552631579$ &   \\ \cline{2-5} 
                          & $x=10^6$  &  $71$    & $204$  &  $0.3480392157$  &   \\ \cline{2-5} 
                          & $x=10^7$  &   $178$   & $531$  &  $0.3352165725$ &   \\ \hline
\multirow{4}{*}{$r=-4$} & $x=10^4$  & $5$ & $14$  & $0.3571428571$  &   \multirow{4}{*}{$\frac{1}{3}$} \\ \cline{2-5}
                          & $x=10^5$  &   $10$   & $30$  & $0.3333333333$  &   \\ \cline{2-5} 
                          & $x=10^6$  &  $30$    & $84$  &  $0.3571428571$  &   \\ \cline{2-5} 
                          & $x=10^7$  &   $74$   & $215$  & $0.3441860465$  &   \\ \hline
\multirow{4}{*}{$r=-1$} & $x=10^4$  & $11$ & $28$  & $0.3928571429$  &   \multirow{4}{*}{$\frac{1}{3}$} \\ \cline{2-5}
                          & $x=10^5$  &   $20$   & $64$  &  $0.3125$ &   \\ \cline{2-5} 
                          & $x=10^6$  &  $52$    & $173$  &  $0.3005780347$  &   \\ \cline{2-5} 
                          & $x=10^7$  &   $145$   & $438$  &  $0.3310502283$ &   \\ \hline
\multirow{4}{*}{$r=2$} & $x=10^4$  & $3$ & $11$  & $0.2727272727$  &   \multirow{4}{*}{$\frac{1}{3}$} \\ \cline{2-5}
                          & $x=10^5$  &   $9$   & $29$   &  $0.3103448276$  &   \\ \cline{2-5} 
                          & $x=10^6$  &  $25$    & $81$  &  $0.3086419753$  &   \\ \cline{2-5} 
                          & $x=10^7$  &   $77$   &  $228$  & $0.3377192982$   &   \\ \hline
\multirow{4}{*}{$r=5$} & $x=10^4$  & $5$ & $20$  & $0.25$  &   \multirow{4}{*}{$\frac{1}{3}$} \\ \cline{2-5}
                          & $x=10^5$  &   $18$   & $55$  &  $0.3272727273$ &   \\ \cline{2-5} 
                          & $x=10^6$  &  $46$    & $137$  &  $0.3357664234$  &   \\ \cline{2-5} 
                          & $x=10^7$  &   $126$   & $361$  &  $0.3490304709$ &   \\ \hline
\multirow{4}{*}{$r=8$} & $x=10^4$  & $4$ & $16$  & $0.25$  &   \multirow{4}{*}{$\frac{1}{3}$} \\ \cline{2-5}
                          & $x=10^5$  &   $12$   & $32$   &  $0.375$  &   \\ \cline{2-5} 
                          & $x=10^6$  &  $31$    & $88$  &  $0.3522727273$  &   \\ \cline{2-5} 
                          & $x=10^7$  &   $77$   &  $234$  & $0.3290598291$   &   \\ \hline
\multirow{4}{*}{$r=11$} & $x=10^4$  & $7$ & $22$  & $0.3181818182$  &   \multirow{4}{*}{$\frac{1}{3}$} \\ \cline{2-5}
                          & $x=10^5$  &   $21$   & $57$   &  $0.3684210526$  &   \\ \cline{2-5} 
                          & $x=10^6$  &  $56$    & $158$  &  $0.3544303797$  &   \\ \cline{2-5} 
                          & $x=10^7$  &   $138$   &  $409$  & $0.3374083130$   &   \\ \bottomrule
\end{tabular}
\\~\\
\vskip 2mm
\caption{$E:y^2=x^3-432.$}\label{table:y2=x3-432}
\end{table}

For other $r$ with $|r|\leqslant10$, the computation by MATLAB finds no primes captured by $\pi_{E,r}(x)$, which coincides with the congruence obstruction.

\section{CM by $\Q(\sqrt{-7})$}

For the case of CM by $\Q(\sqrt{-D})$ with $D\neq1,3$, the treatments  are relatively easy since there are only two roots of unity in $R_D$. We now choose
\begin{align*}
y^2=4x^3-4\cdot 7\cdot 125x-8\cdot 49\cdot 125
\end{align*}
as a particular example, which is with CM by $\Q(\sqrt{-7})$,
and computations for other CM cases can be done similarly.

\begin{table}[hp]
\renewcommand{\arraystretch}{1.1}
\renewcommand{\tabcolsep}{2.2mm}
\centering
\begin{tabular}{cl|cccc}
\toprule
\multicolumn{2}{l}{}            & $\pi_{E,r}(x)$ & $\pi_{1,r}(x)$& $\kappa(x;E,r)$ &  $\kappa(E,r)$ \\ \hline
\multirow{4}{*}{$r=-6$} & $x=10^4$  & $4$ & $7$  & $0.5714285714$  &  \multirow{4}{*}{$\frac{3}{5}$} \\ \cline{2-5}
                          & $x=10^5$  &   $12$   & $18$  &  $0.6666666667$ &   \\ \cline{2-5} 
                          & $x=10^6$  &  $19$    & $27$  &  $0.7037037037$  &   \\ \cline{2-5} 
                          & $x=10^7$  &   $27$   & $39$  &  $0.6923076923$ &      \\ \hline
\multirow{4}{*}{$r=-4$} & $x=10^4$  & $6$ & $11$  & $0.5454545454$  &   \multirow{4}{*}{$\frac{2}{5}$} \\ \cline{2-5}
                          & $x=10^5$  &   $12$   & $29$  & $0.4137931034$  &   \\ \cline{2-5} 
                          & $x=10^6$  &  $16$    & $40$  &  $0.4$  &   \\ \cline{2-5} 
                          & $x=10^7$  &   $21$   & $60$  & $0.35$  &      \\ \hline
\multirow{4}{*}{$r=-2$} & $x=10^4$  & $4$ & $12$  & $0.3333333333$  &   \multirow{4}{*}{$\frac{2}{5}$} \\ \cline{2-5}
                          & $x=10^5$  &   $12$   & $28$  &  $0.4285714286$ &   \\ \cline{2-5} 
                          & $x=10^6$  &  $17$    & $40$  &  $0.425$  &   \\ \cline{2-5} 
                          & $x=10^7$  &   $28$   & $64$  &  $0.4375$ &    \\ \hline
\multirow{4}{*}{$r=2$} & $x=10^4$  & $8$ & $12$  & $0.6666666667$  &   \multirow{4}{*}{$\frac{3}{5}$} \\ \cline{2-5}
                          & $x=10^5$  &   $16$   & $28$   &  $0.5714285714$  &   \\ \cline{2-5} 
                          & $x=10^6$  &  $23$    & $40$  &  $0.575$  &   \\ \cline{2-5} 
                          & $x=10^7$  &   $36$   &  $64$  & $0.5625$   &   \\ \hline
\multirow{4}{*}{$r=4$} & $x=10^4$  & $5$ & $11$  & $0.4545454545$  &   \multirow{4}{*}{$\frac{3}{5}$} \\ \cline{2-5}
                          & $x=10^5$  &   $17$   & $29$  &  $0.5862068966$ &   \\ \cline{2-5} 
                          & $x=10^6$  &  $24$    & $40$  &  $0.6$  &   \\ \cline{2-5} 
                          & $x=10^7$  &   $39$   & $60$  &  $0.65$ &    \\ \hline
\multirow{4}{*}{$r=6$} & $x=10^4$  & $3$ & $7$  & $0.4285714286$  &   \multirow{4}{*}{$\frac{2}{5}$} \\ \cline{2-5}
                          & $x=10^5$  &   $6$   & $18$   &  $0.3333333333$  &   \\ \cline{2-5} 
                          & $x=10^6$  &  $8$    & $27$  &  $0.2962962963$  &   \\ \cline{2-5} 
                          & $x=10^7$  &   $12$   &  $39$  & $ 0.3076923077$   &   \\ \bottomrule
\end{tabular}
\\~\\
\vskip 2mm
\caption{$E:y^2=4x^3-4\cdot 7\cdot 125x-8\cdot 49\cdot 125.$}\label{table:4x3-4*7*125x-8*49*125}
\end{table}

MATLAB finds no primes captured by $\pi_{E,r}(x)$ for $r=-5$, agreeing with the congruence obstruction, and $\kappa(x;E,r)$ is always $1$ for $r=5$.

%% file: appendix.tex
%
%
%

\appendix

\chapter{Evaluation of $\varpi_{E,r}$, with CM by $\Q(\sqrt{-3})$}
\label{ch:appendix}

In this appendix, we give the details on the computation of $\varpi_{E,r}$ for $E$ with CM by $\Q(\sqrt{-3}).$ We keep all conventions and notation as in Chapter \ref{ch:D=3}.
Recall that 
\begin{align*}
\varpi_{E,r}=\frac{h_{3,r}}{216}\cdot C,
\end{align*}
where
\begin{align*}
C&=\sum_{0\leqslant d\leqslant2}C_d
\end{align*}
with
\begin{align*}
C_d
&=\mathop{\sum\sum}_{\substack{0\leqslant a,b\leqslant 1\\(a,b)\neq (0,0)}}\mathop{\sum\sum}_{0\leqslant s,t\leqslant 2}\sum_{0\leqslant k\leqslant 3}\sum_{\pm}\Delta(\mp r,k,a+b\omega,\omega^d,s+t\omega)\\
&\ \ \ \times E_{d-d_0,k,\pm\epsilon}(r,g_1,\omega^d).
\end{align*}
We simplify the above expression by employing the definitions of 
$\Delta(r,k,\gamma,\beta,\eta)$ and $E_{d,k,\epsilon}(r,\alpha,\beta)$ as given in Conjecture \ref{conj:Ed(x;r,alpha,beta,gamma)}.

The case $r\equiv\pm16\bmod{24}$ has been treated in Section \ref{sec:proofofTheorem{thm:HL-LT:D=3}}, and the computation here would proceed along with the following cases:

\begin{itemize}
\item $r\equiv\pm4\bmod{24};$
\item $r\equiv\pm10\bmod{12};$
\item $r\equiv\pm1\bmod{6}.$
\end{itemize}
Only $r\equiv0\bmod3$ has not been covered, in which case we have $h_{3,r}=0$.

Recall that $\widehat{\mu}$ vanishes unless $3\mid\mu$, in which case it is equal to $\frac{2}{3}.$ We also write $\Omega_j=\Omega_j(3;g_1,r).$

\smallskip

\underline{\bf I. The case $r\equiv\pm4\bmod{24}$:}

We first consider the case $r\equiv4\bmod{72},$ for which it follows from \eqref{eq:DeltaE(r,beta)1} and \eqref{eq:Ed,k,epsilon(r,alpha,beta)}
that
\begin{align*}
C_0
&=\mathop{\sum\sum}_{\substack{0\leqslant s,t\leqslant 2\\2s-t\equiv0\bmod3}}\Delta(r,3,1,1,s+t\omega)E_{\mu t,3,-(-1)^{\delta+\mu}}(r,g_1,1),\\
C_1
&=\mathop{\sum\sum}_{\substack{0\leqslant s,t\leqslant 2\\s+t\equiv1\bmod3}}\Delta(r,3,1+\omega,\omega,s+t\omega)E_{1+\mu t,3,-(-1)^{\delta+\mu}}(r,g_1,\omega),\\
C_2
&=\mathop{\sum\sum}_{\substack{0\leqslant s,t\leqslant 2\\2t-s\equiv-1\bmod3}}\Delta(r,3,\omega,\omega^2,s+t\omega)E_{2+\mu t,3,-(-1)^{\delta+\mu}}(r,g_1,\omega^2).
\end{align*}
Note that
\begin{align*}
C_0
&=8\mathop{\sum\sum}_{\substack{0\leqslant s,t\leqslant 2\\2s-t\equiv0\bmod3}}\Big\{1-(-1)^{\delta+\mu+\frac{g_1-1}{2}}\Big(\frac{3}{g_1}\Big)\Omega_2\\
&\ \ \ \ +2\Re(\omega^{\mu t})\cdot\Big(\Omega_3-(-1)^{\delta+\mu+\frac{g_1-1}{2}}\Big(\frac{3}{g_1}\Big)\Omega_6\Big)\Big\}\\
&=24\Big\{1-(-1)^{\delta+\mu+\frac{g_1-1}{2}}\Big(\frac{3}{g_1}\Big)\Omega_2\Big\}\\
&\ \ \ \ +16\sum_{0\leqslant t\leqslant 2}\Re(\omega^{\mu t})\cdot\Big\{\Omega_3-(-1)^{\delta+\mu+\frac{g_1-1}{2}}\Big(\frac{3}{g_1}\Big)\Omega_6\Big\}\\
&=24\Big\{1-(-1)^{\delta+\mu+\frac{g_1-1}{2}}\Big(\frac{3}{g_1}\Big)\Omega_2\Big\}+24\widehat{\mu}\cdot\Big\{\Omega_3-(-1)^{\delta+\mu+\frac{g_1-1}{2}}\Big(\frac{3}{g_1}\Big)\Omega_6\Big\},
\end{align*}
and
\begin{align*}
C_1
&=8\mathop{\sum\sum}_{\substack{0\leqslant s,t\leqslant 2\\s+t\equiv1\bmod3}}\Big\{1-(-1)^{\delta+\mu+\frac{g_1-1}{2}}\Big(\frac{3}{g_1}\Big)\Omega_2\\
&\ \ \ \ +2\Re\Big(\omega^{\lambda+\mu t}\Big(\frac{\omega}{g_1}\Big)_3\Big)\cdot\Big(\Omega_3-(-1)^{\delta+\mu+\frac{g_1-1}{2}}\Big(\frac{3}{g_1}\Big)\Omega_6\Big)\Big\}\\
&=24\Big\{1-(-1)^{\delta+\mu+\frac{g_1-1}{2}}\Big(\frac{3}{g_1}\Big)\Omega_2\Big\}\\
&\ \ \ \ +16\sum_{0\leqslant t\leqslant 2}\Re\Big(\omega^{\lambda+\mu t}\Big(\frac{\omega}{g_1}\Big)_3\Big)\cdot\Big\{\Omega_3-(-1)^{\delta+\mu+\frac{g_1-1}{2}}\Big(\frac{3}{g_1}\Big)\Omega_6\Big\}\\
&=24\Big\{1-(-1)^{\delta+\mu+\frac{g_1-1}{2}}\Big(\frac{3}{g_1}\Big)\Omega_2\Big\}\\
&\ \ \ \ +72\widehat{\mu}\cdot \Re\big(\omega^{\lambda+\frac{g_1^2-1}{3}}\big)\cdot\Big\{\Omega_3-(-1)^{\delta+\mu+\frac{g_1-1}{2}}\Big(\frac{3}{g_1}\Big)\Omega_6\Big\}.
\end{align*}
Similarly, we also have
\begin{align*}
C_2
&=24\Big\{1-(-1)^{\delta+\mu+\frac{g_1-1}{2}}\Big(\frac{3}{g_1}\Big)\Omega_2\Big\}\\
&\ \ \ \ +72\widehat{\mu}\cdot \Re\big(\omega^{\lambda+\frac{g_1^2-1}{3}}\big)\cdot\Big\{\Omega_3-(-1)^{\delta+\mu+\frac{g_1-1}{2}}\Big(\frac{3}{g_1}\Big)\Omega_6\Big\}.
\end{align*}
These yield
\begin{align*}
C&=72\Big\{1-(-1)^{\delta+\mu+\frac{g_1-1}{2}}\Big(\frac{3}{g_1}\Big)\Omega_2\Big\}\\
&\ \ \ \ +72\widehat{\mu}\cdot\Big(1+2\Re\big(\omega^{\lambda+\frac{g_1^2-1}{3}}\big)\Big)\cdot\Big\{\Omega_3-(-1)^{\delta+\mu+\frac{g_1-1}{2}}\Big(\frac{3}{g_1}\Big)\Omega_6\Big\}.
\end{align*}
Following similar arguments, we also find the above expression of $C$ is also valid for $r\equiv 28,52\bmod{72},$ i.e., for all $r\equiv4\bmod{24}$.

In a similar manner, we also obtain
\begin{align*}
C&=72\Big\{1+(-1)^{\delta+\mu+\frac{g_1-1}{2}}\Big(\frac{3}{g_1}\Big)\Omega_2\Big\}\\
&\ \ \ \ +72\widehat{\mu}\cdot\Big(1+2\Re\big(\omega^{\lambda+\frac{g_1^2-1}{3}}\big)\Big)\cdot\Big\{\Omega_3+(-1)^{\delta+\mu+\frac{g_1-1}{2}}\Big(\frac{3}{g_1}\Big)\Omega_6\Big\}
\end{align*}
for all $r\equiv-4\bmod{24}$.

\underline{\bf II. The case $r\equiv\pm10\bmod{12}$:}
We first consider $r\equiv22\bmod{36},$ for which it follows that
\begin{align*}
C_0
&=\mathop{\sum\sum}_{\substack{0\leqslant s,t\leqslant 2\\2s\equiv t\bmod3}}
\big\{\Delta(r,0,1,1,s+t\omega)E_{\mu t,0,-1}(r,g_1,1)\\
&\ \ \ \ \ \ \ \ \ \ \ \ \ \ \ \ \ \ \ \ +\Delta(r,2,1,1,s+t\omega)E_{\mu t,2,-(-1)^{\lambda}}(r,g_1,1)\big\}\\
&=4\mathop{\sum\sum}_{\substack{0\leqslant s,t\leqslant 2\\2s\equiv t\bmod3}}\Big\{1-\Big(\frac{3}{g_1}\Big)\Omega_2+2\Re(\omega^{\mu t})\cdot\Big(\Omega_3-\Big(\frac{3}{g_1}\Big)\Omega_6\Big)\Big\}\\
&\ \ \ \ +4\mathop{\sum\sum}_{\substack{0\leqslant s,t\leqslant 2\\2s\equiv t\bmod3}}\Big\{1-(-1)^{\lambda}\Big(\frac{3}{g_1}\Big)\Omega_2+2\Re(\omega^{\mu t})\cdot\Big(\Omega_3-(-1)^{\lambda}\Big(\frac{3}{g_1}\Big)\Omega_6\Big)\Big\}\\
&=24\Big\{1-\frac{1+(-1)^{\lambda}}{2}\Big(\frac{3}{g_1}\Big)\Omega_2\Big\}+72\widehat{\mu}\cdot\Big\{\Omega_3-\frac{1+(-1)^{\lambda}}{2}\Big(\frac{3}{g_1}\Big)\Omega_6\Big\}.
\end{align*}
Moreover,
\begin{align*}
C_1
&=\mathop{\sum\sum}_{\substack{0\leqslant s,t\leqslant 2\\s+t\equiv1\bmod3}}
\big\{\Delta(r,0,1+\omega,\omega,s+t\omega)E_{\mu t,0,-1}(r,g_1,\omega)\\
&\ \ \ \ \ \ \ \ \ \ \ \ \ \ \ \ \ \ \ \ +\Delta(r,2,1+\omega,\omega,s+t\omega)E_{\mu t,2,-(-1)^{\lambda}}(r,g_1,\omega)\big\}\\
&=4\mathop{\sum\sum}_{\substack{0\leqslant s,t\leqslant 2\\s+t\equiv1\bmod3}}\Big\{1-\Big(\frac{3}{g_1}\Big)\Omega_2+2\Re\Big(\omega^{\lambda+\mu t}\Big(\frac{\omega}{g_1}\Big)_3\Big)\cdot\Big(\Omega_3-\Big(\frac{3}{g_1}\Big)\Omega_6\Big)\Big\}\\
&\ \ \ \ +4\mathop{\sum\sum}_{\substack{0\leqslant s,t\leqslant 2\\s+t\equiv1\bmod3}}\Big\{1-(-1)^{\lambda}\Big(\frac{3}{g_1}\Big)\Omega_2+2\Re\Big(\omega^{\lambda+\mu t}\Big(\frac{\omega}{g_1}\Big)_3\Big)\\
&\ \ \ \ \ \ \ \ \ \ \ \ \ \ \ \ \ \ \ \ \ \ \ \ \ \ \ \ \ \ \ \ \ \ \ \ \ \ \ \ \times\Big(\Omega_3-(-1)^{\lambda}\Big(\frac{3}{g_1}\Big)\Omega_6\Big)\Big\}\\
&=24\Big\{1-\frac{1+(-1)^{\lambda}}{2}\Big(\frac{3}{g_1}\Big)\Omega_2\Big\}\\
&\ \ \ \ +72\widehat{\mu}\cdot\Re\Big(\omega^{\lambda}\Big(\frac{\omega}{g_1}\Big)_3\Big)\cdot\Big\{\Omega_3-\frac{1+(-1)^{\lambda}}{2}\Big(\frac{3}{g_1}\Big)\Omega_6\Big\}.
\end{align*}
Similarly, we also have
\begin{align*}
C_2
&=24\Big\{1-\frac{1+(-1)^{\lambda}}{2}\Big(\frac{3}{g_1}\Big)\Omega_2\Big\}\\
&\ \ \ \ +72\widehat{\mu}\cdot\Re\Big(\omega^{-\lambda}\Big(\frac{\omega^2}{g_1}\Big)_3\Big)\cdot\Big\{\Omega_3-\frac{1+(-1)^{\lambda}}{2}\Big(\frac{3}{g_1}\Big)\Omega_6\Big\}.
\end{align*}
It then follows that
\begin{align*}
C&=72\Big\{1-\frac{1+(-1)^{\lambda}}{2}\Big(\frac{3}{g_1}\Big)\Omega_2\Big\}\\
&\ \ \ \ +72\widehat{\mu}\cdot\Big(1+2\Re\big(\omega^{\lambda+\frac{g_1^2-1}{3}}\big)\Big)\cdot\Big\{\Omega_3-\frac{1+(-1)^{\lambda}}{2}\Big(\frac{3}{g_1}\Big)\Omega_6\Big\}.
\end{align*}

Note that the above expression of $C$ is valid for all $r\equiv10\bmod{12}$, and we may also deduce that
\begin{align*}
C&=72\Big\{1+\frac{1+(-1)^{\lambda}}{2}\Big(\frac{3}{g_1}\Big)\Omega_2\Big\}\\
&\ \ \ \ +72\widehat{\mu}\cdot\Big(1+2\Re\big(\omega^{\lambda+\frac{g_1^2-1}{3}}\Big)\big)\cdot\Big\{\Omega_3+\frac{1+(-1)^{\lambda}}{2}\Big(\frac{3}{g_1}\Big)\Omega_6\Big\}
\end{align*}
for all
$r\equiv2\bmod{12}$.

\underline{\bf III. The case $r\equiv\pm1\bmod{6}$:}
We consider $r\equiv13\bmod{18},$
for which it follows that
\begin{align*}
C_0
&=\mathop{\sum\sum}_{\substack{0\leqslant s,t\leqslant 2\\2s-t\equiv0\bmod3}}
\big\{\Delta(r,0,\omega,1,s+t\omega)E_{1-\lambda+\mu t,0,-1}(r,g_1,1)\\
&\ \ \ \ \ \ \ \ \ \ \ \ \ \ \ \ \ \ \ \ +\Delta(r,1,\omega,1,s+t\omega)E_{1-\lambda+\mu t,1,-(-1)^{\delta+\lambda+\mu}}(r,g_1,1)\\
&\ \ \ \ \ \ \ \ \ \ \ \ \ \ \ \ \ \ \ \ +\Delta(r,2,\omega,1,s+t\omega)E_{1-\lambda+\mu t,2,-(-1)^{\lambda}}(r,g_1,1)\\
&\ \ \ \ \ \ \ \ \ \ \ \ \ \ \ \ \ \ \ \ +\Delta(r,3,\omega,1,s+t\omega)E_{1-\lambda+\mu t,3,-(-1)^{\delta+\mu}}(r,g_1,1)\\
&\ \ \ \ \ \ \ \ \ \ \ \ \ \ \ \ \ \ \ \ +\Delta(r,0,1+\omega,1,s+t\omega)E_{2+\lambda+\mu t,0,-1}(r,g_1,1)\\
&\ \ \ \ \ \ \ \ \ \ \ \ \ \ \ \ \ \ \ \ +\Delta(r,1,1+\omega,1,s+t\omega)E_{2+\lambda+\mu t,1,-(-1)^{\delta+\lambda+\mu}}(r,g_1,1)\\
&\ \ \ \ \ \ \ \ \ \ \ \ \ \ \ \ \ \ \ \ +\Delta(r,2,1+\omega,1,s+t\omega)E_{2+\lambda+\mu t,2,-(-1)^{\lambda}}(r,g_1,1)\\
&\ \ \ \ \ \ \ \ \ \ \ \ \ \ \ \ \ \ \ \ +\Delta(r,3,1+\omega,1,s+t\omega)E_{2+\lambda+\mu t,3,-(-1)^{\delta+\mu}}(r,g_1,1)\big\}\\
&=\sum_{1\leqslant j\leqslant 8}C_0^{(j)},
\end{align*}
say.

Note that
\begin{align*}
C_0^{(1)}
&=\mathop{\sum\sum}_{\substack{0\leqslant s,t\leqslant 2\\2s-t\equiv0\bmod3}}\Big\{1-\Big(\frac{3}{g_1}\Big)\Omega_2+2\Re(\omega^{1-\lambda+\mu t})\cdot\Big(\Omega_3-\Big(\frac{3}{g_1}\Big)\Omega_6\Big)\Big\}\\
&=3\Big\{1-\Big(\frac{3}{g_1}\Big)\Omega_2\Big\}+9\widehat{\mu}\cdot\Re(\omega^{1-\lambda})\cdot\Big\{\Omega_3-\Big(\frac{3}{g_1}\Big)\Omega_6\Big\}.
\end{align*}
Similarly, we have
\begin{align*}
C_0^{(2)}
&=\mathop{\sum\sum}_{\substack{0\leqslant s,t\leqslant 2\\2s-t\equiv0\bmod3}}\Big\{1-(-1)^{\delta+\lambda+\mu+\frac{g_1-1}{2}}\Big(\frac{3}{g_1}\Big)\Omega_2+2\Re(\omega^{1-\lambda+\mu t})\\
&\qquad\qquad\qquad\qquad\qquad\times\Big(\Omega_3-(-1)^{\delta+\lambda+\mu+\frac{g_1-1}{2}}\Big(\frac{3}{g_1}\Big)\Omega_6\Big)\Big\}\\
&=3\Big\{1-(-1)^{\delta+\lambda+\mu+\frac{g_1-1}{2}}\Big(\frac{3}{g_1}\Big)\Omega_2\Big\}\\
&\qquad+9\widehat{\mu}\cdot\Re(\omega^{1-\lambda})\cdot\Big\{\Omega_3-(-1)^{\delta+\lambda+\mu+\frac{g_1-1}{2}}\Big(\frac{3}{g_1}\Big)\Omega_6\Big\},
\end{align*}
and
\begin{align*}
C_0^{(3)}
&=3\Big\{1-(-1)^{\lambda}\Big(\frac{3}{g_1}\Big)\Omega_2\Big\}+9\widehat{\mu}\cdot\Re(\omega^{1-\lambda})\cdot\Big\{\Omega_3-(-1)^{\lambda}\Big(\frac{3}{g_1}\Big)\Omega_6\Big\}\\
C_0^{(4)}
&=3\Big\{1-(-1)^{\delta+\mu+\frac{g_1-1}{2}}\Big(\frac{3}{g_1}\Big)\Omega_2\Big\}\\
&\qquad+9\widehat{\mu}\cdot\Re(\omega^{1-\lambda})\cdot\Big\{\Omega_3-(-1)^{\delta+\mu+\frac{g_1-1}{2}}\Big(\frac{3}{g_1}\Big)\Omega_6\Big\}\\
C_0^{(5)}
&=3\Big\{1-\Big(\frac{3}{g_1}\Big)\Omega_2\Big\}+9\widehat{\mu}\cdot\Re(\omega^{2+\lambda})\cdot\Big\{\Omega_3-\Big(\frac{3}{g_1}\Big)\Omega_6\Big\}\\
C_0^{(6)}
&=3\Big\{1-(-1)^{\delta+\lambda+\mu+\frac{g_1-1}{2}}\Big(\frac{3}{g_1}\Big)\Omega_2\Big\}\\
&\qquad+9\widehat{\mu}\cdot\Re(\omega^{2+\lambda})\cdot\Big\{\Omega_3-(-1)^{\delta+\lambda+\mu+\frac{g_1-1}{2}}\Big(\frac{3}{g_1}\Big)\Omega_6\Big\}\\
C_0^{(7)}
&=3\Big\{1-(-1)^{\lambda}\Big(\frac{3}{g_1}\Big)\Omega_2\Big\}+9\widehat{\mu}\cdot\Re(\omega^{2+\lambda})\cdot\Big\{\Omega_3-(-1)^{\lambda}\Big(\frac{3}{g_1}\Big)\Omega_6\Big\}\\
C_0^{(8)}
&=3\Big\{1-(-1)^{\delta+\mu+\frac{g_1-1}{2}}\Big(\frac{3}{g_1}\Big)\Omega_2\Big\}\\
&\qquad+9\widehat{\mu}\cdot\Re(\omega^{2+\lambda})\cdot\Big\{\Omega_3-(-1)^{\delta+\mu+\frac{g_1-1}{2}}\Big(\frac{3}{g_1}\Big)\Omega_6\Big\}.
\end{align*}
We then find
\begin{align*}
C_0&=24\Big\{1-\frac{(1+(-1)^{\lambda})(1+(-1)^{\delta+\mu+\frac{g_1-1}{2}})}{4}\Big(\frac{3}{g_1}\Big)\Omega_2\Big\}\\
&\ \ \ \ +72\widehat{\mu}\cdot\Re(\omega^{1-\lambda})\cdot\Big\{\Omega_3-\frac{(1+(-1)^{\lambda})(1+(-1)^{\delta+\mu+\frac{g_1-1}{2}})}{4}\Big(\frac{3}{g_1}\Big)\Omega_6\Big\}.
\end{align*}

We now move to $C_1$. In fact,
\begin{align*}
C_1
&=\mathop{\sum\sum}_{\substack{0\leqslant s,t\leqslant 2\\2s-t\equiv0\bmod3}}
\big\{\Delta(r,0,1,\omega,s+t\omega)E_{1+\mu t,0,-1}(r,g_1,\omega)\\
&\ \ \ \ \ \ \ \ \ \ \ \ \ \ \ \ \ \ \ \ +\Delta(r,1,1,\omega,s+t\omega)E_{1+\mu t,1,-(-1)^{\delta+\lambda+\mu}}(r,g_1,\omega)\\
&\ \ \ \ \ \ \ \ \ \ \ \ \ \ \ \ \ \ \ \ +\Delta(r,2,1,\omega,s+t\omega)E_{1+\mu t,2,-(-1)^{\lambda}}(r,g_1,\omega)\\
&\ \ \ \ \ \ \ \ \ \ \ \ \ \ \ \ \ \ \ \ +\Delta(r,3,1,\omega,s+t\omega)E_{1+\mu t,3,-(-1)^{\delta+\mu}}(r,g_1,\omega)\}\\
&\ \ \ \ \ \ \ \ \ \ \ \ \ \ \ \ \ \ \ \ +\Delta(r,0,\omega,\omega,s+t\omega)E_{2-\lambda+\mu t,0,-1}(r,g_1,\omega)\\
&\ \ \ \ \ \ \ \ \ \ \ \ \ \ \ \ \ \ \ \ +\Delta(r,1,\omega,\omega,s+t\omega)E_{2-\lambda+\mu t,1,-(-1)^{\delta+\lambda+\mu}}(r,g_1,\omega)\\
&\ \ \ \ \ \ \ \ \ \ \ \ \ \ \ \ \ \ \ \ +\Delta(r,2,\omega,\omega,s+t\omega)E_{2-\lambda+\mu t,2,-(-1)^{\lambda}}(r,g_1,\omega)\\
&\ \ \ \ \ \ \ \ \ \ \ \ \ \ \ \ \ \ \ \ +\Delta(r,3,\omega,\omega,s+t\omega)E_{2-\lambda+\mu t,3,-(-1)^{\delta+\mu}}(r,g_1,\omega)\big\}\\
&=\sum_{1\leqslant j\leqslant 8}C_1^{(j)},
\end{align*}
say. It follows that
\begin{align*}
C_1^{(1)}
&=\mathop{\sum\sum}_{\substack{0\leqslant s,t\leqslant 2\\2s-t\equiv0\bmod3}}\Big\{1-\Big(\frac{3}{g_1}\Big)\Omega_2+2\Re\Big(\omega^{1+\mu t}\Big(\frac{\omega}{g_1}\Big)_3\Big)\cdot\Big(\Omega_3-\Big(\frac{3}{g_1}\Big)\Omega_6\Big)\Big\}\\
&=3\Big\{1-\Big(\frac{3}{g_1}\Big)\Omega_2\Big\}+9\widehat{\mu}\cdot\Re\big(\omega^{1+\frac{g_1^2-1}{3}}\big)\cdot\Big\{\Omega_3-\Big(\frac{3}{g_1}\Big)\Omega_6\Big\}.
\end{align*}
Similarly, we find
\begin{align*}
C_1^{(2)}
&=3\Big\{1-(-1)^{\delta+\lambda+\mu+\frac{g_1-1}{2}}\Big(\frac{3}{g_1}\Big)\Omega_2\Big\}\\
&\qquad+9\widehat{\mu}\cdot\Re\big(\omega^{1+\frac{g_1^2-1}{3}}\big)\Big\{\Omega_3-(-1)^{\delta+\lambda+\mu+\frac{g_1-1}{2}}\Big(\frac{3}{g_1}\Big)\Omega_6\Big\}\\
C_1^{(3)}
&=3\Big\{1-(-1)^{\lambda}\Big(\frac{3}{g_1}\Big)\Omega_2\Big\}\\
&\qquad+9\widehat{\mu}\cdot\Re\big(\omega^{1+\frac{g_1^2-1}{3}}\big)\Big\{\Omega_3-(-1)^{\lambda}\Big(\frac{3}{g_1}\Big)\Omega_6\Big\}\\
C_1^{(4)}
&=3\Big\{1-(-1)^{\delta+\mu+\frac{g_1-1}{2}}\Big(\frac{3}{g_1}\Big)\Omega_2\Big\}\\
&\qquad+9\widehat{\mu}\cdot\Re\big(\omega^{1+\frac{g_1^2-1}{3}}\big)\Big\{\Omega_3-(-1)^{\delta+\mu+\frac{g_1-1}{2}}\Big(\frac{3}{g_1}\Big)\Omega_6\Big\}\\
C_1^{(5)}
&=3\Big\{1-\Big(\frac{3}{g_1}\Big)\Omega_2\Big\}+9\widehat{\mu}\cdot\Re\big(\omega^{2-\lambda+\frac{g_1^2-1}{3}}\big)\cdot\Big\{\Omega_3-\Big(\frac{3}{g_1}\Big)\Omega_6\Big\}\\
C_1^{(6)}
&=3\Big\{1-(-1)^{\delta+\lambda+\mu+\frac{g_1-1}{2}}\Big(\frac{3}{g_1}\Big)\Omega_2\Big\}\\
&\qquad+9\widehat{\mu}\cdot\Re\big(\omega^{2-\lambda+\frac{g_1^2-1}{3}}\big)\Big\{\Omega_3-(-1)^{\delta+\lambda+\mu+\frac{g_1-1}{2}}\Big(\frac{3}{g_1}\Big)\Omega_6\Big\}\\
C_1^{(7)}
&=3\Big\{1-(-1)^{\lambda}\Big(\frac{3}{g_1}\Big)\Omega_2\Big\}\\
&\qquad+9\widehat{\mu}\cdot\Re\big(\omega^{2-\lambda+\frac{g_1^2-1}{3}}\big)\Big\{\Omega_3-(-1)^{\lambda}\Big(\frac{3}{g_1}\Big)\Omega_6\Big\}\\
C_1^{(8)}
&=3\Big\{1-(-1)^{\delta+\mu+\frac{g_1-1}{2}}\Big(\frac{3}{g_1}\Big)\Omega_2\Big\}\\
&\qquad+9\widehat{\mu}\cdot\Re\big(\omega^{2-\lambda+\frac{g_1^2-1}{3}}\big)\Big\{\Omega_3-(-1)^{\delta+\mu+\frac{g_1-1}{2}}\Big(\frac{3}{g_1}\Big)\Omega_6\Big\}.
\end{align*}
Hence
\begin{align*}
C_1
&=24\Big\{1-\frac{(1+(-1)^{\lambda})(1+(-1)^{\delta+\mu+\frac{g_1-1}{2}})}{4}\Big(\frac{3}{g_1}\Big)\Omega_2\Big\}\\
&\qquad+36\widehat{\mu}\cdot\Re\big(\omega^{1+\frac{g_1^2-1}{3}}\big)\Big\{\Omega_3-(-1)^{\delta+\mu+\frac{g_1-1}{2}}\Big(\frac{3}{g_1}\Big)\Omega_6\Big\}\\
&\qquad+36\widehat{\mu}\cdot\Re\big(\omega^{2-\lambda+\frac{g_1^2-1}{3}}\big)\Big\{\Omega_3-(-1)^{\delta+\mu+\frac{g_1-1}{2}}\Big(\frac{3}{g_1}\Big)\Omega_6\Big\}.
\end{align*}
Moreover, we may also deduce that
\begin{align*}
C_2
&=24\Big\{1-\frac{(1+(-1)^{\lambda})(1+(-1)^{\delta+\mu+\frac{g_1-1}{2}})}{4}\Big(\frac{3}{g_1}\Big)\Omega_2\Big\}\\
&\qquad+36\widehat{\mu}\cdot\Re\big(\omega^{2-\frac{g_1^2-1}{3}}\big)\Big\{\Omega_3-(-1)^{\delta+\mu+\frac{g_1-1}{2}}\Big(\frac{3}{g_1}\Big)\Omega_6\Big\}\\
&\qquad+36\widehat{\mu}\cdot\Re\big(\omega^{1+\lambda-\frac{g_1^2-1}{3}}\big)\Big\{\Omega_3-(-1)^{\delta+\mu+\frac{g_1-1}{2}}\Big(\frac{3}{g_1}\Big)\Omega_6\Big\}.
\end{align*}

In total we have
\begin{align*}
C
&=72\Big\{1-\frac{(1+(-1)^{\lambda})(1+(-1)^{\delta+\mu+\frac{g_1-1}{2}})}{4}\Big(\frac{3}{g_1}\Big)\Omega_2\Big\}\\
&\qquad+72\widehat{\mu}\cdot \Re\big(\omega^{1-\lambda}+\omega^{1+\frac{g_1^2-1}{3}}+\omega^{2-\lambda+\frac{g_1^2-1}{3}}\big)\\
&\qquad\ \ \times\Big\{\Omega_3-\frac{(1+(-1)^{\lambda})(1+(-1)^{\delta+\mu+\frac{g_1-1}{2}})}{4}\Big(\frac{3}{g_1}\Big)\Omega_6\Big\},
\end{align*}
which is proven for $r\equiv13\bmod{18},$ and is in fact valid for all $r\equiv1\bmod{6}$.

As argued before, the above computation can be modified to give
\begin{align*}
C
&=72\Big\{1+\frac{(1+(-1)^{\lambda})(1+(-1)^{\delta+\mu+\frac{g_1-1}{2}})}{4}\Big(\frac{3}{g_1}\Big)\Omega_2\Big\}\\
&\qquad+72\widehat{\mu}\cdot \Re\big(\omega^{1-\lambda}+\omega^{1+\frac{g_1^2-1}{3}}+\omega^{2-\lambda+\frac{g_1^2-1}{3}}\big)\\
&\qquad\ \ \times\Big\{\Omega_3+\frac{(1+(-1)^{\lambda})(1+(-1)^{\delta+\mu+\frac{g_1-1}{2}})}{4}\Big(\frac{3}{g_1}\Big)\Omega_6\Big\}
\end{align*}
for all $r\equiv-1\bmod{6}$.

We are done!

%% file: LT.bbl
\begin{bibsection}

\begin{biblist}



\bibitem[An52]{An52}
N. C. Ankeny, 
Representations of primes by quadratic forms,
\emph{Amer. J. Math.} \textbf{74} (1952), 913--919. 


\bibitem[ADJ04]{ADJ04}
A. Akbary, C. David \& R. Juricevic, 
Average distributions and products of special values of L-series,
\emph{Acta Arith.} \textbf{111} (2004), 239--268.

\bibitem[AP19]{AP19}
A. Akbary \& J. Parks, 
On the Lang--Trotter conjecture for two elliptic curves,
\emph{Ramanujan J.} \textbf{49} (2019), 585--623.


\bibitem[Ba66]{Ba66}
A. Baker, 
Linear forms in logarithms of algebraic numbers, 
\emph{Mathematika} \textbf{13} (1966), 204--216.



\bibitem[BGHT11]{BGHT11}
T. Barnet-Lamb, D. Geraghty, M. Harris \& R. Taylor, 
A family of Calabi--Yau varieties
and potential automorphy II,
\emph{Publ. Res. Inst. Math. Sci.} \textbf{47} (2011), 29--98.


\bibitem[BH62]{BH62}
P. T. Bateman, \& R. A. Horn, 
A heuristic asymptotic formula concerning the distribution of prime numbers,
\emph{Math. Comp.} \textbf{16} (1962), 363--367. 



\bibitem[BE81]{BE81}
B. C. Berndt \& R. J. Evans,
The determination of Gauss sums,
\emph{Bull. Amer. Math. Soc.} \textbf{5} (1981), 107--129.



\bibitem[BEW98]{BEW98}
B. C. Berndt, R. J. Evans \& K. S. Williams, 
Gauss and Jacobi Sums,
Wiley, New York, 1998.





\bibitem[BD20]{BD20} 
R. de la Bret\'eche \& S. Drappeau,
Niveau de r\'epartition des polyn\^omes quadratiques et crible majorant pour les entiers friables,
\emph{J. Eur. Math. Soc.} \textbf{22} (2020), 1577--1624.



\bibitem[BCDT01]{BCDT01} 
C. Breuil, B. Conrad, F. Diamond \& R. Taylor, 
On the modularity of elliptic curves over $\Q$: wild 3-adic exercises, 
\emph{J. Amer. Math. Soc.} \textbf{14} (2001), 843--939.


\bibitem[CP22]{CP22} 
F. Campagna \& R. Pengo,
Entanglement in the family of division fields of elliptic curves with complex multiplication,
\emph{Pacific J. Math.} \textbf{317} (2022), 21--66.


\bibitem[CJS20]{CJS20} 
H. Chen, N. Jones \& V. Serban,
The Lang--Trotter Conjecture for products of non-CM elliptic curves,
arXiv:2006.11269 [math.NT].



\bibitem[Ch49]{Ch49} 
S. Chowla, 
The last entry in Gauss's diary, 
\emph{Proc. Nat. Acad. Sci. U.S.A.} \textbf{35} (1949), 244--246. 


\bibitem[CHT08]{CHT08}
L. Clozel, M. Harris \& R. Taylor,
Automorphy for some $\ell$-adic lifts of automorphic mod $\ell$ Galois representations, With Appendix A, summarizing unpublished work of Russ Mann, and Appendix B by Marie-France Vignéras,
\emph{Publ. Math. IH\'ES} \textbf{108} (2008), 1--181.

\bibitem[CDSS17]{CDSS17} 
A. C. Cojocaru, R. Davis, A. Silverberg \& K. E. Stange, 
Arithmetic properties of the Frobenius traces defined by a rational abelian variety (with two appendices by J-P. Serre),
\emph{IMRN} 2017, no. 12, 3557--3602.


\bibitem[CFM05]{CFM05} 
A. C. Cojocaru, \'E. Fouvry \& M. R. Murty,
The square sieve and the Lang--Trotter conjecture,
\emph{Canad. J. Math.} \textbf{57} (2005), 1155--1177.


\bibitem[DH34]{DH34}
H. Davenport \& H. Hasse, 
Die Nullstellen der Kongruenzzetafunktionen in gewissen zyklischen Fallen, 
\emph{J. Reine Angew. Math.} \textbf{172} (1934), 151--182. 



\bibitem[DKS17]{DKS17}
C. David, D. Koukoulopoulos \& E. Smith,
Sums of Euler products and statistics of elliptic curves,
\emph{Math. Ann.} \textbf{368} (2017), 685--752. 


\bibitem[DP99]{DP99}
C. David \& F. Pappalardi,
Average Frobenius distributions of elliptic curves,
\emph{IMRN} 1999, 165--183.


\bibitem[De80]{De80}
P. Deligne,
La conjecture de Weil II,
\emph{Publ. Math. IH\'ES} \textbf{52} (1980), 137--252.



\bibitem[De41]{De41} 
M. Deuring, 
Die Typen der Multiplikatorenringe elliptischer Functionenk\"orper, 
\emph{Abh. Math. Sem. Hamburg} \textbf{14} (1941), 197--272. 








\bibitem[El87]{El87} 
N. D. Elkies, The existence of infinitely many supersingular primes for every elliptic curve over $\Q$, 
\emph{Invent. math.} \textbf{89} (1987), 561--568. 


\bibitem[Elthesis]{Elthesis} 
N. D. Elkies,
Supersingular primes of a given elliptic curve over a number field,
Ph.D. thesis in Harvard University, 1987. 41 pp. UMI Number: 8800772.



\bibitem[El92]{El92} 
N. D. Elkies,
Distribution of supersingular primes, Journ\'ees Arithm\'etiques, 1989 (Luminy, 1989). 
\emph{Ast\'erisque} No. \textbf{198-200} (1991), 127--132 (1992).


\bibitem[FM95]{FM95}
\'E. Fouvry \& M. R. Murty,
Supersingular primes common to two elliptic curves,
Number theory (Paris, 1992–1993), 91–102, 
London Math. Soc. Lecture Note Ser., 215, Cambridge Univ. Press, Cambridge, 1995. 

\bibitem[FM96]{FM96}
\'E. Fouvry \& M. R. Murty,
On the distribution of supersingular primes,
\emph{Canad. J. Math.} \textbf{48} (1996), 81--104. 




\bibitem[GaDA]{GaDA} 
C. F. Gau{\ss}, 
Disquisitiones Arithmeticae,
Fleischer, Leipzig 1801. 
Translation by A. Clarke: Disquisitiones Aritmeticae. Yale University Press, 1965.

\bibitem[GaII]{GaII} 
C. F. Gau{\ss}, 
Werke. Band II. Reprint of the 1863 original. Georg Olms Verlag, Hildesheim, 1973. viii+504 pp.

\bibitem[GaX1]{GaX1} 
C. F. Gauss,
Werke, Band $X_1$. K\"on. Geselsch. Wissensch. G\"ottingen, Teubner, Leipzig 1917. Abdruch Tagebuch, pp. 483--572.

\bibitem[Gr84]{Gr84} 
J. Gray, 
A commentary on Gauss’s mathematical diary, 1796–1814, with an English translation. 
\emph{Expositiones Math.} \textbf{2} (1984), 97--130.

\bibitem[Gr80]{Gr80} 
B. H. Gross, Arithmetic on elliptic curves with complex multiplication, 
Lect. Notes in Math. 776, Springer, Berlin, 1980.



\bibitem[HL23]{HL23} 
G. H. Hardy \& J. E. Littlewood,
Some problems of ``Partitio Numerorum": III. On the expression of a number as a sum of primes, 
\emph{Acta Math.} \textbf{44} (1923), 1--70.


\bibitem[HL01]{HL01} 
G. Harman \& P. A. Lewis, 
Gaussian primes in narrow sectors, 
\emph{Mathematika} \textbf{48} (2001), 119--135.


\bibitem[Ha09]{Ha09}
M. Harris, 
Potential automorphy of odd-dimensional symmetric powers of elliptic curves and applications,
Algebra, arithmetic, and geometry: in honor of Yu. I. Manin. Vol. II, 1–21, Progr. Math., 270, Birkh\"auser Boston, Boston, MA, 2009.

\bibitem[HST10]{HST10}
M. Harris, N. Shepherd-Barron \& R. Taylor,
A family of Calabi--Yau varieties and potential automorphy,
\emph{Ann. of Math. (2)} \textbf{171} (2010), 779--813. 



\bibitem[He18]{He18}
E. Hecke, 
Eine neue Art von Zetafunktionen und ihre Beziehungen zur Verteilung der Primzahlen,
\emph{Math. Z.} \textbf{1} (1918), 357--376.


\bibitem[He21]{He21}
G. Herglotz,
Zur letzten Eintragung im Gau{\ss}schen Tagebuch,
\emph{Berichte \"uer die Verhandlungen der S\"ahsischen Akademie der Wissenschaften Math.-phys. Klasse} \textbf{73} (1921), 271--276. Reprinted in \emph{Gesammelte Schriften}, ed. H. Schwerdtfeger, pp. 415--420. G\"otingen: Vandenhoeck 1979.




\bibitem[IR90]{IR90} 
K. Ireland \& M. Rosen,
A Classical Introduction to Modern Number Theory, 2nd edition, 
Graduate Texts in Mathematics 84. New York: Springer--Verlag, 1990.

\bibitem[Iw78]{Iw78} 
H. Iwaniec, Almost-primes represented by quadratic polynomials, 
\emph{Invent. math.} \textbf{47} (1978), 171--188.


\bibitem[Iw80]{Iw80} 
H. Iwaniec, 
A new form of the error term in the linear sieve, 
\emph{Acta Arith.} \textbf{37} (1980), 307--320.


\bibitem[Ja16]{Ja16} 
K. James,
Variants of the Sato--Tate and Lang--Trotter Conjectures, Frobenius Distributions: Lang--Trotter and Sato--Tate Conjectures,
\emph{Contemp. Math.} \textbf{663}, 2016.


\bibitem[JP17]{JP17} 
K. James \& P. Pollack, Extremal primes for elliptic 
curves with complex multiplication, 
\emph{J. Number Theory} \textbf{172} (2017), 383--391. 

\bibitem[JS11]{JS11} 
K. James \& E. Smith, 
Average Frobenius distribution for elliptic curves defined over finite Galois extensions of the rationals,
\emph{Math. Proc. Cambridge Philos. Soc.} \textbf{150} (2011), 439--458. 

\bibitem[JS13]{JS13} 
K. James \& E. Smith,
Average Frobenius distribution for the degree two primes of a number field. 
\emph{Math. Proc. Cambridge Philos. Soc.} \textbf{154} (2013), 499--525. 


\bibitem[JQ14]{JQ14}
 Q. Ji \& H. Qin, 
CM elliptic curves and prime captured by quadratic polynomials, 
\emph{Asian J. Math.} \textbf{18} (2014), 707--726. 


\bibitem[JU08]{JU08}
J. Jim\'enez Urroz, 
Almost prime orders of CM elliptic curves modulo $p$,
Algorithmic Number Theory, pp. 74--87, Lecture Notes in Computer Science 5011, Springer-Verlag, Berlin, Heidelberg, 2008.


\bibitem[Jo09]{Jo09} 
N. Jones, 
Average of elliptic curve constants, 
\emph{Math. Ann.} \textbf{345} (2009), 685--710. 

\bibitem[Ka09]{Ka09} 
N. Katz, 
Lang--Trotter revisited, 
\emph{Bull. Amer. Math. Soc.} \textbf{46} (2009), 413--457. 


\bibitem[Kl03]{Kl03} 
F. Klein, 
Gau{\ss}’ wissenschafftliches Tagebuch 1796--1814, \emph{Math. Annalen} \textbf{57} (1903), 1--34. 


\bibitem[Ku55]{Ku55} 
J. Kubilius,
On a problem in the $n$-dimensional analytic theory of numbers (Lithuanian. Russian summary),
\emph{Vilniaus Valst. Univ. Mokslo Darbai. Mat. Fiz. Chem. Mokslu Ser.} \textbf{4} (1955), 5--43.



\bibitem[LT76]{LT76}
S. Lang \& H. Trotter, 
Frobenius Distributions in $\GL_2$-Extensions, 
Lecture Notes in Mathematics 504. Berlin: Springer, 1976.


\bibitem[Le00]{Le00}
F. Lemmermeyer,
Reciprocity Laws: From Euler to Eisenstein, 
Springer Monographs in Mathematics. Berlin: Springer--Verlag, 2000.


\bibitem[Lo17]{Lo17}
D. Lombardo, 
Galois representations attached to abelian varieties of CM type,
\emph{Bull. Soc. Math. France} \textbf{145} (2017), 469--501.


\bibitem[LR22]{LR22}
\'A. Lozano-Robledo,
Galois representations attached to elliptic curves with complex multiplication,
\emph{Algebra Number Theory} \textbf{16} (2022), 777--837.

\bibitem[Ma72]{Ma72} 
B. Mazur, 
Rational points of abelian varieties with values in towers of number fields, 
\emph{Invent. math.} \textbf{18} (1972), 183--266. 


\bibitem[Mi98]{Mi98} 
W. Miller, 
Counting points on certain CM elliptic curves modulo primes, 
Ph.D. thesis in UCSD, 1998. 74pp. UMI Number: 9906470.



\bibitem[MMS88]{MMS88}
M. R. Murty, V. K. Murty \& N. Saradha, 
Modular forms and the Chebotarev density theorem, 
\emph{Amer. J. Math.} \textbf{110} (1988), 253--281. 




\bibitem[Mu85]{Mu85} 
V. K. Murty, 
Explicit formulae and the Lang--Trotter conjecture, 
\emph{Rocky Mountain J. Math.} \textbf{15} (1985), 535--551. 

\bibitem[Mu97]{Mu97} 
V. K. Murty, 
Modular forms and the Chebotarev density theorem. II.  
In Analytic Number Theory (Kyoto, 1996), 287--306, 
\emph{London Math. Soc. Lecture Note Ser.} \textbf{247}, 
Cambridge Univ. Press, 1997. 




\bibitem[Oo16]{Oo16} 
F. Oort, 
Gauss: the last entry,
\emph{ICCM Not.} \textbf{4} (2016), 29--39. 



\bibitem[Po16]{Po16}
P. Pollack, 
A Titchmarsh divisor problem for elliptic curves, 
\emph{Math. Proc. Cambridge Philos. Soc.} \textbf{160} (2016), 167--189.  


\bibitem[Qi16]{Qi16} 
H. Qin, 
Anomalous primes of the elliptic curve $E_D: y^2=x^3+D$, 
\emph{Proc. London Math. Soc.} \textbf{112} (2016), 415--453.  

\bibitem[Qi21]{Qi21} 
H. Qin, 
The Lang--Trotter Conjecture for the elliptic curve $y^2=x^3+Dx$, 
arXiv:2108.06292 [math.NT].



\bibitem[Ra68]{Ra68} 
A. R. Rajwade, 
Arithmetic on curves with complex multiplication by $\sqrt{-2}$,
\emph{Proc. Camb. Phil. Soc.} \textbf{64} (1968), 659--672.


\bibitem[Ra70]{Ra70} 
A. R. Rajwade, 
A note on the number of solutions $N_p$ of the congruence $y^2\equiv x^3-Dx\bmod p$, 
\emph{Proc. Cambridge Philos. Soc.} \textbf{67} (1970), 603--605.

\bibitem[Ra77]{Ra77} 
A. R. Rajwade, 
The Diophantine equation $y^2=x(x^2+21Dx+112D^2)$ and the conjectures
of Birch and Swinnerton-Dyer, 
\emph{J. Austral. Math. Soc. Ser. A} \textbf{24} (1977), 286--295.

\bibitem[RP81]{RP81} 
A. R. Rajwade \& J. C. Parnami, 
A new cubic character sum, 
\emph{Acta Arith.} \textbf{40} (1981/82), 347--356.

\bibitem[Ra23]{Ra23} 
A. Ray,
On the constants of the Lang-Trotter conjecture for CM elliptic curves, 
arXiv:2309.09938 [math.NT].

\bibitem[RPR84]{RPR84} 
D. B. Rishi, J. C. Parnami \& A. R. Rajwade, 
Evaluation of a cubic character sum using the $\sqrt{-19}$ division points of the curve $Y^2=X^3-2^3\cdot 19X+2\cdot 19^2$, 
\emph{J. Number Theory} \textbf{19} (1984), 184--194.


\bibitem[RV01]{RV01} 
A. Reverter \& N. Vila,
Images of mod $p$ Galois representations associated to elliptic curves,
\emph{Canad. Math. Bull.} \textbf{44} (2001), 313--322.


\bibitem[RZ15]{RZ15} 
J. Rouse, D. Zureick-Brown, 
Elliptic curves over $\Q$ and 2-adic images of Galois,
\emph{Res. Number Theory} \textbf{1} (2015), Paper No. 12, 34 pp.


\bibitem[RS09]{RS09} 
K. Rubin \& A. Silverberg, 
Point counting on reductions of CM elliptic curves,
\emph{J. Number Theory} \textbf{129} (2009), 2903--2923.

\bibitem[SS58]{SS58} 
A. Schinzel \& W. Sierpi\'nski, 
Sur certaines hypoth\`eses concernant les nombres premiers,
\emph{Acta Arith.} \textbf{4} (1958), 185--208. 



\bibitem[Se72]{Se72}
J.-P. Serre, 
Propri\'et\'es galoisiennes des points d'ordre fini des courbes elliptiques,
\emph{Invent. math.} \textbf{15} (1972), 259--331.

\bibitem[Se81]{Se81}
J.-P. Serre, 
Quelques applications du th\'eor\`eme de densit\'e de Chebotarev,
\emph{Publ. Math. IH\'ES} \textbf{54} (1981), 323--401.


\bibitem[Si10]{Si10} 
A. Silverberg, 
Group order formulas for reductions of CM elliptic curves,
Arithmetic, geometry, cryptography and coding theory 2009, 107–120, Contemp. Math., 521, Amer. Math. Soc., Providence, RI, 2010.



\bibitem[Si94]{Si94} 
J. Silverman, 
Advanced Topics in the Arithmetic of Elliptic Curves, 
Graduate Texts in Mathematics 151. New York: Springer--Verlag, 1990.



\bibitem[St67]{St67}
H. M. Stark, 
A complete determination of the complex quadratic fields of class number one, 
\emph{Michigan Math. J.} \textbf{14} (1967), 1--27. 

\bibitem[St96]{St96} 
H. M. Stark, 
Counting points on CM elliptic curves,
\emph{Rocky Mountain J. Math.} \textbf{26} (1996), 1115--1138.



\bibitem[Su01]{Su01} 
Z.-H. Sun, Supplements to the theory of quartic residues,
\emph{Acta Arith.} \textbf{XCVII.4} (2001), 361--377.

\bibitem[Su16]{Su16} 
A. Sutherland, 
Computing images of Galois representations attached to elliptic curves,
\emph{Forum Math. Sigma} \textbf{4} (2016), Paper No. e4, 79 pp. 


\bibitem[Ta08]{Ta08}
R. Taylor,
Automorphy for some $\ell$-adic lifts of automorphic mod $\ell$ Galois representations. II,
\emph{Publ. Math. IH\'ES} \textbf{108} (2008), 183--239. 



\bibitem[TW95]{TW95} 
R. Taylor \& A. Wiles, 
Ring-theoretic properties of certain Hecke algebras,
\emph{Ann. of Math. (2)} \textbf{141} (1995), 
553--572.

\bibitem[TZ18]{TZ18}
J. Thorner \& A. Zaman, 
A Chebotarev variant of the Brun--Titchmarsh theorem and bounds for the Lang--Trotter conjectures,
\emph{Int. Math. Res. Not.} 2018, no. 16, 4991--5027. 





\bibitem[Wa90]{Wa90}
D. Wan,
On the Lang--Trotter conjecture, 
\emph{J. Number Theory} \textbf{35} (1990), 247--268.

\bibitem[Wi95]{Wi95} 
A. Wiles, 
Modular elliptic curves and Fermat's last theorem,
\emph{Ann. of Math. (2)} \textbf{141} (1995), 443--551.


\bibitem[Zy15a]{Zy15a}
D. Zywina,
On the possible images of the mod ell representations associated to elliptic curves over $\Q$, 
arXiv:1508.07660 [math.NT].


\bibitem[Zy15b]{Zy15b} 
D. Zywina,
Bounds for the Lang--Trotter conjecture, SCHOLAR--a Scientific Celebration Highlighting Open Lines of Arithmetic Research, \emph{Contemp. Math.} \textbf{655}, 2015.



\end{biblist}
\end{bibsection}
